# Introduction into Geometry over Division Ring

Aleks Kleyn

*E-mail address*: Aleks_Kleyn@MailAPS.org
*URL*: http://sites.google.com/site/alekskleyn/
*URL*: http://arxiv.org/a/kleyn_a_1
*URL*: http://AleksKleyn.blogspot.com/

ABSTRACT. Theory of representations of universal algebra is a natural development of the theory of universal algebra. Exploring of morphisms of the representation leads to the concepts of generating set and basis of representation. In the book I considered the notion of tower of $T\star$-representations of $\Omega_i$-algebras, $i = 1, ..., n$, as the set of coordinated $T\star$-representations of $\Omega_i$-algebras.

I explore the geometry of affine space as example of tower of representations. Exploration of curvilinear coordinates allows us to see how look at the main structures of the manifold with affine connection. I explore Euclidean space over division ring.

# Contents

CHAPTER 1

# Preface

## 1.1. Tower of Representations

The main goal of this book is considering of simple geometry over division ring.

However I decided not to bind this book by description of affine and Euclidean geometry. During description of affine geometry I discovered interesting structure. At first I define $D_*{}^*$-vector space $\overline{V}$. This is $T\star$-representation of the division ring $D$ in additive group. Then I consider a set of points where the representation of $D_*{}^*$-vector space $\overline{V}$ is defined. This many-tier structure stimulate my interest and I returned to exploration of $T\star$-representations of universal algebra.

Thus the concept of tower of $T\star$-representations emerged. As soon as I began draw diagrams related to the tower of $T\star$-representations, I discovered pattern similar to diagrams in paper [11]. Although I did not succeed to carry this analogy to its logical conclusion, I hope that this analogy can bring to interesting results.

## 1.2. Morphism and Basis of Representation

Exploring of the theory of representations of $\mathfrak{F}$-algebra shows that this theory has a lot of common with theory of $\mathfrak{F}$-algebra. In [6, 7, 8], I explored morphisms of representation. They are maps that preserve the structure of representation of $\mathfrak{F}$-algebra (or of tower of representations).

Linear maps of vector spaces are examples of morphism of representations. The definition of basis is important structure in the theory of linear spaces. Concept of basis is directly linked to statement that the group of morphisms of linear space has two single transitive representations. The main function of basis is to form linear space.

The natural question arises. Can we generalize this structure to arbitrary representation? The basis is not the only set that forms the vector space. This statement is initial point where I started exploring of generating set of representation. Generating set of representation is one more interesting parallel between theory of representations and theory of universal algebra.

For some representations, the possibility to find generating set is not equivalent to the possibility to build basis. The problem to find the set of representations that allow to construct basis is interesting and important problem. I wonder is it possible to find effective algorithm of construction of basis of quantum geometry.

## 1.3. $D$-Affine space

Like in the case of $D$-vector space, I can start from considering of $D_*{}^*$-affine space and then move to exploring of $D$-affine space. However I do not expect something new on this way. Because I explore operations related to differentiation, I start immediately from exploring of $D$-affine space.

From the moment when I realized that affine space is tower of representations and linear maps of affine space are morphisms of tower of representations, I realized what form the connection has in manifold with $D$-affine connection. However, I decided to explore curvilinear coordinates in the affine space before I explore the manifold with affine connections.

Although it was clear in my mind what happens when I change coordinates, the picture that I saw appears more interesting then I expected. I have met two very important ideas.

I started from construction of affine space as tower of representations. it turned out, that there exists one more model where we represent vector fields as sets of maps.[1.1] However these two models do not contradict one to another, but rather are complementary.

Second observation did not come as a surprise for me. When I started to explore $D$-vector space, I was ready to consider linear combination of vectors as polylinear form; however I got the conclusion that this does not change dimension of space, and for this reason I did not go beyond exploration of $D_*{}^*$-basis.

It is interesting to observe how new concept appears during writing of book. Construction that I made in chapter 8 again returned me to a linear combination as sum of linear forms. However at this time I have linear dependence of 1-forms. And coefficients of linear dependence are not outside of 1-form, but inside. If we consider 1-forms over field, then 1-form is transparent for scalars of the field, and there is no difference for us whether coefficients are outside or inside the form.

Linear dependence is not subject of this book. I will explore this question in following papers. In this book I use this definition assuming that its properties are clear from written equations.

Exploring of curvilinear coordinates also gave opportunity to see transformation rule of coordinates of $D$-affine connection when we change coordinates. However, I think it was luck that I was able to write the explicit transformation law, and this luck was associated with the fact that the connection has only one contravariant index. I have reason to believe that such an operation, like raising or lowering the index of the tensor in Riemannian space can have only an implicit form of record.

For years affine geometry was symbol for me of the most simple geometry. I am glad to see that I was wrong. At the right moment affine geometry turned out to be a source of inspiration. Affine geometry is the gate of differential geometry. Although, in order to plunge deeply into the differential geometry, I need good knowledge about differential equations, I conclude the book with very brief excursion into the manifold with affine connections. It can give to reader the taste of the new geometry.

## 1.4. Conventions

(1) Function and mapping are synonyms. However according to tradition, correspondence between either rings or vector spaces is called mapping and a mapping of either real field or quaternion algebra is called function.
(2) We can consider division ring $D$ as $D$-vector space of dimension 1. According to this statement, we can explore not only homomorphisms of division ring $D_1$ into division ring $D_2$, but also linear maps of division rings. This

[1.1] At this time it is not clear whether this model leads me to geometry that Alain Connes explores in [14]. Answer on this question requires further research.

means that map is multiplicative over maximum possible field. In particular, linear map of division ring $D$ is multiplicative over center $Z(D)$. This statement does not contradict with definition of linear map of field because for field $F$ is true $Z(F) = F$. When field $F$ is different from maximum possible, I explicit tell about this in text.

(3) In spite of noncommutativity of product a lot of statements remain to be true if we substitute, for instance, right representation by left representation or right vector space by left vector space. To keep this symmetry in statements of theorems I use symmetric notation. For instance, I consider $D\star$-vector space and $\star D$-vector space. We can read notation $D\star$-vector space as either D-star-vector space or left vector space. We can read notation $D\star$-linear dependent vectors as either D-star-linear dependent vectors or vectors that are linearly dependent from left.

(4) Let $A$ be free finite dimensional algebra. Considering expansion of element of algebra $A$ relative basis $\overline{\overline{e}}$ we use the same root letter to denote this element and its coordinates. However we do not use vector notation in algebra. In expression $a^2$, it is not clear whether this is component of expansion of element $a$ relative basis, or this is operation $a^2 = aa$. To make text clearer we use separate color for index of element of algebra. For instance,

$$a = a^{i}\overline{e}_{i}$$

(5) If free finite dimensional algebra has unit, then we identify the vector of basis $\overline{e}_0$ with unit of algebra.

(6) In [12], an arbitrary operation of algebra is denoted by letter $\omega$, and $\Omega$ is the set of operations of some universal algebra. Correspondingly, the universal algebra with the set of operations $\Omega$ is denoted as $\Omega$-algebra. Similar notations we see in [2] with small difference that an operation in the algebra is denoted by letter $f$ and $\mathcal{F}$ is the set of operations. I preferred first case of notations because in this case it is easier to see where I use operation.

(7) Since the number of universal algebras in the tower of representations is varying, then we use vector notation for a tower of representations. We denote the set $(A_1, ..., A_n)$ of $\Omega_i$-algebras $A_i$, $i = 1, ..., n$ as $\overline{A}$. We denote the set of representations $(f_{1,2}, ..., f_{n-1,n})$ of these algebras as $\overline{f}$. Since different algebras have different type, we also talk about the set of $\overline{\Omega}$-algebras. In relation to the set $\overline{A}$, we also use matrix notations that we discussed in section [6]-2.1. For instance, we use the symbol $\overline{A}_{[1]}$ to denote the set of $\overline{\Omega}$-algebras $(A_2, ..., A_n)$. In the corresponding notation $(\overline{A}_{[1]}, \overline{f})$ of tower of representation, we assume that $\overline{f} = (f_{2,3}, ..., f_{n-1,n})$.

(8) Since we use vector notation for elements of the tower of representations, we need convention about notation of operation. We assume that we get result of operation componentwise. For instance,

$$\overline{r}(\overline{a}) = (r_1(a_1), ..., r_n(a_n))$$

(9) Let $A$ be $\Omega_1$-algebra. Let $B$ be $\Omega_2$-algebra. Notation

$$A \xrightarrow{*} B$$

means that there is representation of $\Omega_1$-algebra $A$ in $\Omega_2$-algebra $B$.

(10) Without a doubt, the reader may have questions, comments, objections. I will appreciate any response.

CHAPTER 2

# Representation of Universal Algebra

## 2.1. Representation of Universal Algebra

**Definition 2.1.1.** Suppose[2.1] we defined the structure of $\Omega_2$-algebra on the set $M$ ([2, 12]). Endomorphism of $\Omega_2$-algebra

$$t : M \to M$$

is called **transformation of universal algebra** $M$.[2.2] □

We denote $\delta$ identical transformation.

**Definition 2.1.2.** If transformation acts from left

$$u' = tu$$

it is called **left-side transformation** or $T\star$**-transformation** We denote ${}^\star M$ the set of $T\star$-transformations of universal algebra $M$. □

**Definition 2.1.3.** If transformation acts from right

$$u' = ut$$

it is called **right-side transformation** or $\star T$**-transformation** We denote $M^\star$ the set of nonsingular $\star T$-transformations of universal algebra $M$. □

**Definition 2.1.4.** Suppose we defined the structure of $\Omega_1$-algebra on the set ${}^\star M$ ([2]). Let $A$ be $\Omega_1$-algebra. We call homomorphism

$$f : A \to {}^\star M \tag{2.1.1}$$

**left-side** or $T\star$**-representation of** $\Omega_1$**-algebra** $A$ **in** $\Omega_2$**-algebra** $M$ □

**Definition 2.1.5.** Suppose we defined the structure of $\Omega_1$-algebra on the set $M^\star$ ([2]). Let $A$ be $\Omega_1$-algebra. We call homomorphism

$$f : A \to M^\star$$

**right-side** or $\star T$**-representation of** $\Omega_1$**-algebra** $A$ **in** $\Omega_2$**-algebra** $M$ □

We extend to representation theory convention described in remark [6]-2.2.15. We can write duality principle in the following form

**Theorem 2.1.6** (duality principle)**.** *Any statement which holds for $T\star$-representation of $\Omega_1$-algebra $A$ holds also for $\star T$-representation of $\Omega_1$-algebra $A$.*

---

[2.1] Definitions and theorems from this chapter are included also into chapters [10]-2, [10]-3.

[2.2] If the set of operations of $\Omega_2$-algebra is empty, then $t$ is a mapping.

**Remark 2.1.7.** There exist two forms of notation for transformation of $\Omega_2$-algebra $M$. In operational notation, we write the transformation $A$ as either $Aa$ which corresponds to the $T\star$-transformation or $aA$ which corresponds to the $\star T$-transformation. In functional notation, we write the transformation $A$ as $A(a)$ regardless of the fact whether this is $T\star$-transformation or this is $\star T$-transformation. This notation is in agreement with duality principle.

This remark serves as a basis for the following convention. When we use functional notation we do not make a distinction whether this is $T\star$-transformation or this is $\star T$-transformation. We denote ${}^*M$ the set of transformations of $\Omega_2$-algebra $M$. Suppose we defined the structure of $\Omega_1$-algebra on the set ${}^*M$. Let $A$ be $\Omega_1$-algebra. We call homomorphism

$$f : A \to {}^*M \tag{2.1.2}$$

**representation of $\Omega_1$-algebra $A$ in $\Omega_2$-algebra $M$**. We also use record

$$f : A \xrightarrow{*} M$$

to denote the representation of $\Omega_1$-algebra $A$ in $\Omega_2$-algebra $M$.

Correspondence between operational notation and functional notation is unambiguous. We can select any form of notation which is convenient for presentation of particular subject. □

There are several ways to describe the representation. We can define the mapping $f$ keeping in mind that the domain is $\Omega_1$-algebra $A$ and range is $\Omega_1$-algebra ${}^*M$. Either we can specify $\Omega_1$-algebra $A$ and $\Omega_2$-algebra $M$ keeping in mind that we know the structure of the mapping $f$.[2.3]

Diagram

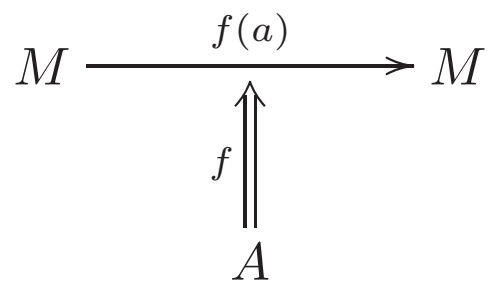

means that we consider the representation of $\Omega_1$-algebra $A$. The mapping $f(a)$ is image of $a \in A$.

**Definition 2.1.8.** Suppose the mapping (2.1.2) is an isomorphism of the $\Omega_1$-algebra $A$ into ${}^*M$. Then the representation of the $\Omega_1$-algebra $A$ is called **effective**. □

**Remark 2.1.9.** If the $T\star$-representation of $\Omega_1$-algebra is effective, then we identify an element of $\Omega_1$-algebra and its image and write $T\star$-transformation caused by element $a \in A$ as

$$v' = av$$

Suppose the $\star T$-representation of $\Omega_1$-algebra is effective. Then we identify an element of $\Omega_1$-algebra and its image and write $\star T$-transformation caused by element $a \in A$ as

$$v' = va$$

□

[2.3]For instance, we consider vector space $\overline{V}$ over field $D$ ( definition [6]-4.1.4 ).

**Definition 2.1.10.** We call a representation of $\Omega_1$-algebra **transitive** if for any $a, b \in V$ exists such $g$ that

$$a = f(g)(b)$$

We call a representation of $\Omega_1$-algebra **single transitive** if it is transitive and effective. □

**Theorem 2.1.11.** *$T\star$-representation is single transitive if and only if for any $a, b \in M$ exists one and only one $g \in A$ such that $a = f(g)(b)$*

PROOF. Corollary of definitions 2.1.8 and 2.1.10. □

## 2.2. Morphism of Representations of Universal Algebra

**Theorem 2.2.1.** *Let $A$ and $B$ be $\Omega_1$-algebras. Representation of $\Omega_1$-algebra $B$*

$$g : B \to {}^\star M$$

*and homomorphism of $\Omega_1$-algebra*

$$h : A \to B \tag{2.2.1}$$

*define representation $f$ of $\Omega_1$-algebra $A$*

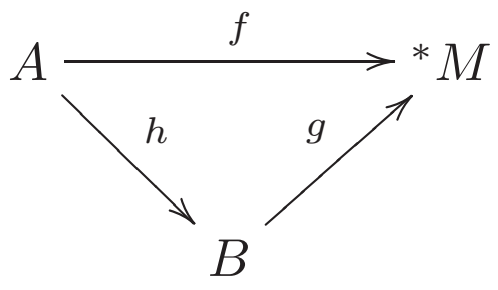

PROOF. Since mapping $g$ is homomorphism of $\Omega_1$-algebra $B$ into $\Omega_1$-algebra ${}^*M$, the mapping $f$ is homomorphism of $\Omega_1$-algebra $A$ into $\Omega_1$-algebra ${}^*M$. □

Considering representations of $\Omega_1$-algebra in $\Omega_2$-algebras $M$ and $N$, we are interested in a mapping that preserves the structure of representation.

**Definition 2.2.2.** Let

$$f : A \to {}^*M$$

be representation of $\Omega_1$-algebra $A$ in $\Omega_2$-algebra $M$ and

$$g : B \to {}^*N$$

be representation of $\Omega_1$-algebra $B$ in $\Omega_2$-algebra $N$. Tuple of mappings

$$(r : A \to B, R : M \to N) \tag{2.2.2}$$

such, that

- $r$ is homomorphism of $\Omega_1$-algebra
- $R$ is homomorphism of $\Omega_2$-algebra
-

$$R \circ f(a) = g(r(a)) \circ R \tag{2.2.3}$$

is called **morphism of representations from $f$ into $g$**. We also say that **morphism of representations of $\Omega_1$-algebra in $\Omega_2$-algebra** is defined. □

**Remark 2.2.3.** We may consider a pair of mappins $r$, $R$ as mapping

$$F : A \cup M \to B \cup N$$

such that

$$F(A) = B \quad F(M) = N$$

Therefore, hereinafter we will say that we have the mapping $(r, R)$. □

**Definition 2.2.4.** If representation $f$ and $g$ coincide, then morphism of representations $(r, R)$ is called **morphism of representation $f$**. □

For any $m \in M$ equation (2.2.3) has form

$$R(f(a)(m)) = g(r(a))(R(m)) \tag{2.2.4}$$

**Remark 2.2.5.** Consider morphism of representations (2.2.2). We denote elements of the set $B$ by letter using pattern $b \in B$. However if we want to show that $b$ is image of element $a \in A$, we use notation $r(a)$. Thus equation

$$r(a) = r(a)$$

means that $r(a)$ (in left part of equation) is image $a \in A$ (in right part of equation). Using such considerations, we denote element of set $N$ as $R(m)$. We will follow this convention when we consider correspondences between homomorphisms of $\Omega_1$-algebra and mappings between sets where we defined corresponding representations.

There are two ways to interpret (2.2.4)

- Let transformation $f(a)$ map $m \in M$ into $f(a)(m)$. Then transformation $g(r(a))$ maps $R(m) \in N$ into $R(f(a)(m))$.
- We represent morphism of representations from $f$ into $g$ using diagram

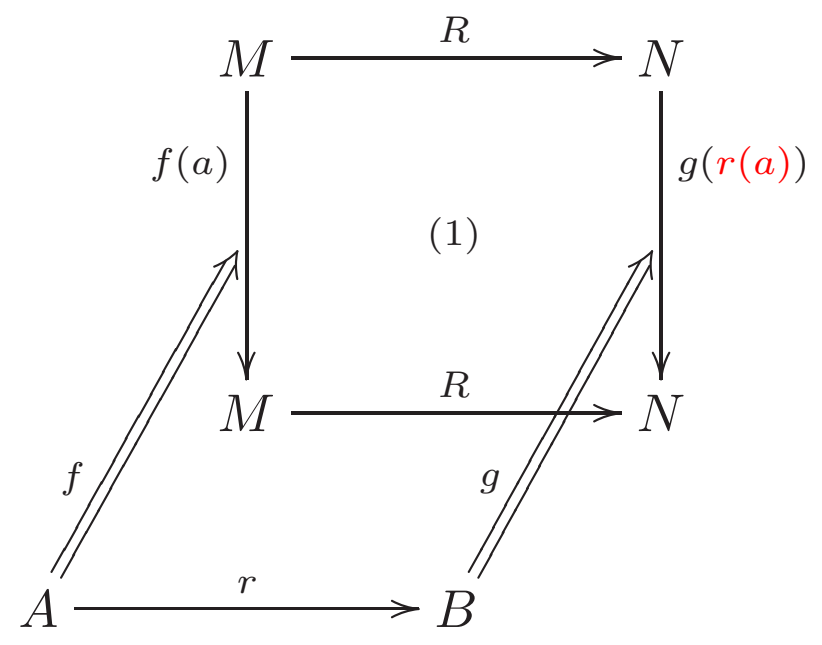

From (2.2.3), it follows that diagram (1) is commutative.

□

**Theorem 2.2.6.** *Consider representation*

$$f : A \to {}^*M$$

*of $\Omega_1$-algebra $A$ and representation*

$$g : B \to {}^*N$$

*of $\Omega_1$-algebra $B$. Morphism*

$$h : A \longrightarrow B \qquad H : M \longrightarrow N$$

*of representations from $f$ into $g$ satisfies equation*

$$H \circ \omega(f(a_1), ..., f(a_n)) = \omega(g(h(a_1)), ..., g(h(a_n))) \circ H \tag{2.2.5}$$

*for any $n$-ary operation $\omega$ of $\Omega_1$-algebra.*

PROOF. Since $f$ is homomorphism, we have

$$H\circ\omega(f(a_1),...,f(a_n))=H\circ f(\omega(a_1,...,a_n)) \tag{2.2.6}$$

From (2.2.3) and (2.2.6) it follows that

$$H\circ\omega(f(a_1),...,f(a_n))=g(h(\omega(a_1,...,a_n)))\circ H \tag{2.2.7}$$

Since $h$ is homomorphism, from (2.2.7) it follows that

$$H\circ\omega(f(a_1),...,f(a_n))=g(\omega(h(a_1),...,h(a_n)))\circ H \tag{2.2.8}$$

Since $g$ is homomorphism, (2.2.5) follows from (2.2.8). □

**Theorem 2.2.7.** *Let the mapping*

$$h:A\longrightarrow B\qquad H:M\longrightarrow N$$

*be morphism from representation*

$$f:A\rightarrow {}^*M$$

*of $\Omega_1$-algebra $A$ into representation*

$$g:B\rightarrow {}^*N$$

*of $\Omega_1$-algebra $B$. If representation $f$ is effective, then the mapping*

$${}^*H:{}^*M\rightarrow {}^*N$$

*defined by equation*

$${}^*H(f(a))=g(h(a)) \tag{2.2.9}$$

*is homomorphism of $\Omega_1$-algebra.*

PROOF. Because representation $f$ is effective, then for given transformation $f(a)$ element $a$ is determined uniquely. Therefore, transformation $g(h(a))$ is properly defined in equation (2.2.9).

Since $f$ is homomorphism, we have

$${}^*H(\omega(f(a_1),...,f(a_n)))={}^*H(f(\omega(a_1,...,a_n))) \tag{2.2.10}$$

From (2.2.9) and (2.2.10) it follows that

$${}^*H(\omega(f(a_1),...,f(a_n)))=g(h(\omega(a_1,...,a_n))) \tag{2.2.11}$$

Since $h$ is homomorphism, from (2.2.11) it follows that

$${}^*H(\omega(f(a_1),...,f(a_n)))=g(\omega(h(a_1),...,h(a_n))) \tag{2.2.12}$$

Since $g$ is homomorphism,

$${}^*H(\omega(f(a_1),...,f(a_n)))=\omega(g(h(a_1)),...,g(h(a_n)))=\omega({}^*H(f(a_1)),...,{}^*H(f(a_n)))$$

follows from (2.2.12). Therefore, the mapping ${}^*H$ is homomorphism of $\Omega_1$-algebra. □

**Theorem 2.2.8.** *Given single transitive representation*

$$f : A \to {}^*M$$

*of $\Omega_1$-algebra $A$ and single transitive $T\star$-representation*

$$g : B \to {}^*N$$

*of $\Omega_1$-algebra $B$, there exists morphism*

$$h : A \to B \quad H : M \to N$$

*of representations from $f$ into $g$.*

Proof. Let us choose homomorphism $h$. Let us choose element $m \in M$ and element $n \in N$. To define mapping $H$, consider following diagram

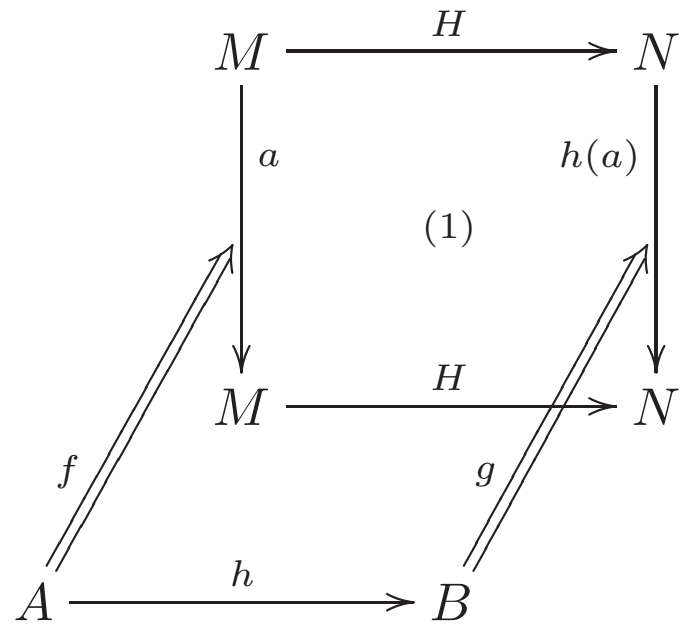

From commutativity of diagram (1), it follows that

$$H(am) = h(a)H(m)$$

For arbitrary $m' \in M$, we defined unambiguously $a \in A$ such that $m' = am$. Therefore, we defined mapping $H$ which satisfies to equation (2.2.3). □

**Theorem 2.2.9.** *Let*

$$f : A \to {}^*M$$

*be single transitive representation of $\Omega_1$-algebra $A$ and*

$$g : B \to {}^*N$$

*be single transitive representation of $\Omega_1$-algebra $B$. Given homomorphism of $\Omega_1$-algebra*

$$h : A \longrightarrow B$$

*consider a mapping*

$$H : M \longrightarrow N$$

*such that $(h, H)$ is morphism of representations from $f$ into $g$. This mapping is unique up to choice of image $n = H(m) \in N$ of given element $m \in M$.*

Proof. From proof of theorem 2.2.8, it follows that choice of homomorphism $h$ and elements $m \in M$, $n \in N$ uniquely defines the mapping $H$. □

**Theorem 2.2.10.** *Given single transitive representation*

$$f : A \to {}^*M$$

*of $\Omega_1$-algebra $A$, for any endomorphism $h$ of $\Omega_1$-algebra $A$ there exists morphism of representation $f$*

$$h : A \to B \quad H : M \to N$$

Proof. Consider following diagram

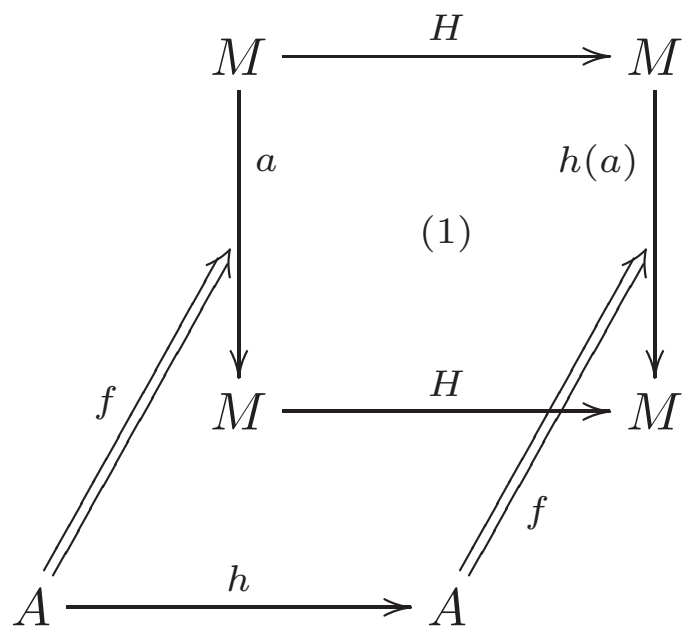

Statement of theorem is corollary of the theorem 2.2.8. □

**Theorem 2.2.11.** *Let*

$$f : A \to {}^*M$$

*be representation of $\Omega_1$-algebra $A$,*

$$g : B \to {}^*N$$

*be representation of $\Omega_1$-algebra $B$,*

$$h : C \to {}^*L$$

*be representation of $\Omega_1$-algebra $C$. Given morphisms of representations of $\Omega_1$-algebra*

$$p : A \longrightarrow B \qquad P : M \longrightarrow N$$

$$q : B \longrightarrow C \qquad Q : N \longrightarrow L$$

*There exists morphism of representations of $\Omega_1$-algebra*

$$r : A \longrightarrow C \qquad R : M \longrightarrow L$$

*where $r = qp$, $R = QP$. We call morphism $(r, R)$ of representations from $f$ into $h$* **product of morphisms** $(p, P)$ **and** $(q, Q)$ **of representations of universal algebra**.

Proof. We represent statement of theorem using diagram

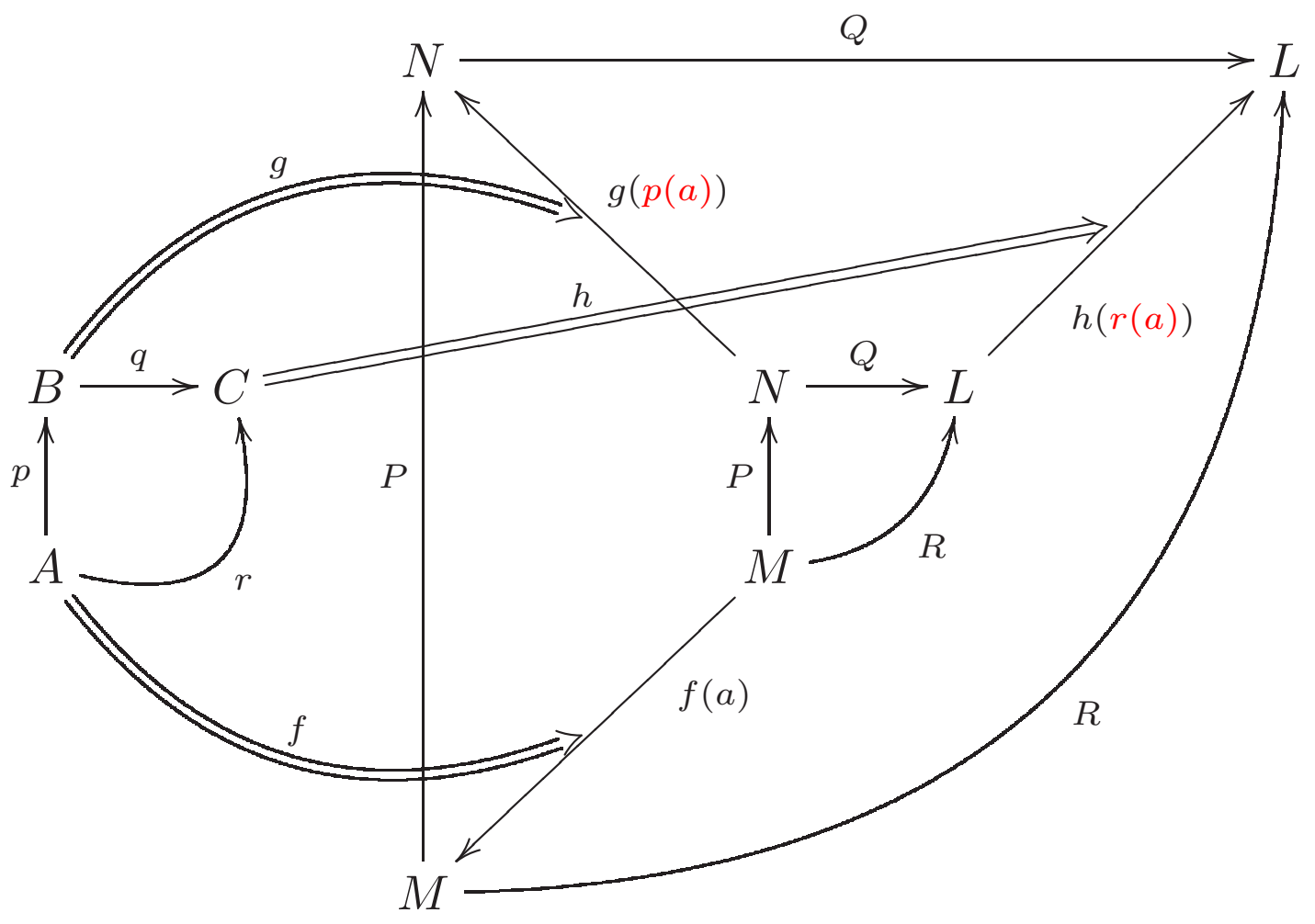

Mapping $r$ is homomorphism of $\Omega_1$-algebra $A$ into $\Omega_1$-algebra $C$. We need to show that tuple of mappings $(r, R)$ satisfies to (2.2.3):

$$\begin{aligned}R(f(a)m) &= QP(f(a)m)\\ &= Q(g(p(a))P(m))\\ &= h(qp(a))QP(m))\\ &= h(r(a))R(m)\end{aligned}$$

□

**Definition 2.2.12.** Let $\mathcal{A}$ be category of $\Omega_1$-algebras. We define **category $T\star\mathcal{A}$ of $T\star$-representations of $\Omega_1$-algebra from category $\mathcal{A}$**. $T\star$-representations of $\Omega_1$-algebra are objects of this category. Morphisms of $T\star$-representations of $\Omega_1$-algebra are morphisms of this category. □

**Definition 2.2.13.** Let us define equivalence $S$ on the set $M$. Transformation $f$ is called **coordinated with equivalence** $S$, when $f(m_1) \equiv f(m_2)(\mathrm{mod}\ S)$ follows from condition $m_1 \equiv m_2(\mathrm{mod} S)$ . □

**Theorem 2.2.14.** *Consider equivalence $S$ on set $M$. Consider $\Omega_1$-algebra on set ${}^*M$. If any transformation $f \in {}^*M$ is coordinated with equivalence $S$, then we can define the structure of $\Omega_1$-algebra on the set ${}^*(M/S)$.*

Proof. Let $h = \mathrm{nat}\ S$. If $m_1 \equiv m_2(\mathrm{mod} S)$ , then $h(m_1) = h(m_2)$. Since $f \in {}^*M$ is coordinated with equivalence $S$, then $h(f(m_1)) = h(f(m_2))$. This allows us to define transformation $F$ according to rule

$$F([m]) = h(f(m)) \tag{2.2.13}$$

Let $\omega$ be n-ary operation of $\Omega_1$-algebra. Suppose $f_1$, ..., $f_n \in {}^*M$ and

$$F_1([m]) = h(f_1(m)) \quad ... \quad F_n([m]) = h(f_n(m))$$

According to condition of theorem, the transformation

$$f = \omega(f_1, ..., f_n) \in {}^*M$$

is coordinated with equivalence $S$. Therefore,

$$\begin{aligned}f(m_1) &\equiv f(m_2)(\mathrm{mod} S)\\ \omega(f_1, ..., f_n)(m_1) &\equiv \omega(f_1, ..., f_n)(m_2)(\mathrm{mod} S)\end{aligned} \tag{2.2.14}$$

follows from condition $m_1 \equiv m_2(\mathrm{mod} S)$ and the definition 2.2.13. Therefore, we can define operation $\omega$ on the set ${}^\star(M/S)$ according to rule

$$\omega(F_1, ..., F_n)[m] = h(\omega(f_1, ..., f_n)(m)) \tag{2.2.15}$$

From the definition (2.2.13) and equation (2.2.14), it follows that we properly defined operation $\omega$ on the set ${}^\star(M/S)$. □

**Theorem 2.2.15.** *Let*

$$f : A \to {}^*M$$

*be representation of $\Omega_1$-algebra $A$,*

$$g : B \to {}^*N$$

*be representation of $\Omega_1$-algebra $B$. Let*

$$r : A \longrightarrow B \qquad R : M \longrightarrow N$$

*be morphism of representations from* $f$ *into* $g$. *Suppose*

$$s = rr^{-1} \qquad S = RR^{-1}$$

*Then there exist decompositions of* $r$ *and* $R$, *which we describe using diagram*

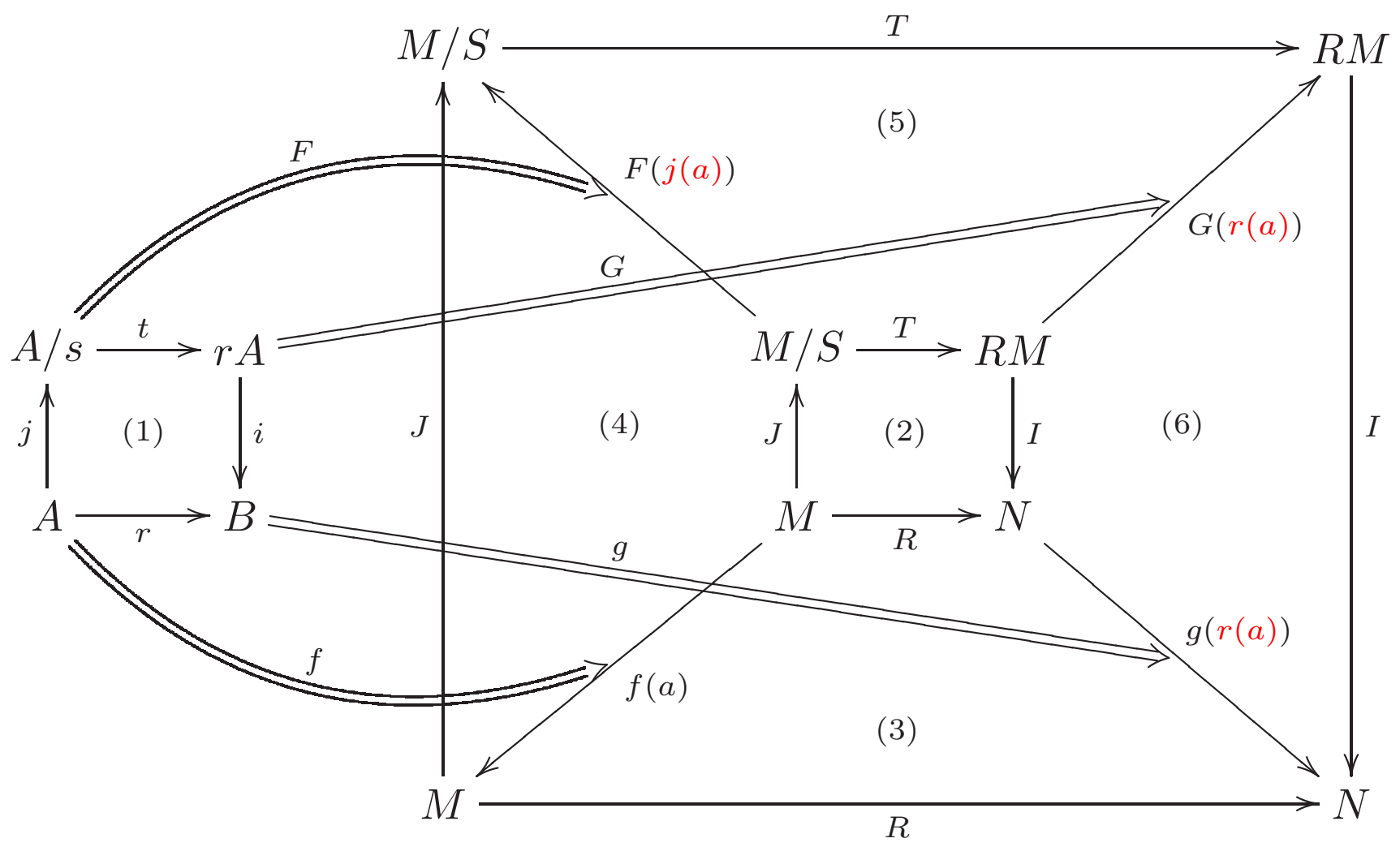

(1) $s = \ker r$ *is a congruence on* $A$. *There exists decompositions of homomorphism* $r$

$$r = itj \tag{2.2.16}$$

$j = \text{nat } s$ *is the natural homomorphism*

$$j(a) = j(a)$$

$t$ *is isomorphism*

$$r(a) = t(j(a)) \tag{2.2.17}$$

$i$ *is the inclusion mapping*

$$r(a) = i(r(a)) \tag{2.2.18}$$

(2) $S = \ker R$ *is an equivalence on* $M$. *There exists decompositions of homomorphism* $R$

$$R = ITJ \tag{2.2.19}$$

$J = \text{nat } S$ *is surjection*

$$J(m) = J(m)$$

$T$ *is bijection*

$$R(m) = T(J(m)) \tag{2.2.20}$$

$I$ *is the inclusion mapping*

$$R(m) = I(R(m)) \tag{2.2.21}$$

(3) $F$ *is* $T\star$*-representation of* $\Omega_1$*-algebra* $A/s$ *in* $M/S$
(4) $G$ *is* $T\star$*-representation of* $\Omega_1$*-algebra* $rA$ *in* $RM$
(5) $(j, J)$ *is morphism of representations* $f$ *and* $F$
(6) $(t, T)$ *is morphism of representations* $F$ *and* $G$

(7) *$(t^{-1}, T^{-1})$ is morphism of representations $G$ and $F$*
(8) *$(i, I)$ is morphism of representations $G$ and $g$*
(9) *There exists decompositions of morphism of representations*

$$(r,R)=(i,I)(t,T)(j,J) \tag{2.2.22}$$

Proof. Existence of diagrams (1) and (2) follows from theorem II.3.7 ([12], p. 60).

We start from diagram (4).

Let $m_1 \equiv m_2 (\text{mod } S)$. Then

$$R(m_1)=R(m_2) \tag{2.2.23}$$

Since $a_1 \equiv a_2 (\text{mod } s)$, then

$$r(a_1)=r(a_2) \tag{2.2.24}$$

Therefore, $j(a_1)=j(a_2)$. Since $(r,R)$ is morphism of representations, then

$$R(f(a_1)(m_1))=g(r(a_1))(R(m_1)) \tag{2.2.25}$$

$$R(f(a_2)(m_2))=g(r(a_2))(R(m_2)) \tag{2.2.26}$$

From (2.2.23), (2.2.24), (2.2.25), (2.2.26), it follows that

$$R(f(a_1)(m_1))=R(f(a_2)(m_2)) \tag{2.2.27}$$

From (2.2.27) it follows

$$f(a_1)(m_1) \equiv f(a_2)(m_2)(\text{mod } S) \tag{2.2.28}$$

and, therefore,

$$J(f(a_1)(m_1))=J(f(a_2)(m_2)) \tag{2.2.29}$$

From (2.2.29) it follows that mapping

$$F(j(a))(J(m))=J(f(a)(m)) \tag{2.2.30}$$

is well defined and this mapping is transformation of set $M/S$.

From equation (2.2.28) (in case $a_1=a_2$) it follows that for any $a$ transformation is coordinated with equivalence $S$. From theorem 2.2.14 it follows that we defined structure of $\Omega_1$-algebra on the set $^*(M/S)$. Consider $n$-ary operation $\omega$ and $n$ transformations

$$F(j(a_i))(J(m))=J(f(a_i)(m)) \quad i=1,...,n$$

of the set $M/S$. We assume

$$\omega(F(j(a_1)),...,F(j(a_n)))(J(m))=J(\omega(f(a_1),...,f(a_n)))(m))$$

Therefore, mapping $F$ is representations of $\Omega_1$-algebra $A/s$.

From (2.2.30) it follows that $(j,J)$ is morphism of representations $f$ and $F$ (the statement (5) of the theorem).

Consider diagram (5).

Since $T$ is bijection, then we identify elements of the set $M/S$ and the set $MR$, and this identification has form

$$T(J(m))=R(m) \tag{2.2.31}$$

We can write transformation $F(j(a))$ of the set $M/S$ as

$$F(j(a)):J(m)\rightarrow F(j(a))(J(m)) \tag{2.2.32}$$

Since $T$ is bijection, we define transformation

$$T(J(m))\rightarrow T(F(j(a))(J(m))) \tag{2.2.33}$$

of the set $RM$. Transformation (2.2.33) depends on $j(a)\in A/s$. Since $t$ is bijection, we identify elements of the set $A/s$ and the set $rA$, and this identification has form

$$t(j(a))=r(a) \tag{2.2.34}$$

Therefore, we defined map

$$G:rA\rightarrow {}^\star RM$$

according to equation

$$G(t(j(a)))(T(J(m)))=T(F(j(a))(J(m))) \tag{2.2.35}$$

Consider $n$-ary operation $\omega$ and $n$ transformations

$$G(r(a_i))(R(m))=T(F(j(a_i))(J(m)))\quad i=1,...,n$$

of space $RM$. We assume

$$\omega(G(r(a_1)),...,G(r(a_n)))(R(m))=T(\omega(F(j(a_1),...,F(j(a_n)))(J(m))) \tag{2.2.36}$$

According to (2.2.35) operation $\omega$ is well defined on the set ${}^\star RM$. Therefore, the mapping $G$ is representations of $\Omega_1$-algebra.

From (2.2.35), it follows that $(t,T)$ is morphism of representations $F$ and $G$ (the statement (6) of the theorem).

Since $T$ is bijection, then from equation (2.2.31), it follows that

$$J(m)=T^{-1}(R(m)) \tag{2.2.37}$$

We can write transformation $G(r(a))$ of the set $RM$ as

$$G(r(a)):R(m)\rightarrow G(r(a))(R(m)) \tag{2.2.38}$$

Since $T$ is bijection, we define transformation

$$T^{-1}(R(m))\rightarrow T^{-1}(G(r(a))(R(m))) \tag{2.2.39}$$

of the set $M/S$. Transformation (2.2.39) depends on $r(a)\in rA$. Since $t$ is bijection, then from equation (2.2.34) it follows that

$$j(a)=t^{-1}(r(a)) \tag{2.2.40}$$

Since, by construction, diagram (5) is commutative, then transformation (2.2.39) coincides with transformation (2.2.32). We can write the equation (2.2.36) as

$$T^{-1}(\omega(G(r(a_1)),...,G(r(a_n)))(R(m)))=\omega(F(j(a_1),...,F(j(a_n)))(J(m)) \tag{2.2.41}$$

Therefore $(t^{-1},T^{-1})$ is morphism of representations $G$ and $F$ (the statement (7) of the theorem).

Diagram (6) is the most simple case in our proof. Since mapping $I$ is immersion and diagram (2) is commutative, we identify $n\in N$ and $R(m)$ when $n\in \mathrm{Im}R$. Similarly, we identify corresponding transformations.

$$g'(i(r(a)))(I(R(m)))=I(G(r(a))(R(m))) \tag{2.2.42}$$

$$\omega(g'(r(a_1)),...,g'(r(a_n)))(R(m))=I(\omega(G(r(a_1),...,G(r(a_n)))(R(m)))$$

Therefore, $(i, I)$ is morphism of representations $G$ and $g$ (the statement (8) of the theorem).

To prove the statement (9) of the theorem we need to show that defined in the proof representation $g'$ is congruent with representation $g$, and operations over transformations are congruent with corresponding operations over ${}^*N$.

$$\begin{aligned} g'(i(r(a)))(I(R(m))) &= I(G(r(a))(R(m))) && \text{by (2.2.42)} \\ &= I(G(t(j(a)))(T(J(m)))) && \text{by (2.2.17), (2.2.20)} \\ &= IT(F(j(a))(J(m))) && \text{by (2.2.35)} \\ &= ITJ(f(a)(m)) && \text{by (2.2.30)} \\ &= R(f(a)(m)) && \text{by (2.2.19)} \\ &= g(r(a))(R(m)) && \text{by (2.2.3)} \end{aligned}$$

$$\begin{aligned} \omega(G(r(a_1)), ..., G(r(a_n)))(R(m)) &= T(\omega(F(j(a_1), ..., F(j(a_n)))(J(m))) \\ &= T(F(\omega(j(a_1), ..., j(a_n)))(J(m))) \\ &= T(F(j(\omega(a_1, ..., a_n)))(J(m))) \\ &= T(J(f(\omega(a_1, ..., a_n))(m))) \end{aligned}$$

□

**Definition 2.2.16.** Let

$$f : A \to {}^*M$$

be representation of $\Omega_1$-algebra $A$,

$$g : B \to {}^*N$$

be representation of $\Omega_1$-algebra $B$. Let

$$r : A \longrightarrow B \qquad R : M \longrightarrow N$$

be morphism of representations from $f$ into $g$ such that $f$ is isomorphism of $\Omega_1$-algebra and $g$ is isomorphism of $\Omega_2$-algebra. Then mapping $(r, R)$ is called **isomorphism of repesentations**. □

**Theorem 2.2.17.** *In the decomposition (2.2.22), the mapping $(t, T)$ is isomorphism of representations $F$ and $G$.*

Proof. The statement of the theorem is corollary of definition 2.2.16 and statements (6) and (7) of the theorem 2.2.15. □

From theorem 2.2.15 it follows that we can reduce the problem of studying of morphism of representations of $\Omega_1$-algebra to the case described by diagram

(2.2.43)

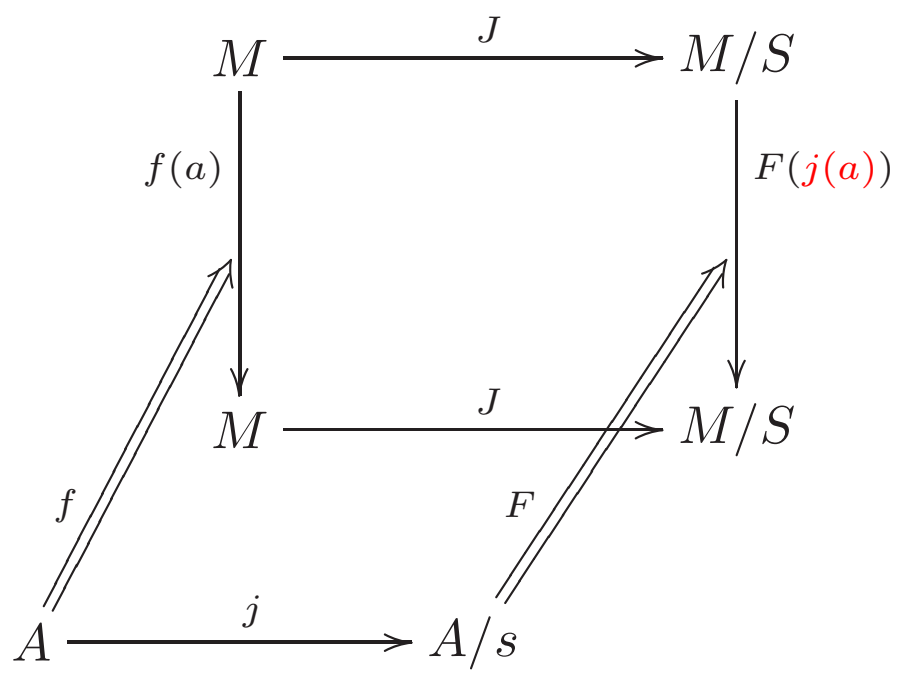

**Theorem 2.2.18.** *We can supplement diagram* (2.2.43) *with representation $F_1$ of $\Omega_1$-algebra $A$ into set $M/S$ such that diagram*

(2.2.44)

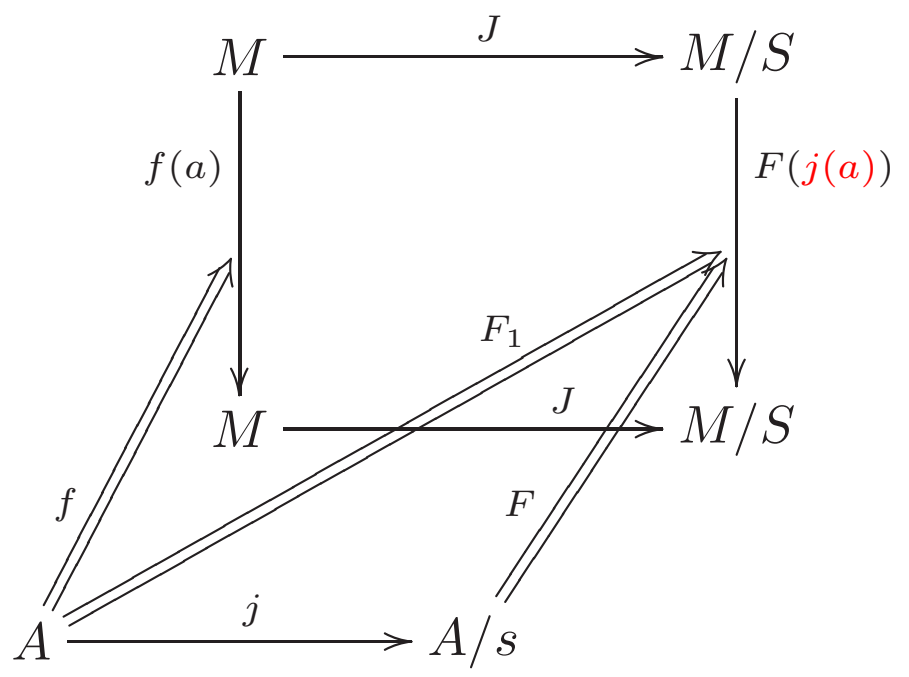

*is commutative. The set of transformations of representation $F$ and the set of transformations of representation $F_1$ coincide.*

Proof. To prove theorem it is enough to assume

$$F_1(a) = F(j(a))$$

Since mapping $j$ is surjection, then $\mathrm{Im}F_1 = \mathrm{Im}F$. Since $j$ and $F$ are homomorphisms of $\Omega_1$-algebra, then $F_1$ is also homomorphism of $\Omega_1$-algebra. □

Theorem 2.2.18 completes the series of theorems dedicated to the structure of morphism of representations $\Omega_1$-algebra. From these theorems it follows that we can simplify task of studying of morphism of representations $\Omega_1$-algebra and not go beyond morphism of representations of form

$$id : A \longrightarrow A \qquad R : M \longrightarrow N$$

In this case we identify morphism of $(\mathrm{id}, R)$ representations of $\Omega_1$-algebra and corresponding homomorphism $R$ of $\Omega_2$-algebra and use the same letter $R$ to denote these mappings. We will use diagram

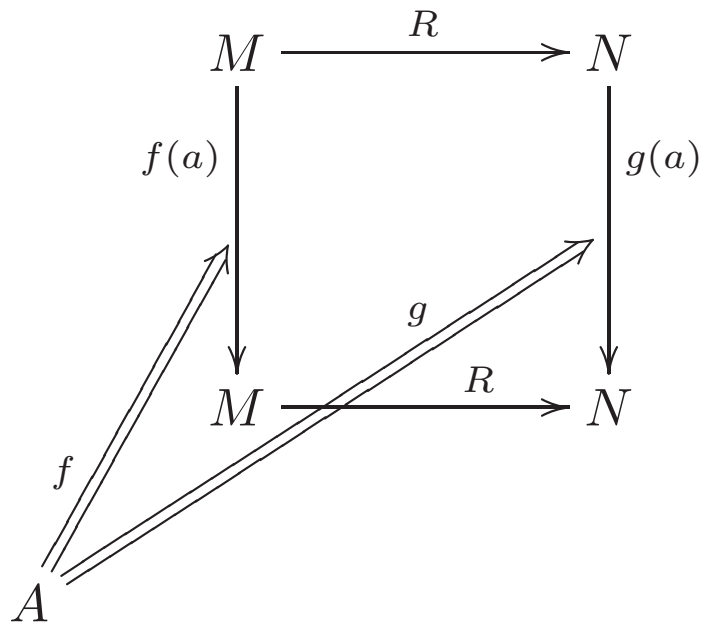

to represent morphism $(\mathrm{id}, R)$ of representations of $\Omega_1$-algebra. From diagram it follows

$$R \circ f(a) = g(a) \circ R \tag{2.2.45}$$

By analogy with definition 2.2.12. we give following definition.

**Definition 2.2.19.** We define **category $T\star A$ $T\star$-representations of $\Omega_1$-algebra** $A$. $T\star$-representations of $\Omega_1$-algebra $A$ are objects of this category. Morphisms $(id, R)$ of $T\star$-representations of $\Omega_1$-algebra $A$ are morphisms of this category. □

## 2.3. Automorphism of Representation of Universal Algebra

**Definition 2.3.1.** Let

$$f : A \to {}^*M$$

be representation of $\Omega_1$-algebra $A$ in $\Omega_2$-algebra $M$. The morphism of representations of $\Omega_1$-algebra

$$(\mathrm{id} : A \to A, R : M \to M)$$

such, that $R$ is endomorphism of $\Omega_2$-algebra is called **endomorphism of representation** $f$. □

**Theorem 2.3.2.** *Given single transitive representation*

$$f : A \to {}^*M$$

*of $\Omega_1$-algebra $A$, for any $p$, $q \in M$ there exists unique endomorphism*

$$H : M \to M$$

*of representation $f$ such that $H(p) = q$.*

Proof. Consider following diagram

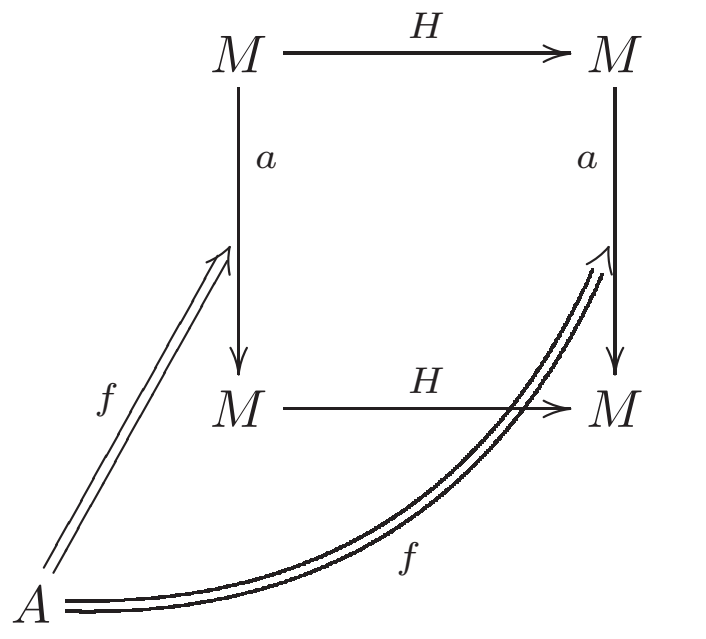

Existence of endomorphism is corollary of the theorem 2.2.8. For given $p, q \in M$, uniqueness of endomorphism follows from the theorem 2.2.9 when $h = \mathrm{id}$. □

**Theorem 2.3.3.** *Endomorphisms of representation $f$ form semigroup.*

Proof. From theorem 2.2.11, it follows that the product of endomorphisms $(p, P)$, $(r, R)$ of the representation $f$ is endomorphism $(pr, PR)$ of the representation $f$. □

**Definition 2.3.4.** Let

$$f : A \to {}^*M$$

be representation of $\Omega_1$-algebra $A$ in $\Omega_2$-algebra $M$. The morphism of representations of $\Omega_1$-algebra

$$(\mathrm{id} : A \to A, R : M \to M)$$

such, that $R$ is automorphism of $\Omega_2$-algebra is called **automorphism of representation** $f$. □

**Theorem 2.3.5.** *Let*

$$f : A \to {}^*M$$

*be representation of $\Omega_1$-algebra $A$ in $\Omega_2$-algebra $M$. The set of automorphisms of the representation $f$ forms* **group** $GA(f)$ .

Proof. Let $R$, $P$ be automorphisms of the representation $f$. According to definition 2.3.4, mappings $R$, $P$ are automorphisms of $\Omega_2$-algebra $M$. According to theorem II.3.2, ([12], p. 57), the mapping $R \circ P$ is automorphism of $\Omega_2$-algebra $M$. From the theorem 2.2.11 and the definition 2.3.4, it follows that product of automorphisms $R \circ P$ of the representation $f$ is automorphism of the representation $f$.

Let $R$, $P$, $Q$ be automorphisms of the representation $f$. The associativity of product of mappings $R$, $P$, $Q$ follows from the chain of equations[2.4]

$$\begin{aligned}((R \circ P) \circ Q)(a) &= (R \circ P)(Q(a)) = R(P(Q(a))) \\ &= R((P \circ Q)(a)) = (R \circ (P \circ Q))(a)\end{aligned}$$

Let $R$ be an automorphism of the representation $f$. According to definition 2.3.4 the mapping $R$ is automorphism of $\Omega_2$-algebra $M$. Therefore, the mapping $R^{-1}$ is automorphism of $\Omega_2$-algebra $M$. The equation (2.2.4) is true for automorphism $R$ of representation. Assume $m' = R(m)$. Since $R$ is automorphism of $\Omega_2$-algebra, then $m = R^{-1}(m')$ and we can write (2.2.4) in the form

$$R(f(a')(R^{-1}(m'))) = f(a')(m') \tag{2.3.1}$$

Since the mapping $R$ is automorphism of $\Omega_2$-algebra $M$, then from the equation (2.3.1) it follows that

$$f(a')(R^{-1}(m')) = R^{-1}(f(a')(m')) \tag{2.3.2}$$

The equation (2.3.2) corresponds to the equation (2.2.4) for the mapping $R^{-1}$. Therefore, mapping $R^{-1}$ of the representation $f$. □

[2.4] To prove the associativity of product I follow to the example of the semigroup from [5], p. 20, 21.

## 2.4. Generating Set of Representation

**Definition 2.4.1.** Let

$$f : A \to {}^*M$$

be representation of $\Omega_1$-algebra $A$ in $\Omega_2$-algebra $M$. The set $N \subset M$ is called **stable set of representation** $f$, if $f(a)(m) \in N$ for each $a \in A$, $m \in N$. □

We also say that the set $M$ is stable with respect to the representation $f$.

**Theorem 2.4.2.** *Let*

$$f : A \to {}^*M$$

*be representation of $\Omega_1$-algebra $A$ in $\Omega_2$-algebra $M$. Let set $N \subset M$ be subalgebra of $\Omega_2$-algebra $M$ and stable set of representation $f$. Then there exists representation*

$$f_N : A \to {}^*N$$

*such that $f_N(a) = f(a)|_N$. Representation $f_N$ is called* **subrepresentation of representation** *$f$.*

Proof. Let $\omega_1$ be $n$-ary operation of $\Omega_1$-algebra $A$. Then for each $a_1$, ..., $a_n \in A$ and each $b \in N$

$$\begin{aligned}(f_N(a_1)...f_N(a_n)\omega_1)(b) &= (f(a_1)...f(a_n)\omega_1)(b)\\ &= f(a_1...a_n\omega_1)(b)\\ &= f_N(a_1...a_n\omega_1)(b)\end{aligned}$$

Let $\omega_2$ be $n$-ary operation of $\Omega_2$-algebra $M$. Then for each $b_1$, ..., $b_n \in N$ and each $a \in A$

$$\begin{aligned}f_N(a)(b_1)...f_N(a)(b_n)\omega_2 &= f(a)(b_1)...f(a)(b_n)\omega_2\\ &= f(a)(b_1...b_n\omega_2)\\ &= f_N(a)(b_1...b_n\omega_2)\end{aligned}$$

We proved the statement of theorem. □

From the theorem 2.4.2, it follows that if $f_N$ is subrepresentation of representation $f$, then the mapping $(id : A \to A, id_N : N \to M)$ is morphism of representations.

**Theorem 2.4.3.** *The set*[2.5] *$\mathcal{B}_f$ of all subrepresentations of representation $f$ generates a closure system on $\Omega_2$-algebra $M$ and therefore is a complete lattice.*

Proof. Let $(K_\lambda)_{\lambda\in\Lambda}$ be the set off subalgebras of $\Omega_2$-algebra $M$ that are stable with respect to representation $f$. We define the operation of intersection on the set $\mathcal{B}_f$ according to rule

$$\bigcap f_{K_\lambda} = f_{\cap K_\lambda}$$

We defined the operation of intersection of subrepresentations properly. $\cap K_\lambda$ is subalgebra of $\Omega_2$-algebra $M$. Let $m \in \cap K_\lambda$. For each $\lambda \in \Lambda$ and for each $a \in A$, $f(a)(m) \in K_\lambda$. Therefore, $f(a)(m) \in \cap K_\lambda$. Therefore, $\cap K_\lambda$ is the stable set of representation $f$. □

[2.5] This definition is similar to definition of the lattice of subalgebras ([12], p. 79, 80)

We denote the corresponding closure operator by $J(f)$. Thus $J(f,X)$ is the intersection of all subalgebras of $\Omega_2$-algebra $M$ containing $X$ and stable with respect to representation $f$.

**Theorem 2.4.4.** *Let*[2.6]

$$f : A \rightarrow {}^*M$$

*be representation of $\Omega_1$-algebra $A$ in $\Omega_2$-algebra $M$. Let $X \subset M$. Define a subset $X_k \subset M$ by induction on $k$.*

$$X_0 = X$$
$$x \in X_k => x \in X_{k+1}$$
$$x_1 \in X_k, ..., x_n \in X_k, \omega \in \Omega_2(n) => x_1...x_n\omega \in X_{k+1}$$
$$x \in X_k, a \in A => f(a)(x) \in X_{k+1}$$

*Then*

$$\bigcup_{k=0}^{\infty} X_k = J(f,X)$$

Proof. If we put $U = \cup X_k$, then by definition of $X_k$, we have $X_0 \subset J(f,X)$, and if $X_k \subset J(f,X)$, then $X_{k+1} \subset J(f,X)$. By induction it follows that $X_k \subset J(f,X)$ for all $k$. Therefore,

$$U \subset J(f,X) \tag{2.4.1}$$

If $a \in U^n$, $a = (a_1, ..., a_n)$, where $a_i \in X_{k_i}$, and if $k = \max\{k_1, ..., k_n\}$, then $a_1...a_n\omega \in X_{k+1} \subset U$. Therefore, $U$ is subalgebra of $\Omega_2$-algebra $M$.

If $m \in U$, then there exists such $k$ that $m \in X_k$. Therefore, $f(a)(m) \in X_{k+1} \subset U$ for any $a \in A$. Therefore, $U$ is stable set of the representation $f$.

Since $U$ is subalgebra of $\Omega_2$-algebra $M$ and is a stable set of the representation $f$, then subrepresentation $f_U$ is defined. Therefore,

$$J(f,X) \subset U \tag{2.4.2}$$

From (2.4.1), (2.4.2), it follows that $J(f,X) = U$. □

**Definition 2.4.5.** $J(f,X)$ is called **subrepresentation generated by set** $X$, and $X$ is a **generating set of subrepresentation** $J(f,X)$. In particular, a **generating set of representation** $f$ is a subset $X \subset M$ such that $J(f,X) = M$. □

It is easy to see that the definition of generating set of representation does not depend on whether representation is effective or not. For this reason hereinafter we will assume that the representation is effective and we will use convention for effective $T\star$-representation in remark 2.1.9. We also will use notation

$$R \circ m = R(m)$$

for image of $m \in M$ under the endomorphism of effective representation. According to the definition of product of mappings, for any endomorphisms $R$, $S$ the following equation is true

$$(R \circ S) \circ m = R \circ (S \circ m) \tag{2.4.3}$$

[2.6] The statement of theorem is similar to the statement of theorem 5.1, [12], p. 79.

The equation (2.4.3) is associative law for $\circ$ and allows us to write expression

$$R\circ S\circ m$$

without brackets.

From theorem 2.4.4, it follows next definition.

**Definition 2.4.6.** Let $X\subset M$. For each $x\in J(f,X)$ there exists $\Omega_2$-word defined according to following rules.

(1) If $m\in X$, then $m$ is $\Omega_2$-word.
(2) If $m_1$, ..., $m_n$ are $\Omega_2$-words and $\omega\in\Omega_2(n)$, then $m_1...m_n\omega$ is $\Omega_2$-word.
(3) If $m$ is $\Omega_2$-word and $a\in A$, then $am$ is $\Omega_2$-word.

**$\Omega_2$-word** $w(f,X,m)$ represents given element $m\in J(f,X)$. We will identify an element $m\in J(f,X)$ and corresponding it $\Omega_2$-word using equation

$$m=w(f,X,m)$$

Similarly, for an arbitrary set $B\subset J(f,X)$ we consider the set of $\Omega_2$-words[2.7]

$$w(f,X,B)=\{w(f,X,m):m\in B\}$$

We also use notation

$$w(f,X,B)=(w(f,X,m),m\in B)$$

Denote $w(f,X)$ **the set of $\Omega_2$-words of representation** $J(f,X)$. □

**Theorem 2.4.7.** *Endomorphism $R$ of representation*

$$f:A\rightarrow {}^*M$$

*generates the mapping of $\Omega_2$-words*

$$w[f,X,R]:w(f,X)\rightarrow w(f,X')\quad X\subset M\quad X'=R\circ X$$

*such that*

(1) *If* $m\in X$, $m'=R\circ m$, *then*

$$w[f,X,R](m)=m'$$

(2) *If*

$$m_1,...,m_n\in w(f,X)$$
$$m_1'=w[f,X,R](m_1)\quad ...\quad m_n'=w[f,X,R](m_n)$$

*then for operation* $\omega\in\Omega_2(n)$ *holds*

$$w[f,X,R](m_1...m_n\omega)=m_1'...m_n'\omega$$

(3) *If*

$$m\in w(f,X)\quad m'=w[f,X,R](m)\quad a\in A$$

*then*

$$w[f,X,R](am)=am'$$

[2.7] The expression $w(f,X,m)$ is a special case of the expression $w(f,X,B)$, namely

$$w(f,X,\{m\})=\{w(f,X,m)\}$$

Proof. Statements (1), (2) of the theorem are true by definition of the endomorpism $R$. Because $r=\mathrm{id}$, the statement (3) of the theorem follows from the equation (2.2.4). □

**Remark 2.4.8.** Let $R$ be endomorphism of representation $f$. Let

$$m\in J(f,X)\quad m'=R\circ m\quad X'=R\circ X$$

The theorem 2.4.7 states that $m'\in J_f(X')$ . The theorem 2.4.7 also states that $\Omega_2$-word representing $m$ relative $X$ and $\Omega_2$-word representing $m'$ relative $X'$ are generated according to the same algorithm. This allows considering of the set of $\Omega_2$-words $w(f,X',m')$ as mapping

$$W(f,X,m):X'\to w(f,X',m')$$
$$W(f,X,m)(X')=W(f,X,m)\circ X'$$

such that, if for certain endomorphism $R$

$$X'=R\circ X\quad m'=R\circ m$$

then

$$W(f,X,m)\circ X'=w(f,X',m')=m'$$

The mapping $W(f,X,m)$ is called **coordinates of element $m$ relative to set $X$**. Similarly, we consider coordinates of a set $B\subset J(f,X)$ relative to the set $X$

$$W(f,X,B)=\{W(f,X,m):m\in B\}=(W(f,X,m),m\in B)$$

Denote $W(f,X)$ **the set of coordinates of representation** $J(f,X)$. □

**Theorem 2.4.9.** *There is a structure of $\Omega_2$-algebra on the set of coordinates $W(f,X)$.*

Proof. Let $\omega\in\Omega_2(n)$. Then for any $m_1$, ..., $m_n\in J(f,X)$ , we assume

$$W(f,X,m_1)...W(f,X,m_n)\omega=W(f,X,m_1...m_n\omega)\tag{2.4.4}$$

According to the remark 2.4.8,

$$\begin{aligned}(W(f,X,m_1)...W(f,X,m_n)\omega)\circ X&=W(f,X,m_1...m_n\omega)\circ X\\&=w(f,X,m_1...m_n\omega)\end{aligned}\tag{2.4.5}$$

follows from the equation (2.4.4). According to rule (2) of the definition 2.4.6, from the equation (2.4.5), it follows that

$$\begin{aligned}(W(f,X,m_1)...W(f,X,m_n)\omega)\circ X&=w(f,X,m_1)...w(f,X,m_n)\omega\\&=(W(f,X,m_1)\circ X)...(W(f,X,m_n)\circ X)\omega\end{aligned}\tag{2.4.6}$$

From the equation (2.4.6), it follows that the operation $\omega$ defined by the equation (2.4.4) on the set of coordinates is defined properly. □

**Theorem 2.4.10.** *There exists the representation of $\Omega_1$-algebra $A$ in $\Omega_2$-algebra $W(f,X)$.*

PROOF. Let $a \in A$. Then for any $m \in J(f,X)$ we assume

$$aW(f,X,m) = W(f,X,am) \tag{2.4.7}$$

According to the remark 2.4.8,

$$(aW(f,X,m)) \circ X = W(f,X,am) \circ X = w(f,X,am) \tag{2.4.8}$$

follows from the equation (2.4.7). According to rule (3) of the definition 2.4.6, from the equation (2.4.8), it follows that

$$(aW(f,X,m)) \circ X = aw(f,X,m) = a(W(f,X,m) \circ X) \tag{2.4.9}$$

From the equation (2.4.9), it follows that the representation (2.4.7) of $\Omega_1$-algebra $A$ in $\Omega_2$-algebra $W(f,X)$ is defined properly. □

**Definition 2.4.11.** Let $X$ be the generating set of the representation $f$. Let $R$ be the endomorphism of the representation $f$. The set of coordinates $W(f,X,R \circ X)$ is called **coordinates of endomorphism of representation**. □

**Definition 2.4.12.** Let $X$ be the generating set of the representation $f$. Let $R$ be the endomorphism of the representation $f$. Let $m \in M$. We define **superposition of coordinates** of the representation $f$ and the element $m$ as coordinates defined according to rule

$$W(f,X,m) \circ W(f,X,R \circ X) = W(f,X,R \circ m) \tag{2.4.10}$$

Let $Y \subset M$. We define superposition of coordinates of the representation $f$ and the set $Y$ according to rule

$$W(f,X,Y) \circ W(f,X,R \circ X) = (W(f,X,m) \circ W(f,X,R \circ X), m \in Y) \tag{2.4.11}$$

□

**Theorem 2.4.13.** *Endomorphism $R$ of representation*

$$f : A \to {}^*M$$

*generates the mapping of coordinates of representation*

$$W[f,X,R] : W(f,X) \to W(f,X) \tag{2.4.12}$$

*such that*

$$\begin{array}{c} W(f,X,m) \to W[f,X,R] \star W(f,X,m) = W(f,X,R \circ m) \\ W[f,X,R] \star W(f,X,m) = W(f,X,m) \circ W(f,X,R \circ X) \end{array} \tag{2.4.13}$$

PROOF. According to the remark 2.4.8, we consider equations (2.4.10), (2.4.12) relative to given generating set $X$. The word

$$W(f,X,m) \circ X = w(f,X,m) \tag{2.4.14}$$

corresponds to coordinates $W(f,X,m)$; the word

$$W(f,X,R \circ m) \circ X = w(f,X,R \circ m) \tag{2.4.15}$$

corresponds to coordinates $W(f,X,R \circ m)$. Therefore, in order to prove the theorem, it is sufficient to show that the mapping $W[f,X,R]$ corresponds to mapping $w[f,X,R]$. We prove this statement by induction over complexity of $\Omega_2$-word.

If $m \in X$, $m' = R \circ m$, then, according to equations (2.4.14), (2.4.15), mappings $W[f,X,R]$ and $w[f,X,R]$ are coordinated.

Let for $m_1$, ..., $m_n \in X$ mappings $W[f,X,R]$ and $w[f,X,R]$ be coordinated. Let $\omega \in \Omega_2(n)$. According to the theorem 2.4.9

$$W(f,X,m_1...m_n\omega) = W(f,X,m_1)...W(f,X,m_n)\omega \tag{2.4.16}$$

Because $R$ is endomorphism of $\Omega_2$-algebra $M$, then from the equation (2.4.16), it follows that

$$\begin{aligned} W(f,X,R\circ(m_1...m_n\omega)) &= W(f,X,(R\circ m_1)...(R\circ m_n)\omega) \\ &= W(f,X,R\circ m_1)...W(f,X,R\circ m_n)\omega \end{aligned} \tag{2.4.17}$$

From equations (2.4.16), (2.4.17) and the statement of induction, it follows that the mappings $W[f,X,R]$ and $w[f,X,R]$ are coordinated for $m = m_1...m_n\omega$.

Let for $m_1 \in M$ mappings $W[f,X,R]$ and $w[f,X,R]$ are coordinated. Let $a \in A$. According to the theorem 2.4.10

$$W(f,X,am_1) = aW(f,X,m_1) \tag{2.4.18}$$

Because $R$ is endomorphism of representation $f$, then, from the equation (2.4.18), it follows that

$$W(f,X,R\circ(am_1)) = W(f,X,a\,R\circ m_1) = aW(f,X,R\circ m_1) \tag{2.4.19}$$

From equations (2.4.18), (2.4.19) and the statement of induction, it follows that mappings $W[f,X,R]$ and $w[f,X,R]$ are coordinated for $m = am_1$. □

**Corollary 2.4.14.** *Let $X$ be the generating set of the representation $f$. Let $R$ be the endomorphism of the representation $f$. The mapping $W[f,X,R]$ is endomorphism of representation of $\Omega_1$-algebra $A$ in $\Omega_2$-algebra $W(f,X)$.* □

Hereinafter we will identify mapping $W[f,X,R]$ and the set of coordinates $W(f,X,R\circ X)$.

**Theorem 2.4.15.** *Let $X$ be the generating set of the representation $f$. Let $R$ be the endomorphism of the representation $f$. Let $Y \subset M$. Then*

$$W(f,X,Y)\circ W(f,X,R\circ X) = W(f,X,R\circ Y) \tag{2.4.20}$$

$$W[f,X,R]\star W(f,X,Y) = W(f,X,R\circ Y) \tag{2.4.21}$$

Proof. The equation (2.4.20) follows from the equation

$$R\circ Y = (R\circ m, m\in Y)$$

as well from equations (2.4.10), (2.4.11). The equation (2.4.21) is corollary of equations (2.4.20), (2.4.13). □

**Theorem 2.4.16.** *Let $X$ be the generating set of the representation $f$. Let $R$, $S$ be the endomorphisms of the representation $f$. Then*

$$W(f,X,S\circ X)\circ W(f,X,R\circ X) = W(f,X,R\circ S\circ X) \tag{2.4.22}$$

$$W[f,X,R]\star W[f,X,S] = W[f,X,R\circ S] \tag{2.4.23}$$

Proof. The equation (2.4.22) follows from the equation (2.4.20), if we assume $Y = S\circ X$. The equation (2.4.23) follows from the equation (2.4.22) and chain of

equations

$$\begin{aligned}
&(W[f,X,R]\star W[f,X,S])\star W(f,X,Y)\\
=&W[f,X,R]\star(W[f,X,S]\star W(f,X,Y))\\
=&(W(f,X,Y)\circ W(f,X,S\circ X))\circ W(f,X,R\circ X)\\
=&W(f,X,Y)\circ(W(f,X,S\circ X)\circ W(f,X,R\circ X))\\
=&W(f,X,Y)\circ W(f,X,R\circ S\circ X)\\
=&W(f,X,R\circ S)\star W(f,X,Y)
\end{aligned}$$

□

The concept of superposition of the coordinates is very simple and resembles a kind of Turing machine. If element $m\in M$ has form either

$$m=m_1...m_n\omega$$

or

$$m=am_1$$

then we are looking for the coordinates of elements $m_i$ to substitute them in an appropriate expression. As soon as an element $m\in M$ belongs to the generating set of $\Omega_2$-algebra $M$, we choose the coordinates of the corresponding element of the second factor. Therefore, we require that the second factor in the superposition has been the set of coordinates of the image of the generating set $X$.

We can generalize the definition of the superposition of coordinates and assume that one of the factors is a set of $\Omega_2$-words. Accordingly, the definition of the superposition of coordinates has the form

$$w(f,X,Y)\circ W(f,X,R\circ X)=W(f,X,Y)\circ w(f,X,R\circ X)=w(f,X,R\circ Y)$$

The following forms of writing an image of the set $Y$ under endomorphism $R$ are equivalent.

$$R\circ Y=W(f,X,Y)\circ(R\circ X)=W(f,X,Y)\circ(W(f,X,R\circ X)\circ X) \tag{2.4.24}$$

From equations (2.4.20), (2.4.24), it follows that

$$(W(f,X,Y)\circ W(f,X,R\circ X))\circ X=W(f,X,Y)\circ(W(f,X,R\circ X)\circ X) \tag{2.4.25}$$

The equation (2.4.25) is associative law for composition and allows us to write expression

$$W(f,X,Y)\circ W(f,X,R\circ X)\circ X$$

without brackets.

Consider equation (2.4.22), where we see change in the order of endomorphisms in a superposition of the coordinates. This equation also follows from the chain of equations, where we can immediately see when order of endomorphisms changes

$$\begin{aligned}
&W(f,X,m)\circ W(f,X,R\circ S\circ X)\circ X\\
=&(R\circ S)\circ m=R\circ(S\circ m)\\
=&R\circ((W(f,X,m)\circ W(f,X,S\circ X))\circ X)\\
=&W(f,X,m)\circ W(f,X,S\circ X)\circ W(f,X,R\circ X)\circ X
\end{aligned} \tag{2.4.26}$$

From the equation (2.4.26), it follows that coordinates of endomorphism act over coordinates of element of $\Omega_2$-algebra $M$ from the right.

**Definition 2.4.17.** Let $X \subset M$ be generating set of representation

$$f : A \to {}^*M$$

Let the mapping

$$H : M \to M$$

be endomorphism of the representation $f$. Let the set $X' = H \circ X$ be the image of the set $X$ under the mapping $H$. Endomorphism $H$ of representation $f$ is called **regular on generating set** $X$, if the set $X'$ is the generating set of representation $f$. Otherwise, endomorphism $H$ of representation $f$ is called **singular on generating set** $X$, □

**Definition 2.4.18.** Endomorphism of representation $f$ is called **regular**, if it is regular on every generating set. □

**Theorem 2.4.19.** *Automorphism $R$ of representation*

$$f : A \to {}^*M$$

*is regular endomorphism.*

Proof. Let $X$ be generating set of representation $f$. Let $X' = R \circ X$.

According to theorem 2.4.7 endomorphism $R$ forms the map of $\Omega_2$-words $w[f, X, R]$.

Let $m' \in M$. Since $R$ is automorphism, then there exists $m \in M$, $R \circ m = m'$. According to definition 2.4.6, $w(f, X, m)$ is $\Omega_2$-word, representing $m$ relative to generating set $X$. According to theorem 2.4.7, $w(f, X', m')$ is $\Omega_2$-word, representing of $m'$ relative to generating set $X'$

$$w(f, X', m') = w[f, X, R](w(f, X, m))$$

Therefore, $X'$ is generating set of representation $f$. According to definition 2.4.18, automorphism $R$ is regular. □

## 2.5. Basis of representation

**Definition 2.5.1.** If the set $X \subset M$ is generating set of representation $f$, then any set $Y$, $X \subset Y \subset M$ also is generating set of representation $f$. If there exists minimal set $X$ generating the representation $f$, then the set $X$ is called **basis of representation** $f$. □

**Theorem 2.5.2.** *The generating set $X$ of representation $f$ is basis iff for any $m \in X$ the set $X \setminus \{m\}$ is not generating set of representation $f$.*

Proof. Let $X$ be generating set of representation $f$. Assume that for some $m \in X$ there exist $\Omega_2$-word

$$w = w(f, X \setminus \{m\}, m) \tag{2.5.1}$$

Consider element $m'$ such that it has $\Omega_2$-word

$$w' = w(f, X, m') \tag{2.5.2}$$

that depends on $m$. According to the definition 2.4.6, any occurrence of $m$ into $\Omega_2$-word $w'$ can be substituted by the $\Omega_2$-word $w$. Therefore, the $\Omega_2$-word $w'$ does not depend on $m$, and the set $X \setminus \{m\}$ is generating set of representation $f$. Therefore, $X$ is not basis of representation $f$. □

**Remark 2.5.3.** The proof of the theorem 2.5.2 gives us effective method for constructing the basis of the representation $f$. Choosing an arbitrary generating set, step by step, we remove from set those elements which have coordinates relative to other elements of the set. If the generating set of the representation is infinite, then this construction may not have the last step. If the representation has finite generating set, then we need a finite number of steps to construct a basis of this representation.

As noted by Paul Cohn in [12], p. 82, 83, the representation may have inequivalent bases. For instance, the cyclic group of order six has bases $\{a\}$ and $\{a^2, a^3\}$ which we cannot map one into another by endomorphism of the representation. □

**Remark 2.5.4.** We write a basis also in following form

$$X = (x, x \in X)$$

If basis is finite, then we also use notation

$$X = (x_i, i \in I) = (x_1, ..., x_n)$$

□

**Remark 2.5.5.** We introduced $\Omega_2$-word of $x \in M$ relative generating set $X$ in the definition 2.4.6. From the theorem 2.5.2, it follows that if the generating set $X$ is not a basis, then a choice of $\Omega_2$-word relative generating set $X$ is ambiguous. However, even if the generating set $X$ is a basis, then a representation of $m \in M$ in form of $\Omega_2$-word is ambiguous. If $m_1$, ..., $m_n$ are $\Omega_2$-words, $\omega \in \Omega_2(n)$ and $a \in A$, then[2.8]

$$a(m_1...m_n\omega) = (am_1)...(am_n)\omega \tag{2.5.3}$$

It is possible that there exist equations related to specific character of representation. For instance, if $\omega$ is operation of $\Omega_1$-algebra $A$ and operation of $\Omega_2$-algebra $M$, then we require that $\Omega_2$-words $a_1...a_n\omega x$ and $a_1x...a_nx\omega$ describe the same element of $\Omega_2$-algebra $M$.[2.9]

In addition to the above equations in $\Omega_2$-algebra there may be relations of the form

$$w_1(f, X, m) = w_2(f, X, m) \quad m \notin X \tag{2.5.4}$$

The feature of the equation (2.5.4) is that this equation cannot be reduced.

On the set of $\Omega_2$-words $w(f, X)$, above equations determine **equivalence** $\rho(f)$ **generated by representation** $f$. It is evident that for any $m \in M$ the choice of appropriate $\Omega_2$-word is unique up to equivalence relations $\rho(f)$. However, if during the construction, we obtain the equality of two $\Omega_2$-word relative to given basis, then we can say without worrying about the equivalence $\rho(f)$ that these $\Omega_2$-words are equal.

[2.8]For instance, let $\{\overline{e}_1, \overline{e}_2\}$ be the basis of vector space over field $k$. The equation (2.5.3) has the form of distributive law

$$a(b^1\overline{e}_1 + b^2\overline{e}_2) = (ab^1)\overline{e}_1 + (ab^2)\overline{e}_2$$

[2.9]For vector space, this requirement has the form of distributive law

$$(a + b)\overline{e}_1 = a\overline{e}_1 + b\overline{e}_1$$

A similar remark concerns the mapping $W(f, X, m)$ defined in the remark 2.4.8.[2.10] □

**Theorem 2.5.6.** *Automorphism of the representation $f$ maps a basis of the representation $f$ into basis.*

Proof. Let the mapping $R$ be automorphism of the representation $f$. Let the set $X$ be a basis of the representation $f$. Let $X' = R \circ X$.

Assume that the set $X'$ is not basis. According to the theorem 2.5.2 there exists such $m' \in X'$ that $X' \setminus \{x'\}$ is generating set of the representation $f$. According to the theorem 2.3.5 the mapping $R^{-1}$ is automorphism of the representation $f$. According to the theorem 2.4.19 and definition 2.4.18, the set $X \setminus \{m\}$ is generating set of the representation $f$. The contradiction completes the proof of the theorem. □

[2.10]If vector space has finite basis, then we represent the basis as matrix

$$\overline{\overline{e}} = \begin{pmatrix} \overline{e}_1 & ... & \overline{e}_2 \end{pmatrix}$$

We present the mapping $W(f, \overline{\overline{e}}, \overline{v})$ as matrix

$$W(f, \overline{\overline{e}}, \overline{v}) = \begin{pmatrix} v^1 \\ ... \\ v^n \end{pmatrix}$$

Then

$$W(f, \overline{\overline{e}}, \overline{v})(\overline{\overline{e}}') = W(f, \overline{\overline{e}}, \overline{v}) \begin{pmatrix} \overline{e}'_1 & ... & \overline{e}'_n \end{pmatrix} = \begin{pmatrix} v^1 \\ ... \\ v^n \end{pmatrix} \begin{pmatrix} \overline{e}'_1 & ... & \overline{e}'_n \end{pmatrix}$$

has form of matrix product.

CHAPTER 3

# Tower of Representations of Universal Algebras

## 3.1. Tower of Representations of Universal Algebras

**Definition 3.1.1.** Consider[3.1] set of $\Omega_k$-algebras $A_k$, $k = 1, ..., n$. Assume $\overline{A} = (A_1, ..., A_n)$. Assume $\overline{f} = (f_{1,2}, ..., f_{n-1,n})$. Set of representations $f_{k,k+1}$, $k = 1, ..., n$, of $\Omega_k$-algebra $A_k$ in $\Omega_{k+1}$-algebra $A_{k+1}$ is called **tower** $(\overline{A}, \overline{f})$ **of representations of** $\overline{\Omega}$**-algebras**. □

Consider following diagram for the purposes of illustration of definition 3.1.1

(3.1.1)

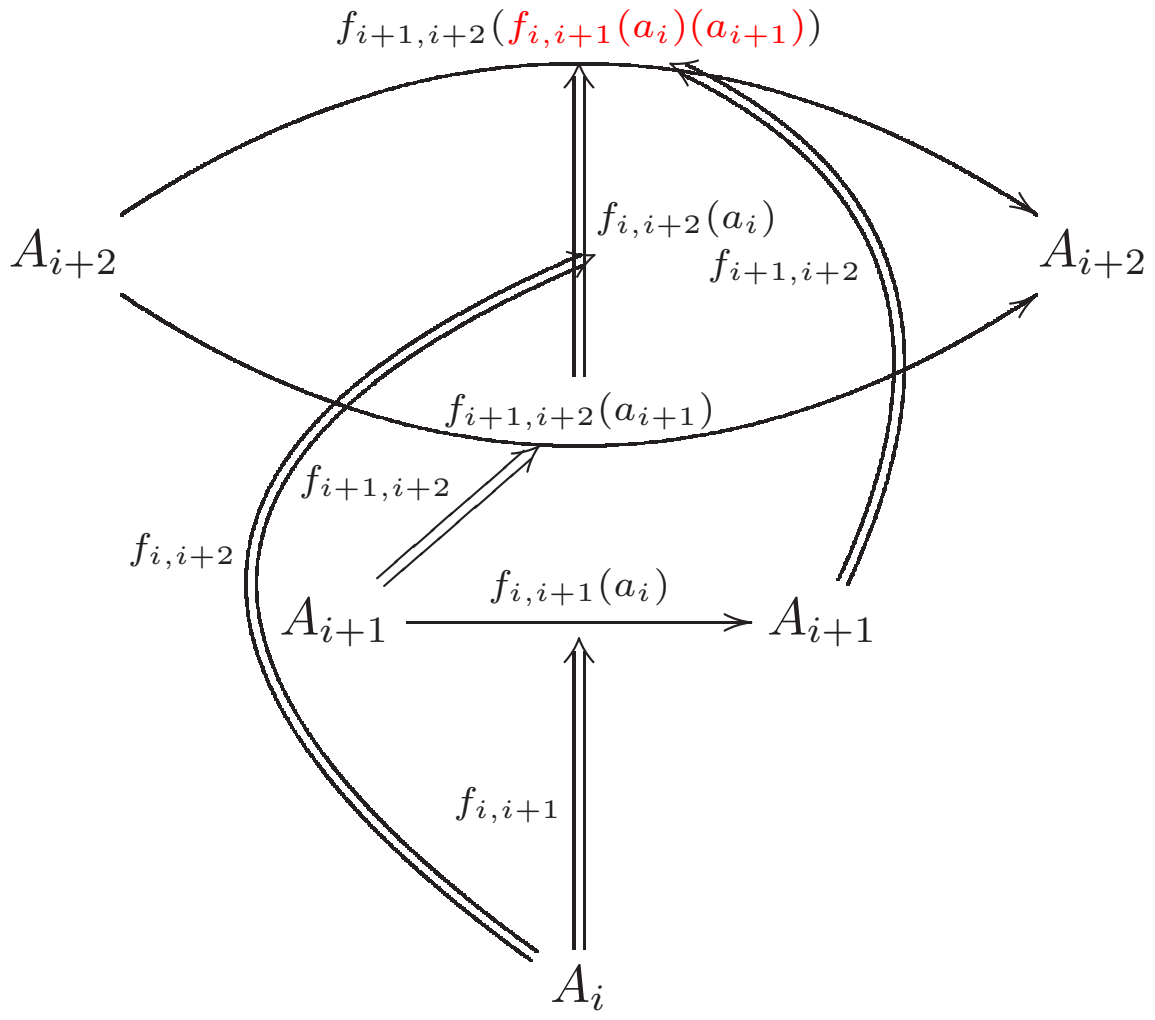

$f_{i,i+1}$ is representation of $\Omega_i$-algebra $A_i$ in $\Omega_{i+1}$-algebra $A_{i+1}$. $f_{i+1,i+2}$ is representation of $\Omega_{i+1}$-algebra $A_{i+1}$ in $\Omega_{i+2}$-algebra $A_{i+2}$.

**Theorem 3.1.2.** *The mapping*

$$f_{i,i+2} : A_i \to {}^{**}A_{i+2}$$

*is representation of* $\Omega_i$*-algebra* $A_i$ *in* $\Omega_{i+1}$*-algebra* ${}^*A_{i+2}$.

Proof. Automorphism $f_{i+1,i+2}(a_{i+1}) \in {}^*A_{i+2}$ corresponds to arbitrary $a_{i+1} \in A_{i+1}$. Automorphism $f_{i,i+1}(a_i) \in {}^*A_{i+1}$ corresponds to arbitrary $a_i \in A_i$. Since $f_{i,i+1}(a_i)(a_{i+1}) \in A_{i+1}$ is image of $a_{i+1} \in A_{i+1}$, then element $a_i \in A_i$ generates transformation of $\Omega_{i+1}$-algebra ${}^*A_{i+2}$ which is defined by equation

$$f_{i,i+2}(a_i)(f_{i+1,i+2}(a_{i+1})) = f_{i+1,i+2}(f_{i,i+1}(a_i)(a_{i+1})) \tag{3.1.2}$$

[3.1] Definitions and theorems from this chapter are included also into chapters [10]-4, [10]-5.

Let $\omega$ be $n$-ary operation of $\Omega_i$-algebra. Because $f_{i,i+1}$ is homomorphism of $\Omega_i$-algebra, then

$$f_{i,i+1}(a_{i,1}...a_{i,n}\omega) = f_{i,i+1}(a_{i,1})...f_{i,i+1}(a_{i,n})\omega \tag{3.1.3}$$

For $a_i$, ..., $a_n \in A_i$ we define operation $\omega$ on the set ${}^{**}A_{i+2}$ using equations

$$\begin{aligned} &f_{i,i+2}(a_{i,1})(f_{i+1,i+2}(a_{i+1}))...f_{i,i+2}(a_{i,n})(f_{i+1,i+2}(a_{i+1}))\omega \\ =&f_{i+1,i+2}(f_{i,i+1}(a_{i,1})(a_{i+1}))...f_{i+1,i+2}(f_{i,i+1}(a_{i,n})(a_{i+1}))\omega \\ =&f_{i+1,i+2}(f_{i,i+1}(a_{i,1})...f_{i,i+1}(a_{i,n})\omega)(a_{i+1}) \end{aligned} \tag{3.1.4}$$

First equation follows from equation (3.1.2). We wrote second equation on the base of demand that map $f_{i+1,i+2}$ is homomorphism of $\Omega_i$-algebra. Therefore, we defined the structure of $\Omega_i$-algebra on the set ${}^{**}A_{i+2}$.

From (3.1.4) and (3.1.3) it follows that

$$\begin{aligned} &f_{i,i+2}(a_{i,1})(f_{i+1,i+2}(a_{i+1}))...f_{i,i+2}(a_{i,n})(f_{i+1,i+2}(a_{i+1}))\omega \\ =&f_{i+1,i+2}(f_{i,i+1}(a_{i,1}...a_{i,n}\omega)(a_{i+1})) \\ =&f_{i,i+2}(a_{i,1}...a_{i,n}\omega)(f_{i+1,i+2}(a_{i+1})) \end{aligned} \tag{3.1.5}$$

From equation (3.1.5) it follows that

$$f_{i,i+2}(a_{i,1})...f_{i,i+2}(a_{i,n})\omega = f_{i,i+2}(a_{i,1}...a_{i,n}\omega)$$

Therefore, the map $f_{i,i+2}$ is homomorphism of $\Omega_i$-algebra. Therefore, the map $f_{i,i+2}$ is $T\star$-representation of $\Omega_i$-algebra $A_i$ in $\Omega_{i+1}$-algebra ${}^*A_{i+2}$. □

**Theorem 3.1.3.** *$(id, f_{i+1,i+2})$ is morphism of $T\star$-representations of $\Omega_i$-algebra from $f_{i,i+1}$ into $f_{i,i+2}$.*

Proof. Consider diagram (3.1.1) in more detail.

(3.1.6)

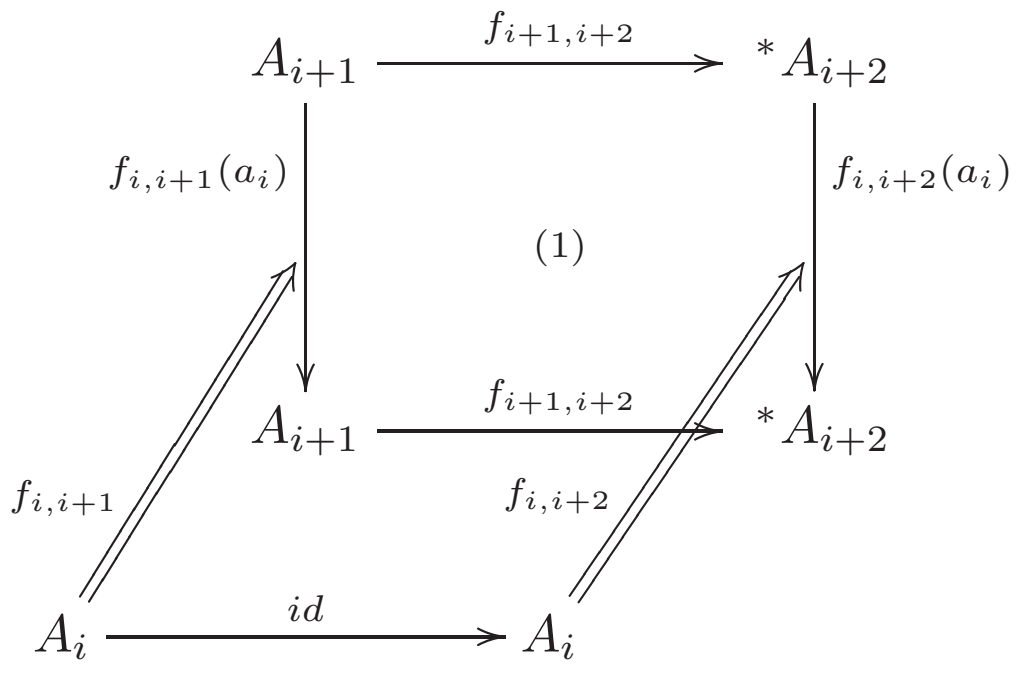

The statement of theorem follows from equation (3.1.2) and definition 2.2.2. □

**Theorem 3.1.4.** *Consider tower $(\overline{A}, \overline{f})$ of representations of $\overline{\Omega}$-algebras. Since identity transformation*

$$\delta_{i+2} : A_{i+2} \to A_{i+2}$$

*of $\Omega_{i+2}$-algebra $A_{i+2}$ belongs to representation $f_{i+1,i+2}$, then the representation $f_{i,i+2}$ of $\Omega_i$-algebra $A_i$ in $\Omega_{i+1}$-algebra ${}^*A_{i+2}$ can be extended to representation*

$$f'_{i,i+2} : A_i \to {}^*A_{i+2}$$

*of $\Omega_i$-algebra $A_i$ in $\Omega_{i+2}$-algebra $A_{i+2}$.*

Proof. According to the equation (3.1.2), for given $a_{i+1}\in A_{i+1}$ we can define a mapping

$$\begin{array}{c} f'_{i,i+2}:A_i\rightarrow {}^*A_{i+2}\\ f'_{i,i+2}(a_i)(f_{i+1,i+2}(a_{i+1})(a_{i+2}))=f_{i+1,i+2}(f_{i,i+1}(a_i)(a_{i+1}))(a_{i+2})\end{array}\tag{3.1.7}$$

Let there exist $a_{i+1}\in A_{i+1}$ such that

$$f_{i+1,i+2}(a_{i+1})=\delta_{i+2}\tag{3.1.8}$$

Then from (3.1.7),(3.1.8), it follows that

$$f'_{i,i+2}(a_i)(a_{i+2})=f_{i+1,i+2}(f_{i,i+1}(a_i)(a_{i+1}))(a_{i+2})\tag{3.1.9}$$

Let $\omega$ be $n$-ari operation of $\Omega_i$-algebra $A_i$. Since $f_{i,i+2}$ is homomorphism of $\Omega_i$-algebra, then, from the equation (3.1.2), it follows that

$$\begin{aligned}&f_{i,i+2}(a_{i\cdot 1}...a_{i\cdot n}\omega)(f_{i+1,i+2}(a_{i+1}))\\ =&f_{i,i+2}(a_{i\cdot 1})(f_{i+1,i+2}(a_{i+1}))...f_{i,i+2}(a_{i\cdot n})(f_{i+1,i+2}(a_{i+1}))\omega\\ =&f_{i+1,i+2}(f_{i,i+1}(a_{i\cdot 1})(a_{i+1})...f_{i,i+1}(a_{i\cdot n})(a_{i+1})\omega)\end{aligned}\tag{3.1.10}$$

From the equation (3.1.10), it follows that
(3.1.11)

$$\begin{aligned}&f_{i,i+2}(a_{i\cdot 1}...a_{i\cdot n}\omega)(f_{i+1,i+2}(a_{i+1})(a_{i+2}))\\ =&f_{i,i+2}(a_{i\cdot 1})(f_{i+1,i+2}(a_{i+1})(a_{i+2}))...f_{i,i+2}(a_{i\cdot n})(f_{i+1,i+2}(a_{i+1})(a_{i+2}))\omega\\ =&f_{i+1,i+2}(f_{i,i+1}(a_{i\cdot 1})(a_{i+1})...f_{i,i+1}(a_{i\cdot n})(a_{i+1})\omega)(a_{i+2})\end{aligned}$$

From equations (3.1.8), (3.1.9), (3.1.11), it follows that

$$\begin{aligned}&f'_{i,i+2}(a_{i\cdot 1}...a_{i\cdot n}\omega)(a_{i+2})\\ =&f'_{i,i+2}(a_{i\cdot 1})(a_{i+2})...f'_{i,i+2}(a_{i\cdot n})(a_{i+2})\omega\\ =&f_{i+1,i+2}(f_{i,i+1}(a_{i\cdot 1})(a_{i+1})...f_{i,i+1}(a_{i\cdot n})(a_{i+1})\omega)(a_{i+2})\end{aligned}\tag{3.1.12}$$

The equation (3.1.12) defines an operation $\omega$ on the set ${}^*A_{i+2}$

$$f'_{i,i+2}(a_{i\cdot 1})...f'_{i,i+2}(a_{i\cdot n})\omega=f'_{i,i+2}(a_{i\cdot 1}...a_{i\cdot n}\omega)\tag{3.1.13}$$

and the mapping $f'_{i,i+2}$ is homomorphism of $\Omega_i$-algebra. Therefore, we have a representation of $\Omega_i$-algebra $A_i$ in $\Omega_{i+2}$-algebra $A_{i+2}$. □

**Definition 3.1.5.** Consider the tower of representations

$$((A_1,A_2,A_3),(f_{1,2},f_{2,3}))$$

The map $f_*=(f_{1,2},f_{1,3})$ is called **representation of $\Omega_1$-algebra $A_1$ in representation $f_{2,3}$**. □

**Definition 3.1.6.** Consider the tower of representations $(\overline{A},\overline{f})$. The map

$$f_*=(f_{1,2},...,f_{1,n})$$

is called **representation of $\Omega_1$-algebra $A_1$ in tower of representations**

$$(\overline{A}_{[1]},\overline{f})=((A_2,...,A_n),(f_{2,3},...,f_{n-1,n}))$$

□

## 3.2. Morphism of Tower of $T*$-Representations

**Definition 3.2.1.** Consider the set of $\Omega_k$-algebr $A_k$, $B_k$, $k = 1, ..., n$. The set of mappings $(h_1, ..., h_n)$ is called **morphism from tower of $T\star$-representations $(\overline{A}, \overline{f})$ into tower of $T\star$-representations** $(\overline{B}, \overline{g})$, if for any $k$, $k = 1, ..., n-1$, the tuple of mappings $(h_k, h_{k+1})$ is morphism of $T\star$-representations from $f_{k,k+1}$ into $g_{k,k+1}$. □

For any $k$, $k = 1, ..., n-1$, we have diagram

(3.2.1)

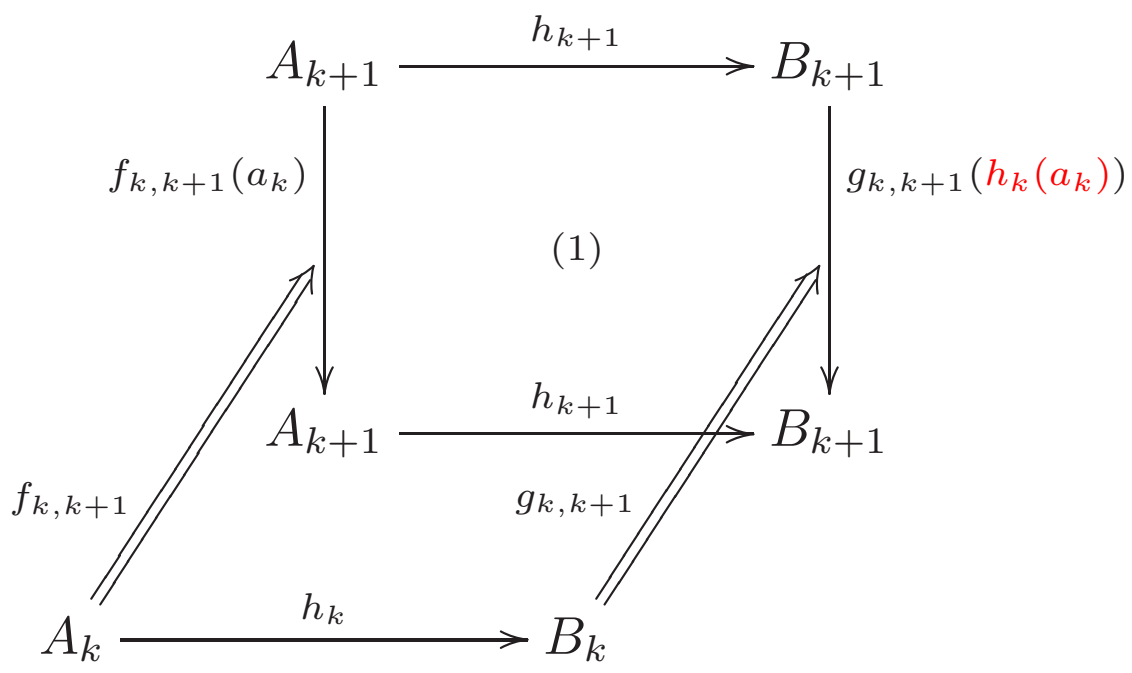

Equations

$$h_{k+1} \circ f_{k,k+1}(a_k) = g_{k,k+1}(h_k(a_k)) \circ h_{k+1} \tag{3.2.2}$$

$$h_{k+1}(f_{k,k+1}(a_k)(a_{k+1})) = g_{k,k+1}(h_k(a_k))(h_{k+1}(a_{k+1})) \tag{3.2.3}$$

express commutativity of diagram (1). However for morphism $(h_k, h_{k+1})$, $k > 1$, diagram (3.2.1) is not complete. Assuming similar diagram for morphism $(h_k, h_{k+1})$, this diagram on the top layer has form

(3.2.4)

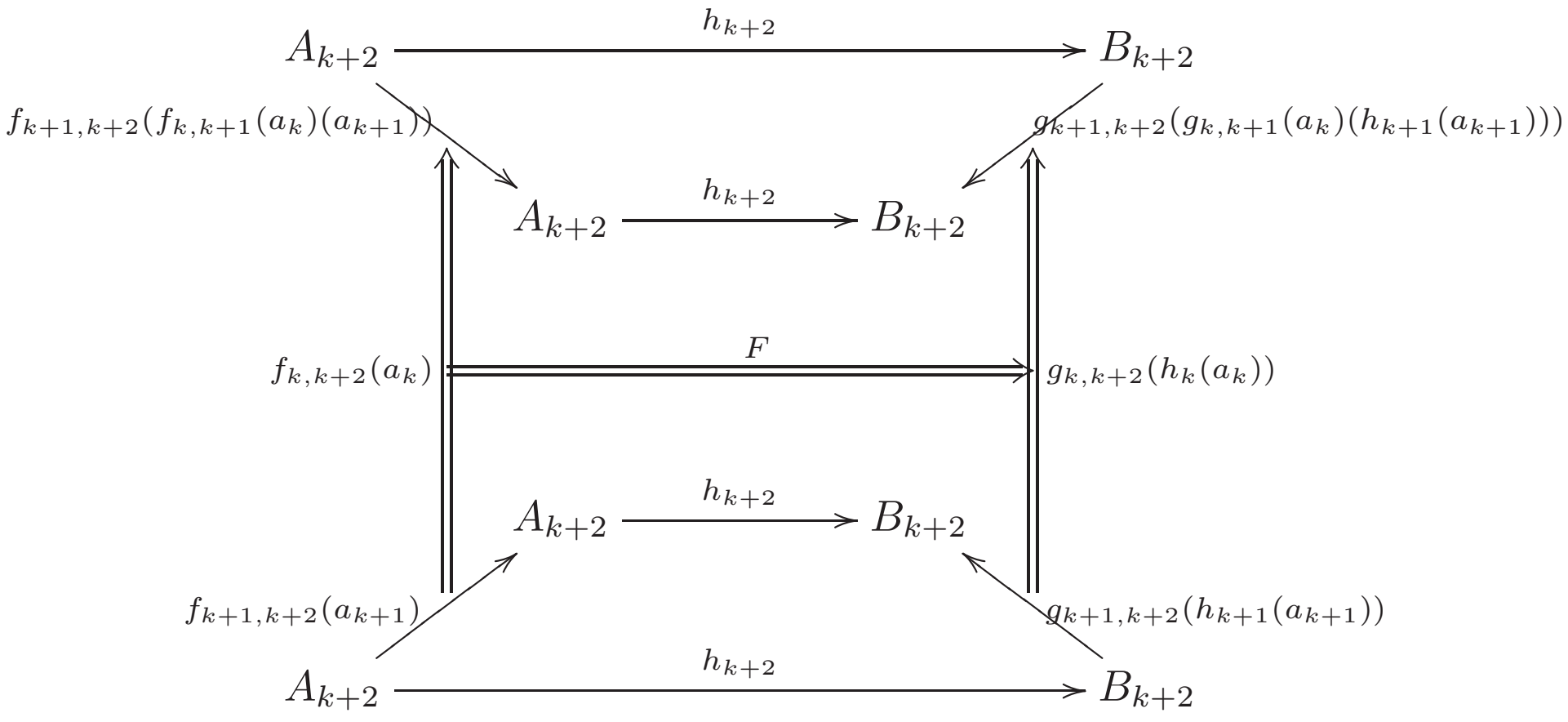

Unfortunately, the diagram (3.2.4) is not too informative. It is evident that there exists morphism from ${}^*A_{k+2}$ into ${}^*B_{k+2}$, mapping $f_{k,k+2}(a_k)$ into $g_{k,k+2}(h_k(a_k))$. However, the structure of this morphism is not clear from the diagram. We need consider mapping from ${}^*A_{k+2}$ into ${}^*B_{k+2}$, like we have done this in theorem 3.1.3.

**Theorem 3.2.2.** *Since $T\star$-representation $f_{i+1,i+2}$ is effective, then $(h_i, {}^*h_{i+2})$ is morphism of $T\star$-representations from $T\star$-representation $f_{i,i+2}$ into $T\star$-representation $g_{i,i+2}$ of $\Omega_i$-algebra.*

Proof. Consider the diagram

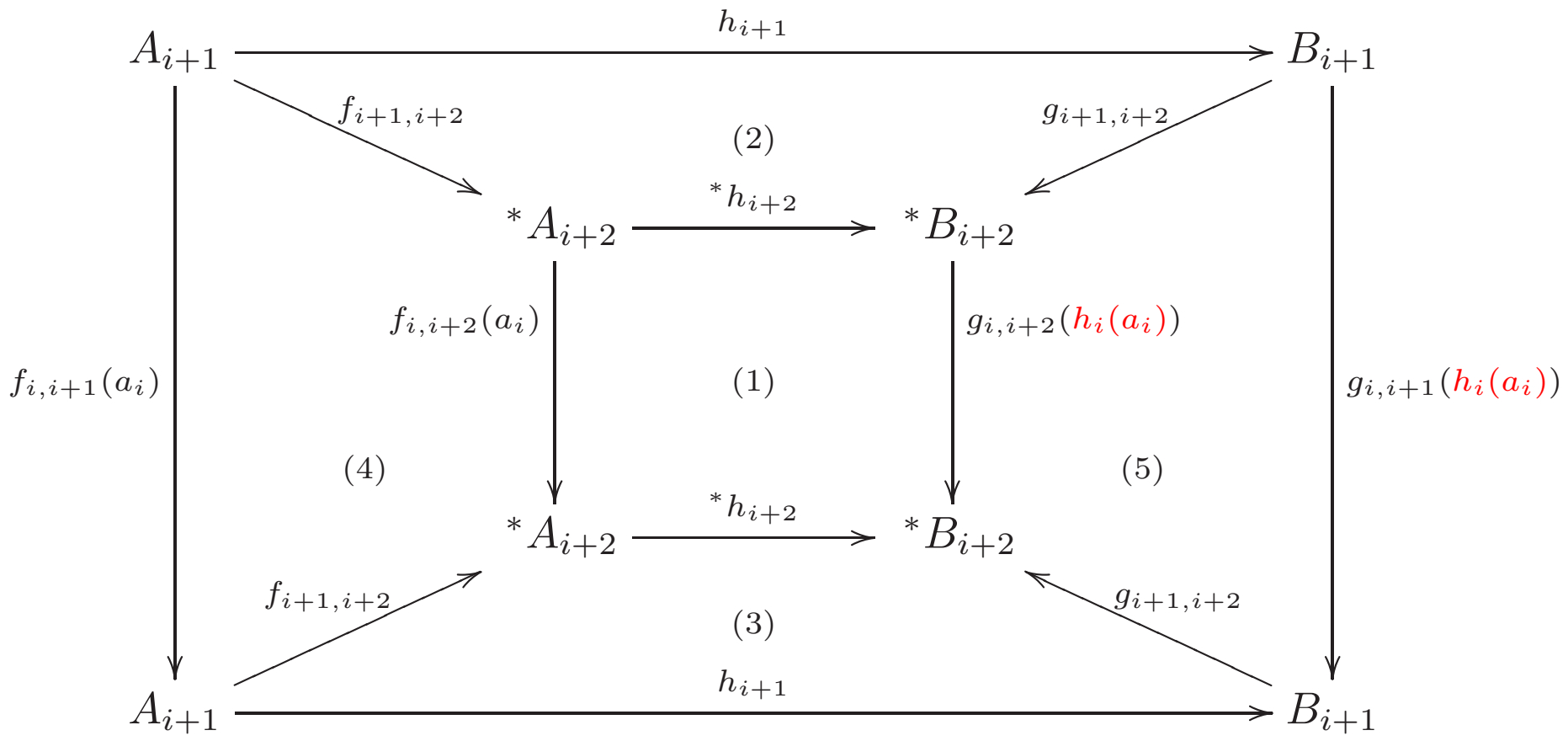

The existence of mapping ${}^*h_{i+2}$ and commutativity of the diagram (2) and (3) follows from effectiveness of mapping $f_{i+1,i+2}$ and theorem 2.2.7. Commutativity of diagrams (4) and (5) follows from theorem 3.1.3.

From commutativity of the diagram (4) it follows that

$$f_{i+1,i+2} \circ f_{i,i+1}(a_i) = f_{i,i+2}(a_i) \circ f_{i+1,i+2} \tag{3.2.5}$$

From equation (3.2.5) it follows that

$${}^*h_{i+2} \circ f_{i+1,i+2} \circ f_{i,i+1}(a_i) = {}^*h_{i+2} \circ f_{i,i+2}(a_i) \circ f_{i+1,i+2} \tag{3.2.6}$$

From commutativity of diagram (3) it follows that

$${}^*h_{i+2} \circ f_{i+1,i+2} = g_{i+1,i+2} \circ h_{i+1} \tag{3.2.7}$$

From equation (3.2.7) it follows

$${}^*h_{i+2} \circ f_{i+1,i+2} \circ f_{i,i+1}(a_i) = g_{i+1,i+2} \circ h_{i+1} \circ f_{i,i+1}(a_i) \tag{3.2.8}$$

From equations (3.2.6) and (3.2.8) it follows that

$${}^*h_{i+2} \circ f_{i,i+2}(a_i) \circ f_{i+1,i+2} = g_{i+1,i+2} \circ h_{i+1} \circ f_{i,i+1}(a_i) \tag{3.2.9}$$

From commutativity of the diagram (5) it follows that

$$g_{i+1,i+2} \circ g_{i,i+1}(h_i(a_i)) = g_{i,i+2}(h_i(a_i)) \circ g_{i+1,i+2} \tag{3.2.10}$$

From equation (3.2.10) it follows that

$$g_{i+1,i+2} \circ g_{i,i+1}(h_i(a_i)) \circ h_{i+1} = g_{i,i+2}(h_i(a_i)) \circ g_{i+1,i+2} \circ h_{i+1} \tag{3.2.11}$$

From commutativity of the diagram (2) it follows that

$${}^*h_{i+2} \circ f_{i+1,i+2} = g_{i+1,i+2} \circ h_{i+1} \tag{3.2.12}$$

From equation (3.2.12) it follows that

$$g_{i,i+2}(h_i(a_i)) \circ {}^*h_{i+2} \circ f_{i+1,i+2} = g_{i,i+2}(h_i(a_i)) \circ g_{i+1,i+2} \circ h_{i+1} \tag{3.2.13}$$

From equations (3.2.11) and (3.2.13) it follows that

$$g_{i+1,i+2} \circ g_{i,i+1}(h_i(a_i)) \circ h_{i+1} = g_{i,i+2}(h_i(a_i)) \circ {}^*h_{i+2} \circ f_{i+1,i+2} \tag{3.2.14}$$

External diagram is diagram (3.2.1) when $i=1$. Therefore, external diagram is commutative

$$h_{i+1}\circ f_{i,i+1}(a_i)=g_{i,i+1}(h_i(a_i))\circ h_{i+1} \tag{3.2.15}$$

From equation (3.2.15) it follows that

$$g_{i+1,i+2}\circ h_{i+1}\circ f_{i,i+1}(a_i)=g_{i+1,i+2}\circ g_{i,i+1}(h_i(a_i))\circ h_{i+1}(a_{i+1}) \tag{3.2.16}$$

From equations (3.2.9), (3.2.14) and (3.2.16) it follows that

$$^*h_{i+2}\circ f_{i,i+2}(a_i)\circ f_{i+1,i+2}=g_{i,i+2}(h_i(a_i))\circ {}^*h_{i+2}\circ f_{i+1,i+2} \tag{3.2.17}$$

Because the mapping $f_{i+1,i+2}$ is injection, then from equation (3.2.17) it follows that

$$^*h_{i+2}\circ f_{i,i+2}(a_i)=g_{i,i+2}(h_i(a_i))\circ {}^*h_{i+2} \tag{3.2.18}$$

From equation (3.2.18) commutativity of the diagram (1) follows. This proves the statement of theorem. □

Theorems 3.1.3 and 3.2.2 are true for any layer of tower of $T\star$-representations. In each particular case we need properly show sets and direction of the map. Meaning of given theorems is that all mappings in tower of $T\star$-representations act coherently.

The theorem 3.2.2 states that unknown mapping on the diagram (3.2.4) is the mapping $^*h_{i+2}$.

**Theorem 3.2.3.** *Consider the set of $\Omega_i$-algebr $A_i$, $B_i$, $C_i$, $i=1,\ldots,n$. Given morphisms of tower of representations*

$$\overline{p}:(\overline{A},\overline{f})\rightarrow(\overline{B},\overline{g})$$
$$\overline{q}:(\overline{B},\overline{g})\rightarrow(\overline{C},\overline{h})$$

*There exists morphism of representations of $\Omega$-algebra*

$$\overline{r}:(\overline{A},\overline{f})\rightarrow(\overline{C},\overline{h})$$

*where $r_k=q_kp_k$, $k=1,\ldots,n$. We call morphism $\overline{r}$ of tower of representations from $\overline{f}$ into $\overline{h}$* **product of morphisms $\overline{p}$ and $\overline{q}$ of tower of representations**.

Proof. For each $k$, $k = 2, ..., n$, we represent statement of theorem using diagram

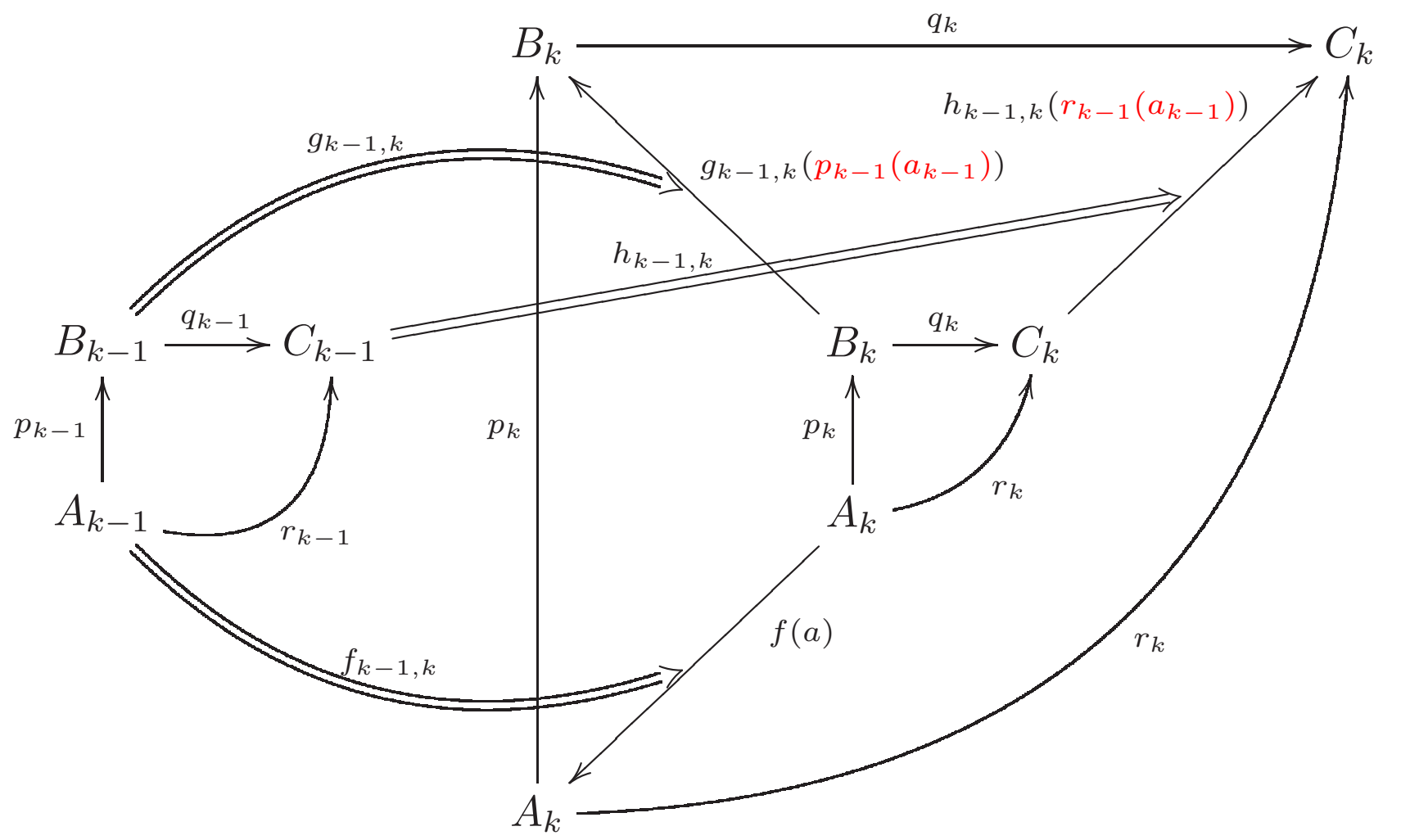

Mapping $r_{k-1}$ is homomorphism of $\Omega_{k-1}$-algebra $A_{k-1}$ into $\Omega_{k-1}$-algebra $C_{k-1}$. We need to show that tuple of mappings $(r_{k-1}, r_k)$ satisfies to (3.2.2):

$$\begin{aligned} r_k(f_{k-1,k}(a_{k-1})a_k) &= q_k p_k(f_{k-1,k}(a_{k-1})a_k) \\ &= q_k(g_{k-1,k}(p_{k-1}(a_{k-1}))p_k(a_k)) \\ &= h_{k-1,k}(q_{k-1}p_{k-1}(a_{k-1}))q_k p_k(a_k)) \\ &= h_{k-1,k}(r(a_{k-1}))r_k(a_k) \end{aligned}$$

□

## 3.3. Automorphism of Tower of Representations

**Definition 3.3.1.** Let $(\overline{A}, \overline{f})$ be tower of representations of $\overline{\Omega}$-algebras. The morphism of tower of representations $(\mathrm{id}, h_2, ..., h_n)$ such, that for each $k$, $k = 2, ..., n$, $h_k$ is endomorphism of $\Omega_k$-algebra $A_k$ is called **endomorphism of tower of representations** $\overline{f}$. □

**Definition 3.3.2.** Let $(\overline{A}, \overline{f})$ be tower of representations of $\overline{\Omega}$-algebras. The morphism of tower of representations $(\mathrm{id}, h_2, ..., h_n)$ such, that for each $k$, $k = 2, ..., n$, $h_k$ is automorphism of $\Omega_k$-algebra $A_k$ is called **automorphism of tower of representations** $f$. □

**Theorem 3.3.3.** *Let $(\overline{A}, \overline{f})$ be tower of representations of $\overline{\Omega}$-algebras. The set of automorphisms of the representation $\overline{f}$ forms group.*

Proof. Let $\overline{r}_{[1]}$, $\overline{p}_{[1]}$ be automorphisms of the tower of representations $\overline{f}$. According to definition 3.3.2, for each $k$, $k = 2, ..., n$, mappings $r_k$, $p_k$ are automorphisms of $\Omega_k$-algebra $A_k$. According to theorem II.3.2 ([12], p. 57), for each $k$, $k = 2, ..., n$, the mapping $r_k p_k$ is automorphism of $\Omega_k$-algebra $A_k$. From the theorem 3.2.3 and the definition 3.3.2, it follows that product of automorphisms $\overline{r}_{[1]}\overline{p}_{[1]}$ of the tower of representations $\overline{f}$ is automorphism of the tower of representations $\overline{f}$.

According to proof of the theorem 2.3.5, for any $k$, $k = 2, ..., n$, the product of automorphisms of $\Omega_i$-algebra is associative. Therefore, the product of automorphisms of tower of representations is associative.

Let $\overline{r}_{[1]}$ be an automorphism of the tower of representations $\overline{f}$. According to definition 3.3.2 for each $k$, $k = 2, ..., n$, the mapping $r_k$ is automorphism of $\Omega_k$-algebra $A_k$. Therefore, for each $k$, $k = 2, ..., n$, the mapping $r_k^{-1}$ is automorphism of $\Omega_k$-algebra $A_k$. The equation (3.2.3) is true for automorphism $\overline{r}_{[1]}$. Assume $a'_k = r_k(a_k)$, $k = 2, ..., n$. Since $r_k$, $k = 2, ..., n$, is automorphism then $a_k = r_k^{-1}(a'_k)$ and we can write (3.2.3) in the form

$$h_{k+1}(f_{k,k+1}(h_k^{-1}(a'_k))(h_{k+1}(a'_{k+1}))) = g_{k,k+1}(a'_k)(a'_{k+1}) \tag{3.3.1}$$

Since the mapping $h_{k+1}$ is automorphism of $\Omega_{k+1}$-algebra $A_{k+1}$, then from the equation (3.3.1) it follows that

$$f_{k,k+1}(h_k^{-1}(a'_k)(h_{k+1}(a'_{k+1}))) = h_{k+1}^{-1}(g_{k,k+1}(a'_k)(a'_{k+1})) \tag{3.3.2}$$

The equation (3.3.2) corresponds to the equation (3.2.3) for the mapping $\overline{r}_{[1]}^{-1}$. Therefore, mapping $\overline{r}_{[1]}^{-1}$ is automorphism of the representation $\overline{f}$. □

## 3.4. Generating Set of Tower of Representations

**Definition 3.4.1.** Tower of representations $(\overline{A}, \overline{f})$. is called **tower of effective representations**, if for any $i$ the representation $f_{i,i+1}$ is effective. □

**Theorem 3.4.2.** *Consider the tower of $T\star$-representations $(\overline{A}, \overline{f})$. Let representations $f_{i,i+1}$, ..., $f_{i+k-1,i+k}$ be effective. Then representation $f_{i,i+k}$ is effective.*

Proof. We will prove the statement of theorem by induction.

Let representations $f_{i,i+1}$, $f_{i+1,i+2}$ be effective. Assume that transformation $f_{i,i+1}(a_i)$ is not identity transformation. Then there exists $a_{i+1} \in A_{i+1}$ such that $f_{i,i+1}(a_i)(a_{i+1}) \neq a_{i+1}$. Because the representation $f_{i+1,i+2}$ is effective, then transformations $f_{i+1,i+2}(a_{i+1})$ and $f_{i+1,i+2}(f_{i,i+1}(a_i)(a_{i+1}))$ do not coincide. According to construction in theorem 3.1.2, the transformation $f_{i,i+2}(a_i)$ is not identity transformation.

Assume the statement of theorem is true for $k-1$ representations. Let $f_{i,i+1}$, ..., $f_{i+k-1,i+k}$ be effective. According assumption of induction, representations $f_{i,i+k-1}$, $f_{i+k-1,i+k}$ are effective. According to proven above, the representation $f_{i,i+k}$ is effective. □

**Theorem 3.4.3.** *Consider tower $(\overline{A}, \overline{f})$ of representations of $\overline{\Omega}$-algebras. Let identity transformation*

$$\delta_{i+2} : A_{i+2} \to A_{i+2}$$

*of $\Omega_{i+2}$-algebra $A_{i+2}$ belong to representation $f_{i+1,i+2}$. Let representations $f_{i,i+1}$, $f_{i+1,i+2}$ be effective. Then representation $f'_{i,i+2}$ defined in the theorem 3.1.4, is effective.*

Proof. Let there exist $a_{i+1} \in A_{i+1}$ such that

$$f_{i+1,i+2}(a_{i+1}) = \delta_{i+2}$$

Assume that the representation $f'_{i,i+2}$ is not effective. Then there exist

$$a_{i\cdot 1}, a_{i\cdot 2} \in A_i \quad a_{i\cdot 1} \neq a_{i\cdot 2} \tag{3.4.1}$$

such that

$$f'_{i,i+2}(a_{i\cdot 1})(a_{i+2}) = f'_{i,i+2}(a_{i\cdot 2})(a_{i+2}) \tag{3.4.2}$$

From equations (3.1.9), (3.4.2), it follows that

$$f_{i+1,i+2}(f_{i,i+1}(a_{i\cdot 1})(a_{i+1}))(a_{i+2}) = f_{i+1,i+2}(f_{i,i+1}(a_{i\cdot 2})(a_{i+1}))(a_{i+2}) \tag{3.4.3}$$

Since $a_{i+2}$ is arbitrary, then, from (3.4.3), it follows that

$$f_{i+1,i+2}(f_{i,i+1}(a_{i\cdot 1})(a_{i+1})) = f_{i+1,i+2}(f_{i,i+1}(a_{i\cdot 2})(a_{i+1})) \tag{3.4.4}$$

Since the representation $f_{i,i+1}$ is effective, then, from the condition (3.4.1), it follows that

$$f_{i,i+1}(a_{i\cdot 1})(a_{i+1}) \neq f_{i,i+1}(a_{i\cdot 2})(a_{i+1}) \tag{3.4.5}$$

From the condition (3.4.5) and the equation (3.4.4), it follows that the representation $f_{i+1,i+2}$ is not effective. The contradiction completes the proof of the theorem. □

We construct the basis of the tower of representations in a similar way that we constructed the basis of representation in the section 2.5.

We will write elements of tower of representations $(\overline{A}, \overline{f})$ as tuple

$$\overline{a} = (a_1, ..., a_3) \quad a_i \in A_i \quad i = 1, ..., n$$

We can interpret this record as

$$(a_2, a_3) = f_{2,3}(a_2)(a_3)$$

At the same time this record has additional meaning. For instance, in affine space, $a_3$ is a point of space, $a_2$ is a vector that determine the transformation of paralel transfer. Thus, we can interpret a tuple $(a_2, a_3)$ either as vector $a_2$ with tail in point $a_3$, or as head of this vector.

**Definition 3.4.4.** Let $(\overline{A}, \overline{f})$ be tower of representations. The tuple of sets

$$\overline{N}_{[1]} = (N_2 \subset A_2, ..., N_n \subset A_n)$$

is called **tuple of stable sets of tower of representations** $\overline{f}$, if

$$f_{i-1,i}(a_{i-1})(a_i) \in N_i \quad i = 2, ..., n$$

for every $a_1 \in A_1$, $a_2 \in N_2$, ..., $a_n \in N_n$. We also will say that tuple of sets

$$\overline{N}_{[1]} = (N_2 \subset A_2, ..., N_n \subset A_n)$$

is stable relative to tower of representations $\overline{f}$. □

**Theorem 3.4.5.** *Let $\overline{f}$ be tower of representations. Let set $N_i \subset A_i$ be subalgebra of $\Omega_i$-algebra $A_i$, $i = 2, ..., n$. Let tuple of sets*

$$\overline{N}_{[1]} = (N_2 \subset A_2, ..., N_n \subset A_n)$$

*be stable relative to tower of representations $\overline{f}$. Then there exists representation*

$$((A_1, N_2, ..., N_n), (f_{N_2,1,2}, ..., f_{N_n,n-1,n})) \tag{3.4.6}$$

*such that*

$$f_{N_i,i-1,i}(a_{i-1}) = f_{i-1,i}(a_{i-1})|_{N_i} \quad i = 2, ..., n$$

*The tower of representations (3.4.6) is called* **tower of subrepresentations**.

Proof. Let $\omega_{i-1,1}$ be $m$-ary operation of $\Omega_{i-1}$-algebra $A_{i-1}$, $i = 2, ..., n$. Then for any $a_{i-1,1}, ..., a_{i-1,m} \in N_{i-1}$ [3.2] and any $a_i \in N_i$

$$\begin{aligned}&(f_{N_i,i-1,i}(a_{i-1,1})...f_{N_i,i-1,i}(a_{i-1,m})\omega_{i-1,1})(a_i)\\ =&(f_{i-1,i}(a_{i-1,1})...f_{i-1,i}(a_{i-1,m})\omega_{i-1,1})(a_i)\\ =&f_{i-1,i}(a_{i-1,1}...a_{i-1,m}\omega_{i-1,1})(a_i)\\ =&f_{N_i,i-1,i}(a_{i-1,1}...a_{i-1,m}\omega_{i-1,1})(a_i)\end{aligned}$$

Let $\omega_{i,2}$ be $m$-ary operation of $\Omega_i$-algebra $A_i$, $i = 2, ..., n$. Then for any $a_{i,1}, ..., a_{i,n} \in N_i$ and any $a_{i-1} \in N_{i-1}$

$$\begin{aligned}&f_{N_i,i-1,i}(a_{i-1})(a_{i,1})...f_{N_i,i-1,i}(a_{i-1})(a_{i,m})\omega_{i,2}\\ =&f_{i-1,i}(a_{i-1})(a_{i,1})...f_{i-1,i}(a_{i-1})(a_{i,m})\omega_{i,2}\\ =&f_{i-1,i}(a_{i-1})(a_{i,1}...a_{i,m}\omega_{i,2})\\ =&f_{N_i,i-1,i}(a_{i-1})(a_{i,1}...a_{i,m}\omega_{i,2})\end{aligned}$$

We proved the statement of theorem. □

From theorem 3.4.5, it follows that if mapping $(f_{N_2,1,2}, ..., f_{N_n,1,n})$ is tower of subrepresentations of tower of representations $\overline{f}$, then mapping

$$(id : A_1 \to A_1, id_2 : N_2 \to A_2, ..., id_n : N_n \to A_n)$$

is morphism of towers of representations.

**Theorem 3.4.6.** *The set*[3.3] $\mathcal{B}_{\overline{f}}$ *of all towers of subrepresentations of tower of representations* $\overline{f}$ *generates a closure system on tower of representations* $\overline{f}$ *and therefore is a complete lattice.*

Proof. Let for given $\lambda \in \Lambda$, $K_{\lambda,i}$, $i = 2, ..., n$, be subalgebra of $\Omega_i$-algebra $A_i$ that is stable relative to representation $f_{i-1,i}$. We determine the operation of intersection on the set $\mathcal{B}_{\overline{f}}$ according to rule

$$\bigcap f_{K_{\lambda,i-1,i}} = f_{\cap K_{\lambda,i-1,i}} \quad i = 2, ..., n$$

$$\bigcap \overline{K}_\lambda = \left(K_1 = A_1, K_2 = \bigcap K_{\lambda,2}, ..., K_n = \bigcap K_{\lambda,n}\right)$$

$\cap K_{\lambda,i}$ is subalgebra of $\Omega_i$-algebra $A_i$. Let $a_i \in \cap K_{\lambda,i}$. For any $\lambda \in \Lambda$ and for any $a_{i-1} \in K_{i-1}$,

$$f_{i-1,i}(a_{i-1})(a_i) \in K_{\lambda,i}$$

Therefore,

$$f_{i-1,i}(a_{i-1})(a_i) \in K_i$$

Repeating this construction in the order of increment $i$, $i = 2, ..., n$, we see that $(K_1, ..., K_n)$ is tuple of stable sets of tower of representations $\overline{f}$. Therefore, we determined the operation of intersection of towers of subrepresentations properly. □

[3.2] Assume $N_1 = A_1$.

[3.3] This definition is similar to definition of the lattice of subalgebras ([12], p. 79, 80)

We denote the corresponding closure operator by $\overline{J}(\overline{f})$. If we denote $\overline{X}_{[1]}$ the tuple of sets $(X_2 \subset A_2, ..., X_n \subset A_n)$ then $\overline{J}(\overline{f}, \overline{X}_{[1]})$ is the intersection of all tuples $(K_1, ..., K_n)$ stable with respect to representation $\overline{f}$ and such that for $i = 2$, ..., $n$, $K_i$ is subalgebra of $\Omega_i$-algebra $A_i$ containing $X_i$.[3.4]

**Theorem 3.4.7.** *Let*[3.5] $\overline{f}$ *be the tower of representations. Let* $X_i \subset A_i$, $i = 2$, ..., $n$. *Assume* $Y_1 = A_1$. *Step by step increasing the value of* $i$, $i = 2$, ..., $n$, *we define a subsets* $X_{i,m} \subset A_i$ *by induction on* $m$.

$$X_{i,0} = X_i$$
$$x \in X_{i,m} => x \in X_{i,m+1}$$
$$x_1 \in X_{i,m}, ..., x_p \in X_{i,m}, \omega \in \Omega_i(p) => x_1...x_p\omega \in X_{i,m+1}$$
$$x_i \in X_{i,m}, x_{i-1} \in Y_{i-1} => f_{i-1,i}(x_{i-1})(x_i) \in X_{i,m+1}$$

*For each value of* $i$, *we assume*

$$Y_i = \bigcup_{m=0}^{\infty} X_{i,m}$$

*Then*

$$\overline{Y} = (Y_1, ..., Y_n) = \overline{J}(\overline{f}, \overline{X}_{[1]})$$

Proof. For each value of $i$ the proof of the theorem coincides with the proof of theorem 2.4.4. Because to define stable subset of $\Omega_i$-algebra $A_i$ we need only certain stable subset of $\Omega_{i-1}$-algebra $A_{i-1}$, we have to find stable subset of $\Omega_{i-1}$-algebra $A_{i-1}$ before we do this in $\Omega_i$-algebra $A_i$. □

$\overline{J}(\overline{f}, \overline{X}_{[1]})$ is called **tower of subrepresentations of tower of representations** $\overline{f}$ **generated by tuple of sets** $\overline{X}_{[1]}$, and $\overline{X}_{[1]}$ is a **tuple of generating sets of tower subrepresentations** $\overline{J}(\overline{f}, \overline{X}_{[1]})$. In particular, a **tuple of generating sets of tower of representations** $\overline{f}$ is a tuple $(X_2 \subset A_2, ..., X_n \subset A_n)$ such that $\overline{J}(\overline{f}, \overline{X}_{[1]}) = \overline{A}$.

It is easy to see that definition of the tuple of generating sets of tower of representations does not depend on whether representations of tower are effective or not. For this reason hereinafter we will assume that representations of the tower are effective and we will use convention for effective $T\star$-representation in remark 2.1.9.

We also will use notation

$$\overline{r}_{[1]} \circ \overline{a}_{[1]} = (r_2 \circ a_2, ..., r_n \circ a_n)$$

for image of tuple of elements $a_2 \in A_2$, ..., $a_n \in A_n$ under the endomorphism of tower of effective representations. According to the definition of product of mappings, for any endomorphisms $\overline{r}_{[1]}$, $\overline{s}_{[1]}$ the following equation is true

$$(\overline{r}_{[1]} \circ \overline{s}_{[1]}) \circ \overline{a}_{[1]} = \overline{r}_{[1]} \circ (\overline{s}_{[1]} \circ \overline{a}_{[1]}) \tag{3.4.7}$$

The equation (3.4.7) is associative law for $\circ$ and allows us to write expression

$$\overline{r}_{[1]} \circ \overline{s}_{[1]} \circ \overline{a}_{[1]}$$

---

[3.4] For $n = 2$, $J_2(f_{1,2}, X_2) = J_{f_{1,2}}(X_2)$. It would be easier to use common notation in sections 2.5 and 3.5. However I think that using of vector notation in section 2.5 is premature.

[3.5] The statement of theorem is similar to the statement of theorem 5.1, [12], p. 79.

without brackets.

From theorem 3.4.7, it follows next definition.

**Definition 3.4.8.** Let $(X_2 \subset A_2, ..., X_n \subset A_n)$ be tuple of sets. For each tuple of elements $\overline{a}$, $\overline{a} \in \overline{J}(\overline{f}, \overline{X}_{[1]})$, there exists tuple of $\overline{\Omega}$-words defined according to following rule.

(1) If $a_1 \in A_1$, then $a_1$ is $\Omega_1$-word.
(2) If $a_i \in X_i$, $i = 2, ..., n$, then $a_i$ is $\Omega_i$-word.
(3) If $a_{i,1}$, ..., $a_{i,p}$ are $\Omega_i$-words, $i = 2, ..., n$, and $\omega \in \Omega_i(p)$, then $a_{i,1}...a_{i,p}\omega$ is $\Omega_i$-word.
(4) If $a_i$ is $\Omega_i$-word, $i = 2, ..., n$, and $a_{i-1}$ is $\Omega_{i-1}$-word, then $a_{i-1}a_i$ is $\Omega_i$-word.

**Tuple of $\overline{\Omega}$-words**[3.6]

$$\overline{w}(\overline{f}, \overline{X}_{[1]}, \overline{a}) = (w_1(\overline{f}, \overline{X}_{[1]}, a_1), ..., w_n(\overline{f}, \overline{X}_{[1]}, a_n))$$

represents given element $\overline{a} \in \overline{J}(\overline{f}, \overline{X}_{[1]})$.[3.7] Denote $\overline{w}(\overline{f}, \overline{X}_{[1]})$ **the set of tuples of $\overline{\Omega}$-words of tower of representations $\overline{J}(\overline{f}, \overline{X}_{[1]})$.** □

Representation of $a_i \in A_i$ as $\Omega_i$-word is ambiguous. If $a_{i,1}$, ..., $a_{i,p}$ are $\Omega_i$-words, $\omega \in \Omega_i(p)$ and $a_{i-1} \in A_{i-1}$, then $\Omega_i$-words $a_{i-1}a_{i,1}...a_{i,p}\omega$ and $a_{i-1}a_{i,1}...a_{i-1}a_{i,p}\omega$ represent the same element of $\Omega_i$-algebra $A_i$. It is possible that there exist equations related with a character of a representation. For instance, if $\omega$ is the operation of $\Omega_{i-1}$-algebra $A_{i-1}$ and the operation of $\Omega_i$-algebra $A_i$, then we can request that $\Omega_i$-words $a_{i-1,1}...a_{i-1,p}\omega a_i$ and $a_{i-1,1}a_i...a_{i-1,p}a_i\omega$ represent the same element of $\Omega_i$-algebra $A_i$. Listed above equations for each value $i$, $i = 2, ..., n$, determine equivalence $r_i$ on the set of $\Omega_i$-words $W_i(\overline{f}, \overline{X}_{[1]})$.According to the construction, equivalence $r_i$ on the set of $\Omega_i$-words $W_i(\overline{f}, \overline{X}_{[1]})$ depends not only on the choice of the set $X_i$, but also on the choice of the set $X_{i-1}$.

**Theorem 3.4.9.** *Endomorphism $\overline{r}_{[1]}$ of tower of representations $\overline{f}$ forms the mapping of $\overline{\Omega}_{[1]}$-words*

$$\overline{w}_{[1]}[\overline{f}, \overline{X}_{[1]}, \overline{r}_{[1]}] : \overline{W}_{[1]}(\overline{f}, \overline{X}_{[1]}) \rightarrow \overline{w}_{[1]}(\overline{f}, \overline{X}'_{[1]}) \quad \overline{X}_{[1]} \subset \overline{A}_{[1]} \quad \overline{X}'_{[1]} = \overline{r}_{[1]}(\overline{X}_{[1]})$$

*such that for any $i$, $i = 2, ..., n$,*

(1) *If* $a_i \in X_i$, $a'_i = r_i \circ a_i$, *then*

$$w_i[\overline{f}, \overline{X}_{[1]}, \overline{r}_{[1]}](a_i) = a'_i$$

---

[3.6] According to our design

$$w_2(\overline{f}, \overline{X}_{[1]}, a_2) = w_2(f_{1,2}, X_2, a_2) \quad w_3(\overline{f}, \overline{X}_{[1]}, a_3) = w_3((f_{1,2}, f_{2,3}), (X_2, X_3), a_3)$$

[3.7] These notations have small difference with notations offered earlier in the theory of representation of universal algebra. Namely, we added a word of $\Omega_1$-algebra into tuple of $\overline{\Omega}$-words. However it is evident that this word is always trivial

$$w_1(\overline{f}, \overline{X}_{[1]}, a_1) = a_1$$

(2) *If*

$$a_{i,1},...,a_{i,n}\in w_i(\overline{f},\overline{X}_{[1]})$$

$$a'_{i,1}=w_i[\overline{f},\overline{X}_{[1]},\overline{r}_{[1]}](a_{i,1})\quad ...\quad a'_{i,p}=w_i[\overline{f},\overline{X}_{[1]},\overline{r}_{[1]}](a_{i,p})$$

*then for operation* $\omega\in\Omega_i(p)$ *holds*

$$w_i[\overline{f},\overline{X}_{[1]},\overline{r}_{[1]}](a_{i,1}...a_{i,p}\omega)=a'_{i,1}...a'_{i,p}\omega$$

(3) *If*

$$a_i\in w_i(\overline{f},\overline{X}_{[1]})\qquad a'_i=w_i[\overline{f},\overline{X}_{[1]},\overline{r}_{[1]}](a_i)$$
$$a_{i-1}\in w_{i-1}(\overline{f},\overline{X}_{[1]})\quad a'_{i-1}=w_{i-1}[\overline{f},\overline{X}_{[1]},\overline{r}_{[1]}](a_{i-1})$$

*then*

$$w_i[\overline{f},\overline{X}_{[1]},\overline{r}_{[1]}](a_{i-1}a_i)=a'_{i-1}a'_i$$

Proof. Statements (1), (2) of the theorem are true by definition of the endomorpism $r_i$. The statement (3) of the theorem follows from the equation (3.2.3). □

**Remark 3.4.10.** Let $\overline{r}_{[1]}$ be endomorphism of tower of representations $\overline{f}$. Let

$$\overline{a}_{[1]}\in\overline{J}_{[1]}(\overline{f},\overline{X}_{[1]})\quad \overline{a}'_{[1]}=\overline{r}_{[1]}\circ\overline{a}_{[1]}\quad \overline{X}'_{[1]}=\overline{r}_{[1]}\circ\overline{X}_{[1]}$$

The theorem 3.4.9 states that $\overline{a}'_{[1]}\in\overline{J}_{[1]}(\overline{f},\overline{X}'_{[1]})$ . The theorem 3.4.9 also states that $\overline{\Omega}_{[1]}$-word representing $\overline{a}_{[1]}$ relative $\overline{X}_{[1]}$ and $\overline{\Omega}_{[1]}$-word representing $\overline{a}'_{[1]}$ relative $\overline{X}'_{[1]}$ are generated according to the same algorithm. This allows considering of the set of $\overline{\Omega}_{[1]}$-words $\overline{w}_{[1]}(\overline{f},\overline{X}'_{[1]},\overline{a}'_{[1]})$ as mapping

$$\overline{W}_{[1]}(\overline{f},\overline{X}_{[1]},\overline{a}_{[1]}):\overline{X}'_{[1]}\rightarrow\overline{w}_{[1]}(\overline{f},\overline{X}'_{[1]},\overline{a}'_{[1]})$$
$$\overline{W}_{[1]}(\overline{f},\overline{X}_{[1]},\overline{a}_{[1]})(\overline{X}'_{[1]})=\overline{W}_{[1]}(\overline{f},\overline{X}_{[1]},\overline{a}_{[1]})\circ\overline{X}'_{[1]}$$

such that, if for certain endomorphism $\overline{r}_{[1]}$

$$\overline{X}'_{[1]}=\overline{r}_{[1]}\circ\overline{X}_{[1]}\quad \overline{a}'_{[1]}=\overline{r}_{[1]}\circ\overline{a}_{[1]}$$

then

$$\overline{W}_{[1]}(\overline{f},\overline{X}_{[1]},\overline{a}_{[1]})\circ\overline{X}'_{[1]}=\overline{w}_{[1]}(\overline{f},\overline{X}'_{[1]},\overline{a}'_{[1]})=\overline{a}'_{[1]}$$

Tuple of mappings[3.8]

$$\overline{W}_{[1]}(\overline{f},\overline{X}_{[1]},\overline{a}_{[1]})=(W_2(\overline{f},\overline{X}_{[1]},a_2),...,W_n(\overline{f},\overline{X}_{[1]},a_n))$$

---

[3.8] According to our design

$$W_2(\overline{f},\overline{X}_{[1]},a_2)=W_2(f_{1,2},X_2,a_2)\quad W_3(\overline{f},\overline{X}_{[1]},a_3)=W_3((f_{1,2},f_{2,3}),(X_2,X_3),a_3)$$

We do not care about the mapping $W_1(\overline{f},\overline{X},a_1)$ because this mapping is trivial

(3.4.8) $$W_1(\overline{f},\overline{X}_{[1]},a_1)\circ\overline{X}'_{[1]}=a_1$$

is called **tuple of coordinates of element $\overline{a}$ relative to tuple of sets $\overline{X}_{[1]}$**. Denote $\overline{W}_{[1]}(\overline{f},\overline{X}_{[1]})$ **the set of tuples of coordinates of tower of representations** $\overline{J}(\overline{f},\overline{X}_{[1]})$. Similarly, we consider coordinates of tuple of sets $\overline{B}_{[1]}\subset\overline{J}_{[1]}(f,X)$ relative to the tuple of sets $\overline{X}_{[1]}$

$$\overline{W}_{[1]}(\overline{f},\overline{X}_{[1]},\overline{B}_{[1]})=(W_2(\overline{f},\overline{X}_{[1]},B_2),...,W_n(\overline{f},\overline{X}_{[1]},B_n))$$

$$W_k(\overline{f},\overline{X}_{[1]},B_k)=\{W_k(\overline{f},\overline{X}_{(2...k)},a_k):a_k\in B_k\}=(W_k(\overline{f},\overline{X}_{(2...k)},a_k),a_k\in B_k)$$

□

**Theorem 3.4.11.** *There is a structure of $\Omega_k$-algebra on the set of coordinates $W_k(\overline{f},\overline{X}_{[1]})$.*

Proof. Let $\omega\in\Omega_k(n)$. Then for any $m_1$, ..., $m_n\in J_k(\overline{f},\overline{X}_{[1]})$ , we assume

$$W_k(\overline{f},\overline{X}_{[1]},m_1)...W_k(\overline{f},\overline{X}_{[1]},m_n)\omega=W_k(\overline{f},\overline{X}_{[1]},m_1...m_n\omega)\tag{3.4.9}$$

According to the remark 3.4.10,

$$\begin{aligned}(W_k(\overline{f},\overline{X}_{[1]},m_1)...W_k(\overline{f},\overline{X}_{[1]},m_n)\omega)\circ\overline{X}_{[1]}&=W_k(\overline{f},\overline{X}_{[1]},m_1...m_n\omega)\circ\overline{X}_{[1]}\\&=w(\overline{f},\overline{X}_{[1]},m_1...m_n\omega)\end{aligned}\tag{3.4.10}$$

follows from the equation (3.4.9). According to rule (3) of the definition 3.4.8, from the equation (3.4.10), it follows that

$$\begin{aligned}&(W_k(\overline{f},\overline{X}_{[1]},m_1)...W_k(\overline{f},\overline{X}_{[1]},m_n)\omega)\circ\overline{X}_{[1]}\\=&w_k(\overline{f},\overline{X}_{[1]},m_1)...w_k(\overline{f},\overline{X}_{[1]},m_n)\omega\\=&(W_k(\overline{f},\overline{X}_{[1]},m_1)\circ\overline{X}_{[1]})...(W_k(\overline{f},\overline{X}_{[1]},m_n)\circ\overline{X}_{[1]})\omega\end{aligned}\tag{3.4.11}$$

From the equation (3.4.11), it follows that the operation $\omega$ defined by the equation (3.4.9) on the set of coordinates $W_k(\overline{f},\overline{X}_{[1]})$ is defined properly. □

**Theorem 3.4.12.** *For $k=2$, ..., $n$, there exists the representation of $\Omega_{k-1}$-algebra $W_{k-1}(\overline{f},\overline{X}_{[1]})$ in $\Omega_k$-algebra $W_k(\overline{f},\overline{X}_{[1]})$.*[3.9]

Proof. Let $a_{k-1}\in J_{k-1}(\overline{f},\overline{X}_{[1]})$. Then for any $a_k\in J_k(\overline{f},\overline{X}_{[1]})$, we assume

$$W_{k-1}(\overline{f},\overline{X}_{[1]},a_{k-1})W_k(\overline{f},\overline{X}_{[1]},a_k)=W_k(\overline{f},\overline{X}_{[1]},a_{k-1}a_k)\tag{3.4.12}$$

According to the remark 3.4.10,

$$\begin{aligned}(W_{k-1}(\overline{f},\overline{X}_{[1]},a_{k-1})W_k(\overline{f},\overline{X}_{[1]},a_k))\circ\overline{X}_{[1]}&=W_k(\overline{f},\overline{X}_{[1]},am)\circ\overline{X}_{[1]}\\&=w_k(\overline{f},\overline{X}_{[1]},a_{k-1}a_k)\end{aligned}\tag{3.4.13}$$

follows from the equation (3.4.12). According to rule (4) of the definition 3.4.8, from the equation (3.4.13), it follows that

$$\begin{aligned}&(W_{k-1}(\overline{f},\overline{X}_{[1]},a_{k-1})W_k(\overline{f},\overline{X}_{[1]},a_k))\circ\overline{X}_{[1]}\\=&w_{k-1}(\overline{f},\overline{X}_{[1]},a_{k-1})w_k(\overline{f},\overline{X}_{[1]},a_k)\\=&(W_{k-1}(\overline{f},\overline{X}_{[1]},a_{k-1})\circ\overline{X}_{[1]})(W_k(\overline{f},\overline{X}_{[1]},a_k)\circ\overline{X}_{[1]})\end{aligned}\tag{3.4.14}$$

From the equation (3.4.14), it follows that the representation (3.4.12) of $\Omega_{k-1}$-algebra $W_{k-1}(\overline{f},\overline{X}_{[1]})$ in $\Omega_k$-algebra $W_k(\overline{f},\overline{X}_{[1]})$ is defined properly. □

[3.9]According to equation (3.4.8), we identify $\Omega_1$-algebras $W_1(\overline{f},\overline{X}_{[1]})$ and $A_1$.

**Corollary 3.4.13.** *Tuple of $\overline{\Omega}$-algebras $\overline{W}(\overline{f},\overline{X}_{[1]})$ forms the tower of representations.* □

**Definition 3.4.14.** Let $\overline{X}_{[1]}$ be the tuple of generating sets of tower of representations $\overline{f}$. Let $\overline{r}_{[1]}$ be the endomorphism of the tower of representations $\overline{f}$. The set of coordinates $\overline{W}_{[1]}(\overline{f},\overline{X}_{[1]},\overline{r}_{[1]}\circ\overline{X}_{[1]})$ is called **coordinates of endomorphism of tower of representations**. □

**Definition 3.4.15.** Let $\overline{X}_{[1]}$ be the tuple of generating sets of tower of representations $\overline{f}$. Let $\overline{r}_{[1]}$ be the endomorphism of the tower of representations $\overline{f}$. Let $\overline{a}_{[1]}\in\overline{A}_{[1]}$. We define **superposition of coordinates** of the tower of representations $\overline{f}$ and the element $\overline{a}_{[1]}$ as coordinates defined according to rule

$$\overline{W}_{[1]}(\overline{f},\overline{X}_{[1]},\overline{a}_{[1]})\circ\overline{W}_{[1]}(\overline{f},\overline{X}_{[1]},\overline{r}_{[1]}\circ\overline{X}_{[1]})=\overline{W}_{[1]}(\overline{f},\overline{X}_{[1]},\overline{r}_{[1]}\circ\overline{a}_{[1]}) \tag{3.4.15}$$

Let $\overline{Y}_{[1]}\subset\overline{A}_{[1]}$. We define superposition of coordinates of the tower of representations $\overline{f}$ and the tuple of sets $\overline{Y}_{[1]}$ according to rule

$$\begin{aligned}&\overline{W}_{[1]}(\overline{f},\overline{X}_{[1]},\overline{Y}_{[1]})\circ\overline{W}_{[1]}(\overline{f},\overline{X}_{[1]},\overline{r}_{[1]}\circ\overline{X}_{[1]})\\ =&(\overline{W}_{[1]}(\overline{f},\overline{X}_{[1]},\overline{a}_{[1]})\circ\overline{W}_{[1]}(\overline{f},\overline{X}_{[1]},\overline{r}_{[1]}\circ\overline{X}_{[1]}),\overline{a}_{[1]}\in\overline{Y}_{[1]})\end{aligned} \tag{3.4.16}$$

□

**Theorem 3.4.16.** *Endomorphism $\overline{r}_{[1]}$ of representation $\overline{f}$ generates the mapping of coordinates of tower of representations*

$$\overline{W}_{[1]}[\overline{f},\overline{X}_{[1]},\overline{r}_{[1]}]:\overline{W}_{[1]}(\overline{f},\overline{X}_{[1]})\rightarrow\overline{W}_{[1]}(\overline{f},\overline{X}_{[1]}) \tag{3.4.17}$$

*such that*

(3.4.18)

$$\overline{W}_{[1]}(\overline{f},\overline{X}_{[1]},\overline{a}_{[1]})\rightarrow\overline{W}_{[1]}[\overline{f},\overline{X}_{[1]},\overline{r}_{[1]}]\star\overline{W}_{[1]}(\overline{f},\overline{X}_{[1]},\overline{a}_{[1]})=\overline{W}_{[1]}(\overline{f},\overline{X}_{[1]},\overline{r}_{[1]}\circ\overline{a}_{[1]})$$
$$\overline{W}_{[1]}[\overline{f},\overline{X}_{[1]},\overline{r}_{[1]}]\star\overline{W}_{[1]}(\overline{f},\overline{X}_{[1]},\overline{a}_{[1]})=\overline{W}_{[1]}(\overline{f},\overline{X}_{[1]},\overline{a}_{[1]})\circ\overline{W}_{[1]}(\overline{f},\overline{X}_{[1]},\overline{r}_{[1]}\circ\overline{X}_{[1]})$$

Proof. According to the remark 3.4.10, we consider equations (3.4.15), (3.4.17) relative to given tuple of generating sets $\overline{X}_{[1]}$. The tuple of words

$$\overline{W}_{[1]}(\overline{f},\overline{X}_{[1]},\overline{a}_{[1]})\circ\overline{X}_{[1]}=\overline{w}_{[1]}(\overline{f},\overline{X}_{[1]},\overline{a}_{[1]}) \tag{3.4.19}$$

corresponds to coordinates $\overline{W}_{[1]}(\overline{f},\overline{X}_{[1]},\overline{a}_{[1]})$; the tuple of words

$$\overline{W}_{[1]}(\overline{f},\overline{X}_{[1]},\overline{r}_{[1]}\circ\overline{a}_{[1]})\circ X=\overline{w}_{[1]}(\overline{f},\overline{X}_{[1]},\overline{r}_{[1]}\circ\overline{a}_{[1]}) \tag{3.4.20}$$

corresponds to coordinates $\overline{W}_{[1]}(\overline{f},\overline{X}_{[1]},\overline{r}_{[1]}\circ\overline{a}_{[1]})$. Therefore, in order to prove the theorem, it is sufficient to show that the mapping $\overline{W}_{[1]}[\overline{f},\overline{X}_{[1]},\overline{r}_{[1]}]$ corresponds to mapping $\overline{w}_{[1]}[\overline{f},\overline{X}_{[1]},\overline{r}_{[1]}]$. We prove this statement by induction over number $n$ of universal algebras in tower of representations.

If $n=2$, then tower of representations $\overline{f}$ is representation of $\Omega_1$-algebra $A_1$ in $\Omega_2$-algebra $A_2$. The statement of theorem is corollary of theorem 2.4.13. Let the statement of theorem be true for $n-1$. We prove the theorem by induction over complexity of $\Omega_n$-word.

If $a_n\in X_n$, $a'_n=r_n\circ a_n$, then, according to equations (3.4.19), (3.4.20), mappings $W_n[\overline{f},\overline{X}_{[1]},\overline{r}_{[1]}]$ and $w_n[\overline{f},\overline{X}_{[1]},\overline{r}_{[1]}]$ are coordinated.

Let for $a_{n,1}$, ..., $a_{n,p} \in X_n$ mappings $W_n[\overline{f}, \overline{X}_{[1]}, \overline{r}_{[1]}]$ and $w_n[\overline{f}, \overline{X}_{[1]}, \overline{r}_{[1]}]$ are coordinated. Let $\omega \in \Omega_n(p)$. According to the theorem 2.4.9

$$W_n(\overline{f}, \overline{X}_{[1]}, a_{n,1}...a_{n,p}\omega) = W_n(\overline{f}, \overline{X}_{[1]}, a_{n,1})...W_n(\overline{f}, \overline{X}_{[1]}, a_{n,p})\omega \tag{3.4.21}$$

Because $r_n$ is endomorphism of $\Omega_n$-algebra $A_n$, then from the equation (3.4.21) it follows that

(3.4.22)
$$\begin{aligned} W_n(\overline{f}, \overline{X}_{[1]}, r_n \circ (a_{n,1}...a_{n,p}\omega)) &= W_n(\overline{f}, \overline{X}_{[1]}, (r_n \circ a_{n,1})...(r_n \circ a_{n,p})\omega) \\ &= W_n(\overline{f}, \overline{X}_{[1]}, r_n \circ a_{n,1})...W(f, X, r_n \circ a_{n,p})\omega \end{aligned}$$

From equations (3.4.21), (3.4.22) and the statement of induction, it follows that the mappings $W_n[\overline{f}, \overline{X}_{[1]}, \overline{r}_{[1]}]$ and $w_n[\overline{f}, \overline{X}_{[1]}, \overline{r}_{[1]}]$ are coordinated for $a_n = a_{n,1}...a_{n,p}\omega$.

Let for $a_{n,1} \in A_n$ mappings $W_n[\overline{f}, \overline{X}_{[1]}, \overline{r}_{[1]}]$ and $w_n[\overline{f}, \overline{X}_{[1]}, \overline{r}_{[1]}]$ are coordinated. Let $a \in A$. According to the theorem 3.4.12

$$W_n(\overline{f}, \overline{X}_{[1]}, a_{n-1}a_{n,1}) = W_{n-1}(\overline{f}, \overline{X}_{[1]}, a_{n-1})W_n(\overline{f}, \overline{X}_{[1]}, a_{n,1}) \tag{3.4.23}$$

Because $\overline{r}_{[1]}$ is endomorphism of tower of representations $\overline{f}$, then from the equation (3.4.23) it follows that

(3.4.24)
$$\begin{aligned} W_n(\overline{f}, \overline{X}_{[1]}, r_n \circ (a_{n-1}a_{n,1})) &= W_n(\overline{f}, \overline{X}_{[1]}, (r_{n-1} \circ a_{n-1})(r_n \circ a_{n,1})) \\ &= W_{n-1}(\overline{f}, \overline{X}_{[1]}, r_{n-1} \circ a_{n-1})W_n(\overline{f}, \overline{X}_{[1]}, r_n \circ a_{n,1}) \end{aligned}$$

From equations (3.4.23), (3.4.24) and the statement of induction, it follows that mappings $W_n[\overline{f}, \overline{X}_{[1]}, \overline{r}_{[1]}]$ and $w_n[\overline{f}, \overline{X}_{[1]}, \overline{r}_{[1]}]$ are coordinated for $a_{n,2} = a_{n-1}a_{n,1}$. □

**Corollary 3.4.17.** *Let $\overline{X}_{[1]}$ be the tuple of generating sets of the tower of representations $\overline{f}$. Let $\overline{r}_{[1]}$ be the endomorphism of the tower of representations $\overline{f}$. The mapping $\overline{W}_{[1]}[\overline{f}, \overline{X}_{[1]}, \overline{r}_{[1]}]$ is endomorphism of tower of representations $\overline{W}(\overline{f}, \overline{X}_{[1]})$.* □

Hereinafter we will identify mapping $\overline{W}_{[1]}[\overline{f}, \overline{X}_{[1]}, \overline{r}_{[1]}]$ and the set of coordinates $\overline{W}_{[1]}(\overline{f}, \overline{X}_{[1]}, \overline{r}_{[1]} \circ \overline{X}_{[1]})$.

**Theorem 3.4.18.** *Let $\overline{X}_{[1]}$ be the tuple of generating sets of the tower of representations $\overline{f}$. Let $\overline{r}_{[1]}$ be the endomorphism of the tower of representations $\overline{f}$. Let $\overline{Y}_{[1]} \subset \overline{A}_{[1]}$. Then*

$$\overline{W}_{[1]}(\overline{f}, \overline{X}_{[1]}, \overline{Y}_{[1]}) \circ \overline{W}_{[1]}(\overline{f}, \overline{X}_{[1]}, \overline{r}_{[1]} \circ \overline{X}_{[1]}) = \overline{W}_{[1]}(\overline{f}, \overline{X}_{[1]}, \overline{r}_{[1]} \circ \overline{Y}_{[1]}) \tag{3.4.25}$$

$$\overline{W}_{[1]}[\overline{f}, \overline{X}_{[1]}, \overline{r}_{[1]}] \star \overline{W}_{[1]}(\overline{f}, \overline{X}_{[1]}, \overline{Y}_{[1]}) = \overline{W}_{[1]}(\overline{f}, \overline{X}_{[1]}, \overline{r}_{[1]} \circ \overline{Y}_{[1]}) \tag{3.4.26}$$

Proof. The equation (3.4.25) follows from the equation

$$\overline{r}_{[1]} \circ \overline{Y}_{[1]} = (\overline{r}_{[1]} \circ \overline{a}_{[1]}, \overline{a}_{[1]} \in \overline{Y}_{[1]})$$

as well from equations (3.4.15), (3.4.16). The equation (3.4.26) is corollary of equations (3.4.25), (3.4.18). □

**Theorem 3.4.19.** *Let $\overline{X}_{[1]}$ be the tuple of generating sets of the tower of representations $\overline{f}$. Let $\overline{r}_{[1]}$, $\overline{s}_{[1]}$ be the endomorphisms of the tower of representations $\overline{f}$.*

*Then*

(3.4.27)

$$\overline{W}_{[1]}(\overline{f},\overline{X}_{[1]},\overline{s}_{[1]}\circ\overline{X}_{[1]})\circ\overline{W}_{[1]}(\overline{f},\overline{X}_{[1]},\overline{r}_{[1]}\circ\overline{X}_{[1]})=\overline{W}_{[1]}(\overline{f},\overline{X}_{[1]},\overline{r}_{[1]}\circ\overline{s}_{[1]}\circ\overline{X}_{[1]})$$

$$\overline{W}_{[1]}[\overline{f},\overline{X}_{[1]},\overline{r}_{[1]}]\star\overline{W}_{[1]}[\overline{f},\overline{X}_{[1]},\overline{s}_{[1]}]=\overline{W}_{[1]}[\overline{f},\overline{X}_{[1]},\overline{r}_{[1]}\circ\overline{s}_{[1]}] \tag{3.4.28}$$

Proof. The equation (3.4.27) follows from the equation (3.4.25), if we assume $\overline{Y}_{[1]}=\overline{s}_{[1]}\circ\overline{X}_{[1]}$. The equation (3.4.28) follows from the equation (3.4.27) and chain of equations

$$\begin{aligned}
&(\overline{W}_{[1]}[\overline{f},\overline{X}_{[1]},\overline{r}_{[1]}]\star\overline{W}_{[1]}[\overline{f},\overline{X}_{[1]},\overline{s}_{[1]}])\star\overline{W}_{[1]}(\overline{f},\overline{X}_{[1]},\overline{Y}_{[1]})\\
=&\overline{W}_{[1]}[\overline{f},\overline{X}_{[1]},\overline{r}_{[1]}]\star(\overline{W}_{[1]}[\overline{f},\overline{X}_{[1]},\overline{s}_{[1]}]\star\overline{W}_{[1]}(\overline{f},\overline{X}_{[1]},\overline{Y}_{[1]}))\\
=&(\overline{W}_{[1]}(\overline{f},\overline{X}_{[1]},\overline{Y}_{[1]})\circ\overline{W}_{[1]}(\overline{f},\overline{X}_{[1]},\overline{s}_{[1]}\circ\overline{X}_{[1]}))\circ\overline{W}_{[1]}(\overline{f},\overline{X}_{[1]},\overline{r}_{[1]}\circ\overline{X}_{[1]})\\
=&\overline{W}_{[1]}(\overline{f},\overline{X}_{[1]},\overline{Y}_{[1]})\circ(\overline{W}_{[1]}(\overline{f},\overline{X}_{[1]},\overline{s}_{[1]}\circ\overline{X}_{[1]})\circ\overline{W}_{[1]}(\overline{f},\overline{X}_{[1]},\overline{r}_{[1]}\circ\overline{X}_{[1]}))\\
=&\overline{W}_{[1]}(\overline{f},\overline{X}_{[1]},\overline{Y}_{[1]})\circ\overline{W}_{[1]}(\overline{f},\overline{X}_{[1]},\overline{r}_{[1]}\circ\overline{s}_{[1]}\circ\overline{X}_{[1]})\\
=&\overline{W}_{[1]}(\overline{f},\overline{X}_{[1]},\overline{r}_{[1]}\circ\overline{s}_{[1]})\star\overline{W}_{[1]}(\overline{f},\overline{X}_{[1]},\overline{Y}_{[1]})
\end{aligned}$$

□

We can generalize the definition of the superposition of coordinates and assume that one of the factors is a tuple of sets of $\overline{\Omega}_{[1]}$-words. Accordingly, the definition of the superposition of coordinates has the form

$$\begin{aligned}
&\overline{w}_{[1]}(\overline{f},\overline{X}_{[1]},\overline{Y}_{[1]})\circ\overline{W}_{[1]}(\overline{f},\overline{X}_{[1]},\overline{r}_{[1]}\circ\overline{X}_{[1]})\\
=&\overline{W}_{[1]}(\overline{f},\overline{X}_{[1]},\overline{Y}_{[1]})\circ\overline{w}_{[1]}(\overline{f},\overline{X}_{[1]},\overline{r}_{[1]}\circ\overline{X}_{[1]})\\
=&\overline{w}_{[1]}(\overline{f},\overline{X}_{[1]},\overline{r}_{[1]}\circ\overline{Y}_{[1]})
\end{aligned}$$

The following forms of writing an image of a tuple of sets $\overline{Y}_{[1]}$ under endomorphism $\overline{r}_{[1]}$ are equivalent

$$\begin{aligned}
\overline{r}_{[1]}\circ\overline{Y}_{[1]}&=\overline{W}_{[1]}(\overline{f},\overline{X}_{[1]},Y)\circ(\overline{r}_{[1]}\circ\overline{X}_{[1]})\\
&=\overline{W}_{[1]}(\overline{f},\overline{X}_{[1]},\overline{Y}_{[1]})\circ(\overline{W}_{[1]}(\overline{f},\overline{X}_{[1]},\overline{r}_{[1]}\circ\overline{X}_{[1]})\circ\overline{X}_{[1]})
\end{aligned} \tag{3.4.29}$$

From equations (3.4.25), (3.4.29), it follows that

$$\begin{aligned}
&(\overline{W}_{[1]}(\overline{f},\overline{X}_{[1]},\overline{Y}_{[1]})\circ\overline{W}_{[1]}(\overline{f},\overline{X}_{[1]},\overline{r}_{[1]}\circ\overline{X}_{[1]}))\circ\overline{X}_{[1]}\\
=&\overline{W}_{[1]}(\overline{f},\overline{X}_{[1]},\overline{Y}_{[1]})\circ(\overline{W}_{[1]}(\overline{f},\overline{X}_{[1]},\overline{r}_{[1]}\circ\overline{X}_{[1]})\circ\overline{X}_{[1]})
\end{aligned} \tag{3.4.30}$$

The equation (3.4.30) is associative law for composition and allows us to write expression

$$\overline{W}_{[1]}(\overline{f},\overline{X}_{[1]},\overline{Y}_{[1]})\circ\overline{W}_{[1]}(\overline{f},\overline{X}_{[1]},\overline{r}_{[1]}\circ\overline{X}_{[1]})\circ X$$

without brackets.

**Definition 3.4.20.** Let $(X_2\subset A_2,...,X_n\subset A_n)$ be tuple of generating sets of tower of representations $\overline{f}$. Let map $\overline{r}_{[1]}$ be endomorphism of tower of representations $\overline{f}$. Let the tuple of sets $\overline{X}'_{[1]}=\overline{r}_{[1]}\circ\overline{X}_{[1]}$ be image of tuple of sets $\overline{X}_{[1]}$ under the mapping $\overline{r}_{[1]}$. Endomorphism $\overline{r}_{[1]}$ of tower of representations $\overline{f}$ is called **regular on tuple of generating sets** $\overline{X}_{[1]}$, if the tuple of sets $\overline{X}'_{[1]}$ is tuple of generating sets of tower of representations $\overline{f}$. Otherwise, endomorphism $\overline{r}_{[1]}$ is called **singular on tuple of generating sets** $\overline{X}_{[1]}$, □

**Definition 3.4.21.** Endomorphism of tower of representations $\overline{f}$ is called **regular**, if it is regular on any tuple of generating sets. □

**Theorem 3.4.22.** *Automorphism $\overline{r}$ of tower of representations $\overline{f}$ is regular endomorphism.*

PROOF. Let $\overline{X}_{[1]}$ be tuple of generating sets of tower of representations $\overline{f}$. Let $\overline{X}'_{[1]} = \overline{r}_{[1]} \circ \overline{X}_{[1]}$.

According to theorem 3.4.16 endomorphism $\overline{r}_{[1]}$ forms the mapping of $\overline{\Omega}_{[1]}$-words $\overline{w}_{[1]}[\overline{f}, \overline{X}_{[1]}, \overline{r}_{[1]}]$.

Let $\overline{a}'_{[1]} \in \overline{A}_{[1]}$. Since $\overline{r}_{[1]}$ is automorphism, then there exists $\overline{a}_{[1]} \in \overline{A}_{[1]}$, $\overline{r}_{[1]} \circ \overline{a}_{[1]} = \overline{a}'_{[1]}$. According to definition 3.4.8, $\overline{w}_{[1]}(\overline{f}, \overline{X}_{[1]}, \overline{a}_{[1]})$ is tuple of words representing $\overline{a}_{[1]}$ relative to tuple of generating sets $\overline{X}_{[1]}$. According to theorem 3.4.16, $\overline{w}_{[1]}(\overline{f}, \overline{X}'_{[1]}, \overline{a}'_{[1]})$ is tuple of words representing $\overline{a}'_{[1]}$ relative to tuple of sets $\overline{X}'_{[1]}$

$$\overline{w}_{[1]}(\overline{f}, \overline{X}'_{[1]}, \overline{a}'_{[1]}) = \overline{w}_{[1]}[\overline{f}, \overline{X}_{[1]}, \overline{r}](\overline{w}_{[1]}(\overline{f}, \overline{X}_{[1]}, \overline{a}_{[1]}))$$

Therefore, $\overline{X}'_{[1]}$ is generating set of representation $\overline{f}$. According to definition 3.4.21, automorphism $\overline{r}_{[1]}$ is regular. □

## 3.5. Basis of Tower of Representations

**Definition 3.5.1.** If the tuple of sets $\overline{X}_{[1]}$ is tuple of generating sets of tower of representations $\overline{f}$, then any tuple of sets $\overline{Y}_{[1]}$, $X_i \subset Y_i \subset A_i$, $i = 2, ..., n$, also is tuple of generating sets of tower of representations $\overline{f}$. If there exists tuple of minimal sets $\overline{X}_{[1]}$ generating the tower of representations $\overline{f}$, then the tuple of sets $\overline{X}_{[1]}$ is called **basis of tower of representations** $\overline{f}$. □

**Theorem 3.5.2.** *We define a basis of tower of representations by induction over $n$. For $n = 2$, the basis of tower of representations is the basis of representation $f_{1,2}$. If tuple of sets $\overline{X}_{[1,n]}$ is basis of tower of representations $\overline{f}_{[1]}$, then the tuple of generating sets $\overline{X}_{[1]}$ of tower of representations $\overline{f}$ is basis iff for any $a_n \in X_n$ the tuple of sets $(X_2, ..., X_{n-1}, X_n \setminus \{a_n\})$ is not tuple of generating sets of tower of representations $\overline{f}$.*

PROOF. For $n = 2$, the statement of the theorem is corollary of the theorem 2.5.2.

Let $n > 2$. Let $\overline{X}_{[1]}$ be tuple of generating sets of tower of representations $\overline{f}$. Let tuple of sets $\overline{X}_{[1,n]}$ be basis of tower of representations $\overline{f}_{[n]}$. Assume that for some $a_n \in X_n$ there exist word

$$w_n = w_n(a_n, \overline{f}, (X_1, ..., X_{n-1}, X_n \setminus \{a_n\})) \tag{3.5.1}$$

Consider $a'_n \in A_n$ such that it has word

$$w'_n = w_n(a'_n, \overline{f}, \overline{X}_{[1]}) \tag{3.5.2}$$

that depends on $a_n$. According to the definition 2.4.6, any occurrence $a_n$ into word $w'_n$ can be substituted by the word $w_n$. Therefore, the word $w'_n$ does not depend on $a_n$, and the tuple of sets $(X_2, ..., X_{n-1}, X_n \setminus \{a_n\})$ is the tuple of generating sets of tower of representations $\overline{f}$. Therefore, $\overline{X}_{[1]}$ is not basis of the tower of representations $\overline{f}$. □

**Remark 3.5.3.** The proof of the theorem 3.5.2 gives us effective method for constructing the basis of tower of representations $\overline{f}$. We start to build a basis in the lowest layer. When the basis is constructed in layer $i$, $i = 2, ..., n-1$, we can proceed to the construction of basis in layer $i+1$. □

**Remark 3.5.4.** We write a basis also in following form

$$X_k = (x_k, x_k \in X_k) \quad k = 2, ..., n$$

If basis is finite, then we also use notation

$$X_k = (x_{k\cdot i}, i \in I_k) = (x_{k\cdot 1}, ..., x_{k\cdot p_k}) \quad k = 2, ..., n$$

□

**Remark 3.5.5.** Coordinates $a_k \in A_k$ relative to basis $\overline{X}_{[1]}$ of the tower of representations $\overline{f}$ are defined ambiguously. Just as in the remark 2.5.5, we consider **tuple of equivalence $\overline{\rho}_{[1]}(\overline{f})$ generated by tower of representations $\overline{f}$.** If during the construction, we obtain the equality of two $\Omega_k$-word relative to given basis, then we can say without worrying about the equivalence $\rho_k(\overline{f})$ that these $\Omega_k$-words are equal. □

**Theorem 3.5.6.** *Automorphism of the tower of representations $\overline{f}$ maps a basis of the tower of representations $\overline{f}$ into basis.*

Proof. For $n = 2$, the statement of theorem is corollary of the theorem 3.5.6.

Let the mapping $\overline{r}$ be automorphism of the tower of representations $\overline{f}$. Let the tuple of sets $\overline{X}_{[1]}$ be a basis of the tower of representations $\overline{f}$. Let $\overline{X}'_{[1]} = \overline{r}_{[1]}(\overline{X}_{[1]})$.

Assume that the tuple of sets $\overline{X}'_{[1]}$ is not basis. According to the theorem 3.5.2 there exist $i$, $i = 2, ..., n$, and $a'_i \in X'_i$ such that the tuple of sets $\overline{X}''_{[1]}$[3.10] is tuple of generating sets of the tower of representations $\overline{f}$. According to the theorem 3.3.3 the mapping $\overline{r}^{-1}$ is automorphism of the tower of representations $\overline{f}$. According to the theorem 3.4.22 and definition 3.4.21, the tuple of sets $\overline{X}'''_{[1]}$[3.11] is tuple of generating sets of the tower of representations $\overline{f}$. The contradiction completes the proof of the theorem. □

[3.10] $X''_j = X'_j$, $j \ne i$, $X''_i = X'_i \setminus \{x'_i\}$
[3.11] $X'''_j = X_j$, $j \ne i$, $X'''_i = X_i \setminus \{x_i\}$

CHAPTER 4

# Geometry of Division Ring

## 4.1. Center of Division Ring

**Definition 4.1.1.** Let $D$ be a ring.[4.1] The set $Z(D)$ of elements $a \in D$ such that

$$ax = xa \tag{4.1.1}$$

for all $x \in D$, is called **center of ring** $D$. □

**Theorem 4.1.2.** *The center $Z(D)$ of ring $D$ is subring of ring $D$.*

Proof. The statement follows immediately from definition 4.1.1. □

**Theorem 4.1.3.** *The center $Z(D)$ of division ring $D$ is subfield of division ring $D$.*

Proof. According to theorem 4.1.2 it is enough to verify that $a^{-1} \in Z(D)$ if $a \in Z(D)$. Let $a \in Z(D)$. Repeatedly using the equation (4.1.1) we get chain of equations

$$aa^{-1}x = x = xaa^{-1} = axa^{-1} \tag{4.1.2}$$

From (4.1.2) it follows

$$a^{-1}x = xa^{-1}$$

Therefore, $a^{-1} \in Z(D)$. □

**Definition 4.1.4.** Let $D$ be a ring with unit element $e$.[4.2] The map

$$l : Z \to D$$

such that $l(n) = ne$ is a homomorphism of rings, and its kernel is an ideal $(n)$, generated by integer $n \geq 0$. We have canonical injective homomorphism

$$Z/nZ \to D$$

which is an isomorphism between $Z/nZ$ and subring of $D$. If $nZ$ is prime ideal, then we have two cases.

- $n = 0$. $D$ contains as subring a ring which isomorphic to $Z$, and which is often identified with $Z$. In that case, we say that $D$ has **characteristic** 0.
- $n = p$ for some prime number $p$. $D$ has **characteristic** $p$, and $D$ contains an isomorphic image of $F_p = Z/pZ$.

□

**Theorem 4.1.5.** *Let $D$ be ring of characteristic 0 and let $d \in D$. Then every integer $n \in Z$ commutes with $d$.*

---

[4.1] [1], page 89.

[4.2] I made definition according to definition from [1], pages 89, 90.

Proof. We prove statement by induction. The statement is evident for $n = 0$ and $n = 1$. Let statement be true for $n = k$. From chain of equation

$$(k+1)d = kd + d = dk + d = d(k+1)$$

evidence of statement for $n = k+1$ follows. □

**Theorem 4.1.6.** *Let $D$ be ring of characteristic 0. Then ring of integers $Z$ is subring of center $Z(D)$ of ring $D$.*

Proof. Corollary of theorem 4.1.5. □

Let $D$ be division ring. If $D$ has characteristic 0, $D$ contains as subfield an isomorphic image of the field $Q$ of rational numbers. If $D$ has characteristic $p$, $D$ contains as subfield an isomorphic image of $F_p$. In either case, this subfield will be called the prime field. Since the prime field is the smallest subfield of $D$ containing 1 and has no automorphism except identity, it is customary to identify it with $Q$ or $F_p$ as the case may be.

**Theorem 4.1.7.** *Let $D$ be division ring of characteristic 0 and let $d \in D$. Then for any integer $n \in Z$*

$$n^{-1}d = dn^{-1} \tag{4.1.3}$$

Proof. According to theorem 4.1.5 following chain of equation is true

$$n^{-1}dn = nn^{-1}d = d \tag{4.1.4}$$

Let us multiply right and left sides of equation (4.1.4) by $n^{-1}$. We get

$$n^{-1}d = n^{-1}dnn^{-1} = dn^{-1} \tag{4.1.5}$$

(4.1.3) follows from (4.1.5). □

**Theorem 4.1.8.** *Let $D$ be division ring of characteristic 0 and let $d \in D$. Then every rational number $p \in Q$ commutes with $d$.*

Proof. Let us represent rational number $p \in Q$ as $p = mn^{-1}$, $m$, $n \in Z$. Statement of theorem follows from chain of equations

$$pd = mn^{-1}d = n^{-1}dm = dmn^{-1} = dp$$

based on the statement of theorem 4.1.5 and equation (4.1.3). □

**Theorem 4.1.9.** *Let $D$ be division ring of characteristic 0. Then field of rational numbers $Q$ is subfield of center $Z(D)$ of division ring $D$.*

Proof. Corollary of theorem 4.1.8. □

## 4.2. Geometry of Division Ring over Field

We may consider division ring $D$ as vector space over field $F \subset Z(D)$. We will use following conventions.

(1) We do not use standatrd notation for vector when we consider element of division ring $D$ considering it as vector over field $F$. However we will use another collor for index for notation of coordinates of element of division ring $D$ as vector over field $F$.
(2) Because $F$ is field, we can write all indexes on right side of root letter.

Let $\overline{\overline{e}}$ be basis of division ring $D$ over field $F$. Then we may present any element $a \in D$ as

$$a = \overline{e}_i a^i \quad a^i \in F \tag{4.2.1}$$

When dimension of division ring $D$ over field $F$ infinite, then basis may be either countable, or its power is not less than power of continuum. If basis is countable, then we put constraints on coefficients $a^i$ of expansion (4.2.1). If power of the set $I$ is continuum, then we assume that there is measure on the set $I$ and sum in expansion (4.2.1) is integral over this measure.

**Remark 4.2.1.** Since we defined product in the division ring $D$, we consider the division ring as algebra over field $F \subset Z(D)$. For elements of basis we assume

$$\overline{e}_i \overline{e}_j = \overline{e}_k C^k_{ij} \tag{4.2.2}$$

Coefficients $C^k_{ij}$ of expansion (4.2.2) are called **structural constants of division ring $D$ over field $F$.**

From equations (4.2.1), (4.2.2), it follows

$$ab = \overline{e}_k C^k_{ij} a^i b^j \tag{4.2.3}$$

From equation (4.2.3) it follows that

$$(ab)c = \overline{e}_k C^k_{ij} (ab)^i c^j = \overline{e}_k C^k_{ij} C^i_{mn} a^m b^n c^j \tag{4.2.4}$$

$$a(bc) = \overline{e}_k C^k_{ij} a^i (bc)^j = \overline{e}_k C^k_{ij} a^i C^j_{mn} b^m c^n \tag{4.2.5}$$

From associativity of product

$$(ab)c = a(bc)$$

and equations (4.2.4) and (4.2.5) it follows that

$$\overline{e}_k C^k_{ij} C^i_{mn} a^m b^n c^j = \overline{e}_k C^k_{ij} a^i C^j_{mn} b^m c^n \tag{4.2.6}$$

Because vectors $a$, $b$, $c$ are arbitrary, and vectors $\overline{e}_k$ are linear independent, then from equation (4.2.6) it follows that

$$C^k_{jn} C^j_{im} = C^k_{ij} C^j_{mn} \tag{4.2.7}$$

**Theorem 4.2.2.** *Coordinates $a^j$ of vector $a$ are tensor*

$$a^j = A^j_i a'^i \tag{4.2.8}$$

Proof. Let $\overline{\overline{e}}'$ be another basis. Let

$$\overline{e}'_i = \overline{e}_j A^j_i \tag{4.2.9}$$

be transformation, mapping basis $\overline{\overline{e}}$ into basis $\overline{\overline{e}}'$. Because vector $a$ does not change, then

$$a = \overline{e}'_i a'^i = \overline{e}_j a^j \tag{4.2.10}$$

From equations (4.2.9) and (4.2.10) it follows that

$$\overline{e}_j a^j = \overline{e}'_i a'^i = \overline{e}_j A^j_i a'^i \tag{4.2.11}$$

Because vectors $\overline{e}_j$ are linear independent, then equation (4.2.8) follows from equation (4.2.11). Therefore, coordinates of vector are tensor. □

**Theorem 4.2.3.** *Structural constants of division ring $D$ over field $F$ are tensor*

$$A^l_k B'^k_{ij} A^{-1}{}^i_n A^{-1}{}^j_m = C^l_{nm} \tag{4.2.12}$$

Proof. Let us consider similarly the transformation of product. Equation (4.2.3) has form

$$ab = \overline{e}'_k B'^k_{ij} a'^i b'^j \tag{4.2.13}$$

relative to basis $\overline{e}'$. Let us substitute (4.2.8) and (4.2.9) into (4.2.13). We get

$$ab = \overline{e}_l A^l_k B'^k_{ij} a^n A^{-1\cdot i}_{\ \ n} b^m A^{-1\cdot j}_{\ \ m} \tag{4.2.14}$$

From (4.2.3) and (4.2.14) it follows that

$$\overline{e}_l A^l_k B'^k_{ij} A^{-1\cdot i}_{\ \ n} a^n A^{-1\cdot j}_{\ \ m} b^m = \overline{e}_l C^l_{nm} a^n b^m \tag{4.2.15}$$

Because vectors $a$ and $b$ are arbitrary, and vectors $\overline{e}_l$ are linear independent, then equation (4.2.12) follows from equation (4.2.15). Therefore, structural constants are tensor. □

CHAPTER 5

# Quadratic Map of Division Ring

## 5.1. Bilinear Map of Division Ring

**Theorem 5.1.1.** *Let $D$ be division ring of characteristic $0$. Let $F$ be field which is subring of center of division ring $D$. Let $\overline{\overline{e}}$ be basis in division ring $D$ over field $F$.* **Standard representation over field $F$ of bilinear map of division ring** *has form*

$$f(a,b) = {}_{1\cdot}f^{ijk}\ {}_ie\ a\ {}_je\ b\ {}_ke + {}_{2\cdot}f^{ijk}\ {}_ie\ b\ {}_je\ a\ {}_ke \tag{5.1.1}$$

*Expression ${}_{t\cdot}f^{ijk}$ in equation (5.1.1) is called* **standard component over field $F$ of bilinear map $f$**.

Proof. The statement of theorem is corollary of theorem [6]-9.2.4. Here we have transpositions

$$\sigma_1 = \begin{pmatrix} a & b \\ a & b \end{pmatrix} \quad \sigma_2 = \begin{pmatrix} a & b \\ b & a \end{pmatrix} \tag{5.1.2}$$

□

**Theorem 5.1.2.** *Let field $F$ be subring of center $Z(D)$ of division ring $D$. Let $\overline{e}_i$ be basis of division ring $D$ over field $F$. We can represent bilinear function*

$$g : D \times D \to D$$

*as $D$-valued bilinear form over field $F$*

$$g(a,b) = a^i\ b^j\ {}_{ij}g \tag{5.1.3}$$

*where*

$$a = a^i\ {}_i\overline{e}$$

$$b = b^j\ {}_j\overline{e}$$

$${}_{ij}g = g({}_i\overline{e}, {}_j\overline{e}) \tag{5.1.4}$$

*and values ${}_{ij}g$ are coordinates of $D$-valued covariant tensor over field $F$.*

Proof. The statement of theorem is corollary of theorem [6]-9.2.5. □

The matrix $G = |{}_{ij}g|$ is called **matrix of bilinear function**. If matrix $G$ is nonsingular, then bilinear function $g$ is called **nonsingular**.

**Theorem 5.1.3.** *Bilinear map*

$$g : D \times D \to D$$

*is symmetric iff the matrix $G$ is symmetric*

$${}_{ij}g = {}_{ji}g$$

PROOF. The statement of theorem is corollary of theorem [6]-9.2.6. □

**Theorem 5.1.4.** *Bilinear map*

$$g : D \times D \to D$$

*is skew symmetric iff the matrix $G$ is skew symmetric*

$$_{ij}g = -_{ji}g$$

PROOF. The statement of theorem is corollary of theorem [6]-9.2.7. □

**Theorem 5.1.5.** *Components of bilinear map*

$$g : D \times D \to D$$

*and its matrix over field $F$ satisfy to equation*

$$_{pq}g = ({}_{1\cdot}g^{ijk}\ C^s_{ip}\ C^l_{sj}\ C^t_{lq}\ C^r_{tk}\ + {}_{2\cdot}g^{ijk}\ C^s_{iq}\ C^l_{sj}\ C^t_{lp}\ C^r_{tk})\ \overline{e}_r \tag{5.1.5}$$

$$_{pq}g^r = {}_{1\cdot}g^{ijk}\ C^s_{ip}\ C^l_{sj}\ C^t_{lq}\ C^r_{tk} + {}_{2\cdot}g^{ijk}\ C^s_{iq}\ C^l_{sj}\ C^t_{lp}\ C^r_{tk} \tag{5.1.6}$$

PROOF. The statement of theorem is corollary of theorem [6]-9.2.8. To prove statement it is enough to substitute (5.1.2) into equations [6]-(9.2.12), [6]-(9.2.13). □

**Theorem 5.1.6.** *Suppose bilinear map*

$$g : D \times D \to D$$

*has matrix $G = \|_{ij}g\|$. Then there exists bilinear map*

$$g' : D \times D \to D$$

*with matrix $G' = \|_{ij}g'\|$, $_{ij}g' = {}_{ji}g$.*

PROOF. Assume

$$_{1\cdot}g'^{ijk} = {}_{2\cdot}g^{ijk} \quad {}_{2\cdot}g'^{ijk} = {}_{1\cdot}g^{ijk} \tag{5.1.7}$$

From equations (5.1.6), (5.1.7) it follows that

$$\begin{aligned}
{}_{qp}g'^r &= {}_{1\cdot}g'^{ijk}\ C^s_{iq}\ C^l_{sj}\ C^t_{lp}\ C^r_{tk} + {}_{2\cdot}g'^{ijk}\ C^s_{ip}\ C^l_{sj}\ C^t_{lq}\ C^r_{tk} \\
&= {}_{2\cdot}g^{ijk}\ C^s_{iq}\ C^l_{sj}\ C^t_{lp}\ C^r_{tk} + {}_{1\cdot}g^{ijk}\ C^s_{ip}\ C^l_{sj}\ C^t_{lq}\ C^r_{tk} \\
&= {}_{pq}g^r
\end{aligned}$$

□

## 5.2. Quadratic Map of Division Ring

**Definition 5.2.1.** Let $D$ be division ring. The map

$$h : D \to R$$

is called **quadratic**, if there exists biliear map

$$g : D \times D \to D$$

such that

$$h(a) = g(a, a)$$

Bilinear map $g$ is called map associated with the map $h$. □

For given map $h$ the map $g$ is defined ambiguously. However according to theorem 5.1.6, we always can assume that bilinear map $g$ is symmetric

$$g(a,b)=g(b,a)$$

Indeed, if bilinear map $g_1$, associated with map $h$, is not symmetric, then we assume

$$g(a,b)=(g_1(a,b)+g_1(b,a))/2$$

Because bilinear map is homogeneous of degree 1 over field $R$ with respect to each variable, quadratic map is homogeneous of degree 2 over field $R$.

**Theorem 5.2.2.** *Let $D$ be division ring of characteristic 0. Let $F$ be field which is subring of center of division ring $D$. Let $\overline{\overline{p}}$ be basis in division ring $D$ over field $F$.* ***Standard representation of quadratic map of division ring over field $F$*** *has form*

$$f(a)=f^{ijk}\ {}_{i}\overline{p}\ a\ {}_{j}\overline{p}\ a\ {}_{k}\overline{p} \tag{5.2.1}$$

*Expression $f^{ijk}$ in equation (5.2.1) is called* **standard component of quadratic map $f$ over field $F$**.

Proof. The statement of theorem is corollary of theorem 5.1.1. □

**Theorem 5.2.3.** *Let field $F$ be subring of center $Z(D)$ of division ring $D$. Let $\overline{\overline{e}}$ be basis of division ring $D$ over field $F$. We can represent quadratic map $f$ as $D$-valued quadratic form over field $F$*

$$f(a)=a^{i}a^{j}\ {}_{ij}f \tag{5.2.2}$$

*where*

$$\overline{a}=a^{i}\ {}_{i}\overline{e}$$

$${}_{ij}f=g({}_{i}\overline{e},{}_{j}\overline{e}) \tag{5.2.3}$$

*and $g$ is associated bilinear map. Values ${}_{ij}f$ are coordinates of $D$-valued covariant tensor over field $F$.*

Proof. For selected bilinear map $g$ associated to map $f$, the statement of theorem is corollary of theorem 5.1.2. We need to prove independence of coordinates of quadratic map from choice of an associated map $g$.

Assume, $g_1$, $g_2$ are bilinear maps associated to the map $f$. Then

$${}_{ii}f_1=g_1({}_{i}\overline{e},{}_{i}\overline{e})=g_2({}_{i}\overline{e},{}_{i}\overline{e})={}_{ii}f_2 \tag{5.2.4}$$

For arbitrary $a\in D$,

$$a^{i}\ a^{j}\ {}_{ij}f_1=a^{i}\ a^{j}\ {}_{ij}f_2$$

If for given $i$, $j$ we consider $a\in D$ such that $a_i=a_j=1$, then we get

$${}_{ii}f_1+{}_{ij}f_1+{}_{ji}f_1+{}_{jj}f_1={}_{ii}f_2+{}_{ij}f_2+{}_{ji}f_2+{}_{jj}f_2 \tag{5.2.5}$$

From equations (5.2.4), (5.2.5) it follows that

$${}_{ij}f_1+{}_{ji}f_1={}_{ij}f_2+{}_{ji}f_2$$

Therefore, if associated bilinear map is symmetric, then this map is determined unique. □

The matrix $F=|{}_{ij}f|$ is called **matrix of quadratic map**. If the matrix $F$ is regular, then the quadratic map $f$ is called **regular**. The rank of the matrix $F$ is called **rank of quadratic map $f$**.

**Theorem 5.2.4.** *The quadratic over field $F$ map*

$$f : D \to D$$

*is defined as map over division ring iff*

$$ {}_{pq}f = f^{ijk}\ C^{s}_{ip}\ C^{l}_{sj}\ C^{t}_{lq}\ C^{r}_{tk}\ {}_{r}\overline{p} \tag{5.2.6}$$

$$ {}_{pq}f^{r} = f^{ijk}\ C^{s}_{ip}\ C^{l}_{sj}\ C^{t}_{lq}\ C^{r}_{tk} \tag{5.2.7}$$

Proof. Let us substitute

$$\overline{a} = a^{j}\ {}_{j}\overline{p}$$

in equation (5.2.1). Then equation (5.2.1) gets form

$$\begin{aligned} f(a) =& f^{ijk}\ {}_{i}\overline{p}\ a^{p}\ {}_{p}\overline{p}\ {}_{j}\overline{p}\ a^{q}\ {}_{q}\overline{p}\ {}_{k}\overline{p} \\ =& a^{p}\ a^{q}\ f^{ijk}\ {}_{i}\overline{p}\ {}_{p}\overline{p}\ {}_{j}\overline{p}\ {}_{q}\overline{p}\ {}_{k}\overline{p} \\ =& a^{p}\ a^{q}\ f^{ijk}\ C^{s}_{ip}\ C^{l}_{sj}\ C^{t}_{lq}\ C^{r}_{tk}\ {}_{r}\overline{p} \end{aligned} \tag{5.2.8}$$

Equation (5.2.6) follows from comparision of equations (5.2.8) and (5.2.2). □

From theorem 5.2.3, it follows that coordinates ${}_{pq}f$ form tensor. We described the structure of acceptable transformations in theorem [6]-9.1.13.

**Theorem 5.2.5.** *Let field $F$ be subring of center $Z(D)$ of division ring $D$. Let $\overline{\overline{p}}$ be basis of division ring $D$ over field $F$. Let $C^{k}_{ij}$ be structural constants of division ring $D$ over field $F$. If*

$$\det(a^{i}\ (C^{k}_{ij}\ +\ C^{k}_{ji})) \neq 0$$

*then equation*

$$ax + xa = b \tag{5.2.9}$$

*has the unique solution.*

*If determinant equal 0, then $F$-linear dependence of vector $b$ from vectors $a^{i}\ (C^{k}_{ij}\ +\ C^{k}_{ji})\ {}_{k}\overline{p}$ is condition of existence of solution. In this case, equation (5.2.9) has infinitely many solutions. Otherwise equation does not have solution.*

Proof. From equations (4.2.3), (5.2.9), it follows that

$$a^{i}\ x^{j}\ C^{k}_{ij}\ {}_{k}\overline{p} + x^{i}\ a^{j}\ C^{k}_{ij}\ {}_{k}\overline{p} = b^{k}\ {}_{k}\overline{p} \tag{5.2.10}$$

Since vectors ${}_{k}\overline{p}$ are linear independent, then from equation (5.2.10) it follows that coordinates $x^{j}$ satisfy to the system of linear equations

$$x^{j}\ a^{i}\ (C^{k}_{ij}\ +\ C^{k}_{ji}) = b^{k} \tag{5.2.11}$$

Therefore, if determinant of the system (5.2.11) does not equal 0, then the system has unique solution.

Similarly, if determinant equal 0, then equality of the rank of the extended matrix and the rank of the matrix of the system of linear equations is the condition of existence of solution. □

**Theorem 5.2.6.** *Let field $F$ be subring of center $Z(D)$ of division ring $D$. There exists a quadratic map $f$ such that we can find a basis $\overline{\overline{p}}$ over field $F$ such that map $f$ is represented as sum of squares.*

Proof. Reduction of quadratic form to diagonal representation is done the same way as is done in [4], p. 169 - 172. Considered transformations of variables are associated with the corresponding transformation of basis. We are not concerned in whether these transformations of basis are acceptable. The reason is that the basis chosen initially may not be connected with the canonical basis by acceptable transformation. But for us it is important that the final presentation of mapping $f$ should be acceptable.

We will prove the statement by induction on number of variables in representation of quadratic map. For this we write the quadratic map $f$ in the form of map of $n$ variables $a^1$, ..., $a^n$

$$f(a)=f(a^1,...,a^n)=a^ia^j\ {}_{ij}f \tag{5.2.12}$$

Representation of quadratic map in the form (5.2.12) is called **quadratic form**.

The statement is evident, when quadratic form depends on one variable $a^1$, because in this case the quadratic form has form

$$f(a^1)=(a^1)^2{}_{11}f$$

Assume, the statement of the theorem is true for quadratic form of $n-1$ variables.

(1) Let the quadratic form holds squares of variables. Without loss of generality, assume ${}_{11}f\neq 0$. Let us consider the map

$$f_1(a^1,...,a^n)=f(a^1,...,a^n)-{}_{11}f^{-1}\ (a^j\ {}_jh)^2 \tag{5.2.13}$$

We select values ${}_jh$ so as to satisfy the equation

$$2\ {}_{11}f\ {}_{1j}f={}_1h\ {}_jh+{}_jh\ {}_1h \tag{5.2.14}$$

When $j=1$, from equation (5.2.14) it follows that

$${}_1h=\ {}_{11}f \tag{5.2.15}$$

When $j>1$, from equations (5.2.14), (5.2.15) it follows that

$$2\ {}_{11}f\ {}_{1j}f={}_{11}f\ {}_jh+{}_jh\ {}_{11}f \tag{5.2.16}$$

From theorem 5.2.5 it follows that if

$$\det({}_{11}f^i\ (C^k_{ij}\ +\ C^k_{ji}))\neq 0 \tag{5.2.17}$$

then equation (5.2.16) has unique solution. We assume that condition (5.2.17) is satisfied and consider change of variables

$$b^1=a^j\ {}_jh$$
$$b^i=a^i \qquad\qquad i>1$$

Therefore, the map $f_1(a^1,...,a^n)$ does not depend on $b^1$. According to definition (5.2.13), we represented the map $f(a^1,...,a^n)$ as sum of square of variable $b^1$ and quadratic map that does not depend on $b^1$.

$$f(a^1,...,a^n)={}_{11}f^{-1}\ (b^1)^2+f_1(b^2,...,b^n)$$

Since form $f_1$ depends on $n-1$ variable, we can find a transformation of variables that leads this form to the sum of squares of variables.

(2) If ${}_{11}f = ... = {}_{nn}f = 0$, then we need an additional linear transformation that lead to appearance of squares of variables. Assume ${}_{12}f \ne 0$. Let us consider the linear transformation

$$a^1 = b^1 - b^2$$
$$a^2 = b^1 + b^2$$
$$a^i = b^i \quad i > 2$$

As a result of this transformation term

$$2\ a^1\ a^2\ {}_{12}f$$

gets form

$$2\ (b^1)^2\ {}_{12}f - 2\ (b^2)^2\ {}_{12}f$$

We got quadratic form considered in case (1). □

If all coefficients of the quadratic form in canonical representation are real, then $f(a) \in R$ for any $a \in D$.

**Definition 5.2.7.** Let $D$ be division ring. Quadratic map

$$f : D \to D$$

is called **positive definite** when

$$f(a) \in R$$

for any $a \in D$ and following condition is true

$$f(a) \ge 0$$
$$f(a) = 0 => a = 0$$

Positive definite quadratic map $f$ is called **Euclidean metric on division ring** $D$. Symmetric associative bilinear map $g$ is called **Euclidean scalar product on division ring** $D$. □

If quadratic map $f$ is not positive definite, then this map is called **pseudo-Euclidean metric on division ring** $D$. Symmetric associative bilinear map $g$ is called **pseudo-Euclidean scalar product on division ring** $D$. Let $\overline{\overline{e}}$ be the basis such that form $f$ is presented as sum of squares in this basis. There exist vectors of basis $\overline{e}_i$ such that $f(\overline{e}_i) < 0$. For arbitrary $a = a^i\ \overline{e}_i$ let us define operation of conjugation using the equation

$$a^\star = sign(f(\overline{e}_i))a^i\ \overline{e}_i$$

Expression $a^\star$ is called **hermitian conjugation**.

Let us consider biadditive map

$$g^\star(a, b) = g(a, b^\star)$$
$$f^\star(a) = g^\star(a, a)$$

The map $f^\star(a)$ is positive definite and is called **hermitian metric on division ring** $D$. Symmetric associative bilinear map $g^\star$ is called **hermitian scalar product on division ring** $D$.

If the map

$$a \to a^\star$$

is linear, then hermitian scalar product coinsides with Euclidean scalar product. Similar statement is true for metric.

CHAPTER 6

# $D$-Affine Space

## 6.1. $D$-Affine Space

**Definition 6.1.1.** Consider a set of points $\overset{\circ}{V}$ and a set of free vectors $\overline{V}$.[6.1] The set $\overset{\circ}{V}$ satisfies to following axioms.

(1) There exists at least one point
(2) One and only one vector is in correspondence to any tuple of points $(A, B)$. We denote this vector as $\overline{AB}$. The vector $\overline{AB}$ has tail in the point $A$ and head in the point B.
(3) For any point $A$ and any vector $\overline{a}$ there exists one and only one point $B$ such, that $\overline{AB} = \overline{a}$ . We will use notation[6.2]

$$B = A + \overline{a} \tag{6.1.1}$$

(4) (Axiom of parallelogram.) If $\overline{AB} = \overline{CD}$, then $\overline{AC} = \overline{BD}$.

□

**Definition 6.1.2.** Let $\overline{a}$ and $\overline{b}$ be vectors. Let $A \in \overset{\circ}{V}$ be arbitrary point. Let $B \in \overset{\circ}{V}$, $B = A + \overline{a}$. Let $C \in \overset{\circ}{V}$, $C = B + \overline{b}$. Vector $\overline{AC}$ is called sum of vectors $\overline{a}$ and $\overline{b}$

$$\overline{a} + \overline{b} = \overline{AC} \tag{6.1.2}$$

□

**Theorem 6.1.3.** *Vector $\overline{AA}$ is zero with respect to addition and does not depend on point $A$. Vector $\overline{AA}$ is called zero-vector and we assume $\overline{AA} = \overline{0}$.*

Proof. We can write rule of addition (6.1.2) in form of equation

$$\overline{AB} + \overline{BC} = \overline{AC} \tag{6.1.3}$$

If $B = C$, then from equation (6.1.3) it follows that

$$\overline{AB} + \overline{BB} = \overline{AB} \tag{6.1.4}$$

If $C = A$, $B = D$, then from axiom (4) of definition 6.1.1 it follows that

$$\overline{AA} = \overline{BB} \tag{6.1.5}$$

The statement of theorem follows from equations (6.1.4) and (6.1.5). □

[6.1] I wrote definitions and theorems in this section according to definition of affine space in [3], pp. 86 - 93.
[6.2] [13], p. 9.

**Theorem 6.1.4.** *Let* $\overline{a}=\overline{A\,B}$ *. Then*

$$\overline{B\,A}=-\overline{a} \tag{6.1.6}$$

*and this equation does not depend on a point* $A$*.*

Proof. From equation (6.1.3) it follows that

$$\overline{A\,B}+\overline{B\,A}=\overline{A\,A}=\overline{0} \tag{6.1.7}$$

Equation (6.1.6) follows from equation (6.1.7). Applying axiom (4) from definition 6.1.1 to equation $\overline{A\,B}=\overline{C\,D}$ we get $\overline{A\,C}=\overline{B\,D}$, or (this is equivalent) $\overline{B\,D}=\overline{A\,C}$. Based on axiom (4) of definition 6.1.1 again $\overline{B\,A}=\overline{D\,C}$ follows. Therefore, equation (6.1.6) does not depend on point $A$. □

**Theorem 6.1.5.** *Sum of vectors* $\overline{a}$ *and* $\overline{b}$ *does not depend on point* $A$*.*

Proof. Let $\overline{a}=\overline{A\,B}=\overline{A'\,B'}$. Let $\overline{b}=\overline{B\,C}=\overline{B'\,C'}$. Let

$$\overline{A\,B}+\overline{B\,C}=\overline{A\,C}$$

$$\overline{A'\,B'}+\overline{B'\,C'}=\overline{A'\,C'}$$

According to axiom (4) of definition 6.1.1

$$\overline{A'\,A}=\overline{B'\,B}=\overline{C'\,C} \tag{6.1.8}$$

Applying axiom (4) of definition 6.1.1 to outermost members of equation (6.1.8), we get

$$\overline{A'\,C'}=\overline{A\,C} \tag{6.1.9}$$

From equation (6.1.9) the statement of theorem follows. □

**Theorem 6.1.6.** *Sum of vectors is associative.*

Proof. Let $\overline{a}=\overline{A\,B}$ , $\overline{b}=\overline{B\,C}$ , $\overline{c}=\overline{C\,D}$ . From equation

$$\begin{array}{ccccc}\overline{a} & + & \overline{b} & = & \overline{A\,C}\\ \overline{A\,B} & + & \overline{B\,C} & = & \overline{A\,C}\end{array}$$

it follows that

$$\begin{array}{ccccc}(\overline{a}+\overline{b}) & + & \overline{c} & = & \overline{A\,D}\\ \overline{A\,C} & + & \overline{C\,D} & = & \overline{A\,D}\end{array} \tag{6.1.10}$$

From equation

$$\begin{array}{ccccc}\overline{b} & + & \overline{c} & = & \overline{B\,D}\\ \overline{B\,C} & + & \overline{C\,D} & = & \overline{B\,D}\end{array}$$

it follows that

$$\begin{array}{ccccc}\overline{a} & + & (\overline{b}+\overline{c}) & = & \overline{A\,D}\\ \overline{A\,B} & + & \overline{B\,D} & = & \overline{A\,D}\end{array} \tag{6.1.11}$$

Associativity of sum follows from comparison of equations (6.1.10) and (6.1.11). □

**Theorem 6.1.7.** *The structure of Abelian group is defined on the set* $\overline{V}$*.*

Proof. From theorems 6.1.3, 6.1.4, 6.1.5, 6.1.6 it follows that sum of vectors determines group.

Let $\overline{a} = \overline{AB}$ , $\overline{b} = \overline{BC}$ .

$$\begin{array}{ccccc} \overline{a} & + & \overline{b} & = \overline{AC} \\ \overline{AB} & + & \overline{BC} & = \overline{AC} \end{array} \tag{6.1.12}$$

According to axiom (3) of definition 6.1.1 there exists the point $D$ such that

$$\overline{b} = \overline{AD} = \overline{BC}$$

According to axiom (4) of definition 6.1.1

$$\overline{AB} = \overline{DC} = \overline{a}$$

According to definition of sum of vectors

$$\begin{array}{ccccc} \overline{AD} & + & \overline{DC} & = \overline{AC} \\ \overline{b} & + & \overline{a} & = \overline{AC} \end{array} \tag{6.1.13}$$

Commutativity of sum follows from comparison of equations (6.1.12) and (6.1.13). □

**Theorem 6.1.8.** *The map*

$$\overline{V} \to \overset{\circ}{V}{}^{\star} \tag{6.1.14}$$

*defined by equation* (6.1.1)*, is a single transitive representation of Abelian group* $\overline{V}$*.*

Proof. Axiom (3) of definition 6.1.1 determines the map (6.1.14). From theorem 6.1.5, it follows that the map (6.1.14) is a representation. Efficiency of the representation follows from theorem 6.1.3 and axiom (2) of definition 6.1.1. From the axiom (2) of definition 6.1.1, it also follows that representation is transitive. Effective and transitive representation is single transitive. □

## 6.2. Affine Structure on the Set

**Definition 6.2.1.** Covariant single transitive $\star T$-representation of the Abelian group $G$ in the set $M$ is called **affine structure on the set** $M$. Elements of the set $M$ are called points. We denote the image of the point $A \in M$ under transformation generated by $a \in G$ as $A + a$. □

**Theorem 6.2.2.** *Affine structure on the set* $M$ *satisfies to axioms of definition 6.1.1.*

Proof. We assume that the set $M$ is not empty; therefore affine structure satisfies the axiom (1). Since $a \in G$ generates the transformation of the set, then, for any $A \in M$, $B \in M$ is defined uniquely such that $B = A + a$. This statement proves the axiom (3). Since the representation is single transitive, then for any $A$, $B \in M$ there exists unique $a \in G$ such that $B = A + a$. This statement allows us to introduce notation $\overline{AB} = a$ , as well this statement proves the axiom (2). So, to prove the theorem we need to prove that the axiom (4) follows from the statement of the theorem.

Let

$$\overline{AB} = \overline{CD} = a \tag{6.2.1}$$

From the equation (6.2.1), it follows that

$$B = A + \overline{A\,B} = A + a \tag{6.2.2}$$

$$D = C + \overline{C\,D} = C + a \tag{6.2.3}$$

Let

$$\overline{A\,C} = b \qquad C = A + b \tag{6.2.4}$$

$$\overline{B\,D} = c \qquad D = B + c \tag{6.2.5}$$

From equations (6.2.4), (6.2.3), it follows that

$$D = A + b + a \tag{6.2.6}$$

From equations (6.2.2), (6.2.5), it follows that

$$D = A + a + c \tag{6.2.7}$$

Since the group $G$ is commutative, then, from the equation (6.2.6), it follows that

$$D = A + a + b \tag{6.2.8}$$

From equations (6.2.7), (6.2.8), it follows that

$$A + a + b = A + a + c \tag{6.2.9}$$

Since the representation of group $G$ is single transitive, then, from the equation (6.2.9), it follows that

$$b = c \tag{6.2.10}$$

From equations (6.2.2), (6.2.3), (6.2.10), it follows that

$$\overline{A\,C} = \overline{B\,D} \tag{6.2.11}$$

Therefore, we have proved the axiom (4). □

Despite the fact that the theorem 6.2.2 states that the affine structure on the set $M$ satisfies to axioms of the definition 6.1.1, the definition 6.2.1 is not equivalent to the definition 6.1.1. In the definition 6.1.1, we assume that there is a structure of vector space on the set $G$. Therefore, if we want to explore affine space, we need to consider tower of representations: a representation of field $F$ in Abelian group $G$ and a representation of Abelian group $G$ in the set $M$. Evidently, we can consider division algebra $D$ instead of the field $F$. Thus, we get following definition of affine space.

**Definition 6.2.3.** Let $D$ be division ring of characteristic 0, $\overline{V}$ be Abelian group, and $\overset{\circ}{V}$ be any set. **$D\star$-affine space** is the tower of representations

$$\overrightarrow{V} : D \xrightarrow[*]{f_{1,2}} \overline{V} \xrightarrow[*]{f_{2,3}} \overset{\circ}{V} \qquad \begin{array}{l} f_{1,2}(d) : \overline{v} \rightarrow d\,\overline{v} \\ f_{2,3}(\overline{v}) : A \rightarrow A + \overline{v} \end{array}$$

where $f_{1,2}$ is effective $T\star$-representation of division ring $D$ in Abelian group $\overline{V}$ and $f_{2,3}$ is single transitive $\star T$-representation of Abelian group $\overline{V}$ in the set $\overset{\circ}{V}$. □

**Definition 6.2.4.** Let $D$ be division ring of characteristic 0, $\overline{V}$ be Abelian group, and $\overset{\circ}{V}$ be any set. **$\star D$-affine space** is the tower of representations

$$\overrightarrow{V}: D \xrightarrow[*]{f_{1,2}} \overline{V} \xrightarrow[*]{f_{2,3}} \overset{\circ}{V} \qquad \begin{array}{l} f_{1,2}(d): \overline{v} \rightarrow \overline{v}\ d \\ f_{2,3}(\overline{v}): A \rightarrow A + \overline{v} \end{array}$$

where $f_{1,2}$ is effective $\star T$-representation of division ring $D$ in Abelian group $\overline{V}$ and $f_{2,3}$ is single transitive $\star T$-representation of Abelian group $\overline{V}$ in the set $\overset{\circ}{V}$. □

## 6.3. Single Transitive Representations

Consider single transitive representation $f$ of $\Omega_1$-algebra $A$ in $\Omega_2$-algebra $M$. We can identify $a \in A$ with the corresponding transformation $f(a)$. In this section, we denote the action of $a \in A$ on $v \in M$ by expression

$$m + a$$

**Theorem 6.3.1.** *Since representation $f$ of $\Omega_1$-algebra $A$ in $\Omega_2$-algebra $M$ is single transitive, then the basis of the representation $f$ is arbitrary point $v_0 \in M$. The cardinality of the set $M$ equals to cardinality of the set $A$.*

Proof. We choose an arbitrary $v_0 \in M$. For any $v \in M$ there exists a unique

(6.3.1) $$a = g(v_0, v)$$

such that

(6.3.2) $$v_0 + a = v$$

Expression (6.3.2) is $\Omega_2$-word. Therefore, $v_0$ is the basis of the representation $f$. The equation (6.3.2) establishes a bijection between the sets $A$ and $M$. Therefore, the cardinality of the set $A$ equals to cardinality of the set $M$. □

We will use notation

(6.3.3) $$a = \overline{v_0\, v}$$

for the mapping $g$ and write down the equation (6.3.2) in the following form

(6.3.4) $$v_0 + \overline{v_0\, v} = v$$

**Theorem 6.3.2.** *Let representation $f$ of $\Omega_1$-algebra $A$ in $\Omega_2$-algebra $M$ be single transitive. If we chose a basis $v_0$ of the representation $f$, then correspondence (6.3.3), (6.3.4) allows us to transfer the structure of $\Omega_1$-algebra into $\Omega_2$-algebra $M$.*[6.3]

Proof. Consider operation $\omega \in \Omega_1(n)$. For arbitrary $v_1$, ..., $v_n \in M$, we assume

(6.3.5) $$v_1 ... v_n\omega = v_0 + \overline{v_0\, v_1} ... \overline{v_0\, v_n}\omega$$

□

[6.3] It is important to keep in mind that the structure of $\Omega_1$-algebra on the set $M$ depends on a choice of point $v_0$. This remark is also valid for all following theorems considered below in this section. To emphasize that considered $\Omega_1$ algebra depends on the choice of point $O$, we shall use the notation $(M, v_0)$ .

**Theorem 6.3.3.** *Let representation $f$ of $\Omega_1$-algebra $A$ in $\Omega_2$-algebra $M$ be single transitive. If we chose a basis $v_0$ of the representation $f$, then endomorphism*

$$f' : v \in M \to v_0 + f(\overline{v_0\, v}) \in M \tag{6.3.6}$$

*of $\Omega_1$-algebra $(M, v_0)$ corresponds to arbitrary endomorphism $f$ of $\Omega_1$-algebra $A$.*

Proof. From the definition of the mapping $f'$ (6.3.6) and equation (6.3.4), it immediately follows that

$$f'(v_0 + \overline{v_0\, v}) = v_0 + f(\overline{v_0\, v}) \tag{6.3.7}$$

$$v_0 + \overline{v_0\, f'(v)} = v_0 + f(\overline{v_0\, v}) \tag{6.3.8}$$

From the equation (6.3.8) it follows that

$$\overline{v_0\, f'(v)} = f(\overline{v_0\, v}) \tag{6.3.9}$$

Let $\omega \in \Omega_1(n)$. From the equation (6.3.5), it follows that

$$f'(v_1...v_n\omega) = f'(v_0 + \overline{v_0\, v_1}...\overline{v_0\, v_n}\omega) \tag{6.3.10}$$

From equations (6.3.10), (6.3.7), it follows that

$$f'(v_1...v_n\omega) = v_0 + f(\overline{v_0\, v_1}...\overline{v_0\, v_n}\omega) \tag{6.3.11}$$

Since the mapping $f$ is endomorphism, then from equation (6.3.11) it follows that

$$f'(v_1...v_n\omega) = v_0 + f(\overline{v_0\, v_1})...f(\overline{v_0\, v_n})\omega \tag{6.3.12}$$

From equations (6.3.12), (6.3.9) it follows that

$$f'(v_1...v_n\omega) = v_0 + \overline{v_0\, f'(v_1)}...\overline{v_0\, f'(v_n)}\omega \tag{6.3.13}$$

From equations (6.3.13), (6.3.5) it follows that

$$f'(v_1...v_n\omega) = f'(v_1)...f'(v_n)\omega \tag{6.3.14}$$

From the equation (6.3.14), it follows that the mapping $f'$ is endomorphism of $\Omega_1$-algebra $(M, v_0)$. □

**Theorem 6.3.4.** *Let representation $f$ of $\Omega_1$-algebra $A$ in $\Omega_2$-algebra $M$ be single transitive. If we chose a basis $v_0$ of the representation $f$, then correspondence (6.3.3), (6.3.4) allows us to transfer the structure of $\Omega_2$-algebra into $\Omega_1$-algebra $A$.*

Proof. Consider an operation $\omega \in \Omega_2(n)$. For arbitrary $a_1, ..., a_n \in A$, there exist $v_1, ..., v_n \in M$ such that

$$a_1 = \overline{v_0\, v_1} \quad ... \quad a_n = \overline{v_0\, v_n}$$

We assume

$$a_1...a_n\omega = \overline{v_0\, v_1}...\overline{v_0\, v_n}\omega = \overline{v_0\, v_1...v_n\omega} \tag{6.3.15}$$

□

**Theorem 6.3.5.** *Consider tower $(\overline{V}, \overline{f})$ of representations of $\overline{\Omega}$-algebras. Let representation $f_{i+1,i+2}$ be single transitive, and $a_{i+2\cdot 0}$ be the basis of the representation $f_{i+1,i+2}$. Then the representation*

$$\begin{aligned} f'_{i,i+2} &: A_i \xrightarrow{\ *\ } A_{i+2} \\ f'_{i,i+2}(a_i) &: a_{i+2} \to a_{i+2\cdot 0} + f_{i,i+1}(a_i)(\overline{a_{i+2\cdot 0}\, a_{i+2}}) \end{aligned} \tag{6.3.16}$$

*of $\Omega_i$-algebra $A_i$ in $\Omega_{i+2}$-algebra $A_{i+2}$ corresponds to the representation $f_{i,i+1}$ of $\Omega_i$-algebra $A_i$ in $\Omega_{i+1}$-algebra $A_{i+1}$.*

Proof. Since the mapping $f_{i,i+1}(a_i)$ is endomorphism of $\Omega_{i+1}$-algebra $A_{i+1}$, then, according to the theorem 6.3.3, the mapping $f'_{i,i+2}(a_i)$ is endomorphism of $\Omega_{i+1}$-algebra $(A_{i+2}, a_{i+2\cdot 0})$. □

From the equation (6.3.16), it follows that

$$f'_{i,i+2}(a_i)(a_{i+2\cdot 0}+\overline{a_{i+2\cdot 0}\,a_{i+2}}) = a_{i+2\cdot 0}+f_{i,i+1}(a_i)(\overline{a_{i+2\cdot 0}\,a_{i+2}}) \tag{6.3.17}$$

$$a_{i+2\cdot 0}+\overline{a_{i+2\cdot 0}\,f'_{i,i+2}(a_i)(a_{i+2})} = a_{i+2\cdot 0}+f_{i,i+1}(a_i)(\overline{a_{i+2\cdot 0}\,a_{i+2}}) \tag{6.3.18}$$

From the equation (6.3.18) it follows that

$$\overline{a_{i+2\cdot 0}\,f'_{i,i+2}(a_i)(a_{i+2})} = f_{i,i+1}(a_i)(\overline{a_{i+2\cdot 0}\,a_{i+2}}) \tag{6.3.19}$$

**Theorem 6.3.6.** *Consider tower $(\overline{V}, \overline{f})$ of representations of $\overline{\Omega}$-algebras. Let representation $f_{i+1,i+2}$ be single transitive, and $a_{i+2\cdot 0}$ be the basis of the representation $f_{i+1,i+2}$. Endomorphism $r_{i+1}$ of the representation $f_{i,i+1}$ generates endomorphism $r'_{i+2}$ of the representation $f'_{i,i+2}$*

$$r'_{i+2}(a_{i+2}) = a_{i+2\cdot 0}+r_{i+1}(\overline{a_{i+2\cdot 0}\,a_{i+2}}) \tag{6.3.20}$$

**Remark 6.3.7.** From the equation (6.3.20) it follows that

$$\overline{a_{i+2\cdot 0}\,r'_{i+2}(a_{i+2})} = r_{i+1}(\overline{a_{i+2\cdot 0}\,a_{i+2}}) \tag{6.3.21}$$

Proof. Since the mapping $r_{i+1}$ is endomorphism of $\Omega_{i+1}$-algebra $A_{i+1}$, then, according to the theorem 6.3.3, the mapping $r'_{i+2}$ is endomorphism of $\Omega_{i+1}$-algebra $(A_{i+2}, a_{i+2\cdot 0})$.

To prove the theorem, it is sufficient to show that following equation is valid

$$r'_{i+2}\circ f'_{i,i+2}(a_i) = f'_{i,i+2}(r_i(a_i))\circ r'_{i+2} \tag{6.3.22}$$

We can write the equation (6.3.22) in the following form

$$r'_{i+2}(f'_{i,i+2}(a_i)(a_{i+2})) = f'_{i,i+2}(r_i(a_i))(r'_{i+2}(a_{i+2})) \tag{6.3.23}$$

From equations (6.3.23), (6.3.16), (6.3.20), it follows that

$$\begin{aligned} & a_{i+2\cdot 0}+r_{i+1}(\overline{a_{i+2\cdot 0}\,f'_{i,i+2}(a_i)(a_{i+2})}) \\ = & a_{i+2\cdot 0}+f_{i,i+1}(r_i(a_i))(\overline{a_{i+2\cdot 0}\,r'_{i+2}(a_{i+2})}) \end{aligned} \tag{6.3.24}$$

From equations (6.3.24), (6.3.19), (6.3.21), it follows that

$$r_{i+1}(f_{i,i+1}(a_i)(\overline{a_{i+2\cdot 0}\,a_{i+2}})) = f_{i,i+1}(r_i(a_i))(r_{i+1}(\overline{a_{i+2\cdot 0}\,a_{i+2}})) \tag{6.3.25}$$

The equation (6.3.25) is true because the mapping $(r_i, r_{i+1})$ is endomorphism of the representation $f_{i,i+1}$. Therefore, the mapping $(r_i, r'_{i+2})$ is endomorphism of the representation $f'_{i,i+2}$. □

## 6.4. Basis in ${}_*{}^*D$-Affine Space

Consider $\star D$-affine space $\overrightarrow{V}$.

The Abelian group $\overline{V}$ acts single transitive on the set $\overset{\circ}{V}$. According to the theorem 6.3.1, the basis of the set $\overset{\circ}{V}$ relative to representation of the Abelian group $\overline{V}$ consists of one point. This point is usually denoted by the letter $O$ and is called **origin of coordinate system of $\star D$-affine space**. Therefore, an arbitrary point $A\in\overset{\circ}{V}$ can be represented using vector $\overline{O\,A}\in\overline{V}$

$$A = O+\overline{O\,A} \tag{6.4.1}$$

According to the theorem 6.3.1, the cardinality of the set $\overset{\circ}{V}$ equal to the cardinality of set $\overline{V}$.

The type of vector space in the tower of representations (section [6]-4.2) determines the type of affine space.

Let $\overline{V}$ be ${}_*{}^*D$-vector space over division ring $D$. $\star D$-affine space $\overrightarrow{V}$ is called **${}_*{}^*D$-affine space**. Let $\overline{\overline{e}}$ be the basis of the ${}_*{}^*D$-vector space $\overline{V}$. Then expansion of the vector $\overrightarrow{OA}$ relative to the ${}_*{}^*D$-basis $\overline{\overline{e}}$ has form

$$\overline{OA} = \overline{e}_i\, {}^iA \tag{6.4.2}$$

From equations (6.4.1), (6.4.2), it follows that

$$A = O + \overline{e}_i\, {}^iA \tag{6.4.3}$$

**Definition 6.4.1.** **${}_*{}^*D$-affine basis** $(\overline{\overline{e}}, O) = (\overline{e}_i, i \in I, O)$ is maximal set of ${}_*{}^*D$-linear independent vectors $\overline{e}_i = \overline{OA_i} = ({}^1e_i, ..., {}^ne_i)$ with common start point $O$. The set of coordinates $({}^iA, i \in I)$ of vector $\overline{OA}$ is called **coordinates of point $A$ of ${}_*{}^*D$-affine space $\overrightarrow{V}$ relative to basis** $(\overline{\overline{e}}, O)$. □

**Definition 6.4.2.** **Basis manifold $\mathcal{B}(\overrightarrow{V})$ of ${}_*{}^*D$-affine space** is set of ${}_*{}^*D$-bases of this space. □

**Theorem 6.4.3.** *Let*

$$\overrightarrow{V} : D \xrightarrow[*]{f_{1,2}} \overline{V} \xrightarrow[*]{f_{2,3}} \overset{\circ}{V} \qquad \begin{array}{l} f_{1,2}(d) : \overline{v} \rightarrow \overline{v}\ d \\ f_{2,3}(\overline{v}) : A \rightarrow A + \overline{v} \end{array}$$

*be ${}_*{}^*D$-affine space.*

(1) *Endomorphism of the representation $f_{2,3}$ has form*

$$A' = A + \overline{a} \tag{6.4.4}$$

(2) *Endomorphism of the representation $f_{1,2}$ has form*

$$\overline{a}' = \overline{a}_*{}^*P \tag{6.4.5}$$

Proof. Let $h$ be endomorphism of the representation $f_{2,3}$. From the definition [10]-2.2.2, it follows that

$$h \circ (A + \overline{b}) = h \circ A + \overline{b} \tag{6.4.6}$$

Let

$$B = A + \overline{b} \quad C = h \circ A = A + \overline{c} \quad D = h \circ B = B + \overline{d} \tag{6.4.7}$$

From equations (6.4.6), (6.4.7), it follows that

$$D = h \circ B = h \circ (A + \overline{b}) = h \circ A + \overline{b} = C + \overline{b} \tag{6.4.8}$$

From equations (6.4.7), (6.4.8), it follows that

$$\overline{AB} = \overline{CD} \tag{6.4.9}$$

From the axiom (4) of the definition 6.1.1 and the equation (6.4.9) it follows that

$$\overline{c} = \overline{AC} = \overline{BD} = \overline{d} \tag{6.4.10}$$

Since $\overline{a}$ is an arbitrary vector, then, from equations (6.4.10), (6.4.7), it follows that

$$D = B + \overline{c}$$

for any point $B\in\overset{\circ}{V}$. Therefore, any endomorphism of representation $f_{2,3}$ has form (6.4.4).

The statement (2) of the theorem was considered in detail in the section [6]-4.4. □

**Theorem 6.4.4.** *Transformation* (6.4.4) *is called* **parallel shift of ${}_*{}^*D$-affine space**. *Coordinates of parallel shift* (6.4.4) *of ${}_*{}^*D$-affine space have form*

$$ {}^iA'={}^iA+{}^ia \tag{6.4.11}$$

Proof. Let $({}^iA,i\in I)$ be coordinates of point $A$ relative to the basis $(\overline{\overline{e}},O)$. From the equation (6.4.3), it follows that

$$ A=O+\overline{e}_i\,{}^iA \tag{6.4.12}$$

Let $({}^iA',i\in I)$ be coordinates of point $A'$ relative to the basis $(\overline{\overline{e}},O)$. From the equation (6.4.3), it follows that

$$ A'=O+\overline{e}_i\,{}^iA' \tag{6.4.13}$$

Let

$$ \overline{a}=\overline{e}_i\,{}^ia \tag{6.4.14}$$

be expansion of vector $\overline{a}$ relative to the basis $\overline{\overline{e}}$. Substituting equation (6.4.12), (6.4.13), (6.4.14) into equation (6.4.4), we obtain

$$ O+\overline{e}_i\,{}^iA'=O+\overline{e}_i\,{}^iA+\overline{e}_i\,{}^ia=O+\overline{e}_i({}^iA+{}^ia) \tag{6.4.15}$$

Since the representation $f_{2,3}$ is single transitive, then, from the equation (6.4.15), it follows that

$$ \overline{e}_i\,{}^iA'=\overline{e}_i({}^iA+{}^ia) \tag{6.4.16}$$

Since vectors $\overline{e}_i$ are ${}_*{}^*D$-linear independent, then the equation (6.4.11) follows from the equation (6.4.16). □

**Theorem 6.4.5.** *Let $(\overline{\overline{e}},O)$ be the basis of ${}_*{}^*D$-affine space*

$$\overrightarrow{V}:D\xrightarrow[*]{f_{1,2}}\overline{V}\xrightarrow[*]{f_{2,3}}\overset{\circ}{V}\qquad\begin{array}{l}f_{1,2}(d):\overline{v}\to\overline{v}\,d\\ f_{2,3}(\overline{v}):A\to A+\overline{v}\end{array}$$

(1) *The representation $f_{2,3}$ generates a sum on the set $\overset{\circ}{V}$*

$$ A+B=O+\overline{OA}+\overline{OB} \tag{6.4.17}$$

(2) *The representation $f_{1,2}$ generates a representation*

$$ \begin{array}{c}f'_{1,3}:D\xrightarrow{*}\overset{\circ}{V}\\ f'_{1,3}(d):A\to O+\overline{OA}\,d\end{array} \tag{6.4.18}$$

(3) *Endomorphism $\overline{g}$ of the representation $f_{1,2}$ generates endomorphism $\overline{g}'$ of the representation $f'_{1,2}$*

$$ \overline{g}':A\in\overset{\circ}{V}\to O+g_*{}^*\overline{OA} \tag{6.4.19}$$

Proof. Since the representation $f_{2,3}$ is single transitive, then

- the statement (1) of the theorem is corollary of the theorem 6.3.2.
- the statement (2) of the theorem is corollary of the theorem 6.3.5.

- the statement (3) of the theorem is corollary of the theorem 6.3.6.

□

**Theorem 6.4.6.** *Transformation* (6.4.19) *is called* **linear transformation of** ${}_*{}^*D$**-affine space**. *Coordinates of linear transformation* (6.4.19) *of* ${}_*{}^*D$*-affine space have form*

$$ {}^iA' = {}^ig_j\ {}^jA \tag{6.4.20}$$

Proof. Let $( {}^iA, i \in I)$ be coordinates of point $A$ relative to the basis $(\overline{\overline{e}}, O)$. From the equation (6.4.2), it follows that

$$\overline{O\,A} = \overline{e}_i\ {}^iA \tag{6.4.21}$$

Let $( {}^iA', i \in I)$ be coordinates of point $A'$ relative to the basis $(\overline{\overline{e}}, O)$. From the equation (6.4.3), it follows that

$$A' = O + \overline{e}_i\ {}^iA' \tag{6.4.22}$$

From the equation (6.4.19), it follows that

$$A' = O + g_*{}^*\overline{O\,A} \tag{6.4.23}$$

From equations (6.4.21), (6.4.22), (6.4.23), it follows that

$$O + \overline{e}_i\ {}^iA' = O + \overline{e}_i\ {}^ig_j\ {}^jA \tag{6.4.24}$$

The equation (6.4.20) follows from the equation (6.4.24), because $(\overline{e}, O)$ is the basis of ${}_*{}^*D$-affine space. □

Composition of parallel shift and linear transformation of ${}_*{}^*D$-affine space, generates ${}_*{}^*D$**-affine transformation**

$$A' = O + a_*{}^*\overline{O\,A} + \overline{v} \tag{6.4.25}$$

We can represent this transformation of ${}_*{}^*D$-affine space as the mapping of basis $(\overline{\overline{e}}, O)$ into basis $(\overline{\overline{e}}', O')$. We can represent the mapping of basis $(\overline{\overline{e}}, O)$ into basis $(\overline{\overline{e}}', O')$ as product of mapping of basis $(\overline{\overline{e}}, O)$ into basis $(\overline{\overline{e}}, O')$ and $D_*{}^*$-linear mapping of basis $(\overline{\overline{e}}, O')$ into basis $(\overline{\overline{e}}', O')$. Therefore, $D_*{}^*$-transformation has form

$$\overline{\overline{e}}' = A_*{}^*\overline{\overline{e}}$$

$$O' = O + \overrightarrow{OO'}$$

Introducing coordinates $A^1, ..., A^n$ of a point $A \in \mathcal{A}_n$ as coordinates of vector $\overline{O\,A}$ relative to basis $\overline{\overline{e}}$ we can write a linear transformation as

$$A'^i = A^j\ {}_jP^i + R^i \qquad \operatorname{rank}_{*}{}^{*} P = n \tag{6.4.26}$$

$$A'{}_*{}^*\overline{e} = (A_*{}^*P + R)_*{}^*\overline{e} \tag{6.4.27}$$

Vector $(R^1, ..., R^n)$ expresses displacement in affine space.

**Theorem 6.4.7.** *Set of transformations* (6.4.27) *is the group Lie which we denote as* $GL(\mathcal{A}_n)$ *and call* **affine transformation group**.

Proof. Let us consider transformations $(P,R)$ and $(Q,S)$. Product of these transformations has the form

$$\begin{aligned} A''^i &= A'^j\ {}_jQ^i + S^i & A'' &= (A'{}_*{}^*Q+S){}_*{}^*\overline{e} \\ &= (A^k\ {}_kP^j + R^j)\ {}_jQ^i + S^i & &= ((A{}_*{}^*P+R){}_*{}^*Q+S){}_*{}^*\overline{e} \\ &= A^k\ {}_kP^j\ {}_jQ^i + R^j\ {}_jQ^i + S^i & &= (A{}_*{}^*P{}_*{}^*Q + R{}_*{}^*Q+S){}_*{}^*\overline{e} \end{aligned}$$

Therefore, we can write product of transformations $(P,R)$ and $(Q,S)$ in form

$$(P,R){}_*{}^*(Q,S) = (P{}_*{}^*Q, R{}_*{}^*Q+S) \tag{6.4.28}$$

□

An active transformation is called **affine transformation**.A passive transformation is called **quasi affine transformation**.

If we do not concern about starting point of a vector we see little different type of space which we call central affine space $CA_n$. If we assume that the start point of vector is origin $O$ of coordinate system in space then we can identify any point $A\in CA_n$ with the vector $a=\overline{OA}$. Now transformation is simply map

$$a'^i = P^i_j a^j \qquad\qquad \det P\neq 0$$

and such transformations build up Lie group $GL_n$.

**Definition 6.4.8. Central affine basis** $\overline{\overline{e}} = (\overline{e}_i)$ is set of linearly independent vectors $\overline{e}_i$. □

**Definition 6.4.9. Basis manifold** $\mathcal{B}(CA_n)$ **of central affine space** is set of bases of this space. □

## 6.5. Plane in ${}_*{}^*D$-affine space

Consider ${}_*{}^*D$-affine space

$$\overrightarrow{V} : D \xrightarrow[*]{f_{1,2}} \overline{V} \xrightarrow[*]{f_{2,3}} \overset{\circ}{V} \qquad \begin{aligned} &f_{1,2}(d) : \overline{v} \to \overline{v}\ d \\ &f_{2,3}(\overline{v}) : A \to A+\overline{v} \end{aligned}$$

Let $\overline{V}_1 \subset \overline{V}$ be ${}_*{}^*D$-vector subspace of ${}_*{}^*D$-vector space $\overline{V}$. Therefore, $\overline{V}_1$ is subgroup of Abelian group $\overline{V}$. In the tower of representations

$$D \xrightarrow[*]{f_{1,2}} \overline{V}_1 \xrightarrow[*]{f_{2,3}} \overset{\circ}{V} \qquad \begin{aligned} &f_{1,2}(d) : \overline{v} \to \overline{v}\ d \\ &f_{2,3}(\overline{v}) : A \to A+\overline{v} \end{aligned}$$

the representation $f_{2,3}$ is not single transitive. According to theorem [6]-3.4.13 , we can represent the set $\overset{\circ}{V}$ as intersection of disjoint orbits of the representation $f_{2,3}$. However since the representation $f_{2,3}$ is effective, then the representation $f_{2,3}$ is single transitive on the orbit of the representation. Therefore, if we choose an arbitrary point $A\in\overset{\circ}{V}$ and orbit of the representation

$$\overset{\circ}{V}_1 = A + f_{2,3}(\overline{V}_1)$$

then we obtain the tower of representations

$$\overrightarrow{V}_1 : D \xrightarrow[*]{f_{1,2}} \overline{V}_1 \xrightarrow[*]{f_{2,3}} \overset{\circ}{V}_1 \qquad \begin{aligned} &f_{1,2}(d) : \overline{v} \to \overline{v}\ d \\ &f_{2,3}(\overline{v}) : A \to A+\overline{v} \end{aligned}$$

describing ${}_*{}^*D$-affine space called ${}_*{}^*D$**-affine plane**, passing through the point $A$.

Let

$$\overline{\overline{e}}=(\overline{e}_1,...,\overline{e}_n)$$

be the basis of ${}_*{}^*D$-vector space $\overline{V}$. Let the basis of ${}_*{}^*D$-vector space $\overline{V}_1$

$$\overline{\overline{e}}_1=(\overline{e}_{1\cdot 1},...,\overline{e}_{1\cdot k})$$

have following expansion relative to the basis $\overline{\overline{e}}$

$$\overline{e}_{1\cdot i}=\overline{e}_j\ {}^j e_{1\cdot i}$$

Let $\overline{v}\in\overline{V}_1$ have following expansion relative to the basis $\overline{\overline{e}}$

$$\overline{v}=\overline{e}_j\ {}^j v$$

According to the theorem [6]-4.6.7

$$\text{rank}\begin{pmatrix}{}^1e_{1\cdot 1}&...&{}^1e_{1\cdot k}&{}^1v\\...&...&...&...\\{}^ne_{1\cdot 1}&...&{}^ne_{1\cdot k}&{}^nv\end{pmatrix}=k\tag{6.5.1}$$

If we assume that first $k$ rows compose nonsingular matrix, then we obtain $n-k$ ${}_*{}^*D$-linear equations

$${}^{k+1}\det({}_*{}^*)_r\begin{pmatrix}{}^1e_{1\cdot 1}&...&{}^1e_{1\cdot k}&{}^1v\\...&...&...&...\\{}^ne_{1\cdot 1}&...&{}^ne_{1\cdot k}&{}^nv\end{pmatrix}=0\quad r=k+1,...,n\tag{6.5.2}$$

satisfied by the coordinates of vector $\overline{v}$.

Let point $A$ have representation

$$A=O+e_i\ {}^iA\tag{6.5.3}$$

relative to the basis $(\overline{\overline{e}},O)$. Let

$$B=O+e_i\ {}^iB\tag{6.5.4}$$

be an arbitrary point of ${}_*{}^*D$-affine plane $\overrightarrow{V}_1$. Then

$$\overline{v}=\overline{A\,B}=e_i\ {}^iv\quad \overline{v}\in\overline{V}_1\tag{6.5.5}$$

From the equation

$$B=A+\overline{v}$$

and equations (6.5.3), (6.5.4), (6.5.5), it follows that

$${}^iB=\ {}^iA+\ {}^iv\tag{6.5.6}$$

From equations (6.5.2), (6.5.6), it follows that coordinates of the point $B$ satisfy to the system of ${}_*{}^*D$-linear equations

$${}^{k+1}\det({}_*{}^*)_r\begin{pmatrix}{}^1e_{1\cdot 1}&...&{}^1e_{1\cdot k}&{}^1B-{}^1A\\...&...&...&...\\{}^ne_{1\cdot 1}&...&{}^ne_{1\cdot k}&{}^nB-{}^nA\end{pmatrix}=0\quad r=k+1,...,n\tag{6.5.7}$$

For given set of vectors

$$\overline{a}_i=\overline{A_0\,A_i}\quad i=1,...,n-1$$

we say that vector $\overline{a}=\overrightarrow{A_0A}$ belongs to the plane spanned by $\overline{a}_1,...,\overline{a}_{n-1}$ , if vector $\overline{a}_n$ is linear combination of vectors $\overline{a}_1,...,\overline{a}_{n-1}$ .

CHAPTER 7

# Euclidean Space

## 7.1. Euclidean Space

Let $\overline{V}$ be $D$-vector space. We can present bilinear map

$$g : \overline{V} \times \overline{V} \to D$$

as

$$g(\overline{v}, \overline{w}) = {}_{ij}g(v^i, w^j)$$

where ${}_{ij}g$ are bilinear maps

$${}_{ij}g : D^2 \to D$$

**Definition 7.1.1.** Bilinear map

$$g : \overline{V} \times \overline{V} \to D$$

is called **symmetric** when

$$g(\overline{w}, \overline{v}) = g(\overline{v}, \overline{w})$$

□

**Definition 7.1.2.** Bilinear map

$$g : \overline{V} \times \overline{V} \to D$$

is called **Euclidean scalar product** in $D$-vector space $\overline{V}$, if ${}_{ii}g$ is Euclidean scalar product in division ring $D$ and ${}_{ij}g = 0$ for $i \neq j$. □

**Definition 7.1.3.** Bilinear map

$$g : \overline{V} \times \overline{V} \to D$$

is called **pseudo-Euclidean scalar product** in $D$-vector space $\overline{V}$, if ${}_{ii}g$ is pseudo-Euclidean scalar product in division ring $D$ and ${}_{ij}g = 0$ for $i \neq j$. □

If there exists hermitian conjugation in the division ring $D$, then we can extend this operation to $D$-vector space $\overline{V}$, namely for given vector $\overline{v} = v^i\ {}_i\overline{e}$ we define **hermitian conjugated vector**

$$\overline{v}^\star = (v^i)^\star\ {}_i\overline{e}$$

**Definition 7.1.4.** Bilinear map

$$g : \overline{V} \times \overline{V} \to D$$

is called **hermitian scalar product** in $D$-vector space $\overline{V}$, if ${}_{ii}g$ is hermitian scalar product in division ring $D$ and ${}_{ij}g = 0$ for $i \neq j$. □

## 7.2. Basis in Euclidean Space

When we introduce a metric in a central affine space we get a new geometry because we can measure a distance and a length of vector. If a metric is positive defined we call the space Euclid $\mathcal{E}_n$ otherwise we call the space pseudo Euclid $\mathcal{E}_{nm}$.

Transformations that preserve length form Lie group $SO(n)$ for Euclid space and Lie group $SO(n,m)$ for pseudo Euclid space where $n$ and $m$ number of positive and negative terms in metrics.

**Definition 7.2.1. Orthonornal basis** $\overline{\overline{e}} = (\overline{e}_i)$ is set of linearly independent vectors $\overline{e}_i$ such that length of each vector is 1 and different vectors are orthogonal. □

**Definition 7.2.2. Basis manifold** $\mathcal{B}(\mathcal{E}_n)$ **of Euclid space** is set of orthonornal bases of this space. □

A active transformation is called **movement**. An passive transformation is called **quasi movement**.

CHAPTER 8

# Calculus in $D$-Affine Space

## 8.1. Curvilinear Coordinates in $D$-Affine Space

Let $\overrightarrow{A}$ be $D$-affine space over continuous division ring $D$.[8.1]

Let $({}_1\overline{e},...,{}_n\overline{e},O)$ be $D_*{}^*$-basis of affine space. For each $i$, $i=1$, ..., $i=n$, we define map[8.2]

$$\overline{x}_i : D \rightarrow \overset{\circ}{A} \quad v^i \rightarrow O + v^i\ {}_i\overline{e} \tag{8.1.1}$$

For arbitrary point $M(v^i)$ we can write vector $\overrightarrow{OM}$ in form of linear combination

$$\overrightarrow{OM} = \overline{x}_i(v^i)$$

We consider the set of coordinates of affine space $\overset{\circ}{A}$ as homeomorphism

$$f : \overset{\circ}{A} \rightarrow D^n$$

If homeomorphism

$$f : \overset{\circ}{A} \rightarrow D^n$$

is an arbitrary map of an open set of $D$-affine space inty an open subset of the set $D^n$, then we will say that image of point $A$ under map $f$ is **curvilinear coordinates of point** $A$, and map $f$ is called the system of coordinates. If there are two systems of coordinates, $f_1$ and $f_2$, then map $f_1 f_2^{-1}$ is homeomorphism

$$g : D^n \rightarrow D^n$$
$$x'^i = x'^i(x^1, ..., x^n)$$

and is called coordinate transformation. We can consider as particular case the coordinate transformation caused by affine transformation.

**Remark 8.1.1.** Let us consider a surface $S$ in affine space defined by equation

$$f : D^n \rightarrow D$$
$$f(x^1, ..., x^n) = 0$$

Let $M(x^i)$ be arbitrary point. A small change of coordinates of the point

$$x'^i = x^i + \Delta x^i$$

leads to a small change of function

$$\Delta f = \partial f(\overline{x})(\Delta\overline{x})$$

[8.1] In this section, I explore curvilinear coordinates the way like it was done in [3], part V.

[8.2] There is no sum over index $i$ in equation (8.1.1).

In particular, if neighboring points $M(x^i)$ and $L(x^i+\Delta x^i)$ belong to the same surface, then vector of increment $\overrightarrow{ML}=\Delta\overline{x}$ satisfies to differential equation

$$\partial f(\overline{x})(\Delta\overline{x})=0$$

and is called vector tangent to surface $S$.

Since the Gâteaux derivative of map $f$ in point $X$ is a linear map, then this derivative is 1-$D$-form. Let us consider the set of 1-$D$-forms

$$dx^i=\frac{\partial x^i}{\partial\overline{x}}$$

It is evident that

$$dx^i(\overline{a})=\frac{\partial x^i}{\partial\overline{x}}(\overline{a})=(x^i+a^i)-x^i=a^i$$

We can represent the Gâteaux derivative of arbitrary function in form

$$\frac{\partial f(\overline{x})}{\partial\overline{x}}=\frac{\partial f(\overline{x})}{\partial x^i}\left(\frac{\partial x^i}{\partial\overline{x}}\right)=\frac{\partial f(\overline{x})}{\partial x^i}\left(dx^i\right)$$

(8.1.2) $$\frac{\partial f(\overline{x})}{\partial\overline{x}}(\overline{a})=\frac{\partial f(\overline{x})}{\partial x^i}\left(\frac{\partial x^i}{\partial\overline{x}}(\overline{a})\right)=\frac{\partial f(\overline{x})}{\partial x^i}\left(dx^i(\overline{a})\right)=\frac{\partial f(\overline{x})}{\partial x^i}(a^i)$$

Therefore, the Gâteaux derivative of function $f$ is linear combination of 1-$D$-forms $dx^i$. Let

$$g:D^n\rightarrow D^n$$
$$x'^i=x'^i(x^1,...,x^n)$$

be coordinate transformation of affine space. From chain rule (theorem [9]-6.2.12), it follows

$$dx'^i=\frac{\partial x'^i}{\partial\overline{x}}=\frac{\partial x'^i}{\partial x^j}\left(\frac{\partial x^j}{\partial\overline{x}}\right)=\frac{\partial x'^i}{\partial x^j}\left(dx^j\right)$$

The coordinate line is a curve, along which changes only one coordinate $x^i$. We consider the coordinate line as a curve with a tangent vector

(8.1.3) $$\overline{x}_i(\overline{x})(v)=\frac{\partial\overline{x}}{\partial x^i}(v)$$

Let the point $O$ has coordinates $\overline{x}$. The set of vectors $\overline{x}_1$, ..., $\overline{x}_n$ is the basis of affine space

$$(\overline{x}_1(\overline{x}),...,\overline{x}_n(\overline{x}),\overline{x})$$

(8.1.4) $$\overline{v}=\overline{x}_i(\overline{x})(v^i)$$

Let

$$g:D^n\rightarrow D^n$$
$$x'^i=x'^i(x^1,...,x^n)$$

be coordinate transformation of affine space. From chain rule (theorem [9]-6.2.12), it follows

(8.1.5) $$\overline{x}'_i(\overline{x}')(v'^i)=\frac{\partial\overline{x}}{\partial x'^i}(v'^i)=\frac{\partial\overline{x}}{\partial x^j}\left(\frac{\partial x^j}{\partial x'^i}(v'^i)\right)=\overline{x}_j(\overline{x})\left(\frac{\partial x^j}{\partial x'^i}(v'^i)\right)$$

The rule of transformation of coordinates of vector

(8.1.6) $$v^j=\frac{\partial x^j}{\partial x'^i}(v'^i)$$

follows from comparison of equations (8.1.4), (8.1.5).

**Example 8.1.2.** Let us consider linear coordinate transformation

$$x'^1 = ax^1b + ax^2c$$
$$x'^2 = x^1 + x^2 \tag{8.1.7}$$

Then

$$x^1 = x'^2 - x^2 \tag{8.1.8}$$
$$x'^1 = a(x'^2 - x^2)b + ax^2c = ax'^2b + ax^2(c-b)$$
$$ax^2(c-b) = x'^1 - ax'^2b \tag{8.1.9}$$

From equations (8.1.8), (8.1.9), it follows inverse transformation

$$x^2 = a^{-1}x'^1(c-b)^{-1} - x'^2b(c-b)^{-1}$$
$$x^1 = x'^2 - a^{-1}x'^1(c-b)^{-1} + x'^2b(c-b)^{-1} \tag{8.1.10}$$
$$= -a^{-1}x'^1(c-b)^{-1} + x'^2(1+b(c-b)^{-1})$$

From equations (8.1.7), it follows rule of transformation of 1-$D$-forms

$$dx'^1 = adx^1b + adx^2c$$
$$dx'^2 = dx^1 + dx^2$$

Rule of transformation of vector

$$v^1 = -a^{-1}v'^1(c-b)^{-1} + v'^2(1+b(c-b)^{-1})$$
$$v^2 = a^{-1}v'^1(c-b)^{-1} - v'^2b(c-b)^{-1}$$

follows from equations (8.1.6), (8.1.10). □

## 8.2. Paralel Transfer

Let a vector $\overline{v}$ have coordinates $v_0^i$ in the point $M_0$. Let the vector $\overrightarrow{M_0M_1}$ be the vector of infinitesimal displacement. Let $N_0 = M_0 + \overline{v}$, $N_1 = M_1 + \overline{v}$. According to axiom 3 of definition 6.1.1 coordinates of vectors $\overrightarrow{M_0N_0}$ and $\overrightarrow{M_1N_1}$ equal. How this construction will look like in curvilinear coordinates?

We assume that functions $v^i(\overline{x})$ are continuously differentiable in the Gâteaux sense. Therefore, vectors of the local basis $\overline{x}_i$ and coordinates of the vector $v^i$ are functions continuously differentiable in the Gâteaux sense. The equation

$$\overline{v} = \overline{x}_i(\overline{x})(v^i(\overline{x})) \tag{8.2.1}$$

is true at the point $M(\overline{x})$. Since $\overline{v} = \text{const}$, then differentiating the equation (8.2.1) with respect to $\overline{x}$ we get

$$0 = \frac{\partial \overline{x}_i(\overline{x})}{\partial \overline{x}}(v^i)(\overline{a}) + \overline{x}_i(\overline{x})\left(\frac{\partial v^i}{\partial \overline{x}}(\overline{a})\right) \tag{8.2.2}$$

From equations (8.1.2), (8.2.2), it follows

$$0 = \frac{\partial \overline{x}_i(\overline{x})}{\partial x^j}(v^i)(a^j) + \overline{x}_i(\overline{x})\left(\frac{\partial v^i}{\partial x^j}(a^j)\right) \tag{8.2.3}$$

The expression $\frac{\partial \overline{x}_i(\overline{x})}{\partial x^j}(v^i)(a^j)$ is a vector of affine space and, in the local basis $\overline{x}_i(\overline{x})$, this expression has expansion

$$\frac{\partial \overline{x}_i(\overline{x})}{\partial x^j}(v^i)(a^j) = \overline{x}_k(\overline{x})(\Gamma^k_{ji}(v^i)(a^j)) \tag{8.2.4}$$

where $\Gamma^k_{ji}$ are bilinear maps that we call **connection coefficients in $D$-affine space**. From equations (8.2.4), (8.2.3), it follows

$$0 = \overline{x}_k(\overline{x})(\Gamma^k_{ji}(v^i)(a^j)) + \overline{x}_i(\overline{x})\left(\frac{\partial v^i}{\partial x^j}(a^j)\right) \tag{8.2.5}$$

Because $\overline{x}_k(\overline{x})$ is the local basis, then, from equation (8.2.5), it follows

$$0 = \overline{x}_k(\overline{x})\left(\Gamma^k_{ji}(v^i)(a^j) + \frac{\partial v^k}{\partial x^j}(a^j)\right)$$

$$0 = \Gamma^k_{ji}(v^i)(a^j) + \frac{\partial v^k}{\partial x^j}(a^j) \tag{8.2.6}$$

Because $\overline{a}$ is an arbitrary vector, then we can write the equation (8.2.6) in the following form

$$0 = \Gamma^k_{ji}(v^i) + \frac{\partial v^k}{\partial x^j} \tag{8.2.7}$$

**Theorem 8.2.1.** *In affine space*

$$\Gamma^k_{ji}(v^i)(a^j) = \Gamma^k_{ij}(v^i)(a^j) \tag{8.2.8}$$

Proof. According to equation (8.1.3),

$$\frac{\partial \overline{x}_i(\overline{x})}{\partial x^j}(v^i)(a^j) = \frac{\partial^2 \overline{x}}{\partial x^j \partial x^i}(v^i; a^j) \tag{8.2.9}$$

According to theorem [9]-9.1.6

$$\frac{\partial^2 \overline{x}}{\partial x^j \partial x^i}(v^i; a^j) = \frac{\partial^2 \overline{x}}{\partial x^i \partial x^j}(v^i; a^j) \tag{8.2.10}$$

Equation (8.2.8) follows from equations (8.2.4), (8.2.10). □

Let

$$g : D^n \to D^n$$
$$x'^i = x'^i(x^1, ..., x^n)$$

be coordinate transformation of affine space. Let us differentiate the equation (8.1.5) with respect to $x'^j$

$$\begin{aligned}&\frac{\partial \overline{x}'_i(\overline{x}')}{\partial x'^j}(v'^i)(a'^j)\\ =&\frac{\overline{x}_k(\overline{x})\left(\frac{\partial x^k}{\partial x'^i}(v'^i)\right)}{\partial x^m}\left(\frac{\partial x^m}{\partial x'^j}(a'^j)\right) + \overline{x}_k(\overline{x})\left(\frac{\partial^2 x^k}{\partial x'^j \partial x'^i}(v'^i; a'^j)\right)\end{aligned} \tag{8.2.11}$$

From equations (8.2.4), (8.2.11), it follows that

$$
\begin{aligned}
&\overline{x}'_p(\overline{x}')(\Gamma'^p_{ji}(v'^i)(a'^j))\\
&=\overline{x}_r(\overline{x})\left(\Gamma^r_{mk}\left(\frac{\partial x^k}{\partial x'^i}(v'^i)\right)\left(\frac{\partial x^m}{\partial x'^j}(a'^j)\right)\right)+\overline{x}_r(\overline{x})\left(\frac{\partial^2 x^r}{\partial x'^j\partial x'^i}(v'^i;a'^j)\right)\\
&=\overline{x}_r(\overline{x})\left(\Gamma^r_{mk}\left(\frac{\partial x^k}{\partial x'^i}(v'^i)\right)\left(\frac{\partial x^m}{\partial x'^j}(a'^j)\right)+\frac{\partial^2 x^r}{\partial x'^j\partial x'^i}(v'^i;a'^j)\right)
\end{aligned}
\tag{8.2.12}
$$

From equations (8.1.5), (8.2.12) it follows that

$$
\begin{aligned}
&\overline{x}'_p(\overline{x}')(\Gamma'^p_{ji}(v'^i)(a'^j))\\
&=\overline{x}'_p(\overline{x})\left(\frac{\partial x'^p}{\partial x^r}\left(\Gamma^r_{mk}\left(\frac{\partial x^k}{\partial x'^i}(v'^i)\right)\left(\frac{\partial x^m}{\partial x'^j}(a'^j)\right)+\frac{\partial^2 x^r}{\partial x'^j\partial x'^i}(v'^i;a'^j)\right)\right)
\end{aligned}
\tag{8.2.13}
$$

From equation (8.2.13) it follows that

$$
\begin{aligned}
&\Gamma'^p_{ji}(v'^i)(a'^j)\\
&=\frac{\partial x'^p}{\partial x^r}\left(\Gamma^r_{mk}\left(\frac{\partial x^k}{\partial x'^i}(v'^i)\right)\left(\frac{\partial x^m}{\partial x'^j}(a'^j)\right)+\frac{\partial^2 x^r}{\partial x'^j\partial x'^i}(v'^i;a'^j)\right)
\end{aligned}
\tag{8.2.14}
$$

CHAPTER 9

# Manifolds with $D$-Affine Connections

The main goal of this chapter is to give general presentation about problems that we need to solve in process of exploration the differential geometry.

## 9.1. Manifolds with $D$-Affine Connections

**Definition 9.1.1.** Fibered central $D$-affine space is called **manifold with $D$-affine connections**. □

We assume that dimension of base equal to dimension of fiber. This allows us to identify fiber with tangent space of base. Let us consider subset $U$ of manifold where homeomorphism

$$\overline{x} : D^n \to U$$

is defined.

We assume that map between fibers is morphism of central $D$-affine space

$$d\overline{v} = \Gamma(\overline{v})(d\overline{x})$$

where bilinear map $\Gamma$ is called **$D$-affine connection**. In the local basis we can represent connection as the set of bilinear functions $\Gamma^k_{ij}$

$$\Gamma(\overline{v})(\overline{a}) = \overline{x}_k(\overline{x})(\Gamma^k_{ji}(v^i)(a^j))$$

that we call **$D$-affine connection coefficients**. When we move from coordinate system $\overline{x}$ to coordinate system $\overline{x}'$ $D$-affine connection coefficients transform according to rule

$$\begin{aligned}&\Gamma'^p_{ji}(v'^i)(a'^j)\\ =&\frac{\partial x'^p}{\partial x^r}\left(\Gamma^r_{mk}\left(\frac{\partial x^k}{\partial x'^i}(v'^i)\right)\left(\frac{\partial x^m}{\partial x'^j}(a'^j)\right)+\frac{\partial^2 x^r}{\partial x'^j\partial x'^i}(v'^i;a'^j)\right)\end{aligned}$$

Covariant derivative of vector field $\overline{v}$ has form

$$D(\overline{v})(d\overline{x}) = \partial(\overline{v})(d\overline{x}) - \Gamma(\overline{v})(d\overline{x})$$

A vector field $\overline{v}$ is being parallel transported in given direction if

$$\partial(\overline{v})(d\overline{x}) = \Gamma(\overline{v})(d\overline{x})$$

If $\overline{x} = \overline{x}(t)$ is map of real field into considered space, then this map determines geodesic if tangent vector is being parallel transported along geodesic

$$\partial(\partial\overline{x}(d\overline{x}))(d\overline{x}) = \Gamma(d\overline{x})(d\overline{x}) \tag{9.1.1}$$

Equation (9.1.1) is equivalent to equation

$$\partial^2\overline{x}(d\overline{x}; d\overline{x}) = \Gamma(d\overline{x})(d\overline{x}) \tag{9.1.2}$$

CHAPTER 10

# References

[1] Serge Lang, Algebra, Springer, 2002

[2] S. Burris, H.P. Sankappanavar, A Course in Universal Algebra, Springer-Verlag (March, 1982),
eprint http://www.math.uwaterloo.ca/ snburris/htdocs/ualg.html
(The Millennium Edition)

[3] P. K. Rashevsky, Riemann Geometry and Tensor Calculus,
Moscow, Nauka, 1967

[4] A. G. Kurosh, High Algebra, Moscow, Nauka, 1968

[5] A. G. Kurosh, Lectures on General Algebra, Chelsea Pub Co, 1965

[6] Aleks Kleyn, Lectures on Linear Algebra over Division Ring,
eprint arXiv:math.GM/0701238 (2010)

[7] Aleks Kleyn, Fibered Correspondence,
eprint arXiv:0707.2246 (2007)

[8] Aleks Kleyn, Morphism of $T\star$-Representations,
eprint arXiv:0803.2620 (2008)

[9] Aleks Kleyn, Introduction into Calculus over Division Ring,
eprint arXiv:0812.4763 (2010)

[10] Aleks Kleyn, Representation Theory: Representation of Universal Algebra,
Lambert Academic Publishing, 2011

[11] John C. Baez, An Introduction to n-Categories,
eprint arXiv:q-alg/9705009 (1997)

[12] Paul M. Cohn, Universal Algebra, Springer, 1981

[13] Paul Bamberg, Shlomo Sternberg, A course in mathematics for students of physics, Cambridge University Press, 1991

[14] Alain Connes, Noncommutative Geometry,
Academic Press, 1994

CHAPTER 11

# Index

CHAPTER 12

# Special Symbols and Notations

# Введение в геометрию над телом

Александр Клейн

*E-mail address*: Aleks_Kleyn@MailAPS.org
*URL*: http://sites.google.com/site/alekskleyn/
*URL*: http://arxiv.org/a/kleyn_a_1
*URL*: http://AleksKleyn.blogspot.com/

Аннотация. Теория представлений универсальной алгебры является естественным развитием теории универсальной алгебры. Изучение морфизмов представлений ведёт к понятиям множества образующих и базиса представления. В книге рассмотрено понятие башни $T\star$-представлений $\Omega_i$-алгебр, $i = 1, ..., n$, как множество согласованных $T\star$-представлений $\Omega_i$-алгебр.

Рассматривается геометрия аффинного пространства как пример башни представлений. Изучение криволинейных координат позволяет понять как выглядят основные структуры в многообразии аффинной связности. Рассматривается геометрия эвклидова пространства над телом.

# Оглавление

Глава 1

# Предисловие

## 1.1. Башня представлений

Основная задача предлагаемой книги - это рассмотрение простейших геометрических конструкций над телом.

Однако я решил не ограничивать эту книгу описанием аффинной и эвклидовой геометрии. В процессе построения аффинной геометрии я обнаруживаю интересную конструкцию. Я сперва определяю $D_*{}^*$-векторное пространство $\overline{V}$. Это $T\star$-представление тела $D$ в аддитивной группе. Затем я рассматриваю множество точек, в котором определено представление $D_*{}^*$-векторного пространства $\overline{V}$. Эта многоярусная конструкция пробудила мой интерес и вернула меня к изучению $T\star$-представлений универсальной алгебры.

Таким образом возникла концепция башни $T\star$-представлений. Как только я начал рисовать диаграммы, связанные с башней $T\star$-представлений, я увидел узор, похожий на диаграммы, построенные в статье [11]. Хотя мне не удалось довести эту аналогию до конца, я надеюсь, что эта аналогия может привести к интересным результатам.

## 1.2. Морфизм и базис представления

Изучение теории представлений $\mathfrak{F}$-алгебры показывает, что эта теория имеет много общего с теорией $\mathfrak{F}$-алгебры. В [6, 7, 8] я изучал морфизмы представления. Это отображения, которые сохраняют структуру представления $\mathfrak{F}$-алгебры (или башни представлений).

Линейные отображения векторных пространств являются примерами морфизма представлений. Понятие базиса является важной конструкцией в теории линейных пространств. Понятие базиса непосредственно связано с тем, что группа морфизмов линейного пространства имеет два однотранзитивных представления. Основная функция базиса - это порождать линейное пространство.

Возникает естественный вопрос. Можно ли обобщить эту конструкцию на произвольное представление? Базис - это не единственное множество, которое порождает векторное пространство. Это утверждение является исходной точкой, от которой я начал изучение множества образующих представления. Множество образующих представления - это ещё одна интересная параллель теории представлений с теорией универсальной алгебры.

Не для всякого представления возможность найти множество образующих равносильна возможности построения базиса. Задача найти класс представлений, где возможно построение базиса является интересной и важной задачей. Мне интересно, можно ли найти эффективный алгоритм построения базиса квантовой геометрии.

## 1.3. $D$-аффинное пространство

Так же как и в случае $D$-векторных пространств, я могу начать с рассмотрения $D_*{}^*$-аффинного пространства и затем перейти к рассмотрению $D$-аффинных пространств. Однако на этом пути я не ожидаю новых явлений. Так как меня интересуют операции, связаные с дифференцированием, я сразу начинаю с рассмотрения $D$-аффинного пространства.

С того момента, как я понял, что аффинное пространство является башней представлений, а линейные отображения аффинных пространств являются морфизмами башни представлений, я понял как будет выглядеть связность в многообразии $D$-аффинной связности. Тем не менее, прежде чем изучать многообразие аффинной связности, я решил исследовать криволинейные координаты в аффинном пространстве.

Хотя я хорошо представлял, что должно происходить при замене координат, увиденное превзошло мои ожидания. Я встретил две очень важные концепции.

Я начал с построения аффинного пространства как башни представлений. Оказалось, что существует ещё одна модель, в которой векторные поля представлены в виде множества отображений.[1.1] Тем не менее эти две модели не противоречат друг другу, а скорее дополняют друг друга.

Второе наблюдение для меня не оказалось неожиданностью. Когда я начал изучать $D$-векторные пространства, я был готов рассмотреть линейную комбинацию векторов как полилинейную форму, однако я очень быстро пришёл к выводу, что это не изменит размерности пространства, и потому ограничился рассмотрением $D_*{}^*$-базиса.

Интересно наблюдать, как в процессе создания книги появляется новая концепция. Построения, выполненные в главе 8 снова вернули меня к линейной комбинации как сумме линейных форм. Однако на этот раз речь идёт о линейной зависимости 1-форм. При этом коэффициенты линейной зависимости находятся не вне 1-формы, а внутри. Если мы рассматриваем 1-формы над полем, то 1-форма прозрачна для скаляров поля и для нас не имеет значение находятся ли коэффициенты вне формы или внутри.

Линейная зависимость не является темой этой книги. Я рассмотрю этот вопрос в последующих работах. В этой книге я буду пользоваться этим понятием, полагая, что его свойства понятны из приводимых равенств.

Изучение криволинейных координат также позволило увидеть как преобразуются координаты $D$-аффинной связности при замене координат. Однако, то что мне удалось записать явный закон преобразования, на мой взгляд удача, связанная с тем, что связность имеет только один контравариантный индекс. У меня есть все основания полагать, что такая операция, как поднятие или опускание индекса тензора в римановом пространстве может иметь только неявную форму записи.

Многие годы для меня аффинная геометрия была символом самой простой геометрии. Я рад увидеть, что я ошибался. В нужную минуту аффинная геометрия оказалась источником вдохновения. Аффинная геометрия - это врата дифференциальной геометрии. И хотя, для того, чтобы погрузиться в дифференциальную геометрию, я должен хорошо освоиться с дифференциальными

[1.1] Пока неясно, приведёт ли меня эта модель к геометрии, которую изучает Ален Конн в [14]. Ответ на этот вопрос требует дальнейшего исследования.

уравнениями, я завершаю книгу очень коротким экскурсом в многообразия аффинной связности, для того, чтобы почувствовать вкус новой геометрии.

## 1.4. Соглашения

(1) Функция и отображение - синонимы. Однако существует традиция соответствие между кольцами или векторными пространствами называть отображением, а отображение поля действительных чисел или алгебры кватернионов называть функцией. Я тоже следую этой традиции, хотя встречается текст, в котором неясно, какому термину надо отдать предпочтение.

(2) Тело $D$ можно рассматривать как $D$-векторное пространство размерности 1. Соответственно этому, мы можем изучать не только гомоморфизм тела $D_1$ в тело $D_2$, но и линейное отображение тел. При этом подразумевается, что отображение мультипликативно над максимально возможным полем. В частности, линейное отображение тела $D$ мультипликативно над центром $Z(D)$. Это не противоречит определению линейного отображения поля, так как для поля $F$ справедливо $Z(F) = F$. Если поле $F$ отлично от максимально возможного, то я это явно указываю в тексте.

(3) Несмотря на некоммутативность произведения многие утверждения сохраняются, если заменить например правое представление на левое представление или правое векторное пространство на левое векторное пространство. Чтобы сохранить эту симметрию в формулировках теорем я пользуюсь симметричными обозначениями. Например, я рассматриваю $D\star$-векторное пространство и $\star D$-векторное пространство. Запись $D\star$-векторное пространство можно прочесть как D-star-векторное пространство либо как левое векторное пространство. Запись $D\star$-линейно зависимые векторы можно прочесть как D-star-линейно зависимые векторы либо как векторы, линейно зависимые слева.

(4) Пусть $A$ - свободная конечно мерная алгебра. При разложении элемента алгебры $A$ относительно базиса $\overline{\overline{e}}$ мы пользуемся одной и той же корневой буквой для обозначения этого элемента и его координат. Однако в алгебре не принято использовать векторные обозначения. В выражении $a^2$ не ясно - это компонента разложения элемента $a$ относительно базиса или это операция возведения в степень. Для облегчения чтения текста мы будем индекс элемента алгебры выделять цветом. Например,

$$a = a^{i}\overline{e}_{i}$$

(5) Если свободная конечномерная алгебра имеет единицу, то мы будем отождествлять вектор базиса $\overline{e}_0$ с единицей алгебры.

(6) В [12] произвольная операция алгебры обозначена буквой $\omega$, и $\Omega$ - множество операций некоторой универсальной алгебры. Соответственно, универсальная алгебра с множеством операций $\Omega$ обозначается $\Omega$-алгебра. Аналогичные обозначения мы видим в [2] с той небольшой разницей, что операция в алгебре обозначена буквой $f$ и $\mathcal{F}$ - множество операций. Я выбрал первый вариант обозначений, так как в этом случае легче видно, где я использую операцию.

(7) Так как число универсальных алгебр в башне представлений переменно, то мы будем пользоваться векторными обозначениями для башни представлений. Множество $(A_1, ..., A_n)$ $\Omega_i$-алгебр $A_i$, $i = 1$, ..., $n$ мы будем обозначать $\overline{A}$. Множество представлений $(f_{1,2}, ..., f_{n-1,n})$ этих алгебр мы будем обозначать $\overline{f}$. Так как разные алгебры имеют разный тип, мы также будем говорить о множестве $\overline{\Omega}$-алгебр. По отношению к множеству $\overline{A}$ мы также будем пользоваться матричными обозначениями, предложенными в разделе [6]-2.1. Например, символом $\overline{A}_{[1]}$ мы будем обозначать множество $\overline{\Omega}$-алгебр $(A_2, ..., A_n)$. В соответствующем обозначении $(\overline{A}_{[1]}, \overline{f})$ башни представлений подразумевается, что $\overline{f} = (f_{2,3}, ..., f_{n-1,n})$.

(8) Так как мы пользуемся векторными обозначениями для элементов башни представлений, необходимо соглашение о записи операций. Предполагается, что операции выполняются покомпонентно. Например,

$$\overline{r}(\overline{a}) = (r_1(a_1), ..., r_n(a_n))$$

(9) Пусть $A$ - $\Omega_1$-алгебра. Пусть $B$ - $\Omega_2$-алгебра. Запись

$$A \xrightarrow{*} B$$

означает, что определено представление $\Omega_1$-алгебры $A$ в $\Omega_2$-алгебре $B$.

(10) Без сомнения, у читателя могут быть вопросы, замечания, возражения. Я буду признателен любому отзыву.

Глава 2

# Представление универсальной алгебры

## 2.1. Представление универсальной алгебры

**Определение 2.1.1.** Пусть[2.1] на множестве $M$ определена структура $\Omega_2$-алгебры ([2, 12]). Эндоморфизм $\Omega_2$-алгебры

$$t : M \to M$$

называется **преобразованием универсальной алгебры** $M$.[2.2] □

Мы будем обозначать $\delta$ тождественное преобразование.

**Определение 2.1.2.** Если преобразование действует слева

$$u' = tu$$

оно называется **левосторонним преобразованием** или $T\star$**-преобразованием**, Мы будем обозначать ${}^\star M$ множество $T\star$-преобразований универсальной алгебры $M$. □

**Определение 2.1.3.** Если преобразование действует справа

$$u' = ut$$

оно называется **правосторонним преобразованием** или $\star T$**-преобразованием**, Мы будем обозначать $M^\star$ множество $\star T$-преобразований универсальной алгебры $M$. □

**Определение 2.1.4.** Пусть на множестве ${}^\star M$ определена структура $\Omega_1$-алгебры ([2]). Пусть $A$ является $\Omega_1$-алгеброй. Мы будем называть гомоморфизм

$$f : A \to {}^\star M \tag{2.1.1}$$

**левосторонним** или $T\star$**-представлением** $\Omega_1$**-алгебры** $A$ **в** $\Omega_2$**-алгебре** $M$ □

**Определение 2.1.5.** Пусть на множестве $M^\star$ определена структура $\Omega_1$-алгебры ([2]). Пусть $A$ является $\Omega_1$-алгеброй. Мы будем называть гомоморфизм

$$f : A \to M^\star$$

**правосторонним** или $\star T$**-представлением** $\Omega_1$**-алгебры** $A$ **в** $\Omega_2$**-алгебре** $M$ □

Мы распространим на теорию представлений соглашение, описанное в замечании [6]-2.2.15. Мы можем записать принцип двойственности в следующей форме

[2.1] Определения и теоремы из этой главы включены также в главы [10]-2, [10]-3.

[2.2] Если множество операций $\Omega_2$-алгебры пусто, то $t$ является отображением.

**Теорема 2.1.6** (принцип двойственности)**.** *Любое утверждение, справедливое для $T\star$-представления $\Omega_1$-алгебры $A$, будет справедливо для $\star T$-представления $\Omega_1$-алгебры $A$.*

**Замечание 2.1.7.** Существует две формы записи преобразования $\Omega_2$-алгебры $M$. Если мы пользуемся операторной записью, то преобразование $A$ записывается в виде $Aa$ или $aA$, что соответствует $T\star$-преобразованию или $\star T$-преобразованию. Если мы пользуемся функциональной записью, то преобразование $A$ записывается в виде $A(a)$ независимо от того, это $T\star$-преобразование или $\star T$-преобразование. Эта запись согласована с принципом двойственности.

Это замечание является основой следующего соглашения. Когда мы пользуемся функциональной записью, мы не различаем $T\star$-преобразование и $\star T$-преобразование. Мы будем обозначать ${}^*M$ множество преобразований $\Omega_2$-алгебры $M$. Пусть на множестве ${}^*M$ определена структура $\Omega_1$-алгебры. Пусть $A$ является $\Omega_1$-алгеброй. Мы будем называть гомоморфизм

$$f : A \to {}^*M \tag{2.1.2}$$

**представлением $\Omega_1$-алгебры $A$ в $\Omega_2$-алгебре** $M$. Мы будем также пользоваться записью

$$f : A \xrightarrow{*} M$$

для обозначения представления $\Omega_1$-алгебры $A$ в $\Omega_2$-алгебре $M$.

Соответствие между операторной записью и функциональной записью однозначно. Мы можем выбирать любую форму записи, которая удобна для изложения конкретной темы. □

Существует несколько способов описать представление. Мы можем указать отображение $f$, имея в виду что область определения - это $\Omega_1$-алгебра $A$ и область значений - это $\Omega_1$-алгебра ${}^*M$. Либо мы можем указать $\Omega_1$-алгебру $A$ и $\Omega_2$-алгебру $M$, имея в виду что нам известна структура отображения $f$.[2.3]

Диаграмма

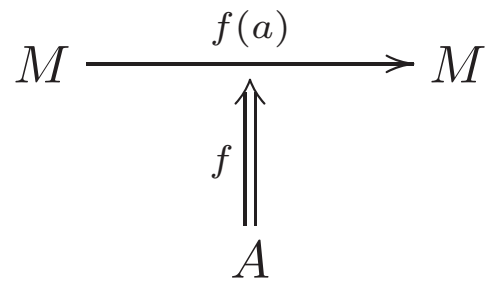

означает, что мы рассматриваем представление $\Omega_1$-алгебры $A$. Отображение $f(a)$ является образом $a \in A$.

**Определение 2.1.8.** Мы будем называть представление $\Omega_1$-алгебры $A$ **эффективным**, если отображение (2.1.2) - изоморфизм $\Omega_1$-алгебры $A$ в ${}^*M$. □

**Замечание 2.1.9.** Если $T\star$-представление $\Omega_1$-алгебры эффективно, мы можем отождествлять элемент $\Omega_1$-алгебры с его образом и записывать $T\star$-преобразование, порождённое элементом $a \in A$, в форме

$$v' = av$$

[2.3]Например, мы рассматриваем векторное пространство $\overline{V}$ над полем $D$ ( определение [6]-4.1.4 ).

Если $\star T$-представление $\Omega_1$-алгебры эффективно, мы можем отождествлять элемент $\Omega_1$-алгебры с его образом и записывать $\star T$-преобразование, порождённое элементом $a \in A$, в форме

$$v' = va$$

□

**Определение 2.1.10.** Мы будем называть представление $\Omega_1$-алгебры **транзитивным**, если для любых $a, b \in V$ существует такое $g$, что

$$a = f(g)(b)$$

Мы будем называть представление $\Omega_1$-алгебры **однотранзитивным**, если оно транзитивно и эффективно. □

**Теорема 2.1.11.** *$T\star$-представление однотранзитивно тогда и только тогда, когда для любых $a, b \in M$ существует одно и только одно $g \in A$ такое, что $a = f(g)(b)$*

Доказательство. Следствие определений 2.1.8 и 2.1.10. □

### 2.2. Морфизм представлений универсальной алгебры

**Теорема 2.2.1.** *Пусть $A$ и $B$ - $\Omega_1$-алгебры. Представление $\Omega_1$-алгебры $B$*

$$g : B \to {}^\star M$$

*и гомоморфизм $\Omega_1$-алгебры*

$$h : A \to B \tag{2.2.1}$$

*определяют представление $f$ $\Omega_1$-алгебры $A$*

$$\begin{array}{ccc} A & \xrightarrow{\quad f \quad} & {}^*M \\ & {\scriptstyle h}\searrow \quad \nearrow{\scriptstyle g} & \\ & B & \end{array}$$

Доказательство. Отображение $f$ является гомоморфизмом $\Omega_1$-алгебры $A$ в $\Omega_1$-алгебру ${}^*M$, так как отображение $g$ является гомоморфизмом $\Omega_1$-алгебры $B$ в $\Omega_1$-алгебру ${}^*M$. □

Если мы изучаем представление $\Omega_1$-алгебры в $\Omega_2$-алгебрах $M$ и $N$, то нас интересуют отображения из $M$ в $N$, сохраняющие структуру представления.

**Определение 2.2.2.** Пусть

$$f : A \to {}^*M$$

представление $\Omega_1$-алгебры $A$ в $\Omega_2$-алгебре $M$ и

$$g : B \to {}^*N$$

представление $\Omega_1$-алгебры $B$ в $\Omega_2$-алгебре $N$. Пара отображений

$$(r : A \to B, R : M \to N) \tag{2.2.2}$$

таких, что

- $r$ - гомоморфизм $\Omega_1$-алгебры
- $R$ - гомоморфизм $\Omega_2$-алгебры

•

$$R \circ f(a) = g(r(a)) \circ R \tag{2.2.3}$$

называется **морфизмом представлений из $f$ в $g$**. Мы также будем говорить, что определён **морфизм представлений $\Omega_1$-алгебры в $\Omega_2$-алгебре**. □

**Замечание 2.2.3.** Мы можем рассматривать пару отображений $r$, $R$ как отображение

$$F : A \cup M \to B \cup N$$

такое, что

$$F(A) = B \quad F(M) = N$$

Поэтому в дальнейшем мы будем говорить, что дано отображение $(r, R)$. □

**Определение 2.2.4.** Если представления $f$ и $g$ совпадают, то морфизм представлений $(r, R)$ называется **морфизмом представления $f$**. □

Для произвольного $m \in M$ равенство (2.2.3) имеет вид

$$R(f(a)(m)) = g(r(a))(R(m)) \tag{2.2.4}$$

**Замечание 2.2.5.** Рассмотрим морфизм представлений (2.2.2). Мы можем обозначать элементы множества $B$, пользуясь буквой по образцу $b \in B$. Но если мы хотим показать, что $b$ является образом элемента $a \in A$, мы будем пользоваться обозначением $r(a)$. Таким образом, равенство

$$r(a) = r(a)$$

означает, что $r(a)$ (в левой части равенства) является образом $a \in A$ (в правой части равенства). Пользуясь подобными соображениями, мы будем обозначать элемент множества $N$ в виде $R(m)$. Мы будем следовать этому соглащению, изучая соотношения между гомоморфизмами $\Omega_1$-алгебр и отображениями между множествами, где определены соответствующие представления.

Мы можем интерпретировать (2.2.4) двумя способами

- Пусть преобразование $f(a)$ отображает $m \in M$ в $f(a)(m)$. Тогда преобразование $g(r(a))$ отображает $R(m) \in N$ в $R(f(a)(m))$.
- Мы можем представить морфизм представлений из $f$ в $g$, пользуясь диаграммой

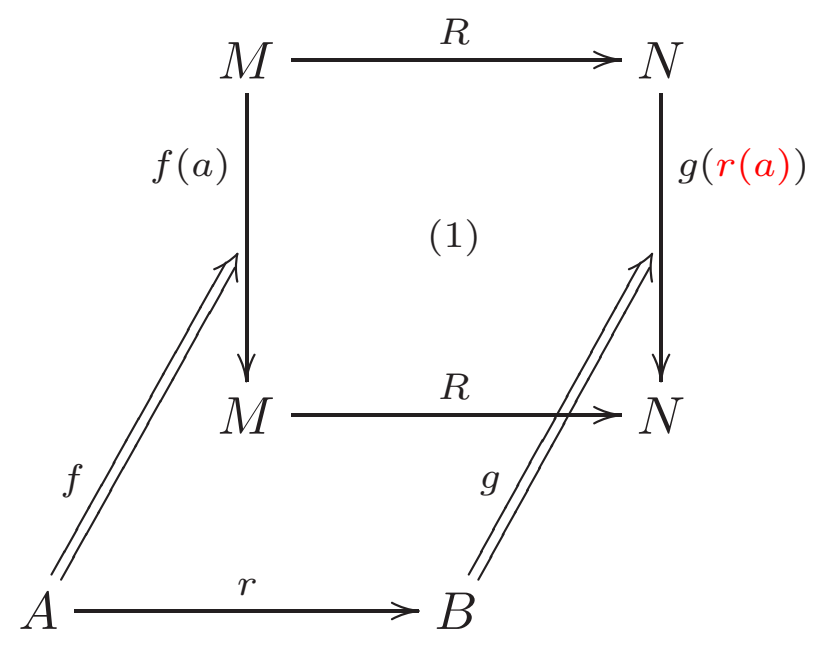

Из (2.2.3) следует, что диаграмма (1) коммутативна.

□

**Теорема 2.2.6.** *Рассмотрим представление*

$$f : A \to {}^*M$$

*$\Omega_1$-алгебры $A$ и представление*

$$g : B \to {}^*N$$

*$\Omega_1$-алгебры $B$. Морфизм*

$$h : A \longrightarrow B \qquad H : M \longrightarrow N$$

*представлений из $f$ в $g$ удовлетворяет соотношению*

(2.2.5) $$H \circ \omega(f(a_1), ..., f(a_n)) = \omega(g(h(a_1)), ..., g(h(a_n))) \circ H$$

*для произвольной $n$-арной операции $\omega$ $\Omega_1$-алгебры.*

Доказательство. Так как $f$ - гомоморфизм, мы имеем

(2.2.6) $$H \circ \omega(f(a_1), ..., f(a_n)) = H \circ f(\omega(a_1, ..., a_n))$$

Из (2.2.3) и (2.2.6) следует

(2.2.7) $$H \circ \omega(f(a_1), ..., f(a_n)) = g(h(\omega(a_1, ..., a_n))) \circ H$$

Так как $h$ - гомоморфизм, из (2.2.7) следует

(2.2.8) $$H \circ \omega(f(a_1), ..., f(a_n)) = g(\omega(h(a_1), ..., h(a_n))) \circ H$$

Так как $g$ - гомоморфизм, из (2.2.8) следует (2.2.5). □

**Теорема 2.2.7.** *Пусть отображение*

$$h : A \longrightarrow B \qquad H : M \longrightarrow N$$

*является морфизмом из представления*

$$f : A \to {}^*M$$

*$\Omega_1$-алгебры $A$ в представление*

$$g : B \to {}^*N$$

*$\Omega_1$-алгебры $B$. Если представление $f$ эффективно, то отображение*

$${}^*H : {}^*M \to {}^*N$$

*определённое равенством*

(2.2.9) $${}^*H(f(a)) = g(h(a))$$

*является гомоморфизмом $\Omega_1$-алгебры.*

Доказательство. Так как представление $f$ эффективно, то для выбранного преобразования $f(a)$ выбор элемента $a$ определён однозначно. Следовательно, преобразование $g(h(a))$ в равенстве (2.2.9) определено корректно.

Так как $f$ - гомоморфизм, мы имеем

(2.2.10) $${}^*H(\omega(f(a_1), ..., f(a_n))) = {}^*H(f(\omega(a_1, ..., a_n)))$$

Из (2.2.9) и (2.2.10) следует

(2.2.11) $${}^*H(\omega(f(a_1), ..., f(a_n))) = g(h(\omega(a_1, ..., a_n)))$$

Так как $h$ - гомоморфизм, из (2.2.11) следует

(2.2.12) $${}^*H(\omega(f(a_1), ..., f(a_n))) = g(\omega(h(a_1), ..., h(a_n)))$$

Так как $g$ - гомоморфизм,

$${}^*H(\omega(f(a_1),...,f(a_n))) = \omega(g(h(a_1)),...,g(h(a_n))) = \omega({}^*H(f(a_1)),...,{}^*H(f(a_n)))$$

следует из (2.2.12). Следовательно, отображение ${}^*H$ является гомоморфизмом $\Omega_1$-алгебры. □

**Теорема 2.2.8.** *Если представление*

$$f : A \to {}^*M$$

*$\Omega_1$-алгебры $A$ однотранзитивно и представление*

$$g : B \to {}^*N$$

*$\Omega_1$-алгебры $B$ однотранзитивно, то существует морфизм*

$$h : A \to B \quad H : M \to N$$

*представлений из $f$ в $g$.*

Доказательство. Выберем гомоморфизм $h$. Выберем элемент $m \in M$ и элемент $n \in N$. Чтобы построить отображение $H$, рассмотрим следующую диаграмму

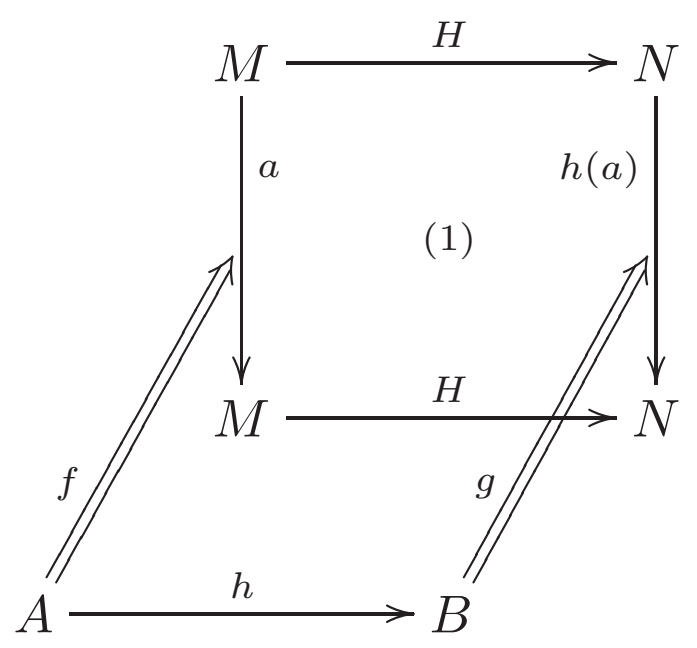

Из коммутативности диаграммы (1) следует

$$H(am) = h(a)H(m)$$

Для произвольного $m' \in M$ однозначно определён $a \in A$ такой, что $m' = am$. Следовательно, мы построили отображении $H$, которое удовлетворяет равенству (2.2.3). □

**Теорема 2.2.9.** *Если представление*

$$f : A \to {}^*M$$

*$\Omega_1$-алгебры $A$ однотранзитивно и представление*

$$g : B \to {}^*N$$

*$\Omega_1$-алгебры $B$ однотранзитивно, то для заданного гомоморфизма $\Omega_1$-алгебры*

$$h : A \longrightarrow B$$

*отображение*

$$H : M \longrightarrow N$$

*такое, что $(h, H)$ является морфизмом представлений из $f$ в $g$, единственно с точностью до выбора образа $n = H(m) \in N$ заданного элемента $m \in M$.*

Доказательство. Из доказательства теоремы 2.2.8 следует, что выбор гомоморфизма $h$ и элементов $m \in M$, $n \in N$ однозначно определяет отображение $H$. □

**Теорема 2.2.10.** *Если представление*

$$f : A \to {}^*M$$

*$\Omega_1$-алгебры $A$ однотранзитивно, то для любого эндоморфизма $h$ $\Omega_1$-алгебры $A$ существует морфизм представления $f$*

$$h : A \to B \quad H : M \to N$$

Доказательство. Рассмотрим следующую диаграмму

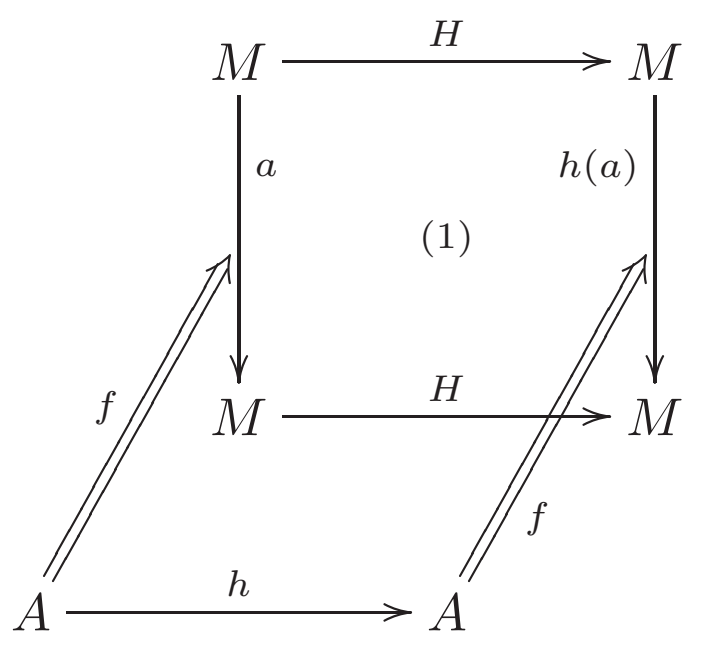

Утверждение теоремы является следствием теоремы 2.2.8. □

**Теорема 2.2.11.** *Пусть*

$$f : A \to {}^*M$$

*представление $\Omega_1$-алгебры $A$,*

$$g : B \to {}^*N$$

*представление $\Omega_1$-алгебры $B$,*

$$h : C \to {}^*L$$

*представление $\Omega_1$-алгебры $C$. Пусть определены морфизмы представлений $\Omega_1$-алгебры*

$$p : A \longrightarrow B \qquad P : M \longrightarrow N$$

$$q : B \longrightarrow C \qquad Q : N \longrightarrow L$$

*Тогда определён морфизм представлений $\Omega_1$-алгебры*

$$r : A \longrightarrow C \qquad R : M \longrightarrow L$$

*где $r = qp$, $R = QP$. Мы будем называть морфизм $(r, R)$ представлений из $f$ в $h$* **произведением морфизмов** $(p, P)$ **и** $(q, Q)$ **представлений универсальной алгебры**.

Доказательство. Мы можем представить утверждение теоремы, пользуясь диаграммой

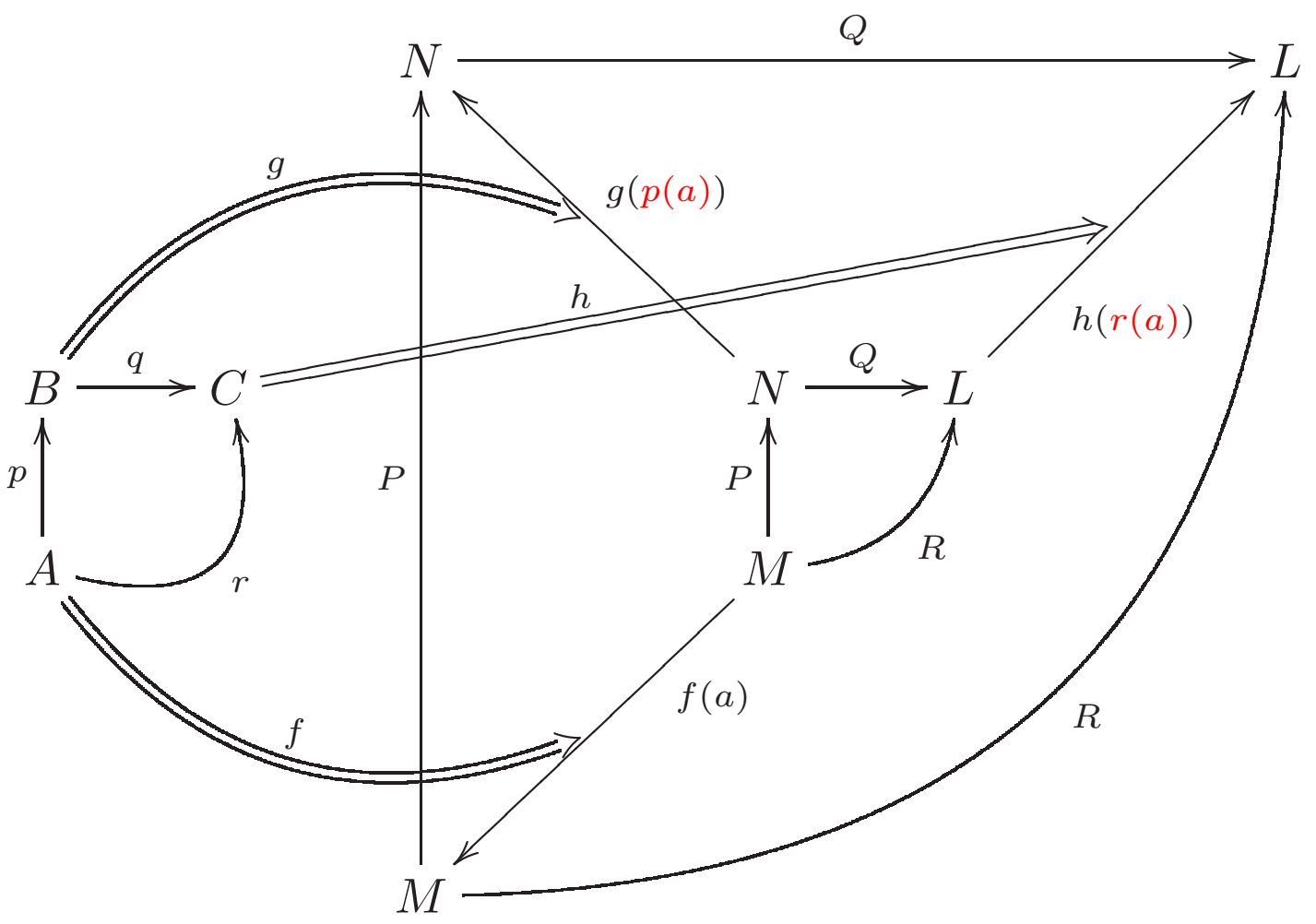

Отображение $r$ является гомоморфизмом $\Omega_1$-алгебры $A$ в $\Omega_1$-алгебру $C$. Нам надо показать, что пара отображений $(r, R)$ удовлетворяет (2.2.3):

$$\begin{aligned} R(f(a)m) &= QP(f(a)m) \\ &= Q(g(p(a))P(m)) \\ &= h(qp(a))QP(m)) \\ &= h(r(a))R(m) \end{aligned}$$

□

**Определение 2.2.12.** Допустим $\mathcal{A}$ категория $\Omega_1$-алгебр. Мы определим **категорию $T\star\mathcal{A}$ $T\star$-представлений универсальной алгебры из категории** $\mathcal{A}$. Объектами этой категории являются $T\star$-представлениями $\Omega_1$-алгебры. Морфизмами этой категории являются морфизмы $T\star$-представлений $\Omega_1$-алгебры. □

**Определение 2.2.13.** Пусть на множестве $M$ определена эквивалентность $S$. Преобразование $f$ называется **согласованным с эквивалентностью** $S$, если из условия $m_1 \equiv m_2 (\mathrm{mod} S)$ следует $f(m_1) \equiv f(m_2) (\mathrm{mod}\ S)$. □

**Теорема 2.2.14.** *Пусть на множестве $M$ определена эквивалентность $S$. Пусть на множестве ${}^*M$ определена $\Omega_1$-алгебра. Если любое преобразование $f \in {}^*M$ согласованно с эквивалентностью $S$, то мы можем определить структуру $\Omega_1$-алгебры на множестве ${}^*(M/S)$.*

Доказательство. Пусть $h = \mathrm{nat}\ S$. Если $m_1 \equiv m_2 (\mathrm{mod} S)$ , то $h(m_1) = h(m_2)$. Поскольку $f \in {}^*M$ согласованно с эквивалентностью $S$, то $h(f(m_1)) = h(f(m_2))$. Это позволяет определить преобразование $F$ согласно правилу

$$F([m]) = h(f(m)) \tag{2.2.13}$$

Пусть $\omega$ - n-арная операция $\Omega_1$-алгебры. Пусть $f_1$, ..., $f_n \in {}^\star M$ и

$$F_1([m]) = h(f_1(m)) \quad ... \quad F_n([m]) = h(f_n(m))$$

Согласно условию теоремы, преобразование

$$f = \omega(f_1, ..., f_n) \in {}^*M$$

согласованно с эквивалентностью $S$. Следовательно, из условия $m_1 \equiv m_2(\mathrm{mod} S)$ и определения 2.2.13 следует

$$f(m_1) \equiv f(m_2)(\mathrm{mod} S) \tag{2.2.14}$$
$$\omega(f_1, ..., f_n)(m_1) \equiv \omega(f_1, ..., f_n)(m_2)(\mathrm{mod} S)$$

Следовательно, мы можем определить операцию $\omega$ на множестве ${}^\star(M/S)$ по правилу

$$\omega(F_1, ..., F_n)[m] = h(\omega(f_1, ..., f_n)(m)) \tag{2.2.15}$$

Из определения (2.2.13) и равенства (2.2.14) следует, что мы корректно определили операцию $\omega$ на множестве ${}^\star(M/S)$. □

**Теорема 2.2.15.** *Пусть*

$$f : A \to {}^*M$$

*представление $\Omega_1$-алгебры $A$,*

$$g : B \to {}^*N$$

*представление $\Omega_1$-алгебры $B$. Пусть*

$$r : A \longrightarrow B \qquad R : M \longrightarrow N$$

*морфизм представлений из $f$ в $g$. Положим*

$$s = rr^{-1} \qquad S = RR^{-1}$$

*Тогда для отображений $r$, $R$ существуют разложения, которые можно описать диаграммой*

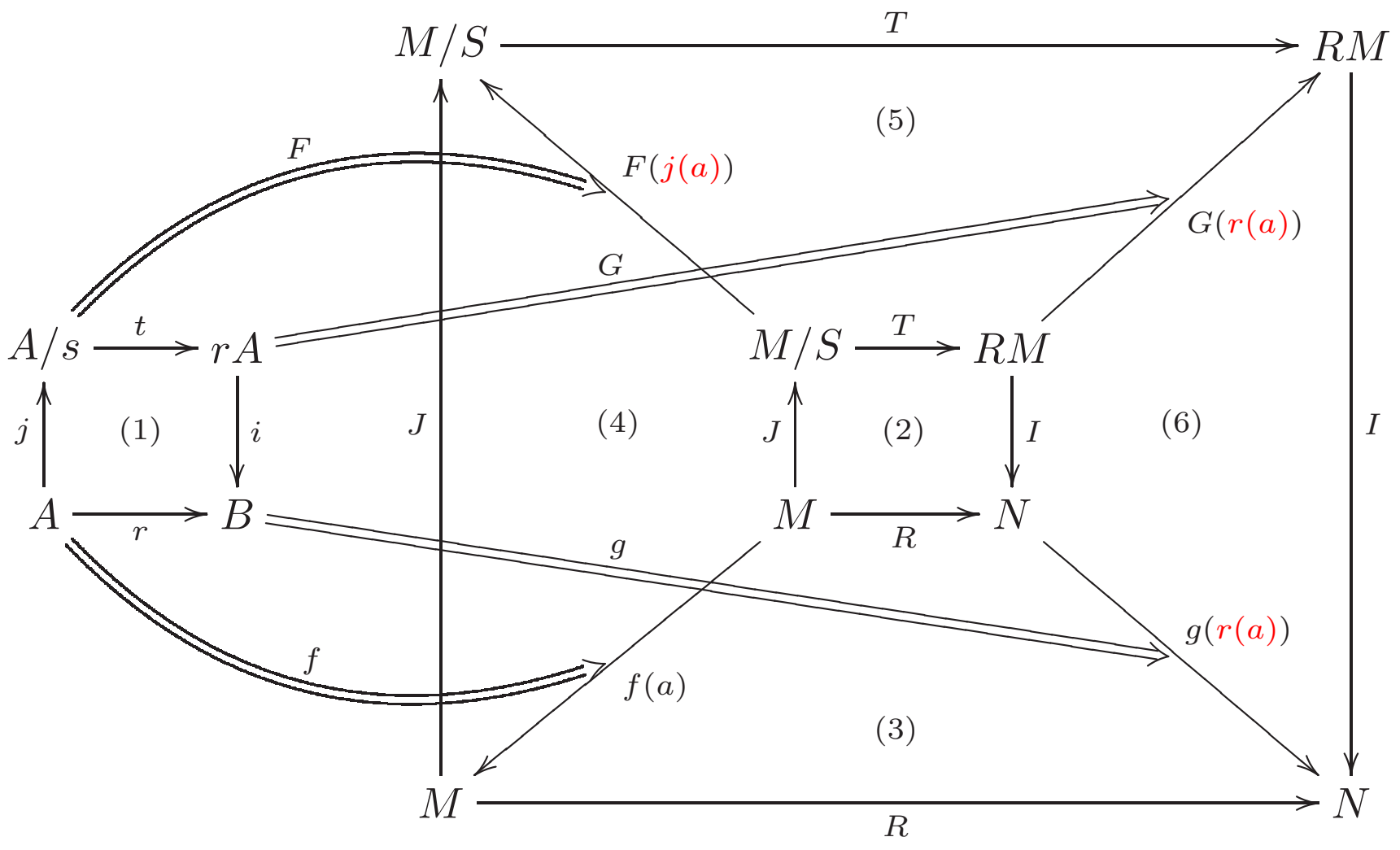

(1) $s = \ker\, r$ *является конгруэнцией на $A$. Существует разложение гомоморфизма $r$*

$$r = itj \tag{2.2.16}$$

$j = \text{nat } s$ - *естественный гомоморфизм*

$$j(a) = j(a)$$

$t$ - *изоморфизм*

$$r(a) = t(j(a)) \tag{2.2.17}$$

$i$ - *вложение*

$$r(a) = i(r(a)) \tag{2.2.18}$$

(2) $S = \ker R$ *является эквивалентностью на* $M$. *Существует разложение отображения* $R$

$$R = ITJ \tag{2.2.19}$$

$J = \text{nat } S$ - *сюръекция*

$$J(m) = J(m)$$

$T$ - *биекция*

$$R(m) = T(J(m)) \tag{2.2.20}$$

$I$ - *вложение*

$$R(m) = I(R(m)) \tag{2.2.21}$$

(3) $F$ - $T\star$-*представление* $\Omega_1$-*алгебры* $A/s$ *в* $M/S$
(4) $G$ - $T\star$-*представление* $\Omega_1$-*алгебры* $rA$ *в* $RM$
(5) $(j,J)$ - *морфизм представлений* $f$ *и* $F$
(6) $(t,T)$ - *морфизм представлений* $F$ *и* $G$
(7) $(t^{-1},T^{-1})$ - *морфизм представлений* $G$ *и* $F$
(8) $(i,I)$ - *морфизм представлений* $G$ *и* $g$
(9) *Существует разложение морфизма представлений*

$$(r,R) = (i,I)(t,T)(j,J) \tag{2.2.22}$$

Доказательство. Существование диаграмм (1) и (2) следует из теоремы II.3.7 ([12], с. 74).

Мы начнём с диаграммы (4).

Пусть $m_1 \equiv m_2 (\text{mod } S)$. Следовательно,

$$R(m_1) = R(m_2) \tag{2.2.23}$$

Если $a_1 \equiv a_2 (\text{mod } s)$, то

$$r(a_1) = r(a_2) \tag{2.2.24}$$

Следовательно, $j(a_1) = j(a_2)$. Так как $(r,R)$ - морфизм представлений, то

$$R(f(a_1)(m_1)) = g(r(a_1))(R(m_1)) \tag{2.2.25}$$

$$R(f(a_2)(m_2)) = g(r(a_2))(R(m_2)) \tag{2.2.26}$$

Из (2.2.23), (2.2.24), (2.2.25), (2.2.26) следует

$$R(f(a_1)(m_1)) = R(f(a_2)(m_2)) \tag{2.2.27}$$

Из (2.2.27) следует

$$f(a_1)(m_1) \equiv f(a_2)(m_2)(\text{mod } S) \tag{2.2.28}$$

и, следовательно,

$$J(f(a_1)(m_1)) = J(f(a_2)(m_2)) \tag{2.2.29}$$

Из (2.2.29) следует, что отображение

$$F(j(a))(J(m)) = J(f(a)(m))) \tag{2.2.30}$$

определено корректно и является преобразованием множества $M/S$.

Из равенства (2.2.28) (в случае $a_1 = a_2$) следует, что для любого $a$ преобразование согласованно с эквивалентностью $S$. Из теоремы 2.2.14 следует, что на множестве $^*(M/S)$ определена структура $\Omega_1$-алгебры. Рассмотрим $n$-арную операцию $\omega$ и $n$ преобразований

$$F(j(a_i))(J(m)) = J(f(a_i)(m))) \quad i = 1, ..., n$$

пространства $M/S$. Мы положим

$$\omega(F(j(a_1)), ..., F(j(a_n)))(J(m)) = J(\omega(f(a_1), ..., f(a_n)))(m))$$

Следовательно, отображение $F$ является представлением $\Omega_1$-алгебры $A/s$.

Из (2.2.30) следует, что $(j, J)$ является морфизмом представлений $f$ и $F$ (утверждение (5) теоремы).

Рассмотрим диаграмму (5).

Так как $T$ - биекция, то мы можем отождествить элементы множества $M/S$ и множества $MR$, причём это отождествление имеет вид

$$T(J(m)) = R(m) \tag{2.2.31}$$

Мы можем записать преобразование $F(j(a))$ множества $M/S$ в виде

$$F(j(a)) : J(m) \to F(j(a))(J(m)) \tag{2.2.32}$$

Так как $T$ - биекция, то мы можем определить преобразование

$$T(J(m)) \to T(F(j(a))(J(m))) \tag{2.2.33}$$

множества $RM$. Преобразование (2.2.33) зависит от $j(a) \in A/s$. Так как $t$ - биекция, то мы можем отождествить элементы множества $A/s$ и множества $rA$, причём это отождествление имеет вид

$$t(j(a)) = r(a) \tag{2.2.34}$$

Следовательно, мы определили отображение

$$G : rA \to {}^\star RM$$

согласно равенству

$$G(t(j(a)))(T(J(m))) = T(F(j(a))(J(m))) \tag{2.2.35}$$

Рассмотрим $n$-арную операцию $\omega$ и $n$ преобразований

$$G(r(a_i))(R(m)) = T(F(j(a_i))(J(m))) \quad i = 1, ..., n$$

пространства $RM$. Мы положим

$$\omega(G(r(a_1)), ..., G(r(a_n)))(R(m)) = T(\omega(F(j(a_1), ..., F(j(a_n)))(J(m))) \tag{2.2.36}$$

Согласно (2.2.35) операция $\omega$ корректно определена на множестве $^\star RM$. Следовательно, отображение $G$ является представлением $\Omega_1$-алгебры.

Из (2.2.35) следует, что $(t, T)$ является морфизмом представлений $F$ и $G$ (утверждение (6) теоремы).

Так как $T$ - биекция, то из равенства (2.2.31) следует

$$J(m)=T^{-1}(R(m)) \tag{2.2.37}$$

Мы можем записать преобразование $G(r(a))$ множества $RM$ в виде

$$G(r(a)):R(m)\rightarrow G(r(a))(R(m)) \tag{2.2.38}$$

Так как $T$ - биекция, то мы можем определить преобразование

$$T^{-1}(R(m))\rightarrow T^{-1}(G(r(a))(R(m))) \tag{2.2.39}$$

множества $M/S$. Преобразование (2.2.39) зависит от $r(a)\in rA$. Так как $t$ - биекция, то из равенства (2.2.34) следует

$$j(a)=t^{-1}(r(a)) \tag{2.2.40}$$

Так как по построению диаграмма (5) коммутативна, то преобразование (2.2.39) совпадает с преобразованием (2.2.32). Равенство (2.2.36) можно записать в виде

$$T^{-1}(\omega(G(r(a_1)),...,G(r(a_n)))(R(m)))=\omega(F(j(a_1),...,F(j(a_n)))(J(m)) \tag{2.2.41}$$

Следовательно, $(t^{-1},T^{-1})$ является морфизмом представлений $G$ и $F$ (утверждение (7) теоремы).

Диаграмма (6) является самым простым случаем в нашем доказательстве. Поскольку отображение $I$ является вложением и диаграмма (2) коммутативна, мы можем отождествить $n\in N$ и $R(m)$, если $n\in \mathrm{Im}R$. Аналогично, мы можем отождествить соответствующие преобразования.

$$g'(i(r(a)))(I(R(m)))=I(G(r(a))(R(m))) \tag{2.2.42}$$

$$\omega(g'(r(a_1)),...,g'(r(a_n)))(R(m))=I(\omega(G(r(a_1),...,G(r(a_n)))(R(m)))$$

Следовательно, $(i,I)$ является морфизмом представлений $G$ и $g$ (утверждение (8) теоремы).

Для доказательства утверждения (9) теоремы осталось показать, что определённое в процессе доказательства представление $g'$ совпадает с представлением $g$, а операции над преобразованиями совпадают с соответствующими операциями на ${}^*N$.

$$\begin{aligned}
g'(i(r(a)))(I(R(m)))&=I(G(r(a))(R(m))) && \text{by } (2.2.42)\\
&=I(G(t(j(a)))(T(J(m)))) && \text{by } (2.2.17), (2.2.20)\\
&=IT(F(j(a))(J(m))) && \text{by } (2.2.35)\\
&=ITJ(f(a)(m)) && \text{by } (2.2.30)\\
&=R(f(a)(m)) && \text{by } (2.2.19)\\
&=g(r(a))(R(m)) && \text{by } (2.2.3)
\end{aligned}$$

$$\begin{aligned}
\omega(G(r(a_1)),...,G(r(a_n)))(R(m))&=T(\omega(F(j(a_1),...,F(j(a_n)))(J(m)))\\
&=T(F(\omega(j(a_1),...,j(a_n)))(J(m)))\\
&=T(F(j(\omega(a_1,...,a_n)))(J(m)))\\
&=T(J(f(\omega(a_1,...,a_n))(m)))
\end{aligned}$$

□

**Определение 2.2.16.** Пусть

$$f : A \to {}^*M$$

представление $\Omega_1$-алгебры $A$,

$$g : B \to {}^*N$$

представление $\Omega_1$-алгебры $B$. Пусть

$$r : A \longrightarrow B \qquad R : M \longrightarrow N$$

морфизм представлений из $r$ в $R$ такой, что $f$ - изоморфизм $\Omega_1$-алгебры и $g$ - изоморфизм $\Omega_2$-алгебры. Тогда отображение $(r, R)$ называется **изоморфизмом представлений**. □

**Теорема 2.2.17.** *В разложении (2.2.22) отображение $(t, T)$ является изоморфизмом представлений $F$ и $G$.*

Доказательство. Следствие определения 2.2.16 и утверждений (6) и (7) теоремы 2.2.15. □

Из теоремы 2.2.15 следует, что мы можем свести задачу изучения морфизма представлений $\Omega_1$-алгебры к случаю, описываемому диаграммой

(2.2.43)

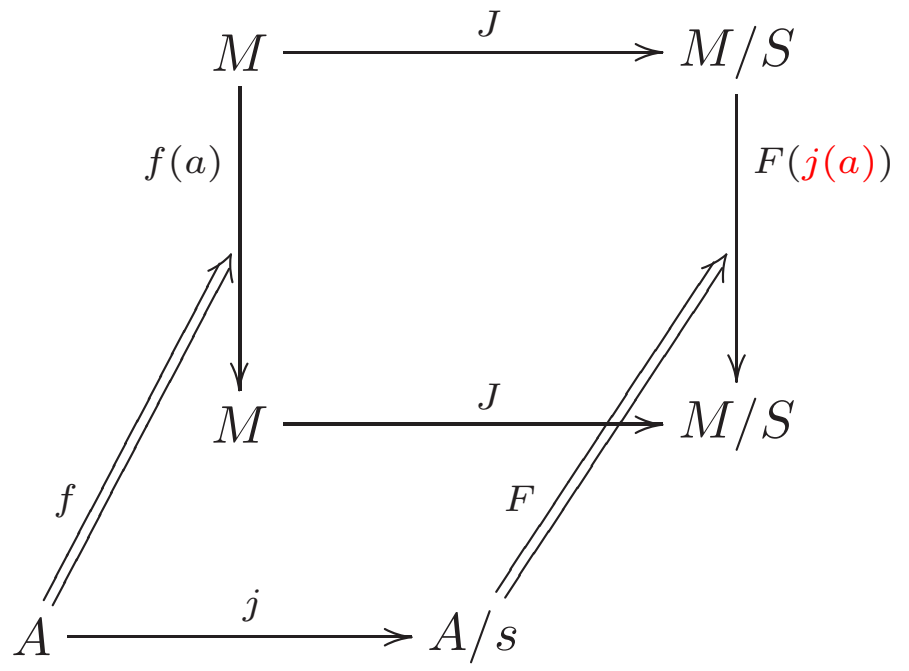

**Теорема 2.2.18.** *Диаграмма (2.2.43) может быть дополнена представлением $F_1$ $\Omega_1$-алгебры $A$ в множестве $M/S$ так, что диаграмма*

(2.2.44)

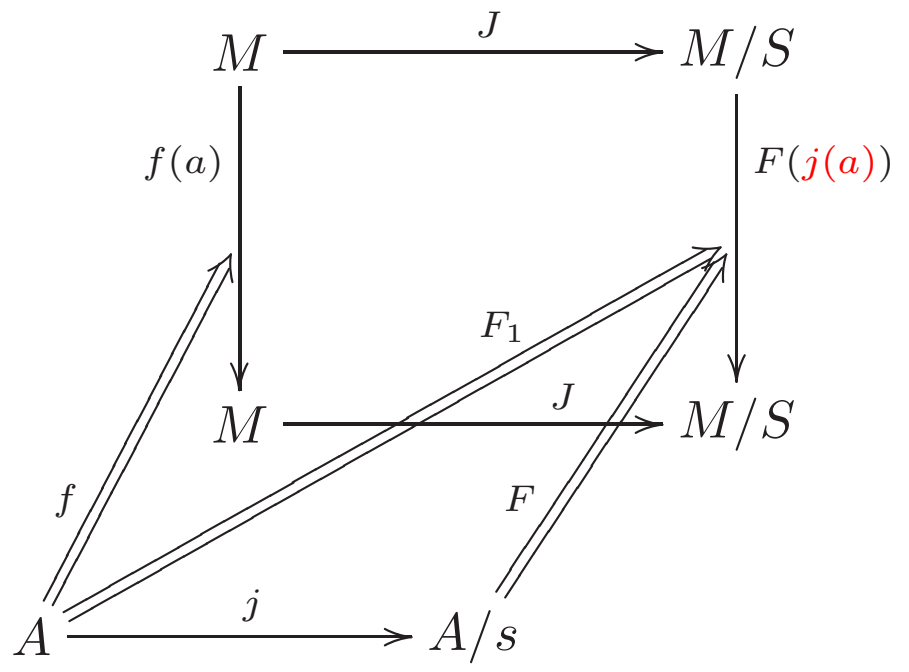

*коммутативна. При этом множество преобразований представления $F$ и множество преобразований представления $F_1$ совпадают.*

Доказательство. Для доказательства теоремы достаточно положить

$$F_1(a) = F(j(a))$$

Так как отображение $j$ - сюрьекция, то $\mathrm{Im}F_1 = \mathrm{Im}F$. Так как $j$ и $F$ - гомоморфизмы $\Omega_1$-алгебры, то $F_1$ - также гомоморфизм $\Omega_1$-алгебры. □

Теорема 2.2.18 завершает цикл теорем, посвящённых структуре морфизма представлений $\Omega_1$-алгебры. Из этих теорем следует, что мы можем упростить задачу изучения морфизма представлений $\Omega_1$-алгебры и ограничиться морфизмом представлений вида

$$id : A \longrightarrow A \qquad R : M \longrightarrow N$$

В этом случае мы можем отождествить морфизм $(\mathrm{id}, R)$ представлений $\Omega_1$-алгебры и соответствующий гомоморфизм $R$ $\Omega_2$-алгебры и пользоваться одной и тойже буквой $R$ для обозначения этих отображений. Мы будем пользоваться диаграммой

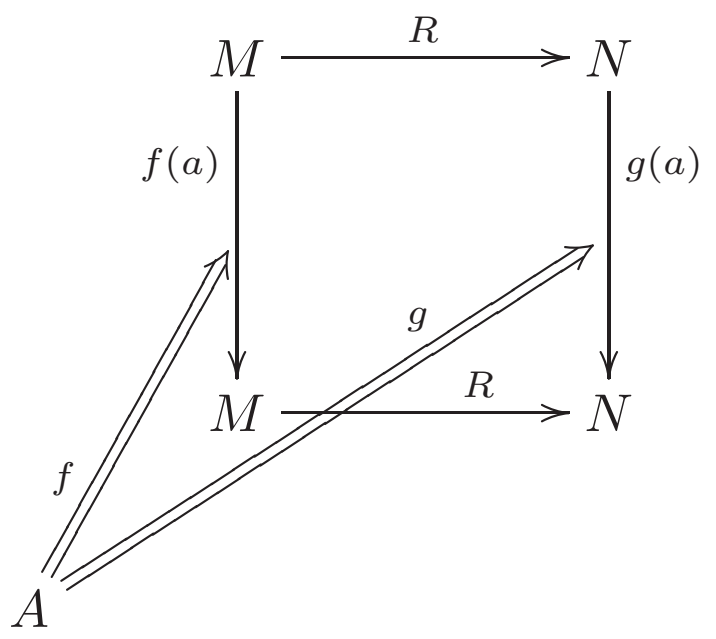

для представления морфизма $(\mathrm{id}, R)$ представлений $\Omega_1$-алгебры. Из диаграммы следует

$$R \circ f(a) = g(a) \circ R \tag{2.2.45}$$

Мы дадим следующее определение по аналогии с определением 2.2.12.

**Определение 2.2.19.** Мы определим **категорию** $T\star A$ $T\star$**-представлений** $\Omega_1$**-алгебры** $A$. Объектами этой категории являются $T\star$-представлениями $\Omega_1$-алгебры $A$. Морфизмами этой категории являются морфизмы $(id, R)$ $T\star$-представлений $\Omega_1$-алгебры $A$. □

### 2.3. Автоморфизм представления универсальной алгебры

**Определение 2.3.1.** Пусть

$$f : A \to {}^*M$$

представление $\Omega_1$-алгебры $A$ в $\Omega_2$-алгебре $M$. Морфизм представлений $\Omega_1$-алгебры

$$(\mathrm{id} : A \to A, R : M \to M)$$

такой, что $R$ - эндоморфизм $\Omega_2$-алгебры называется **эндоморфизмом представления** $f$. □

**Теорема 2.3.2.** *Если представление*

$$f : A \to {}^*M$$

*$\Omega_1$-алгебры $A$ однотранзитивно, то для любых $p$, $q \in M$ существует единственный эндоморфизм*

$$H : M \to M$$

*представления $f$ такой, что $H(p) = q$.*

Доказательство. Рассмотрим следующую диаграмму

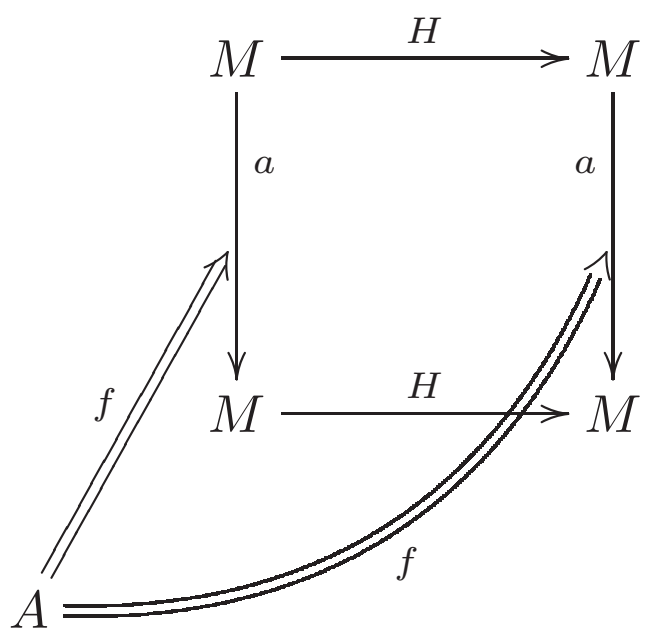

Существование эндоморфизма является следствием теоремы 2.2.8. Единственность эндоморфизма для заданных $p$, $q \in M$ является следствием теоремы 2.2.9, когда $h = \mathrm{id}$. □

**Теорема 2.3.3.** *Эндоморфизмы представления $f$ порождают полугруппу.*

Доказательство. Из теоремы 2.2.11 следует, что произведение эндоморфизмов $(p, P)$, $(r, R)$ представления $f$ является эндоморфизмом $(pr, PR)$ представления $f$. □

**Определение 2.3.4.** Пусть

$$f : A \to {}^*M$$

представление $\Omega_1$-алгебры $A$ в $\Omega_2$-алгебре $M$. Морфизм представлений $\Omega_1$-алгебры

$$(\mathrm{id} : A \to A, R : M \to M)$$

такой, что $R$ - автоморфизм $\Omega_2$-алгебры называется **автоморфизмом представления** $f$. □

**Теорема 2.3.5.** *Пусть*

$$f : A \to {}^*M$$

*представление $\Omega_1$-алгебры $A$ в $\Omega_2$-алгебре $M$. Множество автоморфизмов представления $f$ порождает* **группу** *$GA(f)$ .*

Доказательство. Пусть $R$, $P$ - автоморфизмы представления $f$. Согласно определению 2.3.4, отображения $R$, $P$ являются автоморфизмами $\Omega_2$-алгебры $M$. Согласно теореме II.3.2, ([12], с. 60), отображение $R \circ P$ является автоморфизмом $\Omega_2$-алгебры $M$. Из теоремы 2.2.11 и определения 2.3.4 следует, что произведение автоморфизмов $R \circ P$ представления $f$ является автоморфизмом представления $f$.

Пусть $R$, $P$, $Q$ - автоморфизмы представления $f$. Из цепочки равенств

$$((R\circ P)\circ Q)(a)=(R\circ P)(Q(a))=R(P(Q(a)))$$
$$=R((P\circ Q)(a))=(R\circ(P\circ Q))(a)$$

следует ассоциативность произведения для отображений $R$, $P$, $Q$.[2.4]

Пусть $R$ - автоморфизм представления $f$. Согласно определению 2.3.4 отображение $R$ является автоморфизмом $\Omega_2$-алгебры $M$. Следовательно, отображение $R^{-1}$ является автоморфизмом $\Omega_2$-алгебры $M$. Для автоморфизма $R$ представления справедливо равенство (2.2.4). Положим $m'=R(m)$. Так как $R$ - автоморфизм $\Omega_2$-алгебры, то $m=R^{-1}(m')$ и равенство (2.2.4) можно записать в виде

$$R(f(a')(R^{-1}(m')))=f(a')(m') \tag{2.3.1}$$

Так как отображение $R$ является автоморфизмом $\Omega_2$-алгебры $M$, то из равенства (2.3.1) следует

$$f(a')(R^{-1}(m'))=R^{-1}(f(a')(m')) \tag{2.3.2}$$

Равенство (2.3.2) соответствует равенству (2.2.4) для отображения $R^{-1}$. Следовательно, отображение $R^{-1}$ является автоморфизмом представления $f$. □

## 2.4. Множество образующих представления

**Определение 2.4.1.** Пусть

$$f:A\rightarrow{}^*M$$

представление $\Omega_1$-алгебры $A$ в $\Omega_2$-алгебре $M$. Множество $N\subset M$ называется **стабильным множеством представления** $f$, если $f(a)(m)\in N$ для любых $a\in A$, $m\in N$. □

Мы также будем говорить, что множество $M$ стабильно относительно представления $f$.

**Теорема 2.4.2.** *Пусть*

$$f:A\rightarrow{}^*M$$

*представление $\Omega_1$-алгебры $A$ в $\Omega_2$-алгебре $M$. Пусть множество $N\subset M$ является подалгеброй $\Omega_2$-алгебры $M$ и стабильным множеством представления $f$. Тогда существует представление*

$$f_N:A\rightarrow{}^*N$$

*такое, что $f_N(a)=f(a)|_N$. Представление $f_N$ называется* **подпредставлением представления** *$f$.*

Доказательство. Пусть $\omega_1$ - $n$-арная операция $\Omega_1$-алгебры $A$. Тогда для любых $a_1$, ..., $a_n\in A$ и любого $b\in N$

$$(f_N(a_1)...f_N(a_n)\omega_1)(b)=(f(a_1)...f(a_n)\omega_1)(b)$$
$$=f(a_1...a_n\omega_1)(b)$$
$$=f_N(a_1...a_n\omega_1)(b)$$

[2.4] При доказательстве ассоциативности произведения я следую примеру полугруппы из [5], с. 20, 21.

Пусть $\omega_2$ - $n$-арная операция $\Omega_2$-алгебры $M$. Тогда для любых $b_1$, ..., $b_n \in N$ и любого $a \in A$

$$\begin{aligned} f_N(a)(b_1)...f_N(a)(b_n)\omega_2 &= f(a)(b_1)...f(a)(b_n)\omega_2 \\ &= f(a)(b_1...b_n\omega_2) \\ &= f_N(a)(b_1...b_n\omega_2) \end{aligned}$$

Утверждение теоремы доказано. □

Из теоремы 2.4.2 следует, что если $f_N$ - подпредставление представления $f$, то отображение $(id : A \to A, id_N : N \to M)$ является морфизмом представлений.

**Теорема 2.4.3.** *Множество*[2.5] $\mathcal{B}_f$ *всех подпредставлений представления* $f$ *порождает систему замыканий на* $\Omega_2$*-алгебре* $M$ *и, следовательно, является полной структурой.*

Доказательство. Пусть $(K_\lambda)_{\lambda\in\Lambda}$ - семейство подалгебр $\Omega_2$-алгебры $M$, стабильных относительно представления $f$. Операцию пересечения на множестве $\mathcal{B}_f$ мы определим согласно правилу

$$\bigcap f_{K_\lambda} = f_{\cap K_\lambda}$$

Операция пересечения подпредставлений определена корректно. $\cap K_\lambda$ - подалгебра $\Omega_2$-алгебры $M$. Пусть $m \in \cap K_\lambda$. Для любого $\lambda \in \Lambda$ и для любого $a \in A$, $f(a)(m) \in K_\lambda$. Следовательно, $f(a)(m) \in \cap K_\lambda$. Следовательно, $\cap K_\lambda$ - стабильное множество представления $f$. □

Обозначим соответствующий оператор замыкания через $J(f)$. Таким образом, $J(f, X)$ является пересечением всех подалгебр $\Omega_2$-алгебры $M$, содержащих $X$ и стабильных относительно представления $f$.

**Теорема 2.4.4.** *Пусть*[2.6]

$$f : A \to {}^*M$$

*представление* $\Omega_1$*-алгебры* $A$ *в* $\Omega_2$*-алгебре* $M$*. Пусть* $X \subset M$*. Определим подмножество* $X_k \subset M$ *индукцией по* $k$*.*

$$\begin{gathered} X_0 = X \\ x \in X_k => x \in X_{k+1} \\ x_1 \in X_k, ..., x_n \in X_k, \omega \in \Omega_2(n) => x_1...x_n\omega \in X_{k+1} \\ x \in X_k, a \in A => f(a)(x) \in X_{k+1} \end{gathered}$$

*Тогда*

$$\bigcup_{k=0}^{\infty} X_k = J(f, X)$$

Доказательство. Если положим $U = \cup X_k$, то по определению $X_k$ имеем $X_0 \subset J(f, X)$, и если $X_k \subset J(f, X)$, то $X_{k+1} \subset J(f, X)$. По индукции следует, что $X_k \subset J(f, X)$ для всех $k$. Следовательно,

$$U \subset J(f, X) \tag{2.4.1}$$

[2.5] Это определение аналогично определению структуры подалгебр ([12], стр. 93, 94)

[2.6] Утверждение теоремы аналогично утверждению теоремы 5.1, [12], стр. 94.

Если $a \in U^n$, $a = (a_1, ..., a_n)$, где $a_i \in X_{k_i}$, и если $k = \max\{k_1, ..., k_n\}$, то $a_1...a_n\omega \in X_{k+1} \subset U$. Следовательно, $U$ является подалгеброй $\Omega_2$-алгебры $M$.

Если $m \in U$, то $m \in X_k$ для некоторого $k$. Следовательно, $f(a)(m) \in X_{k+1} \subset U$ для любого $a \in A$. Следовательно, $U$ - стабильное множество представления $f$.

Так как $U$ - подалгеброй $\Omega_2$-алгебры $M$ и стабильное множество представления $f$, то определено подпредставление $f_U$. Следовательно,

$$J(f, X) \subset U \tag{2.4.2}$$

Из (2.4.1), (2.4.2), следует $J(f, X) = U$. □

**Определение 2.4.5.** $J(f, X)$ называется **подпредставлением, порождённым множеством** $X$, а $X$ - **множеством образующих подпредставления** $J(f, X)$. В частности, **множеством образующих представления** $f$ будет такое подмножество $X \subset M$, что $J(f, X) = M$. □

Нетрудно видеть, что определение множества образующих представления не зависит от того, эффективно представление или нет. Поэтому в дальнейшем мы будем предполагать, что представление эффективно и будем опираться на соглашение для эффективного $T\star$-представления в замечании 2.1.9. Мы также будем пользоваться записью

$$R \circ m = R(m)$$

для образа $m \in M$ при эндоморфизме эффективного представления. Согласно определению произведения отображений, для любых эндоморфизмов $R$, $S$ верно равенство

$$(R \circ S) \circ m = R \circ (S \circ m) \tag{2.4.3}$$

Равенство (2.4.3) является законом ассоциативности для $\circ$ и позволяет записать выражение

$$R \circ S \circ m$$

без использования скобок.

Из теоремы 2.4.4 следует следующее определение.

**Определение 2.4.6.** Пусть $X \subset M$. Для любого $x \in J(f, X)$ существует $\Omega_2$-слово, определённое согласно следующему правилам.

(1) Если $m \in X$, то $m$ - $\Omega_2$-слово.
(2) Если $m_1$, ..., $m_n$ - $\Omega_2$-слова и $\omega \in \Omega_2(n)$, то $m_1...m_n\omega$ - $\Omega_2$-слово.
(3) Если $m$ - $\Omega_2$-слово и $a \in A$, то $am$ - $\Omega_2$-слово.

$\Omega_2$**-слово** $w(f, X, m)$ представляет данный элемент $m \in J(f, X)$. Мы будем отождествлять элемент $m \in J(f, X)$ и соответствующее ему $\Omega_2$-слово, выражая это равенством

$$m = w(f, X, m)$$

Аналогично, для произвольного множества $B \subset J(f, X)$ рассмотрим множество $\Omega_2$-слов[2.7]

$$w(f, X, B) = \{w(f, X, m) : m \in B\}$$

[2.7] Выражение $w(f, X, m)$ является частным случаем выражения $w(f, X, B)$, а именно

$$w(f, X, \{m\}) = \{w(f, X, m)\}$$

Мы будем также пользоваться записью

$$w(f,X,B)=(w(f,X,m),m\in B)$$

Обозначим $w(f,X)$ **множество $\Omega_2$-слов представления** $J(f,X)$. □

**Теорема 2.4.7.** *Эндоморфизм $R$ представления*

$$f:A\rightarrow {}^*M$$

*порождает отображение $\Omega_2$-слов*

$$w[f,X,R]:w(f,X)\rightarrow w(f,X')\quad X\subset M\quad X'=R\circ X$$

*такое, что*

(1) *Если* $m\in X$, $m'=R\circ m$, *то*

$$w[f,X,R](m)=m'$$

(2) *Если*

$$m_1,...,m_n\in w(f,X)$$

$$m_1'=w[f,X,R](m_1)\quad ...\quad m_n'=w[f,X,R](m_n)$$

*то для операции $\omega\in\Omega_2(n)$ справедливо*

$$w[f,X,R](m_1...m_n\omega)=m_1'...m_n'\omega$$

(3) *Если*

$$m\in w(f,X)\quad m'=w[f,X,R](m)\quad a\in A$$

*то*

$$w[f,X,R](am)=am'$$

Доказательство. Утверждения (1), (2) теоремы справедливы в силу определения эндоморфизма $R$. Утверждение (3) теоремы следует из равенства (2.2.4), так как $r=\mathrm{id}$. □

**Замечание 2.4.8.** Пусть $R$ - эндоморфизм представления $f$. Пусть

$$m\in J(f,X)\quad m'=R\circ m\quad X'=R\circ X$$

Теорема 2.4.7 утверждает, что $m'\in J_f(X')$ . Теорема 2.4.7 также утверждает, что $\Omega_2$-слово, представляющее $m$, относительно $X$ и $\Omega_2$-слово, представляющее $m'$, относительно $X'$ формируются согласно одному и тому же алгоритму. Это позволяет рассматривать множество $\Omega_2$-слов $w(f,X',m')$ как отображение

$$W(f,X,m):X'\rightarrow w(f,X',m')$$

$$W(f,X,m)(X')=W(f,X,m)\circ X'$$

такое, что, если для некоторого эндоморфизма $R$

$$X'=R\circ X\quad m'=R\circ m$$

то

$$W(f,X,m)\circ X'=w(f,X',m')=m'$$

Отображение $W(f,X,m)$ называется **координатами элемента** $m$ **относительно множества** $X$. Аналогично, мы можем рассмотреть координаты множества $B\subset J(f,X)$ относительно множества $X$

$$W(f,X,B)=\{W(f,X,m):m\in B\}=(W(f,X,m),m\in B)$$

Обозначим $W(f,X)$ **множество координат представления** $J(f,X)$. □

**Теорема 2.4.9.** *На множестве координат* $W(f,X)$ *определена структура* $\Omega_2$*-алгебры.*

Доказательство. Пусть $\omega\in\Omega_2(n)$. Тогда для любых $m_1$, ..., $m_n\in J(f,X)$ положим

$$W(f,X,m_1)...W(f,X,m_n)\omega=W(f,X,m_1...m_n\omega)\tag{2.4.4}$$

Согласно замечанию 2.4.8, из равенства (2.4.4) следует

$$\begin{aligned}(W(f,X,m_1)...W(f,X,m_n)\omega)\circ X&=W(f,X,m_1...m_n\omega)\circ X\\&=w(f,X,m_1...m_n\omega)\end{aligned}\tag{2.4.5}$$

Согласно правилу (2) определения 2.4.6, из равенства (2.4.5) следует

$$\begin{aligned}(W(f,X,m_1)...W(f,X,m_n)\omega)\circ X&=w(f,X,m_1)...w(f,X,m_n)\omega\\&=(W(f,X,m_1)\circ X)...(W(f,X,m_n)\circ X)\omega\end{aligned}\tag{2.4.6}$$

Из равенства (2.4.6) следует корректность определения (2.4.4) операции $\omega$ на множестве координат. □

**Теорема 2.4.10.** *Определено представление* $\Omega_1$*-алгебры* $A$ *в* $\Omega_2$*-алгебре* $W(f,X)$.

Доказательство. Пусть $a\in A$. Тогда для любого $m\in J(f,X)$ положим

$$aW(f,X,m)=W(f,X,am)\tag{2.4.7}$$

Согласно замечанию 2.4.8, из равенства (2.4.7) следует

$$(aW(f,X,m))\circ X=W(f,X,am)\circ X=w(f,X,am)\tag{2.4.8}$$

Согласно правилу (3) определения 2.4.6, из равенства (2.4.8) следует

$$(aW(f,X,m))\circ X=aw(f,X,m)=a(W(f,X,m)\circ X)\tag{2.4.9}$$

Из равенства (2.4.9) следует корректность определения (2.4.7) представления $\Omega_1$-алгебры $A$ в $\Omega_2$-алгебре $W(f,X)$. □

**Определение 2.4.11.** Пусть $X$ - множество образующих представления $f$. Пусть $R$ - эндоморфизм представления $f$. Множество координат $W(f,X,R\circ X)$ называется **координатами эндоморфизма представления**. □

**Определение 2.4.12.** Пусть $X$ - множество образующих представления $f$. Пусть $R$ - эндоморфизм представления $f$. Пусть $m\in M$. Мы определим **суперпозицию координат** эндоморфизма представления $f$ и элемента $m$ как координаты, определённые согласно правилу

$$W(f,X,m)\circ W(f,X,R\circ X)=W(f,X,R\circ m)\tag{2.4.10}$$

Пусть $Y\subset M$. Мы определим суперпозицию координат эндоморфизма представления $f$ и множества $Y$ согласно правилу

$$W(f,X,Y)\circ W(f,X,R\circ X)=(W(f,X,m)\circ W(f,X,R\circ X),m\in Y)\tag{2.4.11}$$

□

**Теорема 2.4.13.** *Эндоморфизм $R$ представления*

$$f : A \to {}^*M$$

*порождает отображение координат представления*

$$W[f,X,R]: W(f,X) \to W(f,X) \tag{2.4.12}$$

*такое, что*

$$\begin{aligned} W(f,X,m) \to W[f,X,R] \star W(f,X,m) = W(f,X,R\circ m) \\ W[f,X,R] \star W(f,X,m) = W(f,X,m) \circ W(f,X,R\circ X) \end{aligned} \tag{2.4.13}$$

Доказательство. Согласно замечанию 2.4.8, мы можем рассматривать равенства (2.4.10), (2.4.12) относительно заданного множества образующих $X$. При этом координатам $W(f,X,m)$ соответствует слово

$$W(f,X,m)\circ X = w(f,X,m) \tag{2.4.14}$$

а координатам $W(f,X,R\circ m)$ соответствует слово

$$W(f,X,R\circ m)\circ X = w(f,X,R\circ m) \tag{2.4.15}$$

Поэтому для того, чтобы доказать теорему, нам достаточно показать, что отображению $W[f,X,R]$ соответствует отображение $w[f,X,R]$. Мы будем доказывать это утверждение индукцией по сложности $\Omega_2$-слова.

Если $m \in X$, $m' = R\circ m$, то, согласно равенствам (2.4.14), (2.4.15), отображения $W[f,X,R]$ и $w[f,X,R]$ согласованы.

Пусть для $m_1$, ..., $m_n \in X$ отображения $W[f,X,R]$ и $w[f,X,R]$ согласованы. Пусть $\omega \in \Omega_2(n)$. Согласно теореме 2.4.9

$$W(f,X,m_1...m_n\omega) = W(f,X,m_1)...W(f,X,m_n)\omega \tag{2.4.16}$$

Так как $R$ - эндоморфизм $\Omega_2$-алгебры $M$, то из равенства (2.4.16) следует

$$\begin{aligned} W(f,X,R\circ(m_1...m_n\omega)) &= W(f,X,(R\circ m_1)...(R\circ m_n)\omega) \\ &= W(f,X,R\circ m_1)...W(f,X,R\circ m_n)\omega \end{aligned} \tag{2.4.17}$$

Из равенств (2.4.16), (2.4.17) и предположения индукции следует, что отображения $W[f,X,R]$ и $w[f,X,R]$ согласованы для $m = m_1...m_n\omega$.

Пусть для $m_1 \in M$ отображения $W[f,X,R]$ и $w[f,X,R]$ согласованы. Пусть $a \in A$. Согласно теореме 2.4.10

$$W(f,X,am_1) = aW(f,X,m_1) \tag{2.4.18}$$

Так как $R$ - эндоморфизм представления $f$, то из равенства (2.4.18) следует

$$W(f,X,R\circ(am_1)) = W(f,X,a\,R\circ m_1) = aW(f,X,R\circ m_1) \tag{2.4.19}$$

Из равенств (2.4.18), (2.4.19) и предположения индукции следует, что отображения $W[f,X,R]$ и $w[f,X,R]$ согласованы для $m = am_1$. □

**Следствие 2.4.14.** *Пусть $X$ - множество образующих представления $f$. Пусть $R$ - эндоморфизм представления $f$. Отображение $W[f,X,R]$ является эндоморфизмом представления $\Omega_1$-алгебры $A$ в $\Omega_2$-алгебре $W(f,X)$.* □

В дальнейшем мы будем отождествлять отображение $W[f,X,R]$ и множество координат $W(f,X,R\circ X)$.

**Теорема 2.4.15.** *Пусть $X$ - множество образующих представления $f$. Пусть $R$ - эндоморфизм представления $f$. Пусть $Y\subset M$. Тогда*

$$W(f,X,Y)\circ W(f,X,R\circ X)=W(f,X,R\circ Y) \tag{2.4.20}$$

$$W[f,X,R]\star W(f,X,Y)=W(f,X,R\circ Y) \tag{2.4.21}$$

Доказательство. Равенство (2.4.20) является следствием равенства

$$R\circ Y=(R\circ m,m\in Y)$$

а также равенств (2.4.10), (2.4.11). Равенство (2.4.21) является следствием равенств (2.4.20), (2.4.13). □

**Теорема 2.4.16.** *Пусть $X$ - множество образующих представления $f$. Пусть $R$, $S$ - эндоморфизмы представления $f$. Тогда*

$$W(f,X,S\circ X)\circ W(f,X,R\circ X)=W(f,X,R\circ S\circ X) \tag{2.4.22}$$

$$W[f,X,R]\star W[f,X,S]=W[f,X,R\circ S] \tag{2.4.23}$$

Доказательство. Равенство (2.4.22) следует из равенства (2.4.20), если положить $Y=S\circ X$. Равенство (2.4.23) следует из равенства (2.4.22) и цепочки равенств

$$\begin{aligned}
&(W[f,X,R]\star W[f,X,S])\star W(f,X,Y)\\
=&W[f,X,R]\star(W[f,X,S]\star W(f,X,Y))\\
=&(W(f,X,Y)\circ W(f,X,S\circ X))\circ W(f,X,R\circ X)\\
=&W(f,X,Y)\circ(W(f,X,S\circ X)\circ W(f,X,R\circ X))\\
=&W(f,X,Y)\circ W(f,X,R\circ S\circ X)\\
=&W(f,X,R\circ S)\star W(f,X,Y)
\end{aligned}$$

□

Концепция суперпозиции координат очень проста и напоминает своеобразную машину Тюринга. Если элемент $m\in M$ имеет вид

$$m=m_1...m_n\omega$$

или

$$m=am_1$$

то мы ищем координаты элементов $m_i$, для того чтобы подставить их в соответствующее выражение. Как только элемент $m\in M$ принадлежит множеству образующих $\Omega_2$-алгебры $M$, мы выбираем координаты соответствующего элемента из второго множителя. Поэтому мы требуем, чтобы второй множитель в суперпозиции был множеством координат образа множества образующих $X$.

Мы можем обобщить определение суперпозиции координат и предположить, что один из множителей является множеством $\Omega_2$-слов. Соответственно, определение суперпозиции координат имеет вид

$$w(f,X,Y)\circ W(f,X,R\circ X)=W(f,X,Y)\circ w(f,X,R\circ X)=w(f,X,R\circ Y)$$

Следующие формы записи образа множества $Y$ при эндоморфизме $R$ эквивалентны.

$$R\circ Y=W(f,X,Y)\circ(R\circ X)=W(f,X,Y)\circ(W(f,X,R\circ X)\circ X) \tag{2.4.24}$$

Из равенств (2.4.20), (2.4.24) следует, что

$$(2.4.25)\quad (W(f,X,Y)\circ W(f,X,R\circ X))\circ X = W(f,X,Y)\circ(W(f,X,R\circ X)\circ X)$$

Равенство (2.4.25) является законом ассоциативности для операции композиции и позволяет записать выражение

$$W(f,X,Y)\circ W(f,X,R\circ X)\circ X$$

без использования скобок.

Рассмотрим равенство (2.4.22) где мы видим изменение порядка эндоморфизмов в суперпозиции координат. Это равенство также следует из цепочки равенств, где мы можем непосредственно видеть, когда изменяется порядок эндоморфизмов

$$(2.4.26)\quad \begin{aligned} &W(f,X,m)\circ W(f,X,R\circ S\circ X)\circ X\\ =&(R\circ S)\circ m = R\circ(S\circ m)\\ =&R\circ((W(f,X,m)\circ W(f,X,S\circ X))\circ X)\\ =&W(f,X,m)\circ W(f,X,S\circ X)\circ W(f,X,R\circ X)\circ X\end{aligned}$$

Из равенства (2.4.26) следует, что координаты эндоморфизма действуют на координаты элементы $\Omega_2$-алгебры $M$ справа.

**Определение 2.4.17.** Пусть $X\subset M$ - множество образующих представления

$$f:A\rightarrow {}^*M$$

Пусть отображение

$$H:M\rightarrow M$$

является эндоморфизмом представления $f$. Пусть множество $X'=H\circ X$ является образом множества $X$ при отображении $H$. Эндоморфизм $H$ представления $f$ называется **невырожденным на множестве образующих** $X$, если множество $X'$ является множеством образующих представления $f$. В противном случае, эндоморфизм $H$ представления $f$ называется **вырожденным на множестве образующих** $X$, □

**Определение 2.4.18.** Эндоморфизм представления $f$ называется **невырожденным**, если он невырожден на любом множестве образующих. □

**Теорема 2.4.19.** *Автоморфизм $R$ представления*

$$f:A\rightarrow {}^*M$$

*является невырожденным эндоморфизмом.*

Доказательство. Пусть $X$ - множество образующих представления $f$. Пусть $X'=R\circ X$.

Согласно теореме 2.4.7 эндоморфизм $R$ порождает отображение $\Omega_2$-слов $w[f,X,R]$.

Пусть $m'\in M$. Так как $R$ - автоморфизм, то существует $m\in M$, $R\circ m=m'$. Согласно определению 2.4.6, $w(f,X,m)$ - $\Omega_2$-слово, представляющее $m$ относительно множества образующих $X$. Согласно теореме 2.4.7, $w(f,X',m')$ - $\Omega_2$-слово, представляющее $m'$ относительно множества образующих $X'$

$$w(f,X',m')=w[f,X,R](w(f,X,m))$$

Следовательно, $X'$ - множество образующих представления $f$. Согласно определению 2.4.18, автоморфизм $R$ - невырожден. □

## 2.5. Базис представления

**Определение 2.5.1.** Если множество $X \subset M$ является множеством образующих представления $f$, то любое множество $Y$, $X \subset Y \subset M$ также является множеством образующих представления $f$. Если существует минимальное множество $X$, порождающее представление $f$, то такое множество $X$ называется **базисом представления** $f$. □

**Теорема 2.5.2.** *Множество образующих $X$ представления $f$ является базисом тогда и только тогда, когда для любого $m \in X$ множество $X \setminus \{m\}$ не является множеством образующих представления $f$.*

Доказательство. Пусть $X$ - множество образующих расслоения $f$. Допустим для некоторого $m \in X$ существует $\Omega_2$-слово

$$w = w(f, X \setminus \{m\}, m) \tag{2.5.1}$$

Рассмотрим элемент $m'$, для которого $\Omega_2$-слово

$$w' = w(f, X, m') \tag{2.5.2}$$

зависит от $m$. Согласно определению 2.4.6, любое вхождение $m$ в $\Omega_2$-слово $w'$ может быть заменено $\Omega_2$-словом $w$. Следовательно, $\Omega_2$-слово $w'$ не зависит от $m$, а множество $X \setminus \{m\}$ является множеством образующих представления $f$. Следовательно, $X$ не является базисом расслоения $f$. □

**Замечание 2.5.3.** Доказательство теоремы 2.5.2 даёт нам эффективный метод построения базиса представления $f$. Выбрав произвольное множество образующих, мы шаг за шагом исключаем те элементы множества, которые имеют координаты относительно остальных элементов множества. Если множество образующих представления бесконечно, то рассмотренная операция может не иметь последнего шага. Если представление имеет конечное множество образующих, то за конечное число шагов мы можем построить базис этого представления.

Как отметил Кон в [12], стр. 96, 97, представление может иметь неэквивалентные базисы. Например, циклическая группа шестого порядка имеет базисы $\{a\}$ и $\{a^2, a^3\}$, которые нельзя отобразить один в другой эндоморфизмом представления. □

**Замечание 2.5.4.** Мы будем записывать базис также в виде

$$X = (x, x \in X)$$

Если базис - конечный, то мы будем также пользоваться записью

$$X = (x_i, i \in I) = (x_1, ..., x_n)$$

□

**Замечание 2.5.5.** Мы ввели $\Omega_2$-слово элемента $x \in M$ относительно множества образующих $X$ в определении 2.4.6. Из теоремы 2.5.2 следует, что если множество образующих $X$ не является базисом, то выбор $\Omega_2$-слова относительно множества образующих $X$ неоднозначен. Но даже если множество образующих $X$ является базисом, то представление $m \in M$ в виде $\Omega_2$-слова неоднозначно. Если $m_1$, ..., $m_n$ - $\Omega_2$-слова, $\omega\in\Omega_2(n)$ и $a\in A$, то[2.8]

$$a(m_1...m_n\omega)=(am_1)...(am_n)\omega \tag{2.5.3}$$

Возможны равенства, связанные со спецификой представления. Например, если $\omega$ является операцией $\Omega_1$-алгебры $A$ и операцией $\Omega_2$-алгебры $M$, то мы можем потребовать, что $\Omega_2$-слова $a_1...a_n\omega x$ и $a_1x...a_nx\omega$ описывают один и тот же элемент $\Omega_2$-алгебры $M$.[2.9]

Помимо перечисленных равенств в $\Omega_2$-алгебре могут существовать соотношения вида

$$w_1(f,X,m)=w_2(f,X,m)\quad m\notin X \tag{2.5.4}$$

Особенностью равенства (2.5.4) является то, что это равенство не может быть приведенно к более простому виду.

На множестве $\Omega_2$-слов $w(f,X)$, перечисленные выше равенства определяют **отношение эквивалентности** $\rho(f)$**, порождённое представлением** $f$. Очевидно, что для произвольного $m\in M$ выбор соответствующего $\Omega_2$-слова однозначен с точностью до отношения эквивалентности $\rho(f)$. Тем не менее, если в процессе построений мы получим равенство двух $\Omega_2$-слов относительно заданного базиса, мы можем утверждать, что эти $\Omega_2$-слова равны, не заботясь об отношении эквивалентности $\rho(f)$.

Аналогичное замечание касается отображения $W(f,X,m)$, определённого в замечании 2.4.8.[2.10] □

**Теорема 2.5.6.** *Автоморфизм представления $f$ отображает базис представления $f$ в базис.*

Доказательство. Пусть отображение $R$ - автоморфизм представления $f$. Пусть множество $X$ - базис представления $f$. Пусть $X'=R\circ X$.

---

[2.8] Например, пусть $\{\overline{e}_1,\overline{e}_2\}$ - базис векторного пространства над полем $k$. Равенство (2.5.3) принимает форму закона дистрибутивности

$$a(b^1\overline{e}_1+b^2\overline{e}_2)=(ab^1)\overline{e}_1+(ab^2)\overline{e}_2$$

[2.9] Для векторного пространства это требование принимает форму закона дистрибутивности

$$(a+b)\overline{e}_1=a\overline{e}_1+b\overline{e}_1$$

[2.10] Если базис векторного пространства - конечен, то мы можем представить базис в виде матрицы строки

$$\overline{\overline{e}}=\begin{pmatrix}\overline{e}_1 & ... & \overline{e}_2\end{pmatrix}$$

Мы можем представить отображение $W(f,\overline{\overline{e}},\overline{v})$ в виде матрицы столбца

$$W(f,\overline{\overline{e}},\overline{v})=\begin{pmatrix}v^1\\ ...\\ v^n\end{pmatrix}$$

Тогда

$$W(f,\overline{\overline{e}},\overline{v})(\overline{\overline{e}}')=W(f,\overline{\overline{e}},\overline{v})\begin{pmatrix}\overline{e}'_1 & ... & \overline{e}'_n\end{pmatrix}=\begin{pmatrix}v^1\\ ...\\ v^n\end{pmatrix}\begin{pmatrix}\overline{e}'_1 & ... & \overline{e}'_n\end{pmatrix}$$

имеет вид произведения матриц.

Допустим множество $X'$ не является базисом. Согласно теореме 2.5.2 существует $m' \in X'$ такое, что $X' \setminus \{x'\}$ является множеством образующих представления $f$. Согласно теореме 2.3.5 отображение $R^{-1}$ является автоморфизмом представления $f$. Согласно теореме 2.4.19 и определению 2.4.18, множество $X \setminus \{m\}$ является множеством образующих представления $f$. Полученное противоречие доказывает теорему. □

Глава 3

# Башня представлений универсальных алгебр

## 3.1. Башня представлений универсальных алгебр

**Определение 3.1.1.** Рассмотрим[3.1] множество $\Omega_k$-алгебр $A_k$, $k = 1, ..., n$. Положим $\overline{A} = (A_1, ..., A_n)$. Положим $\overline{f} = (f_{1,2}, ..., f_{n-1,n})$. Множество представлений $f_{k,k+1}$, $k = 1, ..., n$, $\Omega_k$-алгебры $A_k$ в $\Omega_{k+1}$-алгебре $A_{k+1}$ называется **башней $(\overline{A}, \overline{f})$ представлений $\overline{\Omega}$-алгебр**. □

Мы можем проиллюстрировать определение 3.1.1 с помощью диаграммы

(3.1.1)

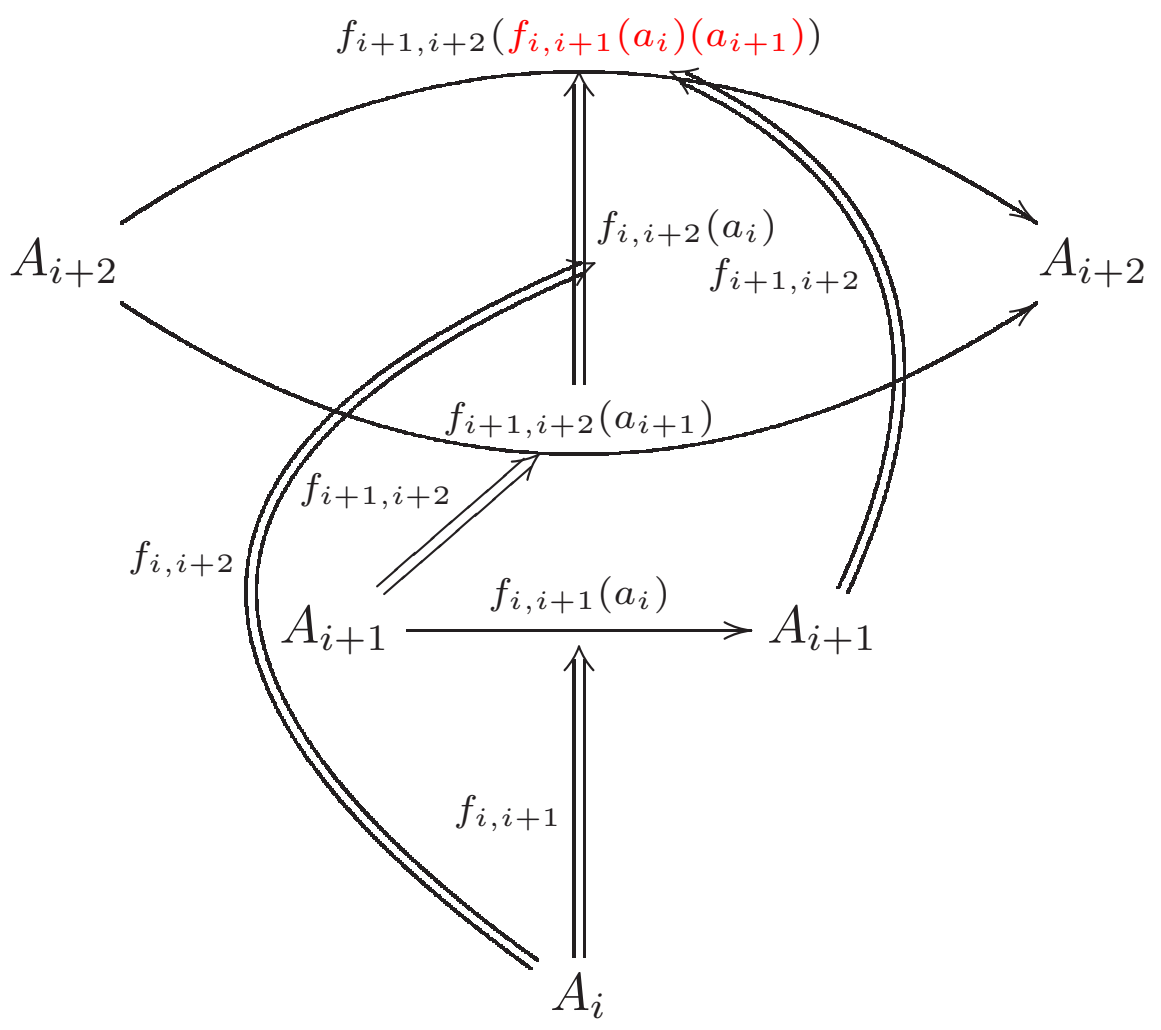

$f_{i,i+1}$ - представление $\Omega_i$-алгебры $A_i$ в $\Omega_{i+1}$-алгебре $A_{i+1}$. $f_{i+1,i+2}$ - представление $\Omega_{i+1}$-алгебры $A_{i+1}$ в $\Omega_{i+2}$-алгебре $A_{i+2}$.

**Теорема 3.1.2.** *Отображение*

$$f_{i,i+2} : A_i \to {}^{**}A_{i+2}$$

*является представлением $\Omega_i$-алгебры $A_i$ в $\Omega_{i+1}$-алгебре ${}^*A_{i+2}$.*

Доказательство. Произвольному $a_{i+1} \in A_{i+1}$ соответствует автоморфизм $f_{i+1,i+2}(a_{i+1}) \in {}^*A_{i+2}$. Произвольному $a_i \in A_i$ соответствует автоморфизм $f_{i,i+1}(a_i) \in {}^*A_{i+1}$. Так как образом элемента $a_{i+1} \in A_{i+1}$ является элемент $f_{i,i+1}(a_i)(a_{i+1}) \in A_{i+1}$, то тем самым элемент $a_i \in A_i$ порождает преобразование $\Omega_{i+1}$-алгебры ${}^*A_{i+2}$, определённое равенством

(3.1.2)
$$f_{i,i+2}(a_i)(f_{i+1,i+2}(a_{i+1})) = f_{i+1,i+2}(f_{i,i+1}(a_i)(a_{i+1}))$$

[3.1] Определения и теоремы из этой главы включены также в главы [10]-4, [10]-5.

Пусть $\omega$ - $n$-арная операция $\Omega_i$-алгебры. Так как $f_{i,i+1}$ - гомоморфизм $\Omega_i$-алгебры, то

$$f_{i,i+1}(a_{i,1}...a_{i,n}\omega) = f_{i,i+1}(a_{i,1})...f_{i,i+1}(a_{i,n})\omega \tag{3.1.3}$$

Для $a_1$, ..., $a_n \in A_i$ мы определим операцию $\omega$ на множестве ${}^{**}A_{i+2}$ с помощью равенства

$$\begin{aligned}&f_{i,i+2}(a_{i,1})(f_{i+1,i+2}(a_{i+1}))...f_{i,i+2}(a_{i,n})(f_{i+1,i+2}(a_{i+1}))\omega\\ =&f_{i+1,i+2}(f_{i,i+1}(a_{i,1})(a_{i+1}))...f_{i+1,i+2}(f_{i,i+1}(a_{i,n})(a_{i+1}))\omega\\ =&f_{i+1,i+2}(f_{i,i+1}(a_{i,1})...f_{i,i+1}(a_{i,n})\omega)(a_{i+1})\end{aligned} \tag{3.1.4}$$

Первое равенство следует из равенства (3.1.2). Второе равенство постулировано на основе требования, что отображение $f_{i+1,i+2}$ является гомоморфизмом $\Omega_i$-алгебры. Следовательно, мы определили структуру $\Omega_i$-алгебры на множестве ${}^{**}A_{i+2}$.

Из (3.1.4) и (3.1.3) следует

$$\begin{aligned}&f_{i,i+2}(a_{i,1})(f_{i+1,i+2}(a_{i+1}))...f_{i,i+2}(a_{i,n})(f_{i+1,i+2}(a_{i+1}))\omega\\ =&f_{i+1,i+2}(f_{i,i+1}(a_{i,1}...a_{i,n}\omega)(a_{i+1}))\\ =&f_{i,i+2}(a_{i,1}...a_{i,n}\omega)(f_{i+1,i+2}(a_{i+1}))\end{aligned} \tag{3.1.5}$$

Из равенства (3.1.5) следует

$$f_{i,i+2}(a_{i,1})...f_{i,i+2}(a_{i,n})\omega = f_{i,i+2}(a_{i,1}...a_{i,n}\omega)$$

Следовательно, отображение $f_{i,i+2}$ является гомоморфизмом $\Omega_i$-алгебры. Следовательно, отображение $f_{i,i+2}$ является $T\star$-представлением $\Omega_i$-алгебры $A_i$ в $\Omega_{i+1}$-алгебре ${}^*A_{i+2}$. □

**Теорема 3.1.3.** *$(id, f_{i+1,i+2})$ является морфизмом $T\star$-представлений $\Omega_i$-алгебры из $f_{i,i+1}$ в $f_{i,i+2}$.*

Доказательство. Рассмотрим более детально диаграмму (3.1.1).

(3.1.6)

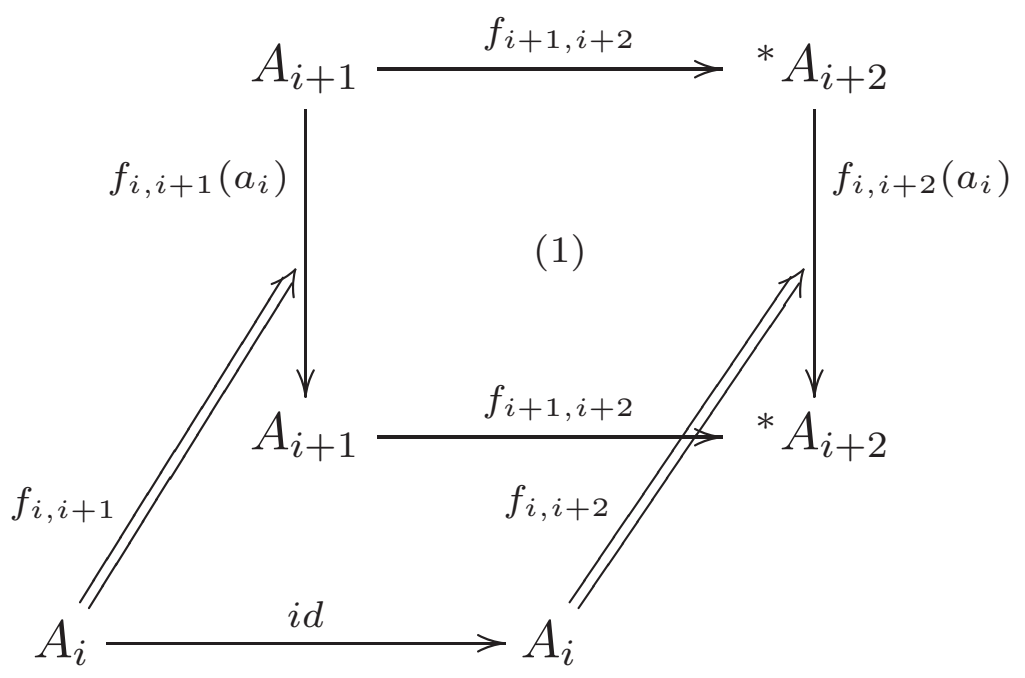

Утверждение теоремы следует из равенства (3.1.2) и определения 2.2.2. □

**Теорема 3.1.4.** *Рассмотрим башню $(\overline{A}, \overline{f})$ представлений $\overline{\Omega}$-алгебр. Если тождественное преобразование*

$$\delta_{i+2} : A_{i+2} \to A_{i+2}$$

*$\Omega_{i+2}$-алгебры $A_{i+2}$ принадлежит представлению $f_{i+1,i+2}$, то представление $f_{i,i+2}$ $\Omega_i$-алгебры $A_i$ в $\Omega_{i+1}$-алгебре ${}^*A_{i+2}$ может быть продолжено до представления*

$$f'_{i,i+2}:A_i\rightarrow {}^*A_{i+2}$$

*$\Omega_i$-алгебры $A_i$ в $\Omega_{i+2}$-алгебре $A_{i+2}$.*

Доказательство. Согласно равенству (3.1.2), для заданного $a_{i+1}\in A_{i+1}$ мы можем определить отображение

$$\begin{gathered}f'_{i,i+2}:A_i\rightarrow {}^*A_{i+2}\\ f'_{i,i+2}(a_i)(f_{i+1,i+2}(a_{i+1})(a_{i+2}))=f_{i+1,i+2}(f_{i,i+1}(a_i)(a_{i+1}))(a_{i+2})\end{gathered}\tag{3.1.7}$$

Пусть существует $a_{i+1}\in A_{i+1}$ такой, что

$$f_{i+1,i+2}(a_{i+1})=\delta_{i+2}\tag{3.1.8}$$

Тогда из (3.1.7),(3.1.8) следует

$$f'_{i,i+2}(a_i)(a_{i+2})=f_{i+1,i+2}(f_{i,i+1}(a_i)(a_{i+1}))(a_{i+2})\tag{3.1.9}$$

Пусть $\omega$ - $n$-арная операция $\Omega_i$-алгебры $A_i$. Поскольку $f_{i,i+2}$ - гомоморфизм $\Omega_i$ -алгебры, то из равенства (3.1.2) следует

$$\begin{aligned}&f_{i,i+2}(a_{i\cdot 1}...a_{i\cdot n}\omega)(f_{i+1,i+2}(a_{i+1}))\\ =&f_{i,i+2}(a_{i\cdot 1})(f_{i+1,i+2}(a_{i+1}))...f_{i,i+2}(a_{i\cdot n})(f_{i+1,i+2}(a_{i+1}))\omega\\ =&f_{i+1,i+2}(f_{i,i+1}(a_{i\cdot 1})(a_{i+1})...f_{i,i+1}(a_{i\cdot n})(a_{i+1})\omega)\end{aligned}\tag{3.1.10}$$

Из равенства (3.1.10) следует

(3.1.11)

$$\begin{aligned}&f_{i,i+2}(a_{i\cdot 1}...a_{i\cdot n}\omega)(f_{i+1,i+2}(a_{i+1})(a_{i+2}))\\ =&f_{i,i+2}(a_{i\cdot 1})(f_{i+1,i+2}(a_{i+1})(a_{i+2}))...f_{i,i+2}(a_{i\cdot n})(f_{i+1,i+2}(a_{i+1})(a_{i+2}))\omega\\ =&f_{i+1,i+2}(f_{i,i+1}(a_{i\cdot 1})(a_{i+1})...f_{i,i+1}(a_{i\cdot n})(a_{i+1})\omega)(a_{i+2})\end{aligned}$$

Из равенств (3.1.8), (3.1.9), (3.1.11) следует

$$\begin{aligned}&f'_{i,i+2}(a_{i\cdot 1}...a_{i\cdot n}\omega)(a_{i+2})\\ =&f'_{i,i+2}(a_{i\cdot 1})(a_{i+2})...f'_{i,i+2}(a_{i\cdot n})(a_{i+2})\omega\\ =&f_{i+1,i+2}(f_{i,i+1}(a_{i\cdot 1})(a_{i+1})...f_{i,i+1}(a_{i\cdot n})(a_{i+1})\omega)(a_{i+2})\end{aligned}\tag{3.1.12}$$

Равенство (3.1.12) определяет операцию $\omega$ на множестве ${}^*A_{i+2}$

$$f'_{i,i+2}(a_{i\cdot 1})...f'_{i,i+2}(a_{i\cdot n})\omega=f'_{i,i+2}(a_{i\cdot 1}...a_{i\cdot n}\omega)\tag{3.1.13}$$

При этом отображение $f'_{i,i+2}$ является гомоморфизмом $\Omega_i$-алгебры. Следовательно, мы построили представление $\Omega_i$-алгебры $A_i$ в $\Omega_{i+2}$-алгебре $A_{i+2}$. □

**Определение 3.1.5.** Рассмотрим башню представлений

$$((A_1,A_2,A_3),(f_{1,2},f_{2,3}))$$

Отображение $f_*=(f_{1,2},f_{1,3})$ называется **представлением $\Omega_1$-алгебры $A_1$ в представлении** $f_{2,3}$. □

**Определение 3.1.6.** Рассмотрим башню представлений $(\overline{A}, \overline{f})$. Отображение

$$f_* = (f_{1,2}, ..., f_{1,n})$$

называется **представлением $\Omega_1$-алгебры $A_1$ в башне представлений**

$$(\overline{A}_{[1]}, \overline{f}) = ((A_2, ..., A_n), (f_{2,3}, ..., f_{n-1,n}))$$

□

## 3.2. Морфизм башни $T*$-представлений

**Определение 3.2.1.** Рассмотрим множество $\Omega_k$-алгебр $A_k$, $B_k$, $k = 1, ..., n$. Множество отображений $(h_1, ..., h_n)$ называется **морфизмом из башни $T\star$-представлений $(\overline{A}, \overline{f})$ в башню $T\star$-представлений $(\overline{B}, \overline{g})$**, если для каждого $k$, $k = 1, ..., n-1$, пара отображений $(h_k, h_{k+1})$ является морфизмом $T\star$-представлений из $f_{k,k+1}$ в $g_{k,k+1}$. □

Для любого $k$, $k = 1, ..., n-1$, мы имеем диаграмму

(3.2.1)

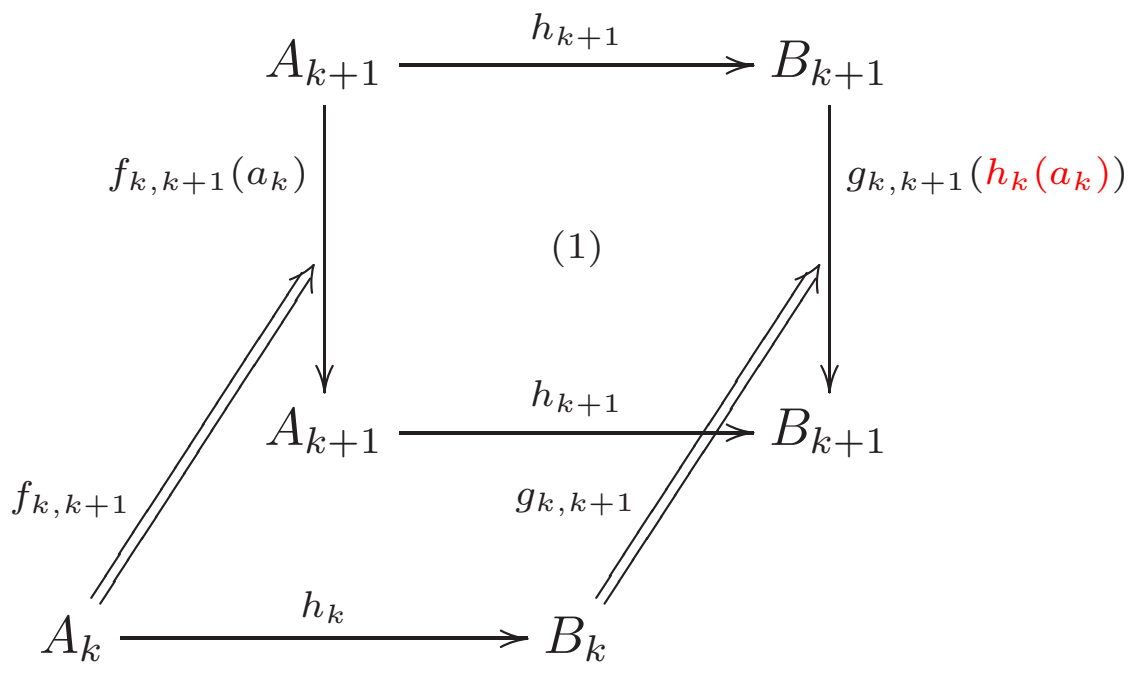

Равенства

(3.2.2) $$h_{k+1} \circ f_{k,k+1}(a_k) = g_{k,k+1}(h_k(a_k)) \circ h_{k+1}$$

(3.2.3) $$h_{k+1}(f_{k,k+1}(a_k)(a_{k+1})) = g_{k,k+1}(h_k(a_k))(h_{k+1}(a_{k+1}))$$

выражают коммутативность диаграммы (1). Однако уже для морфизма $(h_k, h_{k+1})$, $k > 1$, диаграмма (3.2.1) неполна. Учитывая аналогичную диаграмму для морфизма $(h_k, h_{k+1})$ эта диаграмма на верхнем уровне приобретает вид
(3.2.4)

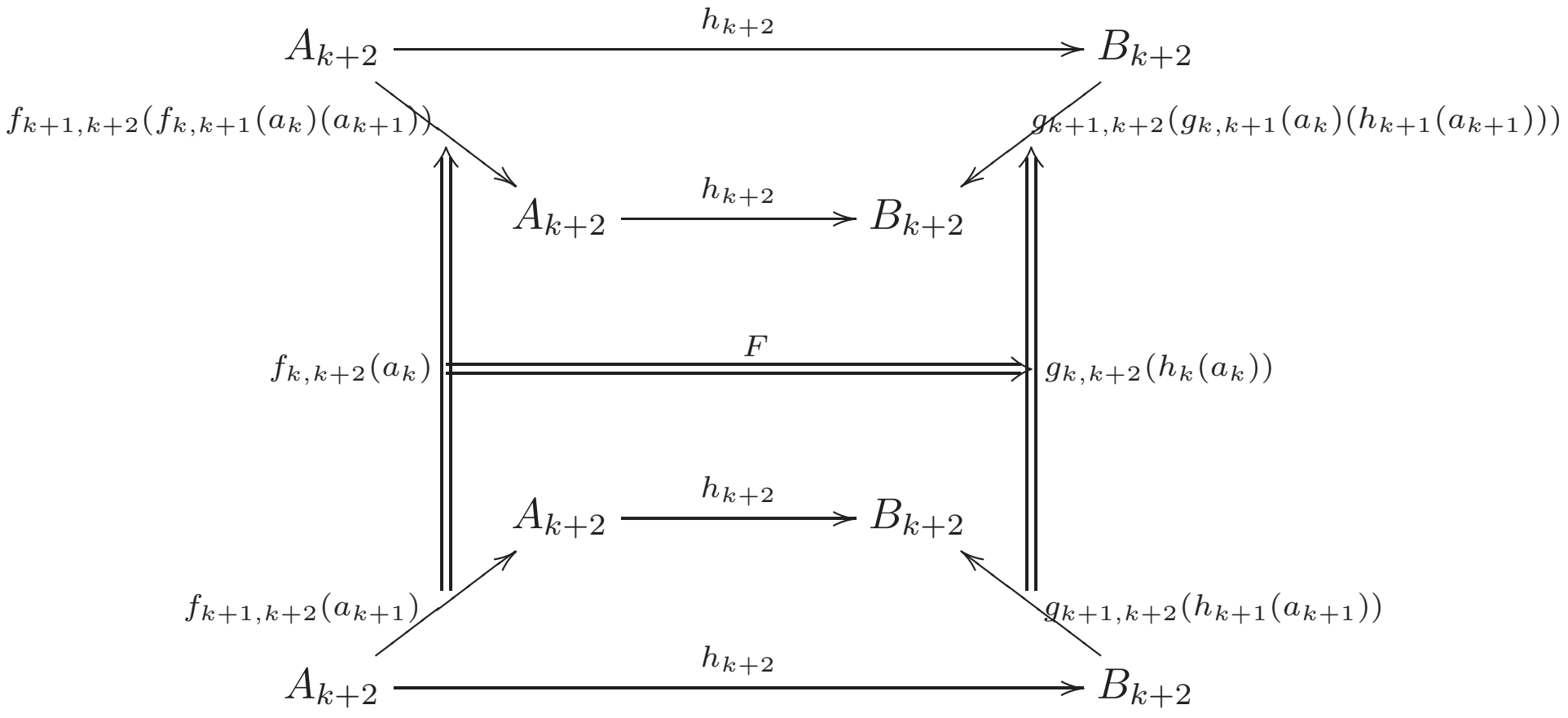

К сожалению, диаграмма (3.2.4) малоинформативна. Очевидно, что существует морфизм из ${}^*A_{k+2}$ в ${}^*B_{k+2}$, отображающий $f_{k,k+2}(a_k)$ в $g_{k,k+2}(h_k(a_k))$. Однако структура этого морфизма из диаграммы неясна. Мы должны рассмотреть отображение из ${}^*A_{k+2}$ в ${}^*B_{k+2}$, так же мы это сделали в теореме 3.1.3.

**Теорема 3.2.2.** *Если $T\star$-представления $f_{i+1,i+2}$ эффективно, то $(h_i, {}^*h_{i+2})$ является морфизмом $T\star$-представлений из $T\star$-представления $f_{i,i+2}$ в $T\star$-представление $g_{i,i+2}$ $\Omega_i$-алгебры.*

Доказательство. Рассмотрим диаграмму

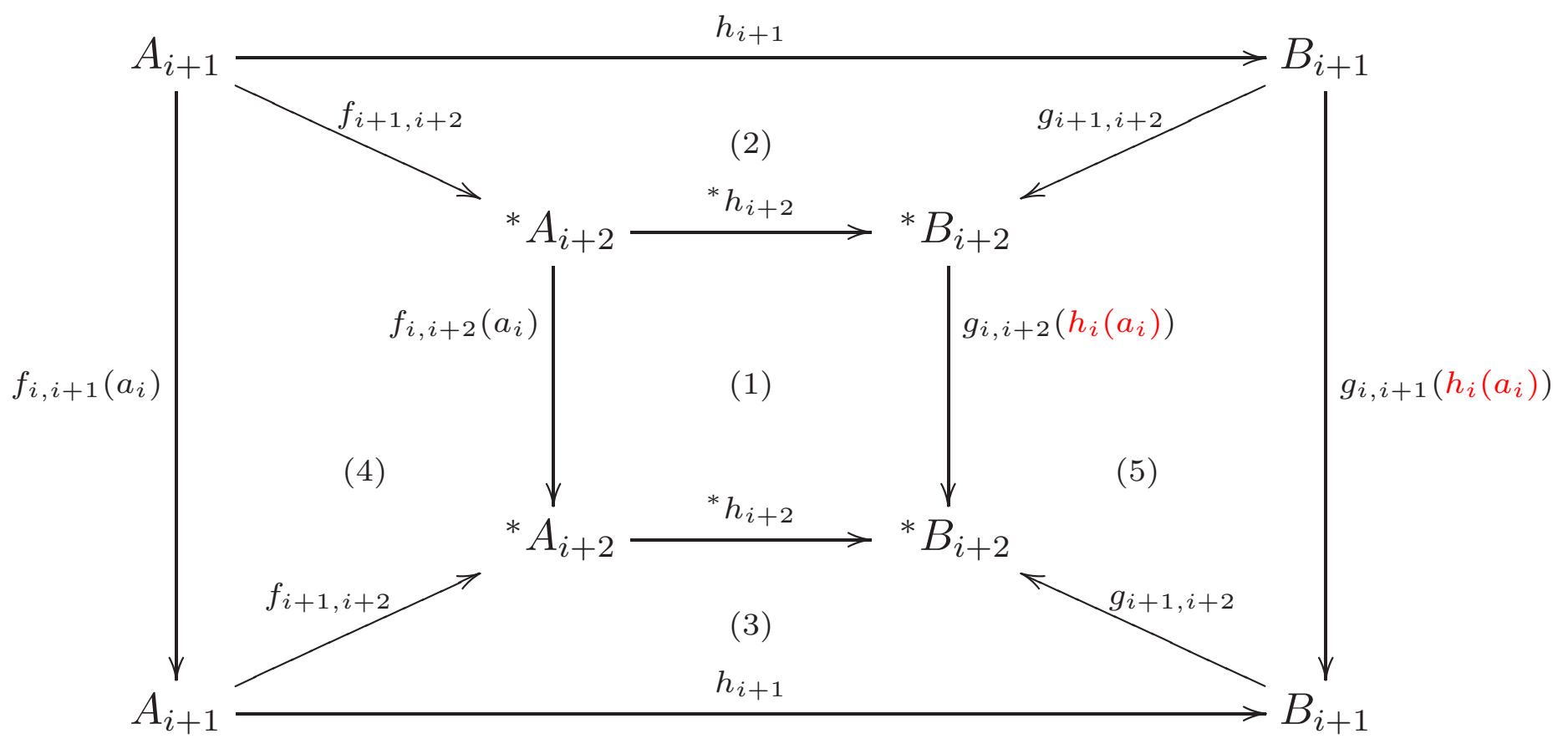

Существование отображения ${}^*h_{i+2}$ и коммутативность диаграмм (2) и (3) следует из эффективности отображения $f_{i+1,i+2}$ и теоремы 2.2.7. Коммутативность диаграмм (4) и (5) следует из теоремы 3.1.3.

Из коммутативности диаграммы (4) следует

$$f_{i+1,i+2} \circ f_{i,i+1}(a_i) = f_{i,i+2}(a_i) \circ f_{i+1,i+2} \tag{3.2.5}$$

Из равенства (3.2.5) следует

$${}^*h_{i+2} \circ f_{i+1,i+2} \circ f_{i,i+1}(a_i) = {}^*h_{i+2} \circ f_{i,i+2}(a_i) \circ f_{i+1,i+2} \tag{3.2.6}$$

Из коммутативности диаграммы (3) следует

$${}^*h_{i+2} \circ f_{i+1,i+2} = g_{i+1,i+2} \circ h_{i+1} \tag{3.2.7}$$

Из равенства (3.2.7) следует

$${}^*h_{i+2} \circ f_{i+1,i+2} \circ f_{i,i+1}(a_i) = g_{i+1,i+2} \circ h_{i+1} \circ f_{i,i+1}(a_i) \tag{3.2.8}$$

Из равенств (3.2.6) и (3.2.8) следует

$${}^*h_{i+2} \circ f_{i,i+2}(a_i) \circ f_{i+1,i+2} = g_{i+1,i+2} \circ h_{i+1} \circ f_{i,i+1}(a_i) \tag{3.2.9}$$

Из коммутативности диаграммы (5) следует

$$g_{i+1,i+2} \circ g_{i,i+1}(h_i(a_i)) = g_{i,i+2}(h_i(a_i)) \circ g_{i+1,i+2} \tag{3.2.10}$$

Из равенства (3.2.10) следует

$$g_{i+1,i+2} \circ g_{i,i+1}(h_i(a_i)) \circ h_{i+1} = g_{i,i+2}(h_i(a_i)) \circ g_{i+1,i+2} \circ h_{i+1} \tag{3.2.11}$$

Из коммутативности диаграммы (2) следует

$${}^*h_{i+2} \circ f_{i+1,i+2} = g_{i+1,i+2} \circ h_{i+1} \tag{3.2.12}$$

Из равенства (3.2.12) следует

$$g_{i,i+2}(h_i(a_i)) \circ {}^*h_{i+2} \circ f_{i+1,i+2} = g_{i,i+2}(h_i(a_i)) \circ g_{i+1,i+2} \circ h_{i+1} \tag{3.2.13}$$

Из равенств (3.2.11) и (3.2.13) следует

$$g_{i+1,i+2} \circ g_{i,i+1}(h_i(a_i)) \circ h_{i+1} = g_{i,i+2}(h_i(a_i)) \circ {}^*h_{i+2} \circ f_{i+1,i+2} \tag{3.2.14}$$

Внешняя диаграмма является диаграммой (3.2.1) при $i = 1$. Следовательно, внешняя диаграмма коммутативна

$$h_{i+1} \circ f_{i,i+1}(a_i) = g_{i,i+1}(h_i(a_i)) \circ h_{i+1} \tag{3.2.15}$$

Из равенства (3.2.15) следует

$$g_{i+1,i+2} \circ h_{i+1} \circ f_{i,i+1}(a_i) = g_{i+1,i+2} \circ g_{i,i+1}(h_i(a_i)) \circ h_{i+1}(a_{i+1}) \tag{3.2.16}$$

Из равенств (3.2.9), (3.2.14) и (3.2.16) следует

$${}^*h_{i+2} \circ f_{i,i+2}(a_i) \circ f_{i+1,i+2} = g_{i,i+2}(h_i(a_i)) \circ {}^*h_{i+2} \circ f_{i+1,i+2} \tag{3.2.17}$$

Так как отображение $f_{i+1,i+2}$ - инъекция, то из равенства (3.2.17) следует

$${}^*h_{i+2} \circ f_{i,i+2}(a_i) = g_{i,i+2}(h_i(a_i)) \circ {}^*h_{i+2} \tag{3.2.18}$$

Из равенства (3.2.18) следует коммутативность диаграммы (1), откуда следует утверждение теоремы. □

Теоремы 3.1.3 и 3.2.2 справедливы для любых уровней башни $T\star$-представлений. В каждом конкретном случае надо правильно указывать множества, откуда и куда направленно отображение. Смысл приведенных теорем состоит в том, что все отображения в башне представлений действуют согласовано.

Теорема 3.2.2 утверждает, что неизвестное отображение на диаграмме (3.2.4) является отображением ${}^*h_{i+2}$.

**Теорема 3.2.3.** *Рассмотрим множество $\Omega_i$-алгебр $A_i$, $B_i$, $C_i$, $i = 1, ..., n$. Пусть определены морфизмы башни представлений*

$$\overline{p} :(\overline{A}, \overline{f}) \to (\overline{B}, \overline{g})$$

$$\overline{q} :(\overline{B}, \overline{g}) \to (\overline{C}, \overline{h})$$

*Тогда определён морфизм представлений $\Omega$-алгебры*

$$\overline{r} : (\overline{A}, \overline{f}) \to (\overline{C}, \overline{h})$$

*где $r_k = q_k p_k$, $k = 1, ..., n$. Мы будем называть морфизм $\overline{r}$ башни представлений из $\overline{f}$ в $\overline{h}$* **произведением морфизмов $\overline{p}$ и $\overline{q}$ башни представлений**.

Доказательство. Для каждого $k$, $k = 2, ..., n$, мы можем представить утверждение теоремы, пользуясь диаграммой

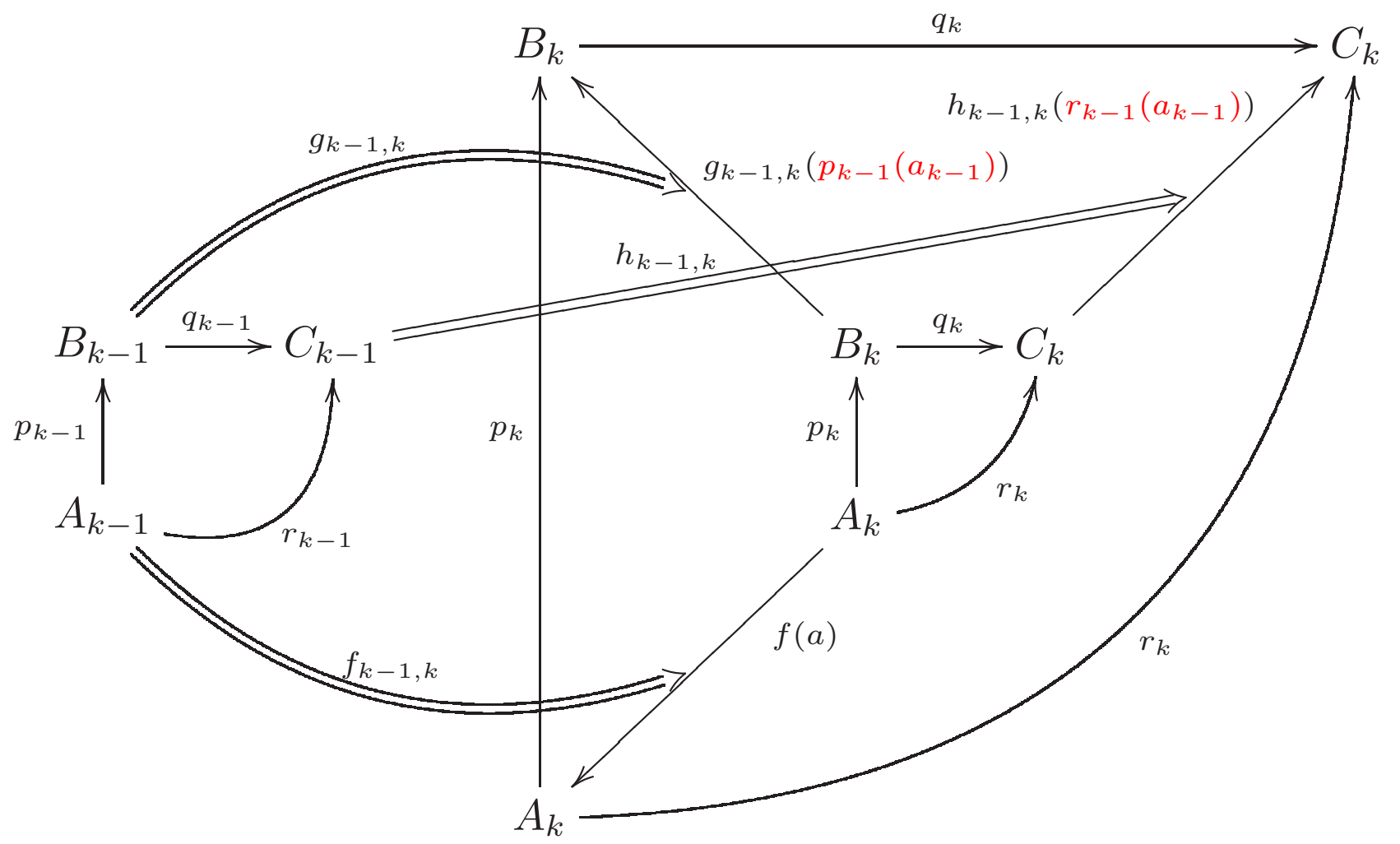

Отображение $r_{k-1}$ является гомоморфизмом $\Omega_{k-1}$-алгебры $A_{k-1}$ в $\Omega_{k-1}$-алгебру $C_{k-1}$. Нам надо показать, что пара отображений $(r_{k-1}, r_k)$ удовлетворяет (3.2.2):

$$\begin{aligned} r_k(f_{k-1,k}(a_{k-1})a_k) &= q_k p_k(f_{k-1,k}(a_{k-1})a_k) \\ &= q_k(g_{k-1,k}(p_{k-1}(a_{k-1}))p_k(a_k)) \\ &= h_{k-1,k}(q_{k-1}p_{k-1}(a_{k-1}))q_k p_k(a_k)) \\ &= h_{k-1,k}(r(a_{k-1}))r_k(a_k) \end{aligned}$$

□

## 3.3. Автоморфизм башни представлений

**Определение 3.3.1.** Пусть $(\overline{A}, \overline{f})$ - башня представлений $\overline{\Omega}$-алгебр. Морфизм башни представлений $(\mathrm{id}, h_2, ..., h_n)$ такой, что для любого $k$, $k = 2, ..., n$, $h_k$ - эндоморфизм $\Omega_k$-алгебры $A_k$, называется **эндоморфизмом башни представлений** $\overline{f}$. □

**Определение 3.3.2.** Пусть $(\overline{A}, \overline{f})$ - башня представлений $\overline{\Omega}$-алгебр. Морфизм башни представлений $(\mathrm{id}, h_2, ..., h_n)$ такой, что для любого $k$, $k = 2, ..., n$, $h_k$ - автоморфизм $\Omega_k$-алгебры $A_k$, называется **автоморфизмом башни представлений** $f$. □

**Теорема 3.3.3.** *Пусть $(\overline{A}, \overline{f})$ - башня представлений $\overline{\Omega}$-алгебр. Множество автоморфизмов представления $\overline{f}$ порождает группу.*

Доказательство. Пусть $\overline{r}_{[1]}$, $\overline{p}_{[1]}$ - автоморфизмы башни представлений $\overline{f}$. Согласно определению 3.3.2 для любого $k$, $k = 2, ..., n$, отображения $r_k$, $p_k$ являются автоморфизмами $\Omega_k$-алгебры $A_k$. Согласно теореме II.3.2 ([12], с. 60) для любого $k$, $k = 2, ..., n$, отображение $r_k p_k$ является автоморфизмом $\Omega_k$-алгебры $A_k$. Из теоремы 3.2.3 и определения 3.3.2 следует, что произведение автоморфизмов $\overline{r}_{[1]}\overline{p}_{[1]}$ башни представлений $\overline{f}$ является автоморфизмом башни представлений $\overline{f}$.

Согласно доказательству теоремы 2.3.5, для любого $k$, $k = 2, ..., n$, произведение автоморфизмов $\Omega_i$-алгебры ассоциативно. Следовательно, ассоциативно произведение автоморфизмов башни представлений.

Пусть $\overline{r}_{[1]}$ - автоморфизм башни представлений $\overline{f}$. Согласно определению 3.3.2 для любого $k$, $k = 2, ..., n$, отображение $r_k$ является автоморфизмом $\Omega_k$-алгебры $A_k$. Следовательно, для любого $k$, $k = 2, ..., n$, отображение $r_k^{-1}$ является автоморфизмом $\Omega_k$-алгебры $A_k$. Для автоморфизма $\overline{r}_{[1]}$ справедливо равенство (3.2.3). Положим $a'_k = r_k(a_k)$, $k = 2, ..., n$. Так как $r_k$, $k = 2, ..., n$, - автоморфизм, то $a_k = r_k^{-1}(a'_k)$ и равенство (3.2.3) можно записать в виде

$$h_{k+1}(f_{k,k+1}(h_k^{-1}(a'_k))(h_{k+1}(a'_{k+1}))) = g_{k,k+1}(a'_k)(a'_{k+1}) \tag{3.3.1}$$

Так как отображение $h_{k+1}$ является автоморфизмом $\Omega_{k+1}$-алгебры $A_{k+1}$, то из равенства (3.3.1) следует

$$f_{k,k+1}(h_k^{-1}(a'_k)(h_{k+1}(a'_{k+1}))) = h_{k+1}^{-1}(g_{k,k+1}(a'_k)(a'_{k+1})) \tag{3.3.2}$$

Равенство (3.3.2) соответствует равенству (3.2.3) для отображения $\overline{r}_{[1]}^{-1}$. Следовательно, отображение $\overline{r}_{[1]}^{-1}$ является автоморфизмом представления $\overline{f}$. □

## 3.4. Множество образующих башни представлений

**Определение 3.4.1.** Башня представлений $(\overline{A}, \overline{f})$. называется **башней эффективных представлений**, если для любого $i$ представление $f_{i,i+1}$ эффективно. □

**Теорема 3.4.2.** *Рассмотрим башню $T\star$-представлений $(\overline{A}, \overline{f})$. Пусть представления $f_{i,i+1}$, ..., $f_{i+k-1,i+k}$ эффективны. Тогда представление $f_{i,i+k}$ эффективно.*

Доказательство. Мы докажем утверждение теоремы по индукции.

Пусть представления $f_{i,i+1}$, $f_{i+1,i+2}$ эффективны. Предположим, что преобразование $f_{i,i+1}(a_i)$ не является тождественным преобразованием. Тогда существует $a_{i+1} \in A_{i+1}$ такой, что $f_{i,i+1}(a_i)(a_{i+1}) \neq a_{i+1}$. Так как представление $f_{i+1,i+2}$ эффективно, то преобразования $f_{i+1,i+2}(a_{i+1})$ и $f_{i+1,i+2}(f_{i,i+1}(a_i)(a_{i+1}))$ не совпадают. Согласно построению, выполненному в теореме 3.1.2, преобразование $f_{i,i+2}(a_i)$ не является тождественным преобразованием. Следовательно, представление $f_{i,i+2}$ эффективно.

Предположим утверждение теоремы верно для $k-1$ представлений и пусть $f_{i,i+1}$, ..., $f_{i+k-1,i+k}$ эффективны. Согласно предположению индукции, представления $f_{i,i+k-1}$, $f_{i+k-1,i+k}$ эффективны. Согласно доказанному выше, представление $f_{i,i+k}$ эффективно. □

**Теорема 3.4.3.** *Рассмотрим башню $(\overline{A}, \overline{f})$ представлений $\overline{\Omega}$-алгебр. Пусть тождественное преобразование*

$$\delta_{i+2} : A_{i+2} \to A_{i+2}$$

*$\Omega_{i+2}$-алгебры $A_{i+2}$ принадлежит представлению $f_{i+1,i+2}$. Пусть представления $f_{i,i+1}$, $f_{i+1,i+2}$ эффективны. Тогда представление $f'_{i,i+2}$, определённое в теореме 3.1.4, эффективно.*

Доказательство. Пусть существует $a_{i+1}\in A_{i+1}$ такой, что

$$f_{i+1,i+2}(a_{i+1})=\delta_{i+2}$$

Допустим, что представление $f'_{i,i+2}$ не является эффективным. Тогда существуют

$$a_{i\cdot 1},a_{i\cdot 2}\in A_i\quad a_{i\cdot 1}\ne a_{i\cdot 2}\tag{3.4.1}$$

такие, что

$$f'_{i,i+2}(a_{i\cdot 1})(a_{i+2})=f'_{i,i+2}(a_{i\cdot 2})(a_{i+2})\tag{3.4.2}$$

Из равенств (3.1.9), (3.4.2), следует

$$f_{i+1,i+2}(f_{i,i+1}(a_{i\cdot 1})(a_{i+1}))(a_{i+2})=f_{i+1,i+2}(f_{i,i+1}(a_{i\cdot 2})(a_{i+1}))(a_{i+2})\tag{3.4.3}$$

Так как $a_{i+2}$ произвольно, то из (3.4.3) следует

$$f_{i+1,i+2}(f_{i,i+1}(a_{i\cdot 1})(a_{i+1}))=f_{i+1,i+2}(f_{i,i+1}(a_{i\cdot 2})(a_{i+1}))\tag{3.4.4}$$

Поскольку представление $f_{i,i+1}$ эффективно, то из условия (3.4.1) следует

$$f_{i,i+1}(a_{i\cdot 1})(a_{i+1})\ne f_{i,i+1}(a_{i\cdot 2})(a_{i+1})\tag{3.4.5}$$

Из условия (3.4.5) и равенства (3.4.4) следует, что представление $f_{i+1,i+2}$ не является эффективным. Полученное противоречие доказывает утверждение теоремы. □

Мы строим базис представления башни представлений по той же схеме, что мы построили базис представления в секции 2.5.

Мы будем записывать элементы башни представлений $(\overline{A},\overline{f})$ в виде кортежа

$$\overline{a}=(a_1,...,a_3)\quad a_i\in A_i\quad i=1,...,n$$

Эту запись можно интерпретировать как

$$(a_2,a_3)=f_{2,3}(a_2)(a_3)$$

В тоже время эта запись имеет дополнительный смысл. Например, в аффинном пространстве $a_3$ - это точка пространства, $a_2$ - вектор, определяющий преобразование параллельного переноса. Таким образом, мы можем интерпретировать кортеж $(a_2,a_3)$ либо как вектор $a_2$ с началом в точке $a_3$, либо как конец этого вектора.

**Определение 3.4.4.** Пусть $(\overline{A},\overline{f})$ - башня представлений. Кортеж множеств

$$\overline{N}_{[1]}=(N_2\subset A_2,...,N_n\subset A_n)$$

называется **кортежем стабильных множеств башни представлений** $\overline{f}$, если

$$f_{i-1,i}(a_{i-1})(a_i)\in N_i\quad i=2,...,n$$

для любых $a_1\in A_1$, $a_2\in N_2$, ..., $a_n\in N_n$. Мы также будем говорить, что кортеж множеств

$$\overline{N}_{[1]}=(N_2\subset A_2,...,N_n\subset A_n)$$

стабилен относительно башни представлений $\overline{f}$. □

**Теорема 3.4.5.** *Пусть $\overline{f}$ - башня представлений. Пусть множество $N_i \subset A_i$ является подалгеброй $\Omega_i$-алгебры $A_i$, $i = 2$, ..., $n$. Пусть кортеж множеств*

$$\overline{N}_{[1]} = (N_2 \subset A_2, ..., N_n \subset A_n)$$

*стабилен относительно башни представлений $\overline{f}$. Тогда существует башня представлений*

$$((A_1, N_2, ..., N_n), (f_{N_2,1,2}, ..., f_{N_n,n-1,n})) \tag{3.4.6}$$

*такая, что*

$$f_{N_i,i-1,i}(a_{i-1}) = f_{i-1,i}(a_{i-1})|_{N_i} \quad i = 2, ..., n$$

*Башня представлений* (3.4.6) *называется* **башней подпредставлений**.

Доказательство. Пусть $\omega_{i-1,1}$ - $m$-арная операция $\Omega_{i-1}$-алгебры $A_{i-1}$, $i = 2$, ..., $n$. Тогда для любых $a_{i-1,1}$, ..., $a_{i-1,m} \in N_{i-1}$ [3.2] и любого $a_i \in N_i$

$$\begin{aligned}&(f_{N_i,i-1,i}(a_{i-1,1})...f_{N_i,i-1,i}(a_{i-1,m})\omega_{i-1,1})(a_i)\\=&(f_{i-1,i}(a_{i-1,1})...f_{i-1,i}(a_{i-1,m})\omega_{i-1,1})(a_i)\\=&f_{i-1,i}(a_{i-1,1}...a_{i-1,m}\omega_{i-1,1})(a_i)\\=&f_{N_i,i-1,i}(a_{i-1,1}...a_{i-1,m}\omega_{i-1,1})(a_i)\end{aligned}$$

Пусть $\omega_{i,2}$ - $m$-арная операция $\Omega_i$-алгебры $A_i$, $i = 2$, ..., $n$. Тогда для любых $a_{i,1}$, ..., $a_{i,n} \in N_i$ и любого $a_{i-1} \in N_{i-1}$

$$\begin{aligned}&f_{N_i,i-1,i}(a_{i-1})(a_{i,1})...f_{N_i,i-1,i}(a_{i-1})(a_{i,m})\omega_{i,2}\\=&f_{i-1,i}(a_{i-1})(a_{i,1})...f_{i-1,i}(a_{i-1})(a_{i,m})\omega_{i,2}\\=&f_{i-1,i}(a_{i-1})(a_{i,1}...a_{i,m}\omega_{i,2})\\=&f_{N_i,i-1,i}(a_{i-1})(a_{i,1}...a_{i,m}\omega_{i,2})\end{aligned}$$

Утверждение теоремы доказано. □

Из теоремы 3.4.5 следует, что если отображение $(f_{N_2,1,2}, ..., f_{N_n,1,n})$ - башня подпредставлений башни представлений $\overline{f}$, то отображение

$$(id : A_1 \to A_1, id_2 : N_2 \to A_2, ..., id_n : N_n \to A_n)$$

является морфизмом башен представлений.

**Теорема 3.4.6.** *Множество*[3.3] $\mathcal{B}_{\overline{f}}$ *всех башен подпредставлений башни представлений $\overline{f}$ порождает систему замыканий на башне представлений $\overline{f}$ и, следовательно, является полной структурой.*

Доказательство. Пусть для данного $\lambda \in \Lambda$, $K_{\lambda,i}$, $i = 2$, ..., $n$, - подалгебра $\Omega_i$-алгебры $A_i$, стабильная относительно представления $f_{i-1,i}$. Операцию пересечения на множестве $\mathcal{B}_{\overline{f}}$ мы определим согласно правилу

$$\bigcap f_{K_{\lambda,i-1,i}} = f_{\cap K_{\lambda,i-1,i}} \quad i = 2, ..., n$$

$$\bigcap \overline{K}_\lambda = \left(K_1 = A_1, K_2 = \bigcap K_{\lambda,2}, ..., K_n = \bigcap K_{\lambda,n}\right)$$

[3.2] Положим $N_1 = A_1$.

[3.3] Это определение аналогично определению структуры подалгебр ([12], стр. 93, 94)

$\cap K_{\lambda,i}$ - подалгебра $\Omega_i$-алгебры $A_i$. Пусть $a_i \in \cap K_{\lambda,i}$. Для любого $\lambda \in \Lambda$ и для любого $a_{i-1} \in K_{i-1}$,

$$f_{i-1,i}(a_{i-1})(a_i) \in K_{\lambda,i}$$

Следовательно,

$$f_{i-1,i}(a_{i-1})(a_i) \in K_i$$

Повторяя приведенное построение в порядке возрастания $i$, $i = 2$, ..., $n$, мы видим, что $(K_1, ..., K_n)$ - кортеж стабильных множеств башни представления $\overline{f}$. Следовательно, операция пересечения башен подпредставлений определена корректно. □

Обозначим соответствующий оператор замыкания через $\overline{J}(\overline{f})$. Если мы обозначим через $\overline{X}_{[1]}$ кортеж множеств $(X_2 \subset A_2, ..., X_n \subset A_n)$ то $\overline{J}(\overline{f}, \overline{X}_{[1]})$ является пересечением всех кортежей $(K_1, ..., K_n)$, стабильных относительно представления $\overline{f}$ и таких, что для $i = 2$, ..., $n$, $K_i$ - подалгебра $\Omega_i$-алгебры $A_i$, содержащая $X_i$.[3.4]

**Теорема 3.4.7.** *Пусть*[3.5] $\overline{f}$ *- башня представлений. Пусть* $X_i \subset A_i$, $i = 2$, ..., $n$. *Положим* $Y_1 = A_1$. *Последовательно увеличивая значение* $i$, $i = 2$, ..., $n$, *определим подмножества* $X_{i,m} \subset A_i$ *индукцией по* $m$.

$$X_{i,0} = X_i$$

$$x \in X_{i,m} => x \in X_{i,m+1}$$

$$x_1 \in X_{i,m}, ..., x_p \in X_{i,m}, \omega \in \Omega_i(p) => x_1...x_p\omega \in X_{i,m+1}$$

$$x_i \in X_{i,m}, x_{i-1} \in Y_{i-1} => f_{i-1,i}(x_{i-1})(x_i) \in X_{i,m+1}$$

*Для каждого значения* $i$ *положим*

$$Y_i = \bigcup_{m=0}^{\infty} X_{i,m}$$

*Тогда*

$$\overline{Y} = (Y_1, ..., Y_n) = \overline{J}(\overline{f}, \overline{X}_{[1]})$$

Доказательство. Для каждого значения $i$ доказательство теоремы совпадает с доказательством теоремы 2.4.4. Так как для определения устойчивого подмножества $\Omega_i$-алгебры $A_i$ нас интересует только некоторое устойчивое подмножество $\Omega_{i-1}$-алгебры $A_{i-1}$, мы должны сперва найти устойчивое подмножество $\Omega_{i-1}$-алгебры $A_{i-1}$. □

$\overline{J}(\overline{f}, \overline{X}_{[1]})$ называется **башней подпредставлений башни представлений $\overline{f}$, порождённой кортежем множеств $\overline{X}_{[1]}$**, а $\overline{X}_{[1]}$ - **кортежем множеств образующих башни подпредставлений** $\overline{J}(\overline{f}, \overline{X}_{[1]})$. В частности, **кортеж множеств образующих башни представлений** $\overline{f}$ будет такой кортеж $(X_2 \subset A_2, ..., X_n \subset A_n)$, что $\overline{J}(\overline{f}, \overline{X}_{[1]}) = \overline{A}$.

Нетрудно видеть, что определение кортежа множеств образующих башни представлений не зависит от того, эффективны представления башни или нет. Поэтому в дальнейшем мы будем предполагать, что представления башни

[3.4] При $n = 2$, $J_2(f_{1,2}, X_2) = J_{f_{1,2}}(X_2)$. Было бы проще использовать единные обозначения в разделах 2.5 и 3.5. Однако использование векторных обозначений в разделе 2.5 мне кажется несвоевременным.

[3.5] Утверждение теоремы аналогично утверждению теоремы 5.1, [12], стр. 94.

эффективны и будем опираться на соглашение для эффективного $T\star$-представления в замечании 2.1.9.

Мы также будем пользоваться записью

$$\overline{r}_{[1]} \circ \overline{a}_{[1]} = (r_2 \circ a_2, ..., r_n \circ a_n)$$

для образа кортежа элементов $a_2 \in A_2$, ..., $a_n \in A_n$ при эндоморфизме башни эффективных представлений. Согласно определению произведения отображений, для любых эндоморфизмов $\overline{r}_{[1]}$, $\overline{s}_{[1]}$ верно равенство

$$(\overline{r}_{[1]} \circ \overline{s}_{[1]}) \circ \overline{a}_{[1]} = \overline{r}_{[1]} \circ (\overline{s}_{[1]} \circ \overline{a}_{[1]}) \tag{3.4.7}$$

Равенство (3.4.7) является законом ассоциативности для $\circ$ и позволяет записать выражение

$$\overline{r}_{[1]} \circ \overline{s}_{[1]} \circ \overline{a}_{[1]}$$

без использования скобок.

Из теоремы 3.4.7 следует следующее определение.

**Определение 3.4.8.** Пусть $(X_2 \subset A_2, ..., X_n \subset A_n)$ - кортеж множеств. Для любого кортежа элементов $\overline{a}$, $\overline{a} \in \overline{J}(\overline{f}, \overline{X}_{[1]})$, существует кортеж $\overline{\Omega}$-слов, определённых согласно следующему правилу.

(1) Если $a_1 \in A_1$, то $a_1$ - $\Omega_1$-слово.
(2) Если $a_i \in X_i$, $i = 2$, ..., $n$, то $a_i$ - $\Omega_i$-слово.
(3) Если $a_{i,1}$, ..., $a_{i,p}$ - $\Omega_i$-слова, $i = 2$, ..., $n$, и $\omega \in \Omega_i(p)$, то $a_{i,1}...a_{i,p}\omega$ - $\Omega_i$-слово.
(4) Если $a_i$ - $\Omega_i$-слово, $i = 2$, ..., 3, и $a_{i-1}$ - $\Omega_{i-1}$-слово, то $a_{i-1}a_i$ - $\Omega_i$-слово.

**Кортеж $\overline{\Omega}$-слов**[3.6]

$$\overline{w}(\overline{f}, \overline{X}_{[1]}, \overline{a}) = (w_1(\overline{f}, \overline{X}_{[1]}, a_1), ..., w_n(\overline{f}, \overline{X}_{[1]}, a_n))$$

представляет данный элемент $\overline{a} \in \overline{J}(\overline{f}, \overline{X}_{[1]})$.[3.7] Обозначим $\overline{w}(\overline{f}, \overline{X}_{[1]})$ **множество кортежей $\overline{\Omega}$-слов башни представлений** $\overline{J}(\overline{f}, \overline{X}_{[1]})$. □

Представление $a_i \in A_i$ в виде $\Omega_i$-слова неоднозначно. Если $a_{i,1}$, ..., $a_{i,p}$ - $\Omega_i$-слова, $\omega \in \Omega_i(p)$ и $a_{i-1} \in A_{i-1}$, то $\Omega_i$-слова $a_{i-1}a_{i,1}...a_{i,p}\omega$ и $a_{i-1}a_{i,1}...a_{i-1}a_{i,p}\omega$ описывают один и тот же элемент $\Omega_i$-алгебры $A_i$. Возможны равенства, связанные со спецификой представления. Например, если $\omega$ является операцией $\Omega_{i-1}$-алгебры $A_{i-1}$ и операцией $\Omega_i$-алгебры $A_i$, то мы можем потребовать, что $\Omega_i$-слова $a_{i-1,1}...a_{i-1,p}\omega a_i$ и $a_{i-1,1}a_i...a_{i-1,p}a_i\omega$ описывают один и тот же элемент $\Omega_i$-алгебры $A_i$. Перечисленные выше равенства для каждого значения $i$, $i = 2$, ..., $n$, определяют отношение эквивалентности $r_i$ на множестве $\Omega_i$-слов $W_i(\overline{f}, \overline{X}_{[1]})$.Согласно построению, отношение эквивалентности $r_i$ на множестве

[3.6] Согласно построению

$$w_2(\overline{f}, \overline{X}_{[1]}, a_2) = w_2(f_{1,2}, X_2, a_2) \quad w_3(\overline{f}, \overline{X}_{[1]}, a_3) = w_3((f_{1,2}, f_{2,3}), (X_2, X_3), a_3)$$

[3.7] Эти обозначения имеют небольшое различие с обозначениями, предложенными ранее в теории представлении универсальной алгебры. А именно, в кортеж $\overline{\Omega}$-слов мы добавили слово $\Omega_1$-алгебры. Однако очевидно, что это слово всегда тривиально

$$w_1(\overline{f}, \overline{X}_{[1]}, a_1) = a_1$$

$\Omega_i$-слов $W_i(\overline{f},\overline{X}_{[1]})$ зависит не только от выбора множества $X_i$, но и от выбора множества $X_{i-1}$.

**Теорема 3.4.9.** *Эндоморфизм $\overline{r}_{[1]}$ башни представлений $\overline{f}$ порождает отображение $\overline{\Omega}_{[1]}$-слов*

$$\overline{w}_{[1]}[\overline{f},\overline{X}_{[1]},\overline{r}_{[1]}]:\overline{W}_{[1]}(\overline{f},\overline{X}_{[1]})\rightarrow\overline{w}_{[1]}(\overline{f},\overline{X}'_{[1]})\quad\overline{X}_{[1]}\subset\overline{A}_{[1]}\quad\overline{X}'_{[1]}=\overline{r}_{[1]}(\overline{X}_{[1]})$$

*такое, что для любого $i$, $i=2$, ..., $n$,*

(1) *Если* $a_i\in X_i$, $a'_i=r_i\circ a_i$, *то*

$$w_i[\overline{f},\overline{X}_{[1]},\overline{r}_{[1]}](a_i)=a'_i$$

(2) *Если*

$$a_{i,1},...,a_{i,n}\in w_i(\overline{f},\overline{X}_{[1]})$$
$$a'_{i,1}=w_i[\overline{f},\overline{X}_{[1]},\overline{r}_{[1]}](a_{i,1})\quad ...\quad a'_{i,p}=w_i[\overline{f},\overline{X}_{[1]},\overline{r}_{[1]}](a_{i,p})$$

*то для операции* $\omega\in\Omega_i(p)$ *справедливо*

$$w_i[\overline{f},\overline{X}_{[1]},\overline{r}_{[1]}](a_{i,1}...a_{i,p}\omega)=a'_{i,1}...a'_{i,p}\omega$$

(3) *Если*

$$\begin{array}{ll}a_i\in w_i(\overline{f},\overline{X}_{[1]}) & a'_i=w_i[\overline{f},\overline{X}_{[1]},\overline{r}_{[1]}](a_i)\\ a_{i-1}\in w_{i-1}(\overline{f},\overline{X}_{[1]}) & a'_{i-1}=w_{i-1}[\overline{f},\overline{X}_{[1]},\overline{r}_{[1]}](a_{i-1})\end{array}$$

*то*

$$w_i[\overline{f},\overline{X}_{[1]},\overline{r}_{[1]}](a_{i-1}a_i)=a'_{i-1}a'_i$$

Доказательство. Утверждения (1), (2) теоремы справедливы в силу определения эндоморфизма $r_i$. Утверждение (3) теоремы следует из равенства (3.2.3). □

**Замечание 3.4.10.** Пусть $\overline{r}_{[1]}$ - эндоморфизм башни представлений $\overline{f}$. Пусть

$$\overline{a}_{[1]}\in\overline{J}_{[1]}(\overline{f},\overline{X}_{[1]})\quad\overline{a}'_{[1]}=\overline{r}_{[1]}\circ\overline{a}_{[1]}\quad\overline{X}'_{[1]}=\overline{r}_{[1]}\circ\overline{X}_{[1]}$$

Теорема 3.4.9 утверждает, что $\overline{a}'_{[1]}\in\overline{J}_{[1]}(\overline{f},\overline{X}'_{[1]})$ . Теорема 3.4.9 также утверждает, что $\overline{\Omega}_{[1]}$-слово, представляющее $\overline{a}_{[1]}$, относительно $\overline{X}_{[1]}$ и $\overline{\Omega}_{[1]}$-слово, представляющее $\overline{a}'_{[1]}$, относительно $\overline{X}'_{[1]}$ формируются согласно одному и тому же алгоритму. Это позволяет рассматривать множество $\overline{\Omega}_{[1]}$-слов $\overline{w}_{[1]}(\overline{f},\overline{X}'_{[1]},\overline{a}'_{[1]})$ как отображение

$$\overline{W}_{[1]}(\overline{f},\overline{X}_{[1]},\overline{a}_{[1]}):\overline{X}'_{[1]}\rightarrow\overline{w}_{[1]}(\overline{f},\overline{X}'_{[1]},\overline{a}'_{[1]})$$
$$\overline{W}_{[1]}(\overline{f},\overline{X}_{[1]},\overline{a}_{[1]})(\overline{X}'_{[1]})=\overline{W}_{[1]}(\overline{f},\overline{X}_{[1]},\overline{a}_{[1]})\circ\overline{X}'_{[1]}$$

такое, что, если для некоторого эндоморфизма $\overline{r}_{[1]}$

$$\overline{X}'_{[1]}=\overline{r}_{[1]}\circ\overline{X}_{[1]}\quad\overline{a}'_{[1]}=\overline{r}_{[1]}\circ\overline{a}_{[1]}$$

то

$$\overline{W}_{[1]}(\overline{f},\overline{X}_{[1]},\overline{a}_{[1]})\circ\overline{X}'_{[1]}=\overline{w}_{[1]}(\overline{f},\overline{X}'_{[1]},\overline{a}'_{[1]})=\overline{a}'_{[1]}$$

Кортеж отображений[3.8]

$$\overline{W}_{[1]}(\overline{f},\overline{X}_{[1]},\overline{a}_{[1]})=(W_2(\overline{f},\overline{X}_{[1]},a_2),...,W_n(\overline{f},\overline{X}_{[1]},a_n))$$

представляющий данный элемент $\overline{a}\in\overline{J}(\overline{f},\overline{X}_{[1]})$, называется **кортежем координат элемента $\overline{a}$ относительно кортежа множеств** $\overline{X}_{[1]}$. Обозначим $\overline{W}_{[1]}(\overline{f},\overline{X}_{[1]})$ **множество кортежей координат башни представлений** $\overline{J}(\overline{f},\overline{X}_{[1]})$. Аналогично, мы можем рассмотреть координаты кортежа множеств $\overline{B}_{[1]}\subset\overline{J}_{[1]}(f,X)$ относительно кортежа множеств $\overline{X}_{[1]}$

$$\overline{W}_{[1]}(\overline{f},\overline{X}_{[1]},\overline{B}_{[1]})=(W_2(\overline{f},\overline{X}_{[1]},B_2),...,W_n(\overline{f},\overline{X}_{[1]},B_n))$$

$$W_k(\overline{f},\overline{X}_{[1]},B_k)=\{W_k(\overline{f},\overline{X}_{(2...k)},a_k):a_k\in B_k\}=(W_k(\overline{f},\overline{X}_{(2...k)},a_k),a_k\in B_k)$$

□

**Теорема 3.4.11.** *На множестве координат $W_k(\overline{f},\overline{X}_{[1]})$ определена структура $\Omega_k$-алгебры.*

Доказательство. Пусть $\omega\in\Omega_k(n)$. Тогда для любых $m_1$, ..., $m_n\in J_k(\overline{f},\overline{X}_{[1]})$ положим

$$W_k(\overline{f},\overline{X}_{[1]},m_1)...W_k(\overline{f},\overline{X}_{[1]},m_n)\omega=W_k(\overline{f},\overline{X}_{[1]},m_1...m_n\omega)\tag{3.4.9}$$

Согласно замечанию 3.4.10, из равенства (3.4.9) следует

$$\begin{aligned}(W_k(\overline{f},\overline{X}_{[1]},m_1)...W_k(\overline{f},\overline{X}_{[1]},m_n)\omega)\circ\overline{X}_{[1]}&=W_k(\overline{f},\overline{X}_{[1]},m_1...m_n\omega)\circ\overline{X}_{[1]}\\&=w(\overline{f},\overline{X}_{[1]},m_1...m_n\omega)\end{aligned}\tag{3.4.10}$$

Согласно правилу (3) определения 3.4.8, из равенства (3.4.10) следует

$$\begin{aligned}&(W_k(\overline{f},\overline{X}_{[1]},m_1)...W_k(\overline{f},\overline{X}_{[1]},m_n)\omega)\circ\overline{X}_{[1]}\\=&w_k(\overline{f},\overline{X}_{[1]},m_1)...w_k(\overline{f},\overline{X}_{[1]},m_n)\omega\\=&(W_k(\overline{f},\overline{X}_{[1]},m_1)\circ\overline{X}_{[1]})...(W_k(\overline{f},\overline{X}_{[1]},m_n)\circ\overline{X}_{[1]})\omega\end{aligned}\tag{3.4.11}$$

Из равенства (3.4.11) следует корректность определения (3.4.9) операции $\omega$ на множестве координат $W_k(\overline{f},\overline{X}_{[1]})$. □

**Теорема 3.4.12.** *Для $k=2$, ..., $n$, определено представление $\Omega_{k-1}$-алгебры $W_{k-1}(\overline{f},\overline{X}_{[1]})$ в $\Omega_k$-алгебре $W_k(\overline{f},\overline{X}_{[1]})$.*[3.9]

Доказательство. Пусть $a_{k-1}\in J_{k-1}(\overline{f},\overline{X}_{[1]})$. Тогда для любого $a_k\in J_k(\overline{f},\overline{X}_{[1]})$, положим

$$W_{k-1}(\overline{f},\overline{X}_{[1]},a_{k-1})W_k(\overline{f},\overline{X}_{[1]},a_k)=W_k(\overline{f},\overline{X}_{[1]},a_{k-1}a_k)\tag{3.4.12}$$

[3.8] Согласно построению

$$W_2(\overline{f},\overline{X}_{[1]},a_2)=W_2(f_{1,2},X_2,a_2)\qquad W_3(\overline{f},\overline{X}_{[1]},a_3)=W_3((f_{1,2},f_{2,3}),(X_2,X_3),a_3)$$

Нас не интересует отображение $W_1(\overline{f},\overline{X},a_1)$ поскольку оно тривиально

$$W_1(\overline{f},\overline{X}_{[1]},a_1)\circ\overline{X}'_{[1]}=a_1\tag{3.4.8}$$

[3.9] Согласно равенству (3.4.8), мы можем отождествить $\Omega_1$-алгебры $W_1(\overline{f},\overline{X}_{[1]})$ и $A_1$.

Согласно замечанию 3.4.10, из равенства (3.4.12) следует

$$(W_{k-1}(\overline{f},\overline{X}_{[1]},a_{k-1})W_k(\overline{f},\overline{X}_{[1]},a_k))\circ\overline{X}_{[1]}=W_k(\overline{f},\overline{X}_{[1]},am)\circ\overline{X}_{[1]} \\ =w_k(\overline{f},\overline{X}_{[1]},a_{k-1}a_k) \tag{3.4.13}$$

Согласно правилу (4) определения 3.4.8, из равенства (3.4.13) следует

$$\begin{aligned}&(W_{k-1}(\overline{f},\overline{X}_{[1]},a_{k-1})W_k(\overline{f},\overline{X}_{[1]},a_k))\circ\overline{X}_{[1]}\\ =&w_{k-1}(\overline{f},\overline{X}_{[1]},a_{k-1})w_k(\overline{f},\overline{X}_{[1]},a_k)\\ =&(W_{k-1}(\overline{f},\overline{X}_{[1]},a_{k-1})\circ\overline{X}_{[1]})(W_k(\overline{f},\overline{X}_{[1]},a_k)\circ\overline{X}_{[1]})\end{aligned} \tag{3.4.14}$$

Из равенства (3.4.14) следует корректность определения (3.4.12) представления $\Omega_{k-1}$-алгебры $W_{k-1}(\overline{f},\overline{X}_{[1]})$ в $\Omega_k$-алгебре $W_k(\overline{f},\overline{X}_{[1]})$. □

**Следствие 3.4.13.** *Кортеж $\overline{\Omega}$-алгебр $\overline{W}(\overline{f},\overline{X}_{[1]})$ порождает башню представлений.* □

**Определение 3.4.14.** Пусть $\overline{X}_{[1]}$ - кортеж множеств образующих башни представлений $\overline{f}$. Пусть $\overline{r}_{[1]}$ - эндоморфизм башни представлений $\overline{f}$. Множество координат $\overline{W}_{[1]}(\overline{f},\overline{X}_{[1]},\overline{r}_{[1]}\circ\overline{X}_{[1]})$ называется **координатами эндоморфизма башни представлений**. □

**Определение 3.4.15.** Пусть $\overline{X}_{[1]}$ - кортеж множеств образующих башни представлений $\overline{f}$. Пусть $\overline{r}_{[1]}$ - эндоморфизм башни представлений $\overline{f}$. Пусть $\overline{a}_{[1]}\in\overline{A}_{[1]}$. Мы определим **суперпозицию координат** эндоморфизма башни представлений $\overline{f}$ и элемента $\overline{a}_{[1]}$ как координаты, определённые согласно правилу

$$\overline{W}_{[1]}(\overline{f},\overline{X}_{[1]},\overline{a}_{[1]})\circ\overline{W}_{[1]}(\overline{f},\overline{X}_{[1]},\overline{r}_{[1]}\circ\overline{X}_{[1]})=\overline{W}_{[1]}(\overline{f},\overline{X}_{[1]},\overline{r}_{[1]}\circ\overline{a}_{[1]}) \tag{3.4.15}$$

Пусть $\overline{Y}_{[1]}\subset\overline{A}_{[1]}$. Мы определим суперпозицию координат эндоморфизма башни представлений $\overline{f}$ и кортежа множеств $\overline{Y}_{[1]}$ согласно правилу

$$\begin{aligned}&\overline{W}_{[1]}(\overline{f},\overline{X}_{[1]},\overline{Y}_{[1]})\circ\overline{W}_{[1]}(\overline{f},\overline{X}_{[1]},\overline{r}_{[1]}\circ\overline{X}_{[1]})\\ =&(\overline{W}_{[1]}(\overline{f},\overline{X}_{[1]},\overline{a}_{[1]})\circ\overline{W}_{[1]}(\overline{f},\overline{X}_{[1]},\overline{r}_{[1]}\circ\overline{X}_{[1]}),\overline{a}_{[1]}\in\overline{Y}_{[1]})\end{aligned} \tag{3.4.16}$$

□

**Теорема 3.4.16.** *Эндоморфизм $\overline{r}_{[1]}$ представления $\overline{f}$ порождает отображение координат башни представлений*

$$\overline{W}_{[1]}[\overline{f},\overline{X}_{[1]},\overline{r}_{[1]}]:\overline{W}_{[1]}(\overline{f},\overline{X}_{[1]})\rightarrow\overline{W}_{[1]}(\overline{f},\overline{X}_{[1]}) \tag{3.4.17}$$

*такое, что*

(3.4.18)

$$\overline{W}_{[1]}(\overline{f},\overline{X}_{[1]},\overline{a}_{[1]})\rightarrow\overline{W}_{[1]}[\overline{f},\overline{X}_{[1]},\overline{r}_{[1]}]\star\overline{W}_{[1]}(\overline{f},\overline{X}_{[1]},\overline{a}_{[1]})=\overline{W}_{[1]}(\overline{f},\overline{X}_{[1]},\overline{r}_{[1]}\circ\overline{a}_{[1]})$$
$$\overline{W}_{[1]}[\overline{f},\overline{X}_{[1]},\overline{r}_{[1]}]\star\overline{W}_{[1]}(\overline{f},\overline{X}_{[1]},\overline{a}_{[1]})=\overline{W}_{[1]}(\overline{f},\overline{X}_{[1]},\overline{a}_{[1]})\circ\overline{W}_{[1]}(\overline{f},\overline{X}_{[1]},\overline{r}_{[1]}\circ\overline{X}_{[1]})$$

Доказательство. Согласно замечанию 3.4.10, мы можем рассматривать равенства (3.4.15), (3.4.17) относительно заданного кортежа множеств образующих $\overline{X}_{[1]}$. При этом координатам $\overline{W}_{[1]}(\overline{f},\overline{X}_{[1]},\overline{a}_{[1]})$ соответствует кортеж слов

$$\overline{W}_{[1]}(\overline{f},\overline{X}_{[1]},\overline{a}_{[1]})\circ\overline{X}_{[1]}=\overline{w}_{[1]}(\overline{f},\overline{X}_{[1]},\overline{a}_{[1]}) \tag{3.4.19}$$

а координатам $\overline{W}_{[1]}(\overline{f}, \overline{X}_{[1]}, \overline{r}_{[1]} \circ \overline{a}_{[1]})$ соответствует кортеж слов

$$\overline{W}_{[1]}(\overline{f}, \overline{X}_{[1]}, \overline{r}_{[1]} \circ \overline{a}_{[1]}) \circ X = \overline{w}_{[1]}(\overline{f}, \overline{X}_{[1]}, \overline{r}_{[1]} \circ \overline{a}_{[1]}) \tag{3.4.20}$$

Поэтому для того, чтобы доказать теорему, нам достаточно показать, что отображению $\overline{W}_{[1]}[\overline{f}, \overline{X}_{[1]}, \overline{r}_{[1]}]$ соответствует отображение $\overline{w}_{[1]}[\overline{f}, \overline{X}_{[1]}, \overline{r}_{[1]}]$. Мы докажем утверждение индукцией по числу $n$ универсальных алгебр в башне представлений.

При $n = 2$ башня представлений $\overline{f}$ является представлением $\Omega_1$-алгебры $A_1$ в $\Omega_2$-алгебре $A_2$. Утверждение теоремы является следствием теоремы 2.4.13. Пусть утверждение теоремы верно для $n-1$. Мы будем доказывать теорему индукцией по сложности $\Omega_n$-слова.

Если $a_n \in X_n$, $a'_n = r_n \circ a_n$, то, согласно равенствам (3.4.19), (3.4.20), отображения $W_n[\overline{f}, \overline{X}_{[1]}, \overline{r}_{[1]}]$ и $w_n[\overline{f}, \overline{X}_{[1]}, \overline{r}_{[1]}]$ согласованы.

Пусть для $a_{n,1}$, ..., $a_{n,p} \in X_n$ отображения $W_n[\overline{f}, \overline{X}_{[1]}, \overline{r}_{[1]}]$ и $w_n[\overline{f}, \overline{X}_{[1]}, \overline{r}_{[1]}]$ согласованы. Пусть $\omega \in \Omega_n(p)$. Согласно теореме 2.4.9

$$W_n(\overline{f}, \overline{X}_{[1]}, a_{n,1}...a_{n,p}\omega) = W_n(\overline{f}, \overline{X}_{[1]}, a_{n,1})...W_n(\overline{f}, \overline{X}_{[1]}, a_{n,p})\omega \tag{3.4.21}$$

Так как $r_n$ - эндоморфизм $\Omega_n$-алгебры $A_n$, то из равенства (3.4.21) следует

(3.4.22)
$$\begin{aligned} W_n(\overline{f}, \overline{X}_{[1]}, r_n \circ (a_{n,1}...a_{n,p}\omega)) &= W_n(\overline{f}, \overline{X}_{[1]}, (r_n \circ a_{n,1})...(r_n \circ a_{n,p})\omega) \\ &= W_n(\overline{f}, \overline{X}_{[1]}, r_n \circ a_{n,1})...W(f, X, r_n \circ a_{n,p})\omega \end{aligned}$$

Из равенств (3.4.21), (3.4.22) и предположения индукции следует, что отображения $W_n[\overline{f}, \overline{X}_{[1]}, \overline{r}_{[1]}]$ и $w_n[\overline{f}, \overline{X}_{[1]}, \overline{r}_{[1]}]$ согласованы для $a_n = a_{n,1}...a_{n,p}\omega$.

Пусть для $a_{n,1} \in A_n$ отображения $W_n[\overline{f}, \overline{X}_{[1]}, \overline{r}_{[1]}]$ и $w_n[\overline{f}, \overline{X}_{[1]}, \overline{r}_{[1]}]$ согласованы. Пусть $a_{n-1} \in A_{n-1}$. Согласно теореме 3.4.12

$$W_n(\overline{f}, \overline{X}_{[1]}, a_{n-1}a_{n,1}) = W_{n-1}(\overline{f}, \overline{X}_{[1]}, a_{n-1})W_n(\overline{f}, \overline{X}_{[1]}, a_{n,1}) \tag{3.4.23}$$

Так как $\overline{r}_{[1]}$ - эндоморфизм башни представлений $\overline{f}$, то из равенства (3.4.23) следует

(3.4.24)
$$\begin{aligned} W_n(\overline{f}, \overline{X}_{[1]}, r_n \circ (a_{n-1}a_{n,1})) &= W_n(\overline{f}, \overline{X}_{[1]}, (r_{n-1} \circ a_{n-1})(r_n \circ a_{n,1})) \\ &= W_{n-1}(\overline{f}, \overline{X}_{[1]}, r_{n-1} \circ a_{n-1})W_n(\overline{f}, \overline{X}_{[1]}, r_n \circ a_{n,1}) \end{aligned}$$

Из равенств (3.4.23), (3.4.24) и предположения индукции следует, что отображения $W_n[\overline{f}, \overline{X}_{[1]}, \overline{r}_{[1]}]$ и $w_n[\overline{f}, \overline{X}_{[1]}, \overline{r}_{[1]}]$ согласованы для $a_{n,2} = a_{n-1}a_{n,1}$. □

**Следствие 3.4.17.** *Пусть $\overline{X}_{[1]}$ - кортеж множеств образующих башни представлений $\overline{f}$. Пусть $\overline{r}_{[1]}$ - эндоморфизм башни представлений $\overline{f}$. Отображение $\overline{W}_{[1]}[\overline{f}, \overline{X}_{[1]}, \overline{r}_{[1]}]$ является эндоморфизмом башни представлений $\overline{W}(\overline{f}, \overline{X}_{[1]})$.* □

В дальнейшем мы будем отождествлять отображение $\overline{W}_{[1]}[\overline{f}, \overline{X}_{[1]}, \overline{r}_{[1]}]$ и множество координат $\overline{W}_{[1]}(\overline{f}, \overline{X}_{[1]}, \overline{r}_{[1]} \circ \overline{X}_{[1]})$.

**Теорема 3.4.18.** *Пусть $\overline{X}_{[1]}$ - кортеж множеств образующих башни представлений $\overline{f}$. Пусть $\overline{r}_{[1]}$ - эндоморфизм башни представлений $\overline{f}$. Пусть $\overline{Y}_{[1]} \subset$*

$\overline{M}_{[1]}$. *Тогда*

(3.4.25) $$\overline{W}_{[1]}(\overline{f},\overline{X}_{[1]},\overline{Y}_{[1]})\circ\overline{W}_{[1]}(\overline{f},\overline{X}_{[1]},\overline{r}_{[1]}\circ\overline{X}_{[1]})=\overline{W}_{[1]}(\overline{f},\overline{X}_{[1]},\overline{r}_{[1]}\circ\overline{Y}_{[1]})$$

(3.4.26) $$\overline{W}_{[1]}[\overline{f},\overline{X}_{[1]},\overline{r}_{[1]}]\star\overline{W}_{[1]}(\overline{f},\overline{X}_{[1]},\overline{Y}_{[1]})=\overline{W}_{[1]}(\overline{f},\overline{X}_{[1]},\overline{r}_{[1]}\circ\overline{Y}_{[1]})$$

Доказательство. Равенство (3.4.25) является следствием равенства

$$\overline{r}_{[1]}\circ\overline{Y}_{[1]}=(\overline{r}_{[1]}\circ\overline{a}_{[1]},\overline{a}_{[1]}\in\overline{Y}_{[1]})$$

а также равенств (3.4.15), (3.4.16). Равенство (3.4.26) является следствием равенств (3.4.25), (3.4.18). □

**Теорема 3.4.19.** *Пусть $\overline{X}_{[1]}$ - кортеж множеств образующих башни представлений $\overline{f}$. Пусть $\overline{r}_{[1]}$, $\overline{s}_{[1]}$ - эндоморфизмы башни представлений $\overline{f}$. Тогда*

(3.4.27)
$$\overline{W}_{[1]}(\overline{f},\overline{X}_{[1]},\overline{s}_{[1]}\circ\overline{X}_{[1]})\circ\overline{W}_{[1]}(\overline{f},\overline{X}_{[1]},\overline{r}_{[1]}\circ\overline{X}_{[1]})=\overline{W}_{[1]}(\overline{f},\overline{X}_{[1]},\overline{r}_{[1]}\circ\overline{s}_{[1]}\circ\overline{X}_{[1]})$$

(3.4.28) $$\overline{W}_{[1]}[\overline{f},\overline{X}_{[1]},\overline{r}_{[1]}]\star\overline{W}_{[1]}[\overline{f},\overline{X}_{[1]},\overline{s}_{[1]}]=\overline{W}_{[1]}[\overline{f},\overline{X}_{[1]},\overline{r}_{[1]}\circ\overline{s}_{[1]}]$$

Доказательство. Равенство (3.4.27) следует из равенства (3.4.25), если положить $\overline{Y}_{[1]}=\overline{s}_{[1]}\circ\overline{X}_{[1]}$. Равенство (3.4.28) следует из равенства (3.4.27) и цепочки равенств

$$\begin{aligned}
&(\overline{W}_{[1]}[\overline{f},\overline{X}_{[1]},\overline{r}_{[1]}]\star\overline{W}_{[1]}[\overline{f},\overline{X}_{[1]},\overline{s}_{[1]}])\star\overline{W}_{[1]}(\overline{f},\overline{X}_{[1]},\overline{Y}_{[1]})\\
=&\overline{W}_{[1]}[\overline{f},\overline{X}_{[1]},\overline{r}_{[1]}]\star(\overline{W}_{[1]}[\overline{f},\overline{X}_{[1]},\overline{s}_{[1]}]\star\overline{W}_{[1]}(\overline{f},\overline{X}_{[1]},\overline{Y}_{[1]}))\\
=&(\overline{W}_{[1]}(\overline{f},\overline{X}_{[1]},\overline{Y}_{[1]})\circ\overline{W}_{[1]}(\overline{f},\overline{X}_{[1]},\overline{s}_{[1]}\circ\overline{X}_{[1]}))\circ\overline{W}_{[1]}(\overline{f},\overline{X}_{[1]},\overline{r}_{[1]}\circ\overline{X}_{[1]})\\
=&\overline{W}_{[1]}(\overline{f},\overline{X}_{[1]},\overline{Y}_{[1]})\circ(\overline{W}_{[1]}(\overline{f},\overline{X}_{[1]},\overline{s}_{[1]}\circ\overline{X}_{[1]})\circ\overline{W}_{[1]}(\overline{f},\overline{X}_{[1]},\overline{r}_{[1]}\circ\overline{X}_{[1]}))\\
=&\overline{W}_{[1]}(\overline{f},\overline{X}_{[1]},\overline{Y}_{[1]})\circ\overline{W}_{[1]}(\overline{f},\overline{X}_{[1]},\overline{r}_{[1]}\circ\overline{s}_{[1]}\circ\overline{X}_{[1]})\\
=&\overline{W}_{[1]}(\overline{f},\overline{X}_{[1]},\overline{r}_{[1]}\circ\overline{s}_{[1]})\star\overline{W}_{[1]}(\overline{f},\overline{X}_{[1]},\overline{Y}_{[1]})
\end{aligned}$$

□

Мы можем обобщить определение суперпозиции координат и предположить, что один из множителей является кортежем множеств $\overline{\Omega}_{[1]}$-слов. Соответственно, определение суперпозиции координат имеет вид

$$\begin{aligned}
&\overline{w}_{[1]}(\overline{f},\overline{X}_{[1]},\overline{Y}_{[1]})\circ\overline{W}_{[1]}(\overline{f},\overline{X}_{[1]},\overline{r}_{[1]}\circ\overline{X}_{[1]})\\
=&\overline{W}_{[1]}(\overline{f},\overline{X}_{[1]},\overline{Y}_{[1]})\circ\overline{w}_{[1]}(\overline{f},\overline{X}_{[1]},\overline{r}_{[1]}\circ\overline{X}_{[1]})\\
=&\overline{w}_{[1]}(\overline{f},\overline{X}_{[1]},\overline{r}_{[1]}\circ\overline{Y}_{[1]})
\end{aligned}$$

Следующие формы записи образа кортежа множеств $\overline{Y}_{[1]}$ при эндоморфизме $\overline{r}_{[1]}$ эквивалентны.

(3.4.29) $$\begin{aligned}
\overline{r}_{[1]}\circ\overline{Y}_{[1]}&=\overline{W}_{[1]}(\overline{f},\overline{X}_{[1]},Y)\circ(\overline{r}_{[1]}\circ\overline{X}_{[1]})\\
&=\overline{W}_{[1]}(\overline{f},\overline{X}_{[1]},\overline{Y}_{[1]})\circ(\overline{W}_{[1]}(\overline{f},\overline{X}_{[1]},\overline{r}_{[1]}\circ\overline{X}_{[1]})\circ\overline{X}_{[1]})
\end{aligned}$$

Из равенств (3.4.25), (3.4.29) следует, что

(3.4.30) $$\begin{aligned}
&(\overline{W}_{[1]}(\overline{f},\overline{X}_{[1]},\overline{Y}_{[1]})\circ\overline{W}_{[1]}(\overline{f},\overline{X}_{[1]},\overline{r}_{[1]}\circ\overline{X}_{[1]}))\circ\overline{X}_{[1]}\\
=&\overline{W}_{[1]}(\overline{f},\overline{X}_{[1]},\overline{Y}_{[1]})\circ(\overline{W}_{[1]}(\overline{f},\overline{X}_{[1]},\overline{r}_{[1]}\circ\overline{X}_{[1]})\circ\overline{X}_{[1]})
\end{aligned}$$

Равенство (3.4.30) является законом ассоциативности для операции композиции и позволяет записать выражение

$$\overline{W}_{[1]}(\overline{f},\overline{X}_{[1]},\overline{Y}_{[1]})\circ\overline{W}_{[1]}(\overline{f},\overline{X}_{[1]},\overline{r}_{[1]}\circ\overline{X}_{[1]})\circ X$$

без использования скобок.

**Определение 3.4.20.** Пусть $(X_2\subset A_2,...,X_n\subset A_n)$ - кортеж множеств образующих башни представлений $\overline{f}$. Пусть отображение $\overline{r}_{[1]}$ является эндоморфизмом башни представления $\overline{f}$. Пусть кортеж множеств $\overline{X}'_{[1]}=\overline{r}_{[1]}\circ\overline{X}_{[1]}$ является образом кортежа множеств $\overline{X}_{[1]}$ при отображении $\overline{r}_{[1]}$. Эндоморфизм $\overline{r}_{[1]}$ башни представлений $\overline{f}$ называется **невырожденным на кортеже множеств образующих** $\overline{X}_{[1]}$, если кортеж множеств $\overline{X}'_{[1]}$ является кортежем множеств образующих башни представлений $\overline{f}$. В противном случае, эндоморфизм $\overline{r}_{[1]}$ называется **вырожденным на кортеже множеств образующих** $\overline{X}_{[1]}$, □

**Определение 3.4.21.** Эндоморфизм башни представлений $\overline{f}$ называется **невырожденным**, если он невырожден на любом кортеже множеств образующих. □

**Теорема 3.4.22.** *Автоморфизм $\overline{r}$ башни представлений $\overline{f}$ является невырожденным эндоморфизмом.*

Доказательство. Пусть $\overline{X}_{[1]}$ - кортеж множеств образующих башни представлений $\overline{f}$. Пусть $\overline{X}'_{[1]}=\overline{r}_{[1]}\circ\overline{X}_{[1]}$.

Согласно теореме 3.4.16 эндоморфизм $\overline{r}_{[1]}$ порождает отображение $\overline{\Omega}_{[1]}$-слов $\overline{w}_{[1]}[\overline{f},\overline{X}_{[1]},\overline{r}_{[1]}]$.

Пусть $\overline{a}'_{[1]}\in\overline{A}_{[1]}$. Так как $\overline{r}_{[1]}$ - автоморфизм, то существует $\overline{a}_{[1]}\in\overline{A}_{[1]}$, $\overline{r}_{[1]}\circ\overline{a}_{[1]}=\overline{a}'_{[1]}$. Согласно определению 3.4.8, $\overline{w}_{[1]}(\overline{f},\overline{X}_{[1]},\overline{a}_{[1]})$ - кортеж $\overline{\Omega}_{[1]}$-слов, представляющих $\overline{a}_{[1]}$ относительно кортежа множеств образующих $\overline{X}_{[1]}$. Согласно теореме 3.4.16, $\overline{w}_{[1]}(\overline{f},\overline{X}'_{[1]},\overline{a}'_{[1]})$ - кортеж $\overline{\Omega}_{[1]}$-слов, представляющих $\overline{a}'_{[1]}$ относительно кортежа множеств $\overline{X}'_{[1]}$

$$\overline{w}_{[1]}(\overline{f},\overline{X}'_{[1]},\overline{a}'_{[1]})=\overline{w}_{[1]}[\overline{f},\overline{X}_{[1]},\overline{r}](\overline{w}_{[1]}(\overline{f},\overline{X}_{[1]},\overline{a}_{[1]}))$$

Следовательно, $\overline{X}'_{[1]}$ - множество образующих представления $\overline{f}$. Согласно определению 3.4.21, автоморфизм $\overline{r}_{[1]}$ - невырожден. □

## 3.5. Базис башни представлений

**Определение 3.5.1.** Если кортеж множеств $\overline{X}_{[1]}$ является кортежем множеств образующих башни представлений $\overline{f}$, то любой кортеж множеств $\overline{Y}_{[1]}$, $X_i\subset Y_i\subset A_i$, $i=2$, ..., $n$, также является кортежем множеств образующих башни представлений $\overline{f}$. Если существует кортеж минимальных множеств $\overline{X}_{[1]}$, порождающих башню представлений $\overline{f}$, то такой кортеж множеств $\overline{X}_{[1]}$ называется **базисом башни представлений** $\overline{f}$. □

**Теорема 3.5.2.** *Базис башни представлений определён индукцией по $n$. При $n=2$ базис башни представлений является базисом представления $f_{1,2}$. Если кортеж множеств $\overline{X}_{[1,n]}$ является базисом башни представлений $\overline{f}_{[1]}$, то* *кортеж множеств образующих* $\overline{X}_{[1]}$ *башни представлений* $\overline{f}$ *является базисом тогда и только тогда, когда для любого* $a_n \in X_n$ *кортеж множеств* $(X_2, ..., X_{n-1}, X_n \backslash \{a_n\})$ *не является кортежем множеств образующих башни представлений* $\overline{f}$.

Доказательство. При $n = 2$ утверждение теоремы является следствием теоремы 2.5.2.

Пусть $n > 2$. Пусть $\overline{X}_{[1]}$ - кортеж множеств образующих башни расслоений $\overline{f}$. Пусть кортеж множеств $\overline{X}_{[1,n]}$ является базисом башни расслоений $\overline{f}_{[n]}$. Допустим для некоторого $a_n \in X_n$ существует слово

$$w_n = w_n(a_n, \overline{f}, (X_1, ..., X_{n-1}, X_n \setminus \{a_n\})) \tag{3.5.1}$$

Рассмотрим $a'_n \in A_n$, для которого слово

$$w'_n = w_n(a'_n, \overline{f}, \overline{X}_{[1]}) \tag{3.5.2}$$

зависит от $a_n$. Согласно определению 2.4.6, любое вхождение $a_n$ в слово $w'_n$ может быть заменено словом $w_n$. Следовательно, слово $w'_n$ не зависит от $a_n$, а кортеж множеств $(X_2, ..., X_{n-1}, X_n \setminus \{a_n\})$ является кортежем множеств образующих башни представлений $\overline{f}$. Следовательно, $\overline{X}_{[1]}$ не является базисом башни представлений $\overline{f}$. □

**Замечание 3.5.3.** Доказательство теоремы 3.5.2 даёт нам эффективный метод построения базиса башни представлений $\overline{f}$. Мы начинаем строить базис в самом нижнем слое. Когда базис построен в слое $i$, $i = 2, ..., n-1$, мы можем перейти к построению базиса в слое $i+1$. □

**Замечание 3.5.4.** Мы будем записывать базис также в виде

$$X_k = (x_k, x_k \in X_k) \quad k = 2, ..., n$$

Если базис - конечный, то мы будем также пользоваться записью

$$X_k = (x_{k\cdot i}, i \in I_k) = (x_{k\cdot 1}, ..., x_{k\cdot p_k}) \quad k = 2, ..., n$$

□

**Замечание 3.5.5.** Координаты $a_k \in A_k$ относительно базиса $\overline{X}_{[1]}$ башни представлений $\overline{f}$ определены неоднозначно. Так же как в замечании 2.5.5, мы рассмотрим **кортеж отношений эквивалентности** $\overline{\rho}_{[1]}(\overline{f})$**, порождённых башней представлений** $\overline{f}$. Если в процессе построений мы получим равенство двух $\Omega_k$-слов относительно заданного базиса, мы можем утверждать, что эти $\Omega_k$-слова равны, не заботясь об отношении эквивалентности $\rho_k(\overline{f})$. □

**Теорема 3.5.6.** *Автоморфизм башни представлений* $\overline{f}$ *отображает базис башни представлений* $\overline{f}$ *в базис.*

Доказательство. При $n = 2$ утверждение теоремы является следствием теоремы 3.5.6.

Пусть отображение $\overline{r}$ - автоморфизм башни представлений $\overline{f}$. Пусть кортеж множеств $\overline{X}_{[1]}$ - базис башни представлений $\overline{f}$. Пусть $\overline{X}'_{[1]} = \overline{r}_{[1]}(\overline{X}_{[1]})$.

Допустим кортеж множеств $\overline{X}'_{[1]}$ не является базисом. Согласно теореме 3.5.2 существуют $i$, $i = 2, ..., n$, и $a'_i \in X'_i$ такие, что кортеж множеств $\overline{X}''_{[1]}$[3.10]

[3.10] $X''_j = X'_j$, $j \neq i$, $X''_i = X'_i \setminus \{x'_i\}$

является кортежем множеств образующих башни представлений $\overline{f}$. Согласно теореме 3.3.3 отображение $\overline{r}^{-1}$ является автоморфизмом башни представлений $\overline{f}$. Согласно теореме 3.4.22 и определению 3.4.21, кортеж множеств $\overline{X}'''_{[1]}$[3.11] является кортежем множеств образующих представления $\overline{f}$. Полученное противоречие доказывает теорему. □

[3.11] $X_j''' = X_j$, $j \ne i$, $X_i''' = X_i \setminus \{x_i\}$

Глава 4

# Геометрия тела

## 4.1. Центр тела

**Определение 4.1.1.** Пусть $D$ - кольцо.[4.1] Множество $Z(D)$ элементов $a \in D$ таких, что

$$ax = xa \tag{4.1.1}$$

для всех $x \in D$, называется **центром кольца** $D$. □

**Теорема 4.1.2.** *Центр $Z(D)$ кольца $D$ является подкольцом кольца $D$.*

Доказательство. Непосредственно следует из определения 4.1.1. □

**Теорема 4.1.3.** *Центр $Z(D)$ тела $D$ является подполем тела $D$.*

Доказательство. Согласно теореме 4.1.2 достаточно проверить, что $a^{-1} \in Z(D)$, если $a \in Z(D)$. Пусть $a \in Z(D)$. Многократно применяя равенство (4.1.1), мы получим цепочку равенств

$$aa^{-1}x = x = xaa^{-1} = axa^{-1} \tag{4.1.2}$$

Из (4.1.2) следует

$$a^{-1}x = xa^{-1}$$

Следовательно, $a^{-1} \in Z(D)$.

□

**Определение 4.1.4.** Пусть $D$ - кольцо с единицей $e$.[4.2] Отображение

$$l : Z \to D$$

для которого $l(n) = ne$ будет гомоморфизмом колец, и его ядро является идеалом $(n)$, порождённым целым числом $n \geq 0$. Канонический инъективный гомоморфизм

$$Z/nZ \to D$$

является изоморфизмом между $Z/nZ$ и подкольцом в $D$. Если $nZ$ - простой идеал, то у нас возникает два случая.

- $n = 0$. $D$ содержит в качестве подкольца кольцо, изоморфное $Z$ и часто отождествляемое с $Z$. В этом случае мы говорим, что $D$ имеет **характеристику** 0.
- $n = p$ для некоторого простого числа $p$. $D$ имеет **характеристику** $p$, и $D$ содержит изоморфный образ $F_p = Z/pZ$.

□

---

[4.1] [1], стр. 84.

[4.2] Определение дано согласно определению из [1], стр. 84, 85.

**Теорема 4.1.5.** *Пусть $D$ - кольцо характеристики 0 и пусть $d \in D$. Тогда любое целое число $n \in Z$ коммутирует с $d$.*

Доказательство. Утверждение теоремы доказывается по индукции. При $n=0$ и $n=1$ утверждение очевидно. Допустим утверждение справедливо при $n=k$. Из цепочки равенств

$$(k+1)d = kd + d = dk + d = d(k+1)$$

следует очевидность утверждения при $n=k+1$. □

**Теорема 4.1.6.** *Пусть $D$ - кольцо характеристики 0. Тогда кольцо целых чисел $Z$ является подкольцом центра $Z(D)$ кольца $D$.*

Доказательство. Следствие теоремы 4.1.5. □

Пусть $D$ - тело. Если $D$ имеет характеристику 0, $D$ содержит в качестве подполя изоморфный образ поля $Q$ рациональных чисел. Если $D$ имеет характеристику $p$, $D$ содержит в качестве подполя изоморфный образ $F_p$. В обоих случаях это подполе будет называться простым полем. Так как простое поле является наименьшим подполем в $D$, содержащим 1 и не имеющим автоморфизмов, кроме тождественного, его обычно отождествляют с $Q$ или $F_p$, в зависимости от того, какой случай имеет место.

**Теорема 4.1.7.** *Пусть $D$ - тело характеристики 0 и пусть $d \in D$. Тогда для любого целого числа $n \in Z$*

$$n^{-1}d = dn^{-1} \tag{4.1.3}$$

Доказательство. Согласно теореме 4.1.5 справедлива цепочка равенств

$$n^{-1}dn = nn^{-1}d = d \tag{4.1.4}$$

Умножив правую и левую части равенства (4.1.4) на $n^{-1}$, получим

$$n^{-1}d = n^{-1}dnn^{-1} = dn^{-1} \tag{4.1.5}$$

Из (4.1.5) следует (4.1.3). □

**Теорема 4.1.8.** *Пусть $D$ - тело характеристики 0 и пусть $d \in D$. Тогда любое рациональное число $p \in Q$ коммутирует с $d$.*

Доказательство. Мы можем представить рациональное число $p \in Q$ в виде $p = mn^{-1}$, $m$, $n \in Z$. Утверждение теоремы следует из цепочки равенств

$$pd = mn^{-1}d = n^{-1}dm = dmn^{-1} = dp$$

основанной на утверждении теоремы 4.1.5 и равенстве (4.1.3). □

**Теорема 4.1.9.** *Пусть $D$ - тело характеристики 0. Тогда поле рациональных чисел $Q$ является подполем центра $Z(D)$ тела $D$.*

Доказательство. Следствие теоремы 4.1.8. □

## 4.2. Геометрия тела над полем

Мы можем рассматривать тело $D$ как векторное пространство над полем $F \subset Z(D)$. При этом мы будем пользоваться следующими соглашениями.

(1) Мы не будем для элемента тела $D$ пользоваться стандартными обозначениями для вектора, если мы будем рассматривать этот элемент как вектор над полем $F$. Однако мы будем пользоваться другим цветом для индекса при записи координат элемента тела $D$ как вектора над полем $F$.

(2) Так как $F$ - поле, то мы можем писать все индексы справа от корневой буквы.

Пусть $\overline{\overline{e}}$ - базис тела $D$ над полем $F$. Тогда произвольный элемент $a \in D$ можно представить в виде

(4.2.1) $$a = \overline{e}_i a^i \quad a^i \in F$$

Если размерность тела $D$ над полем $F$ бесконечна, то базис может быть либо счётным, либо его мощность может быть не меньше, чем мощность континуума. Если базис счётный, то на коэффициенты $a^i$ разложения (4.2.1) накладываются определённые ограничения. Если мощность множества $I$ континуум, то предполагается, что на множестве $I$ определена мера и сумма в разложении (4.2.1) является интегралом по этой мере.

**Замечание 4.2.1.** Поскольку в теле $D$ определена операция произведения, то мы можем рассматривать тело как алгебру над полем $F \subset Z(D)$. Для элементов базиса мы положим

(4.2.2) $$\overline{e}_i \overline{e}_j = \overline{e}_k C^k_{ij}$$

Коэффициенты $C^k_{ij}$ разложения (4.2.2) называются **структурными константами тела $D$ над полем $F$.**

Из равенств (4.2.1), (4.2.2) следует

(4.2.3) $$ab = \overline{e}_k C^k_{ij} a^i b^j$$

Из равенства (4.2.3) следует

(4.2.4) $$(ab)c = \overline{e}_k C^k_{ij} (ab)^i c^j = \overline{e}_k C^k_{ij} C^i_{mn} a^m b^n c^j$$

(4.2.5) $$a(bc) = \overline{e}_k C^k_{ij} a^i (bc)^j = \overline{e}_k C^k_{ij} a^i C^j_{mn} b^m c^n$$

Из ассоциативности произведения

$$(ab)c = a(bc)$$

и равенств (4.2.4) и (4.2.5) следует

(4.2.6) $$\overline{e}_k C^k_{ij} C^i_{mn} a^m b^n c^j = \overline{e}_k C^k_{ij} a^i C^j_{mn} b^m c^n$$

Так как векторы $a$, $b$, $c$ произвольны, а векторы $\overline{e}_k$ линейно независимы, то из равенства (4.2.6) следует

(4.2.7) $$C^k_{jn} C^j_{im} = C^k_{ij} C^j_{mn}$$

**Теорема 4.2.2.** *Координаты $a^j$ вектора $a$ являются тензором*

(4.2.8) $$a^j = A^j_i a'^i$$

Доказательство. Пусть $\overline{\overline{e}}'$ - другой базис. Пусть

$$\overline{e}'_i = \overline{e}_j A^j_i \tag{4.2.9}$$

преобразование, отображающее базис $\overline{\overline{e}}$ в базис $\overline{\overline{e}}'$. Так как вектор $a$ не меняется, то

$$a = \overline{e}'_i a'^i = \overline{e}_j a^j \tag{4.2.10}$$

Из равенств (4.2.9) и (4.2.10) следует

$$\overline{e}_j a^j = \overline{e}'_i a'^i = \overline{e}_j A^j_i a'^i \tag{4.2.11}$$

Так как векторы $\overline{e}_j$ линейно независимы, то равенство (4.2.8) следует из равенства (4.2.11). Следовательно, компоненты вектора являются тензором. □

**Теорема 4.2.3.** *Структурные константы тела $D$ над полем $F$ являются тензором*

$$A^l_k B'^k_{ij} A^{-1\cdot i}{}_n A^{-1\cdot j}{}_m = C^l_{nm} \tag{4.2.12}$$

Доказательство. Рассмотрим аналогичным образом преобразование произведения. Равенство (4.2.3) в базисе $\overline{e}'$ имеет вид

$$ab = \overline{e}'_k B'^k_{ij} a'^i b'^j \tag{4.2.13}$$

Подставив (4.2.8) и (4.2.9) в (4.2.13), получим

$$ab = \overline{e}_l A^l_k B'^k_{ij} a^n A^{-1\cdot i}{}_n b^m A^{-1\cdot j}{}_m \tag{4.2.14}$$

Из (4.2.3) и (4.2.14) следует

$$\overline{e}_l A^l_k B'^k_{ij} A^{-1\cdot i}{}_n a^n A^{-1\cdot j}{}_m b^m = \overline{e}_l C^l_{nm} a^n b^m \tag{4.2.15}$$

Так как векторы $a$ и $b$ произвольны, а векторы $\overline{e}_l$ линейно независимы, то равенство (4.2.12) следует из равенства (4.2.15). Следовательно, структурные константы являются тензором. □

Глава 5

# Квадратичное отображение тела

## 5.1. Билинейное отображение тела

**Теорема 5.1.1.** *Пусть $D$ - тело характеристики 0. Пусть $F$ - поле, которое является подкольцом центра тела $D$. Допустим $\overline{\overline{e}}$ - базис тела $D$ над полем $F$.* **Стандартное представление над полем $F$ билинейного отображения тела** *имеет вид*

(5.1.1) $$f(a,b) = {}_{1\cdot}f^{ijk}\ {}_{i}e\ a\ {}_{j}e\ b\ {}_{k}e + {}_{2\cdot}f^{ijk}\ {}_{i}e\ b\ {}_{j}e\ a\ {}_{k}e$$

*Выражение ${}_{t\cdot}f^{ijk}$ в равенстве (5.1.1) называется* **стандартной компонентой над полем $F$ билинейного отображения $f$**.

Доказательство. Следствие теоремы [6]-9.2.4. Здесь даны перестановки

(5.1.2) $$\sigma_1 = \begin{pmatrix} a & b \\ a & b \end{pmatrix} \quad \sigma_2 = \begin{pmatrix} a & b \\ b & a \end{pmatrix}$$

□

**Теорема 5.1.2.** *Пусть поле $F$ является подкольцом центра $Z(D)$ тела $D$. Пусть $\overline{e}_i$ - базис тела $D$ над полем $F$. Билинейное отображение*

$$g : D \times D \to D$$

*можно представить в виде $D$-значной билинейной формы над полем $F$*

(5.1.3) $$g(a,b) = a^{i}\ b^{j}\ {}_{ij}g$$

*где*

$$a = a^{i}\ {}_{i}\overline{e}$$

$$b = b^{j}\ {}_{j}\overline{e}$$

(5.1.4) $${}_{ij}g = g({}_{i}\overline{e}, {}_{j}\overline{e})$$

*и величины ${}_{ij}g$ являются координатами $D$-значного ковариантного тензора над полем $F$.*

Доказательство. Следствие теоремы [6]-9.2.5. □

Матрица $G = |{}_{ij}g|$ называется **матрицей билинейной функции**. Если матрица $G$ невырождена, то билинейная функция $g$ называется **невырожденной**.

**Теорема 5.1.3.** *Билинейное отображение*

$$g : D \times D \to D$$

*симметрично тогда и только тогда, когда матрица $G$ симметрична*

$$_{ij}g = {}_{ji}g$$

Доказательство. Следствие теоремы [6]-9.2.6. □

**Теорема 5.1.4.** *Билинейное отображение*

$$g : D \times D \to D$$

*косо симметрично тогда и только тогда, когда матрица $G$ косо симметрична*

$$_{ij}g = -{}_{ji}g$$

Доказательство. Следствие теоремы [6]-9.2.7. □

**Теорема 5.1.5.** *Компоненты билинейного отображения*

$$g : D \times D \to D$$

*и его матрица над полем $F$ связаны равенством*

(5.1.5) $$_{pq}g = ({}_{1\cdot}g^{ijk}\ C^s_{ip}\ C^l_{sj}\ C^t_{lq}\ C^r_{tk}\ + {}_{2\cdot}g^{ijk}\ C^s_{iq}\ C^l_{sj}\ C^t_{lp}\ C^r_{tk})\ \overline{e}_r$$

(5.1.6) $$_{pq}g^r = {}_{1\cdot}g^{ijk}\ C^s_{ip}\ C^l_{sj}\ C^t_{lq}\ C^r_{tk} + {}_{2\cdot}g^{ijk}\ C^s_{iq}\ C^l_{sj}\ C^t_{lp}\ C^r_{tk}$$

Доказательство. Следствие теоремы [6]-9.2.8. Для доказательства утверждения достаточно подставить (5.1.2) в равенства [6]-(9.2.13), [6]-(9.2.14). □

**Теорема 5.1.6.** *Если билинейное отображение*

$$g : D \times D \to D$$

*имеет матрицу $G = \|{}_{ij}g\|$, то существует билинейное отображение*

$$g' : D \times D \to D$$

*которое имеет матрицу $G' = \|{}_{ij}g'\|$, $_{ij}g' = {}_{ji}g$.*

Доказательство. Положим

(5.1.7) $$_{1\cdot}g'^{ijk} = {}_{2\cdot}g^{ijk} \quad {}_{2\cdot}g'^{ijk} = {}_{1\cdot}g^{ijk}$$

Из равенств (5.1.6), (5.1.7) следует

$$\begin{aligned} {}_{qp}g'^r &= {}_{1\cdot}g'^{ijk}\ C^s_{iq}\ C^l_{sj}\ C^t_{lp}\ C^r_{tk} + {}_{2\cdot}g'^{ijk}\ C^s_{ip}\ C^l_{sj}\ C^t_{lq}\ C^r_{tk} \\ &= {}_{2\cdot}g^{ijk}\ C^s_{iq}\ C^l_{sj}\ C^t_{lp}\ C^r_{tk} + {}_{1\cdot}g^{ijk}\ C^s_{ip}\ C^l_{sj}\ C^t_{lq}\ C^r_{tk} \\ &= {}_{pq}g^r \end{aligned}$$

□

## 5.2. Квадратичное отображение тела

**Определение 5.2.1.** Пусть $D$ - тело. Отображение

$$h : D \to R$$

называется **квадратичным**, если существует билинейное отображение

$$g : D \times D \to D$$

такое, что

$$h(a) = g(a, a)$$

Билинейное отображение $g$ называется отображением, ассоциированным отображению $h$. □

Отображение $g$ определено неоднозначно для заданного отображения $h$. Однако согласно теореме 5.1.6, мы всегда можем положить, что билинейное отображение $g$ симметрично

$$g(a,b) = g(b,a)$$

Действительно, если билинейное отображение $g_1$, ассоциированное отображению $h$, не симметрично, то мы положим

$$g(a,b) = (g_1(a,b) + g_1(b,a))/2$$

Так как билинейное отображение однородно степени 1 над полем $R$ по каждой переменной, квадратичное отображение однородно степени 2 над полем $R$.

**Теорема 5.2.2.** *Пусть $D$ - тело характеристики* 0*. Пусть $F$ - поле, которое является подкольцом центра тела $D$. Допустим $\overline{\overline{p}}$ - базис тела $D$ над полем $F$.* **Стандартное представление квадратичного отображения тела над полем** *$F$ имеет вид*

$$f(a) = f^{ijk}\ {}_{i}\overline{p}\ a\ {}_{j}\overline{p}\ a\ {}_{k}\overline{p} \tag{5.2.1}$$

*Выражение $f^{ijk}$ в равенстве* (5.2.1) *называется* **стандартной компонентой квадратичного отображения** *$f$* **над полем** *$F$.*

Доказательство. Следствие теоремы 5.1.1. □

**Теорема 5.2.3.** *Пусть поле $F$ является подкольцом центра $Z(D)$ тела $D$. Пусть $\overline{\overline{e}}$ - базис тела $D$ над полем $F$. Квадратичного отображение $f$ можно представить в виде $D$-значной квадратичной формы над полем $F$*

$$f(a) = a^i a^j\ {}_{ij}f \tag{5.2.2}$$

*где*

$$\overline{a} = a^i\ {}_{i}\overline{e}$$

$${}_{ij}f = g({}_{i}\overline{e}, {}_{j}\overline{e}) \tag{5.2.3}$$

*и $g$ - ассоциированное билинейное отображение. Величины ${}_{ij}f$ являются координатами $D$-значного ковариантного тензора над полем $F$.*

Доказательство. Для выбранного билинейного отображения $g$, ассоциированного отображению $f$, утверждение теоремы является следствием теоремы 5.1.2. Нам надо доказать независимость координат квадратичного отображения от выбора ассоциированного отображения $g$.

Допустим, $g_1$, $g_2$ - билинейные отображения, ассоциированные отображению $f$. Тогда

$${}_{ii}f_1 = g_1({}_{i}\overline{e}, {}_{i}\overline{e}) = g_2({}_{i}\overline{e}, {}_{i}\overline{e}) = {}_{ii}f_2 \tag{5.2.4}$$

Для произвольного $a \in D$

$$a^i\ a^j\ {}_{ij}f_1 = a^i\ a^j\ {}_{ij}f_2$$

Если для данных $i$, $j$ мы рассмотрим $a \in D$, для которого $a_i = a_j = 1$, то мы получим

$${}_{ii}f_1 + {}_{ij}f_1 + {}_{ji}f_1 + {}_{jj}f_1 = {}_{ii}f_2 + {}_{ij}f_2 + {}_{ji}f_2 + {}_{jj}f_2 \tag{5.2.5}$$

Из равенств (5.2.4), (5.2.5) следует

$$_{ij}f_1 + {}_{ji}f_1 = {}_{ij}f_2 + {}_{ji}f_2$$

Следовательно, если ассоциированное билинейное отображение симметрично, то оно определенно однозначно. □

Матрица $F = |_{ij}f|$ называется **матрицей квадратичного отображения**. Если матрица $F$ невырождена, то квадратичное отображение $f$ называется **невырожденным**. Ранг матрицы $F$ называется **рангом квадратичного отображения** $f$.

**Теорема 5.2.4.** *Отображение*

$$f : D \to D$$

*квадратичное над полем $F$, определено как отображение над телом тогда и только тогда, когда*

$$_{pq}f = f^{ijk}\ C^s_{ip}\ C^l_{sj}\ C^t_{lq}\ C^r_{tk}\ {}_r\overline{p} \tag{5.2.6}$$

$$_{pq}f^r = f^{ijk}\ C^s_{ip}\ C^l_{sj}\ C^t_{lq}\ C^r_{tk} \tag{5.2.7}$$

Доказательство. Подставим

$$\overline{a} = a^j\ {}_j\overline{p}$$

в равенство (5.2.1). Тогда равенство (5.2.1) примет вид

$$\begin{aligned} f(a) &= f^{ijk}\ {}_i\overline{p}\ a^p\ {}_p\overline{p}\ {}_j\overline{p}\ a^q\ {}_q\overline{p}\ {}_k\overline{p} \\ &= a^p\ a^q\ f^{ijk}\ {}_i\overline{p}\ {}_p\overline{p}\ {}_j\overline{p}\ {}_q\overline{p}\ {}_k\overline{p} \\ &= a^p\ a^q\ f^{ijk}\ C^s_{ip}\ C^l_{sj}\ C^t_{lq}\ C^r_{tk}\ {}_r\overline{p} \end{aligned} \tag{5.2.8}$$

Равенство (5.2.6) следует из сравнения равенств (5.2.8) и (5.2.2). □

Из теоремы 5.2.3 следует, что координаты $_{pq}f$ образуют тензор. В теореме [6]-9.1.13 описана структура допустимых преобразований.

**Теорема 5.2.5.** *Пусть поле $F$ является подкольцом центра $Z(D)$ тела $D$. Пусть $\overline{\overline{p}}$ - базис тела $D$ над полем $F$. Пусть $C^k_{ij}$ - структурные константы тела $D$ над полем $F$. Если*

$$\det(a^i\ (C^k_{ij}\ +\ C^k_{ji})) \neq 0$$

*то уравнение*

$$ax + xa = b \tag{5.2.9}$$

*имеет единственного решение.*

*Если определитель равен 0, то условием существования решения является $F$-линейная зависимость вектора $b$ от векторов $a^i\ (C^k_{ij}\ +\ C^k_{ji})\ {}_k\overline{p}$ . В этом случае, уравнение (5.2.9) имеет бесконечно много решений. В противном случае уравнение не имеет решений.*

Доказательство. Из равенств (4.2.3), (5.2.9) следует

$$a^i\ x^j\ C^k_{ij}\ {}_k\overline{p} + x^i\ a^j\ C^k_{ij}\ {}_k\overline{p} = b^k\ {}_k\overline{p} \tag{5.2.10}$$

Так как векторы $_k\overline{p}$ линейно независимы, то из равенства (5.2.10) следует, что координаты $x^j$ удовлетворяют системе линейных уравнений

$$x^j\ a^i\ (C^k_{ij}\ +\ C^k_{ji}) = b^k \tag{5.2.11}$$

Следовательно, если определитель системы (5.2.11) не равен 0, то система имеет единственное решение.

Аналогично, если определитель равен 0, то условием существования решения является равенство ранга расширенной матрицы и ранга матрицы системы линейных уравнений. □

**Теорема 5.2.6.** *Пусть поле $F$ является подкольцом центра $Z(D)$ тела $D$. Существует квадратичное отображение $f$ такое, что можно найти базис $\overline{\overline{p}}$ над полем $F$, в котором отображение $f$ представлено как сумма квадратов.*

Доказательство. Приведение квадратичной формы к диагональному виду выполнено аналогично тому как это сделано в [4], с. 169 - 172. Рассматриваемые преобразования переменных связаны с соответствующими преобразованиями базиса. Нас не интересует, является ли эти преобразования базиса допустимыми. Причина состоит в том, что базис, выбранный вначале может быть не связан с каноническим базисом допустимым преобразованием. Однако для нас важно, чтобы окончательное представление отображения $f$ было допустимым.

Мы докажем утверждение индукцией по числу переменных в представлении квадратичного отображения. Для этого мы будем записывать квадратичное отображение $f$ в виде отображения от $n$ переменных $a^1$, ..., $a^n$

$$f(a) = f(a^1, ..., a^n) = a^i a^j \ {}_{ij}f \tag{5.2.12}$$

Представление квадратичного отображения в виде (5.2.12) называется **квадратичной формой**.

Утверждение очевидно, если квадратичная форма зависит от одной переменной $a^1$, так как в этом случае квадратичная форма имеет вид

$$f(a^1) = (a^1)^2 {}_{11}f$$

Допустим утверждение теоремы справедливо для квадратичной формы от $n-1$ переменной.

(1) Пусть квадратичная форма содержит квадраты переменных. Не нарушая общности, положим ${}_{11}f \neq 0$. Рассмотрим отображение

$$f_1(a^1, ..., a^n) = f(a^1, ..., a^n) - {}_{11}f^{-1} \ (a^j \ {}_jh)^2 \tag{5.2.13}$$

Мы выбираем значения ${}_jh$ таким образом, чтобы удовлетворить равенство

$$2 \ {}_{11}f \ {}_{1j}f = {}_1h \ {}_jh + {}_jh \ {}_1h \tag{5.2.14}$$

Когда $j = 1$, из равенства (5.2.14) следует

$${}_1h = \ {}_{11}f \tag{5.2.15}$$

Когда $j > 1$, из равенств (5.2.13), (5.2.15) следует

$$2 \ {}_{11}f \ {}_{1j}f = {}_{11}f \ {}_jh + {}_jh \ {}_{11}f \tag{5.2.16}$$

Из теоремы 5.2.5 следует, что если

$$\det({}_{11}f^i \ (C^k_{ij} \ + \ C^k_{ji})) \neq 0 \tag{5.2.17}$$

то уравнение (5.2.16) имеет единственное решение. Мы предположим, что условие (5.2.17) выполнено, и рассмотрим замену переменных

$$b^1 = a^j\ {}_j h$$
$$b^i = a^i \qquad i > 1$$

Следовательно, отображение $f_1(a^1, ..., a^n)$ не зависит от $b^1$. Согласно определению (5.2.13), мы представили отображение $f(a^1, ..., a^n)$ в виде суммы квадрата переменной $b^1$ и квадратичного отображения от $b^1$ не зависящего.

$$f(a^1, ..., a^n) = {}_{11}f^{-1}\ (b^1)^2 + f_1(b^2, ..., b^n)$$

Так как форма $f_1$ зависит от $n-1$ переменной, то мы можем найти преобразование переменных, приводящее эту форму к сумме квадратов переменных.

(2) Если ${}_{11}f = ... = {}_{nn}f = 0$, то нам нужно дополнительное линейное преобразование, приводящее к появлению квадратов переменных. Положим ${}_{12}f \neq 0$. Совершим линейное преобразование

$$a^1 = b^1 - b^2$$
$$a^2 = b^1 + b^2$$
$$a^i = b^i \quad i > 2$$

В результате этого преобразования слагаемое

$$2\ a^1\ a^2\ {}_{12}f$$

примет вид

$$2\ (b^1)^2\ {}_{12}f - 2\ (b^2)^2\ {}_{12}f$$

Мы получили квадратичную форму, рассмотренную в случае (1).

□

Если все коэффициенты квадратичной формы в каноническом представлении действительны, то $f(a) \in R$ для любого $a \in D$.

Допустим, мы выбрали квадратичную форму $f$, принимающую действительные значения.

**Определение 5.2.7.** Пусть $D$ - тело. Квадратичное отображение

$$f : D \times D \to D$$

называется **положительно определённым**, если

$$f(a) \in R$$

для любого $a \in D$ и при этом выполнены условия

$$f(a) \geq 0$$
$$f(a) = 0 => a = 0$$

Положительно определённое квадратичное отображение называется **эвклидовой метрикой на теле** $D$. Симметричное ассоциированное билинейное отображение $g$ называется **эвклидовым скалярным произведением на теле** $D$. □

Если квадратичное отображение $f$ не является положительно определённым, то это отображение называется **псевдоэвклидовой метрикой на теле** $D$. Симметричное ассоциированное билинейное отображение $g$ называется **псевдоэвклидовым скалярным произведением на теле** $D$. Пусть $\overline{\overline{e}}$ - базис, относительно которого форма $f$ представлена в виде суммы квадратов. Существуют вектора базиса $\overline{e}_i$, для которых $f(\overline{e}_i) < 0$. Для произвольного $a = a^i\ \overline{e}_i$ определим операцию сопряжения равенством

$$a^\star = sign(f(\overline{e}_i))a^i\ \overline{e}_i$$

Выражение $a^\star$ называется **эрмитовым сопряжением**.

Рассмотрим биаддитивное отображение

$$g^\star(a, b) = g(a, b^\star)$$

$$f^\star(a) = g^\star(a, a)$$

Отображение $f^\star(a)$ положительно определено и называется **эрмитовой метрикой на теле** $D$. Симметричное ассоциированное отображение $g^\star$ называется **эрмитовым скалярным произведением на теле** $D$.

Если отображение

$$a \to a^\star$$

линейно, то эрмитово скалярное произведение совпадает с эвклидовым скалярнып произведением. Аналогичное утверждение верно для метрики.

Глава 6

# $D$-аффинное пространство

## 6.1. $D$-аффинное пространство

**Определение 6.1.1.** Рассмотрим множество точек $\overset{\circ}{V}$ и множество свободных векторов $\overline{V}$.[6.1] Множество $\overset{\circ}{V}$ удовлетворяет следующим аксиомам

(1) Существует по крайней мере одна точка
(2) Каждой паре точек $(A,B)$ поставлен в соответствие один и только один вектор. Этот вектор мы будем обозначать $\overline{AB}$. Вектор $\overline{AB}$ имеет начало в точке $A$ и конец в точке B.
(3) Для каждой точки $A$ и любого вектора $\overline{a}$ существует одна и только одна точка $B$ такая, что $\overline{AB}=\overline{a}$ . Мы будем также пользоваться записью[6.2]

$$B=A+\overline{a} \tag{6.1.1}$$

(4) (Аксиома параллелограмма.) Если $\overline{AB}=\overline{CD}$, то $\overline{AC}=\overline{BD}$.

□

**Определение 6.1.2.** Пусть даны векторы $\overline{a}$ и $\overline{b}$. Пусть $A\in\overset{\circ}{V}$ - произвольная точка. Пусть $B\in\overset{\circ}{V}$, $B=A+\overline{a}$. Пусть $C\in\overset{\circ}{V}$, $C=B+\overline{b}$. Вектор $\overline{AC}$ называется суммой векторов $\overline{a}$ и $\overline{b}$

$$\overline{a}+\overline{b}=\overline{AC} \tag{6.1.2}$$

□

**Теорема 6.1.3.** *Вектор $\overline{AA}$ является нулём по отношению к операции сложения и не зависит от выбора точки $A$. Вектор $\overline{AA}$ называется нуль-вектором и мы полагаем $\overline{AA}=\overline{0}$.*

Доказательство. Мы можем записать правило сложения (6.1.2) в виде

$$\overline{AB}+\overline{BC}=\overline{AC} \tag{6.1.3}$$

Если $B=C$, то из равенства (6.1.3) следует

$$\overline{AB}+\overline{BB}=\overline{AB} \tag{6.1.4}$$

Если $C=A$, $B=D$, то из аксиомы (4) определения 6.1.1 следует

$$\overline{AA}=\overline{BB} \tag{6.1.5}$$

Утверждение теоремы следует из равенств (6.1.4) и (6.1.5). □

[6.1] Я написал определения и теоремы в этом разделе согласно определению аффинного пространства в [3], с. 86 - 93.
[6.2] [13], с. 9.

**Теорема 6.1.4.** *Пусть* $\overline{a}=\overline{AB}$ *. Тогда*

$$\overline{BA}=-\overline{a}\tag{6.1.6}$$

*и это равенство не зависит от выбора точки* $A$*.*

Доказательство. Из равенства (6.1.3) следует

$$\overline{AB}+\overline{BA}=\overline{AA}=\overline{0}\tag{6.1.7}$$

Равенство (6.1.6) следует из равенства (6.1.7). Применяя аксиому (4) определения 6.1.1 к равенству $\overline{AB}=\overline{CD}$ получим $\overline{AC}=\overline{BD}$, или, что то же, $\overline{BD}=\overline{AC}$. В силу аксиомы (4) определения 6.1.1 снова следует $\overline{BA}=\overline{DC}$ . Следовательно, равенство (6.1.6) не зависит от выбора точки $A$. □

**Теорема 6.1.5.** *Сумма векторов* $\overline{a}$ *и* $\overline{b}$ *не зависит от выбора точки* $A$*.*

Доказательство. Пусть $\overline{a}=\overline{AB}=\overline{A'B'}$. Пусть $\overline{b}=\overline{BC}=\overline{B'C'}$. Пусть

$$\overline{AB}+\overline{BC}=\overline{AC}$$
$$\overline{A'B'}+\overline{B'C'}=\overline{A'C'}$$

Согласно аксиоме (4) определения 6.1.1

$$\overline{A'A}=\overline{B'B}=\overline{C'C}\tag{6.1.8}$$

Применяя аксиому (4) определения 6.1.1 к крайним членам равенства (6.1.8), получаем

$$\overline{A'C'}=\overline{AC}\tag{6.1.9}$$

Из равенства (6.1.9) следует утверждение теоремы. □

**Теорема 6.1.6.** *Сложение векторов ассоциативно.*

Доказательство. Пусть $\overline{a}=\overline{AB}$ , $\overline{b}=\overline{BC}$ , $\overline{c}=\overline{CD}$ . Из равенства

$$\begin{array}{ccccc}\overline{a}&+&\overline{b}&=&\overline{AC}\\ \overline{AB}&+&\overline{BC}&=&\overline{AC}\end{array}$$

следует

$$\begin{array}{ccccc}(\overline{a}+\overline{b})&+&\overline{c}&=&\overline{AD}\\ \overline{AC}&+&\overline{CD}&=&\overline{AD}\end{array}\tag{6.1.10}$$

Из равенства

$$\begin{array}{ccccc}\overline{b}&+&\overline{c}&=&\overline{BD}\\ \overline{BC}&+&\overline{CD}&=&\overline{BD}\end{array}$$

следует

$$\begin{array}{ccccc}\overline{a}&+&(\overline{b}+\overline{c})&=&\overline{AD}\\ \overline{AB}&+&\overline{BD}&=&\overline{AD}\end{array}\tag{6.1.11}$$

Из сравнения равенств (6.1.10) и (6.1.11) следует ассоциативность сложения. □

**Теорема 6.1.7.** *На множестве* $\overline{V}$ *определена структура абелевой группы.*

Доказательство. Из теорем 6.1.3, 6.1.4, 6.1.5, 6.1.6 следует, что операция сложение векторов определяет группу.

Пусть $\overline{a} = \overline{AB}$ , $\overline{b} = \overline{BC}$ .

$$\begin{array}{cccc} \overline{a} & + & \overline{b} & = \overline{AC} \\ \overline{AB} & + & \overline{BC} & = \overline{AC} \end{array} \tag{6.1.12}$$

Согласно аксиоме (3) определения 6.1.1 существует точка $D$ такая, что

$$\overline{b} = \overline{AD} = \overline{BC}$$

Согласно аксиоме (4) определения 6.1.1

$$\overline{AB} = \overline{DC} = \overline{a}$$

Согласно определению суммы векторов

$$\begin{array}{cccc} \overline{AD} & + & \overline{DC} & = \overline{AC} \\ \overline{b} & + & \overline{a} & = \overline{AC} \end{array} \tag{6.1.13}$$

Из сравнения равенств (6.1.12) и (6.1.13) следует коммутативность сложения. □

**Теорема 6.1.8.** *Отображение*

$$\overline{V} \to \overset{\circ}{V}{}^{\star} \tag{6.1.14}$$

*определённое равенством* (6.1.1)*, является однотранзитивным представлением абелевой группы* $\overline{V}$.

Доказательство. Аксиома (3) определения 6.1.1 определяет отображение (6.1.14). Из теоремы 6.1.5 следует, что отображение (6.1.14) является представлением. Эффективность представления следует из теоремы 6.1.3 и аксиомы (2) определения 6.1.1. Из аксиомы (2) определения 6.1.1, следует также, что представление транзитивно. Эффективное и однотранзитивное представление однотранзитивно. □

## 6.2. Аффинная структура на множестве

**Определение 6.2.1.** Ковариантное однотранзитивное $\star T$-представление абелевой группы $G$ в множестве $M$ называется **аффинной структурой на множестве** $M$. Элементы множества $M$ называются точками. Образ точки $A \in M$ при преобразовании, порождённом $a \in G$, мы обозначим $A + a$. □

**Теорема 6.2.2.** *Аффинная структура на множестве $M$ удовлетворяет аксиомам определения 6.1.1.*

Доказательство. Мы полагаем, что множество $M$ не пусто, поэтому аксиома (1) выполнена. Поскольку $a \in G$ порождает преобразование множества, то для любого $A \in M$ однозначно определён $B \in M$ такое, что $B = A + a$. Это утверждение доказывает аксиому (3). Поскольку представление однотранзитивно, то для любых $A$, $B \in M$ существует единственное $a \in G$ такое, что $B = A + a$. Это утверждение позволяет ввести обозначение $\overline{AB} = a$ , а также доказывает аксиому (2). Чтобы доказать теорему, нам осталось доказать, что из утверждения теоремы следует аксиома (4).

Пусть

$$\overline{AB}=\overline{CD}=a \tag{6.2.1}$$

Из равенства (6.2.1) следует

$$B=A+\overline{AB}=A+a \tag{6.2.2}$$

$$D=C+\overline{CD}=C+a \tag{6.2.3}$$

Пусть

$$\overline{AC}=b \quad C=A+b \tag{6.2.4}$$

$$\overline{BD}=c \quad D=B+c \tag{6.2.5}$$

Из равенств (6.2.4), (6.2.3) следует, что

$$D=A+b+a \tag{6.2.6}$$

Из равенств (6.2.2), (6.2.5) следует, что

$$D=A+a+c \tag{6.2.7}$$

Поскольку группа $G$ коммутативна, то из равенства (6.2.6) следует

$$D=A+a+b \tag{6.2.8}$$

Из равенств (6.2.7), (6.2.8) следует, что

$$A+a+b=A+a+c \tag{6.2.9}$$

Поскольку представление группы $G$ однотранзитивно, то из равенства (6.2.9) следует

$$b=c \tag{6.2.10}$$

Из равенств (6.2.2), (6.2.3), (6.2.10) следует, что

$$\overline{AC}=\overline{BD} \tag{6.2.11}$$

Следовательно, мы доказали аксиому (4). □

Несмотря на то, что теорема 6.2.2 утверждает, что аффинная структура на множестве $M$ удовлетворяет аксиомам определения 6.1.1, определение 6.2.1 не эквивалентно определению 6.1.1. В определении 6.1.1 предполагается, что на множестве $G$ определена структура векторного пространства. Поэтому, если мы хотим изучать аффинное пространство, мы должны рассмотреть башню представлений: представление поля $F$ в абелевой группе $G$ и представление абелевой группы $G$ в множестве $M$. Очевидно, что мы можем рассмотреть алгебру с делением $D$ вместо поля $F$. Таким образом, мы приходим к следующему определению аффинного пространства.

**Определение 6.2.3.** Пусть $D$ - тело характеристики 0, $\overline{V}$ - абелева группа и $\overset{\circ}{V}$ - множество. **$D\star$-аффинное пространство** - это башня представлений

$$\overrightarrow{V}:D\xrightarrow[*]{f_{1,2}}\overline{V}\xrightarrow[*]{f_{2,3}}\overset{\circ}{V} \qquad \begin{array}{l} f_{1,2}(d):\overline{v}\rightarrow d\,\overline{v}\\ f_{2,3}(\overline{v}):A\rightarrow A+\overline{v}\end{array}$$

где $f_{1,2}$ - эффективное $T\star$-представление тела $D$ в абелевой группе $\overline{V}$ и $f_{2,3}$ - однотранзитивное $\star T$-представление абелевой группы $\overline{V}$ в множестве $\overset{\circ}{V}$. □

**Определение 6.2.4.** Пусть $D$ - тело характеристики 0, $\overline{V}$ - абелева группа и $\overset{\circ}{V}$ - множество. $\star D$**-аффинное пространство** - это башня представлений

$$\overrightarrow{V}: D \xrightarrow[*]{f_{1,2}} \overline{V} \xrightarrow[*]{f_{2,3}} \overset{\circ}{V} \qquad \begin{array}{l} f_{1,2}(d): \overline{v} \to \overline{v}\, d \\ f_{2,3}(\overline{v}): A \to A + \overline{v} \end{array}$$

где $f_{1,2}$ - эффективное $\star T$-представление тела $D$ в абелевой группе $\overline{V}$ и $f_{2,3}$ - однотранзитивное $\star T$-представление абелевой группы $\overline{V}$ в множестве $\overset{\circ}{V}$. □

## 6.3. Однотранзитивное представление

Рассмотрим однотранзитивное представление $f$ $\Omega_1$-алгебры $A$ в $\Omega_2$-алгебре $M$. Мы можем отождествить $a \in A$ с соответствующим преобразованием $f(a)$. В этой разделе мы будем обозначать действие $a \in A$ на $v \in M$ посредством выражения

$$m + a$$

**Теорема 6.3.1.** *Если представление $f$ $\Omega_1$-алгебры $A$ в $\Omega_2$-алгебре $M$ однотранзитивно, то базисом представления $f$ является произвольная точка $v_0 \in M$. Мощность множества $M$ равна мощности множества $A$.*

Доказательство. Выберем произвольный $v_0 \in M$. Для любого $v \in M$ существует единственный

$$a = g(v_0, v) \tag{6.3.1}$$

такой, что

$$v_0 + a = v \tag{6.3.2}$$

Выражение (6.3.2) является $\Omega_2$-словом. Следовательно, $v_0$ является базисом представления $f$. Равенство (6.3.2) устанавливает взаимно однозначное соответствие между множествами $A$ и $M$. Следовательно, мощности множеств $A$ и $M$ равны. □

Для отображения $g$ мы будем пользоваться обозначением

$$a = \overline{v_0\, v} \tag{6.3.3}$$

и запишем равенство (6.3.2) в следующем виде

$$v_0 + \overline{v_0\, v} = v \tag{6.3.4}$$

**Теорема 6.3.2.** *Пусть представление $f$ $\Omega_1$-алгебры $A$ в $\Omega_2$-алгебре $M$ однотранзитивно. Если мы выбрали базис $v_0$ представления $f$, то соответствие (6.3.3), (6.3.4) позволяет нам перенести структуру $\Omega_1$-алгебры в $\Omega_2$-алгебру $M$.*[6.3]

Доказательство. Рассмотрим операцию $\omega \in \Omega_1(n)$. Для произвольных $v_1$, ..., $v_n \in M$ мы положим

$$v_1 ... v_n \omega = v_0 + \overline{v_0\, v_1} ... \overline{v_0\, v_n} \omega \tag{6.3.5}$$

□

[6.3] Важно иметь в виду, что структура $\Omega_1$-алгебры на множестве $M$ зависит от выбора точки $v_0$. Это замечание также справедливо для всех последующих теорем, которые мы рассмотрим ниже в этом разделе. Чтобы подчеркнуть, что рассматриваемая $\Omega_1$ алгебра зависит от выбора точки $O$, мы будем пользоваться обозначением $(M, v_0)$ .

**Теорема 6.3.3.** *Пусть представление $f$ $\Omega_1$-алгебры $A$ в $\Omega_2$-алгебре $M$ однотранзитивно. Если мы выбрали базис $v_0$ представления $f$, то произвольному эндоморфизму $f$ $\Omega_1$-алгебры $A$ соответствует эндомофизм*

$$f' : v \in M \to v_0 + f(\overline{v_0\, v}) \in M \tag{6.3.6}$$

*$\Omega_1$-алгебры $(M, v_0)$.*

Доказательство. Из определения отображения $f'$ (6.3.6) и равенства (6.3.4) непосредственно следует

$$f'(v_0 + \overline{v_0\, v}) = v_0 + f(\overline{v_0\, v}) \tag{6.3.7}$$

$$v_0 + \overline{v_0\, f'(v)} = v_0 + f(\overline{v_0\, v}) \tag{6.3.8}$$

Из равенства (6.3.8) следует

$$\overline{v_0\, f'(v)} = f(\overline{v_0\, v}) \tag{6.3.9}$$

Пусть $\omega \in \Omega_1(n)$. Из равенства (6.3.5) следует

$$f'(v_1 ... v_n \omega) = f'(v_0 + \overline{v_0\, v_1} ... \overline{v_0\, v_n} \omega) \tag{6.3.10}$$

Из равенств (6.3.10), (6.3.7) следует

$$f'(v_1 ... v_n \omega) = v_0 + f(\overline{v_0\, v_1} ... \overline{v_0\, v_n} \omega) \tag{6.3.11}$$

Так как отображение $f$ является эндоморфизмом, то из равенства (6.3.11) следует

$$f'(v_1 ... v_n \omega) = v_0 + f(\overline{v_0\, v_1}) ... f(\overline{v_0\, v_n}) \omega \tag{6.3.12}$$

Из равенств (6.3.12), (6.3.9) следует

$$f'(v_1 ... v_n \omega) = v_0 + \overline{v_0\, f'(v_1)} ... \overline{v_0\, f'(v_n)} \omega \tag{6.3.13}$$

Из равенств (6.3.13), (6.3.5) следует

$$f'(v_1 ... v_n \omega) = f'(v_1) ... f'(v_n) \omega \tag{6.3.14}$$

Из равенства (6.3.14) следует, что отображение $f'$ является эндоморфизмом $\Omega_1$-алгебры $(M, v_0)$. □

**Теорема 6.3.4.** *Пусть представление $f$ $\Omega_1$-алгебры $A$ в $\Omega_2$-алгебре $M$ однотранзитивно. Если мы выбрали базис $v_0$ представления $f$, то соответствие (6.3.3), (6.3.4) позволяет нам перенести структуру $\Omega_2$-алгебры в $\Omega_1$-алгебру $A$.*

Доказательство. Рассмотрим операцию $\omega \in \Omega_2(n)$. Для произвольных $a_1$, ..., $a_n \in A$ существуют $v_1$, ..., $v_n \in M$ такие, что

$$a_1 = \overline{v_0\, v_1} \quad ... \quad a_n = \overline{v_0\, v_n}$$

Мы положим

$$a_1 ... a_n \omega = \overline{v_0\, v_1} ... \overline{v_0\, v_n} \omega = \overline{v_0\, v_1 ... v_n \omega} \tag{6.3.15}$$

□

**Теорема 6.3.5.** *Рассмотрим башню $(\overline{V},\overline{f})$ представлений $\overline{\Omega}$-алгебр. Пусть представление $f_{i+1,i+2}$ однотранзитивно, и $a_{i+2\cdot 0}$ - базис представления $f_{i+1,i+2}$. Тогда представлению $f_{i,i+1}$ $\Omega_i$-алгебры $A_i$ в $\Omega_{i+1}$-алгебре $A_{i+1}$ соответствует представление*

$$\begin{array}{c}f'_{i,i+2}:A_i\xrightarrow{*}A_{i+2}\\ f'_{i,i+2}(a_i):a_{i+2}\to a_{i+2\cdot 0}+f_{i,i+1}(a_i)(\overline{a_{i+2\cdot 0}\,a_{i+2}})\end{array}\tag{6.3.16}$$

*$\Omega_i$-алгебры $A_i$ в $\Omega_{i+2}$-алгебре $A_{i+2}$.*

Доказательство. Поскольку отображение $f_{i,i+1}(a_i)$ является эндоморфизмом $\Omega_{i+1}$-алгебры $A_{i+1}$, то, согласно теореме 6.3.3, отображение $f'_{i,i+2}(a_i)$ является эндоморфизмом $\Omega_{i+1}$-алгебры $(A_{i+2},a_{i+2\cdot 0})$. □

Из равенства (6.3.16) следует

$$f'_{i,i+2}(a_i)(a_{i+2\cdot 0}+\overline{a_{i+2\cdot 0}\,a_{i+2}})=a_{i+2\cdot 0}+f_{i,i+1}(a_i)(\overline{a_{i+2\cdot 0}\,a_{i+2}})\tag{6.3.17}$$

$$a_{i+2\cdot 0}+\overline{a_{i+2\cdot 0}\,f'_{i,i+2}(a_i)(a_{i+2})}=a_{i+2\cdot 0}+f_{i,i+1}(a_i)(\overline{a_{i+2\cdot 0}\,a_{i+2}})\tag{6.3.18}$$

Из равенства (6.3.18) следует

$$\overline{a_{i+2\cdot 0}\,f'_{i,i+2}(a_i)(a_{i+2})}=f_{i,i+1}(a_i)(\overline{a_{i+2\cdot 0}\,a_{i+2}})\tag{6.3.19}$$

**Теорема 6.3.6.** *Рассмотрим башню $(\overline{V},\overline{f})$ представлений $\overline{\Omega}$-алгебр. Пусть представление $f_{i+1,i+2}$ однотранзитивно, и $a_{i+2\cdot 0}$ - базис представления $f_{i+1,i+2}$. Эндоморфизм $r_{i+1}$ представления $f_{i,i+1}$ порождает эндоморфизм $r'_{i+2}$ представления $f'_{i,i+2}$*

$$r'_{i+2}(a_{i+2})=a_{i+2\cdot 0}+r_{i+1}(\overline{a_{i+2\cdot 0}\,a_{i+2}})\tag{6.3.20}$$

**Замечание 6.3.7.** Из равенства (6.3.20) следует

$$\overline{a_{i+2\cdot 0}\,r'_{i+2}(a_{i+2})}=r_{i+1}(\overline{a_{i+2\cdot 0}\,a_{i+2}})\tag{6.3.21}$$

Доказательство. Поскольку отображение $r_{i+1}$ является эндоморфизмом $\Omega_{i+1}$-алгебры $A_{i+1}$, то, согласно теореме 6.3.3, отображение $r'_{i+2}$ является эндоморфизмом $\Omega_{i+1}$-алгебры $(A_{i+2},a_{i+2\cdot 0})$.

Для доказательства теоремы достаточно показать, что верно следующее равенство

$$r'_{i+2}\circ f'_{i,i+2}(a_i)=f'_{i,i+2}(r_i(a_i))\circ r'_{i+2}\tag{6.3.22}$$

Равенство (6.3.22) можно записать в виде

$$r'_{i+2}(f'_{i,i+2}(a_i)(a_{i+2}))=f'_{i,i+2}(r_i(a_i))(r'_{i+2}(a_{i+2}))\tag{6.3.23}$$

Из равенств (6.3.23), (6.3.16), (6.3.20) следует

$$\begin{array}{c}a_{i+2\cdot 0}+r_{i+1}(\overline{a_{i+2\cdot 0}\,f'_{i,i+2}(a_i)(a_{i+2})})\\ =a_{i+2\cdot 0}+f_{i,i+1}(r_i(a_i))(\overline{a_{i+2\cdot 0}\,r'_{i+2}(a_{i+2})})\end{array}\tag{6.3.24}$$

Из равенств (6.3.24), (6.3.19), (6.3.21) следует

$$r_{i+1}(f_{i,i+1}(a_i)(\overline{a_{i+2\cdot 0}\,a_{i+2}}))=f_{i,i+1}(r_i(a_i))(r_{i+1}(\overline{a_{i+2\cdot 0}\,a_{i+2}}))\tag{6.3.25}$$

Равенство (6.3.25) верно, так как отображение $(r_i,r_{i+1})$ является эндоморфизмом представления $f_{i,i+1}$. Следовательно, отображение $(r_i,r'_{i+2})$ является эндоморфизмом представления $f'_{i,i+2}$. □

## 6.4. Базис в ${}_*{}^*D$-аффинном пространстве

Рассмотрим $\star D$-аффинное пространство $\overrightarrow{V}$.

Абелевая группа $\overline{V}$ действует однотранзитивно на множестве $\overset{\circ}{V}$. Согласно теореме 6.3.1, базис множества $\overset{\circ}{V}$ относительно представления абелевой группы $\overline{V}$ состоит из одной точки. Эту точку обычно обозначают буквой $O$ и называют **началом системы координат $\star D$-аффинного пространства**. Следовательно, произвольную точку $A \in \overset{\circ}{V}$ можно представить с помощью вектора $\overline{OA} \in \overline{V}$

$$A = O + \overline{OA} \tag{6.4.1}$$

Согласно теореме 6.3.1, мощности множеств $\overset{\circ}{V}$ и $\overline{V}$ равны.

Тип аффинного пространства определяется типом векторного пространства в башне представлений (раздел [6]-4.2).

Пусть $\overline{V}$ - ${}_*{}^*D$-векторное пространство над телом $D$. $\star D$-аффинное пространство $\overrightarrow{V}$ называется **${}_*{}^*D$-аффинным пространством**. Пусть $\overline{\overline{e}}$ - базис ${}_*{}^*D$-векторного пространства $\overline{V}$. Тогда разложение вектора $\overrightarrow{OA}$ относительно ${}_*{}^*D$-базиса $\overline{\overline{e}}$ имеет вид

$$\overline{OA} = \overline{e}_i\, {}^iA \tag{6.4.2}$$

Из равенств (6.4.1), (6.4.2) следует

$$A = O + \overline{e}_i\, {}^iA \tag{6.4.3}$$

**Определение 6.4.1.** **${}_*{}^*D$-аффинный базис** $(\overline{\overline{e}}, O) = (\overline{e}_i, i \in I, O)$ - это максимальное множество ${}_*{}^*D$-линейно независимых векторов $\overline{e}_i = \overline{OA_i} = ({}^1e_i, ..., {}^ne_i)$ с общей начальной точкой $O$. Множество координат $({}^iA, i \in I)$ вектора $\overline{OA}$ называется **координатами точки $A$ ${}_*{}^*D$-аффинного пространства $\overrightarrow{V}$ относительно базиса** $(\overline{\overline{e}}, O)$. □

**Определение 6.4.2.** **Многообразие базисов** $\mathcal{B}(\overrightarrow{V})$ **${}_*{}^*D$-аффинного пространства** - это множество ${}_*{}^*D$-базисов этого пространства. □

**Теорема 6.4.3.** *Пусть*

$$\overrightarrow{V} : D \xrightarrow[*]{f_{1,2}} \overline{V} \xrightarrow[*]{f_{2,3}} \overset{\circ}{V} \qquad \begin{array}{l} f_{1,2}(d) : \overline{v} \to \overline{v}\, d \\ f_{2,3}(\overline{v}) : A \to A + \overline{v} \end{array}$$

*${}_*{}^*D$-аффинное пространство.*

(1) *Эндоморфизм представления $f_{2,3}$ имеет вид*

$$A' = A + \overline{a} \tag{6.4.4}$$

(2) *Эндоморфизм представления $f_{1,2}$ имеет вид*

$$\overline{a}' = \overline{a}_*{}^*P \tag{6.4.5}$$

Доказательство. Пусть $h$ - эндоморфизм представления $f_{2,3}$. Из определения [10]-2.2.2 следует, что

$$h \circ (A + \overline{b}) = h \circ A + \overline{b} \tag{6.4.6}$$

Пусть

$$B = A + \overline{b} \quad C = h \circ A = A + \overline{c} \quad D = h \circ B = B + \overline{d} \tag{6.4.7}$$

Из равенств (6.4.6), (6.4.7) следует

$$D = h \circ B = h \circ (A + \overline{b}) = h \circ A + \overline{b} = C + \overline{b} \tag{6.4.8}$$

Из равенств (6.4.7), (6.4.8) следует

$$\overline{AB} = \overline{CD} \tag{6.4.9}$$

Из аксиомы (4) определения 6.1.1 и равенства (6.4.9) следует

$$\overline{c} = \overline{AC} = \overline{BD} = \overline{d} \tag{6.4.10}$$

Поскольку $\overline{a}$ - произвольный вектор, то из равенств (6.4.10), (6.4.7) следует

$$D = B + \overline{c}$$

для любой точки $B \in \overset{\circ}{V}$. Следовательно, любой эндоморфизм представления $f_{2,3}$ имеет вид (6.4.4).

Утверждение (2) теоремы подробно рассмотрено в разделе [6]-4.4. □

**Теорема 6.4.4.** *Преобразование* (6.4.4) *называется* **параллельным сдвигом** ${}_*{}^*$**D-аффинного пространства**. *Координаты параллельного сдвига* (6.4.4) ${}_*{}^*D$*-аффинного пространства имеют вид*

$$ {}^iA' = {}^iA + {}^ia \tag{6.4.11}$$

Доказательство. Пусть $( {}^iA, i \in I)$ - координаты точки $A$ относительно базиса $(\overline{\overline{e}}, O)$. Из равенства (6.4.3) следует

$$A = O + \overline{e}_i\ {}^iA \tag{6.4.12}$$

Пусть $( {}^iA', i \in I)$ - координаты точки $A'$ относительно базиса $(\overline{\overline{e}}, O)$. Из равенства (6.4.3) следует

$$A' = O + \overline{e}_i\ {}^iA' \tag{6.4.13}$$

Пусть

$$\overline{a} = \overline{e}_i\ {}^ia \tag{6.4.14}$$

разложение вектора $\overline{a}$ относительно базиса $\overline{\overline{e}}$. Подставив равенства (6.4.12), (6.4.13), (6.4.14) into the equation (6.4.4), мы получим

$$O + \overline{e}_i\ {}^iA' = O + \overline{e}_i\ {}^iA + \overline{e}_i\ {}^ia = O + \overline{e}_i( {}^iA + {}^ia) \tag{6.4.15}$$

Поскольку представление $f_{2,3}$ однотранзитивно, то из равенства (6.4.15) следует

$$\overline{e}_i\ {}^iA' = \overline{e}_i( {}^iA + {}^ia) \tag{6.4.16}$$

Поскольку векторы $\overline{e}_i$ ${}_*{}^*D$-линейно независимы, то равенство (6.4.11) следует из равенства (6.4.16). □

**Теорема 6.4.5.** *Пусть* $(\overline{\overline{e}}, O)$ *- базис* ${}_*{}^*D$*-аффинного пространства*

$$\overrightarrow{V} : D \xrightarrow[*]{f_{1,2}} \overline{V} \xrightarrow[*]{f_{2,3}} \overset{\circ}{V} \qquad \begin{array}{l} f_{1,2}(d) : \overline{v} \rightarrow \overline{v}\ d \\ f_{2,3}(\overline{v}) : A \rightarrow A + \overline{v} \end{array}$$

(1) *Представление $f_{2,3}$ порождает операцию сложения на множестве* $\overset{\circ}{V}$

$$A+B=O+\overline{OA}+\overline{OB} \tag{6.4.17}$$

(2) *Представление $f_{1,2}$ порождает представление*

$$\begin{aligned}&f'_{1,3}:D\ {-}{*}{\longrightarrow}\ \overset{\circ}{V}\\&f'_{1,3}(d):A\to O+\overline{OA}\ d\end{aligned} \tag{6.4.18}$$

(3) *Эндоморфизм $\overline{g}$ представления $f_{1,2}$ порождает эндоморфизм $\overline{g}'$ представления $f'_{1,2}$*

$$\overline{g}':A\in\overset{\circ}{V}\to O+g_*{}^*\overline{OA} \tag{6.4.19}$$

Доказательство. Так как представление $f_{2,3}$ однотранзитивно, то

- утверждение (1) теоремы является следствием теоремы 6.3.2.
- утверждение (2) теоремы является следствием теоремы 6.3.5.
- утверждение (3) теоремы является следствием теоремы 6.3.6.

□

**Теорема 6.4.6.** *Преобразование* (6.4.19) *называется* **линейным преобразованием ${}_*{}^*D$-аффинного пространства**. *Координаты линейного преобразования* (6.4.19) *${}_*{}^*D$-аффинного пространства имеют вид*

$${}^iA'={}^ig_j\ {}^jA \tag{6.4.20}$$

Доказательство. Пусть $({}^iA, i\in I)$ - координаты точки $A$ относительно базиса $(\overline{\overline{e}}, O)$. Из равенства (6.4.2) следует

$$\overline{OA}=\overline{e}_i\ {}^iA \tag{6.4.21}$$

Пусть $({}^iA', i\in I)$ - координаты точки $A'$ относительно базиса $(\overline{\overline{e}}, O)$. Из равенства (6.4.3) следует

$$A'=O+\overline{e}_i\ {}^iA' \tag{6.4.22}$$

Из равенства (6.4.19) следует

$$A'=O+g_*{}^*\overline{OA} \tag{6.4.23}$$

Из равенств (6.4.21), (6.4.22), (6.4.23) следует

$$O+\overline{e}_i\ {}^iA'=O+\overline{e}_i\ {}^ig_j\ {}^jA \tag{6.4.24}$$

Равенство (6.4.20) следует из равенства (6.4.24), так как $(\overline{e}, O)$ - базис ${}_*{}^*D$-аффинного пространства. □

Последовательное выполнение параллельного сдвига и линейного преобразования ${}_*{}^*D$-аффинного пространства, порождает **${}_*{}^*D$-аффинное преобразование**

$$A'=O+a_*{}^*\overline{OA}+\overline{v} \tag{6.4.25}$$

Это преобразование ${}_*{}^*D$-аффинного пространства можно представить как отображение базиса $(\overline{\overline{e}}, O)$ в базис $(\overline{\overline{e}}', O')$. Отображение базиса $(\overline{\overline{e}}, O)$ в базис $(\overline{\overline{e}}', O')$ можно представить как композицию отображения базиса $(\overline{\overline{e}}, O)$ в базис $(\overline{\overline{e}}, O')$ и $D_*{}^*$-линейного отображения базиса $(\overline{\overline{e}}, O')$ в базис $(\overline{\overline{e}}', O')$. Следовательно, $D_*{}^*$-преобразование имеет вид

$$\overline{\overline{e}}' = A_*{}^*\overline{\overline{e}}$$
$$O' = O + \overrightarrow{OO'}$$

Вводя координаты $A^1, ..., A^n$ точки $A \in \mathcal{A}_n$ как координаты вектора $\overline{O\,A}$ относительно базиса $\overline{\overline{e}}$, мы можем записать линейное преобразование как

(6.4.26) $$A'^i = A^j\ {}_jP^i + R^i \qquad \operatorname{rank}_{*^*} P = n$$

(6.4.27) $$A'_*{}^*\overline{e} = (A_*{}^*P + R)_*{}^*\overline{e}$$

Вектор $(R^1, ..., R^n)$ выражает смещение в аффинном пространстве.

**Теорема 6.4.7.** *Множество преобразований* (6.4.27) *- это группа Ли, которую мы обозначим* $GL(\mathcal{A}_n)$ *и будем называть* **группой аффинных преобразований**.

Доказательство. Рассмотрим преобразования $(P, R)$ и $(Q, S)$. Последовательное выполнение этих преобразований имеет вид

$$\begin{aligned} A''^i &= A'^j\ {}_jQ^i + S^i & A'' &= (A'_*{}^*Q + S)_*{}^*\overline{e} \\ &= (A^k\ {}_kP^j + R^j)\ {}_jQ^i + S^i & &= ((A_*{}^*P + R)_*{}^*Q + S)_*{}^*\overline{e} \\ &= A^k\ {}_kP^j\ {}_jQ^i + R^j\ {}_jQ^i + S^i & &= (A_*{}^*P_*{}^*Q + R_*{}^*Q + S)_*{}^*\overline{e} \end{aligned}$$

Следовательно, мы можем записать произведение преобразований $(P, R)$ и $(Q, S)$ в виде

(6.4.28) $$(P, R)_*{}^*(Q, S) = (P_*{}^*Q, R_*{}^*Q + S)$$

□

Мы будем называть активное преобразование **аффинным преобразованием**.Мы будем называть пассивное преобразование **квазиаффинным преобразованием**.

Если мы не заботимся о начальной точке вектора, мы получим несколько отличный тип пространства, которое мы будем называть центро-аффинным пространством $CA_n$. Если мы предположим, что начальная точка вектора - это начало $O$ координатной системы в пространстве, то мы можем отождествить любую точку $A \in CA_n$ с вектором $a = \overline{O\,A}$. Теперь преобразование - это просто отображение

$$a'^i = P^i_j a^j \qquad \det P \neq 0$$

и такие преобразования порождают группу Ли $GL_n$.

**Определение 6.4.8. Центро-аффинный базис** $\overline{\overline{e}} = (\overline{e}_i)$ - это множество линейно независимых векторов $\overline{e}_i$. □

**Определение 6.4.9. Многообразие базисов** $\mathcal{B}(CA_n)$ **центро-аффинного пространства** - это множество базисов этого пространства. □

## 6.5. Плоскость в ${}_*{}^*D$-аффинном пространстве

Рассмотрим ${}_*{}^*D$-аффинное пространство

$$\overset{\to}{V}:D\xrightarrow[*]{f_{1,2}}\overline{V}\xrightarrow[*]{f_{2,3}}\overset{\circ}{V}\qquad \begin{array}{l}f_{1,2}(d):\overline{v}\to\overline{v}\ d\\ f_{2,3}(\overline{v}):A\to A+\overline{v}\end{array}$$

Пусть $\overline{V}_1\subset\overline{V}$ - ${}_*{}^*D$-векторное подпространство ${}_*{}^*D$-векторного пространства $\overline{V}$. Следовательно, $\overline{V}_1$ является подгруппой абелевой группы $\overline{V}$. В башне подпредставлений

$$D\xrightarrow[*]{f_{1,2}}\overline{V}_1\xrightarrow[*]{f_{2,3}}\overset{\circ}{V}\qquad \begin{array}{l}f_{1,2}(d):\overline{v}\to\overline{v}\ d\\ f_{2,3}(\overline{v}):A\to A+\overline{v}\end{array}$$

представление $f_{2,3}$ не является однотранзитивным. Согласно теореме [6]-3.4.13 множество $\overset{\circ}{V}$ может быть представленно как пересечение не пересекающихся орбит представления $f_{2,3}$. Однако, так как представление $f_{2,3}$ эффективно, то представление $f_{2,3}$ однотранзитивно на орбите представления. Следовательно, если мы выберем произвольную точку $A\in\overset{\circ}{V}$ и орбиту представления

$$\overset{\circ}{V}_1=A+f_{2,3}(\overline{V}_1)$$

то мы получим башню представлений

$$\overset{\to}{V}_1:D\xrightarrow[*]{f_{1,2}}\overline{V}_1\xrightarrow[*]{f_{2,3}}\overset{\circ}{V}_1\qquad \begin{array}{l}f_{1,2}(d):\overline{v}\to\overline{v}\ d\\ f_{2,3}(\overline{v}):A\to A+\overline{v}\end{array}$$

описывающую ${}_*{}^*D$-аффинное пространство, называемое **${}_*{}^*D$-аффинной плоскостью**, проходящей через точку $A$.

Пусть

$$\overline{\overline{e}}=(\overline{e}_1,...,\overline{e}_n)$$

базис ${}_*{}^*D$-векторного пространства $\overline{V}$. Пусть базис ${}_*{}^*D$-векторного пространства $\overline{V}_1$

$$\overline{\overline{e}}_1=(\overline{e}_{1\cdot 1},...,\overline{e}_{1\cdot k})$$

имеет следующее разложение относительно базиса $\overline{\overline{e}}$

$$\overline{e}_{1\cdot i}=\overline{e}_j\ {}^je_{1\cdot i}$$

Пусть $\overline{v}\in\overline{V}_1$ имеет следующее разложение относительно базиса $\overline{\overline{e}}$

$$\overline{v}=\overline{e}_j\ {}^jv$$

Согласно теореме [6]-4.6.7

$$\text{rank}\begin{pmatrix}{}^1e_{1\cdot 1}&...&{}^1e_{1\cdot k}&{}^1v\\ ...&...&...&...\\ {}^ne_{1\cdot 1}&...&{}^ne_{1\cdot k}&{}^nv\end{pmatrix}=k\tag{6.5.1}$$

Если мы предположим, что первые $k$ строк составляют невырожденную матрицу, то мы получим $n-k$ ${}_*{}^*D$-линейных уравнений

$${}^{k+1}\det({}_*{}^*)_r\begin{pmatrix}{}^1e_{1\cdot 1}&...&{}^1e_{1\cdot k}&{}^1v\\ ...&...&...&...\\ {}^ne_{1\cdot 1}&...&{}^ne_{1\cdot k}&{}^nv\end{pmatrix}=0\qquad r=k+1,...,n\tag{6.5.2}$$

которым удовлетворяют координаты вектора $\overline{v}$.

Пусть точка $A$ имеет представление

(6.5.3) $$A = O + e_i\ {}^iA$$

относительно базиса $(\overline{\overline{e}}, O)$. Пусть

(6.5.4) $$B = O + e_i\ {}^iB$$

произвольная точка ${}_*{}^*D$-аффинной плоскости $\overrightarrow{V}_1$. Тогда

(6.5.5) $$\overline{v} = \overline{A\,B} = e_i\ {}^iv \quad \overline{v} \in \overline{V}_1$$

Из равенства

$$B = A + \overline{v}$$

и равенств (6.5.3), (6.5.4), (6.5.5) следует

(6.5.6) $${}^iB = {}^iA + {}^iv$$

Из равенств (6.5.2), (6.5.6) следует, что координаты точки $B$ удовлетворяют системе ${}_*{}^*D$-линейных уравнений

(6.5.7) $${}^{k+1}\det({}_*{}^*)_r \begin{pmatrix} {}^1e_{1\cdot 1} & ... & {}^1e_{1\cdot k} & {}^1B - {}^1A \\ ... & ... & ... & ... \\ {}^ne_{1\cdot 1} & ... & {}^ne_{1\cdot k} & {}^nB - {}^nA \end{pmatrix} = 0 \quad r = k+1, ..., n$$

Пусть задано семейство векторов

$$\overline{a}_i = \overline{A_0\,A_i} \quad i = 1, ..., n-1$$

Мы говорим, что вектор $\overline{a} = \overrightarrow{A_0A}$ лежит в плоскости, натянутой на векторы $\overline{a}_1, ..., \overline{a}_{n-1}$ , если вектор $\overline{a}_n$ является линейной комбинацией векторов $\overline{a}_1, ..., \overline{a}_{n-1}$ .

Глава 7

# Евклидово пространство

## 7.1. Евклидово пространство

Пусть $\overline{V}$ - $D$-векторное пространство. Билинейное отображение

$$g:\overline{V}\times\overline{V}\rightarrow D$$

можно представить в виде

$$g(\overline{v},\overline{w})={}_{ij}g(v^i,w^j)$$

где ${}_{ij}g$ - билинейное отображение

$${}_{ij}g:D^2\rightarrow D$$

**Определение 7.1.1.** Билинейное отображение

$$g:\overline{V}\times\overline{V}\rightarrow D$$

называется **симметричным**, если

$$g(\overline{w},\overline{v})=g(\overline{v},\overline{w})$$

□

**Определение 7.1.2.** Билинейное отображение

$$g:\overline{V}\times\overline{V}\rightarrow D$$

называется **эвклидовым скалярным произведением** в $D$-векторном пространстве $\overline{V}$, если ${}_{ii}g$ - эвклидово скалярное произведение в теле $D$ и ${}_{ij}g=0$ для $i\neq j$. □

**Определение 7.1.3.** Билинейное отображение

$$g:\overline{V}\times\overline{V}\rightarrow D$$

называется **псевдоэвклидовым скалярным произведением** в $D$-векторном пространстве $\overline{V}$, если ${}_{ii}g$ - псевдоэвклидово скалярное произведение в теле $D$ и ${}_{ij}g=0$ для $i\neq j$. □

Если в теле $D$ определена операция эрмитова сопряжения, то мы можем эту операцию распространить на $D$-векторное пространство $\overline{V}$, а именно для данного вектора $\overline{v}=v^i\ {}_i\overline{e}$ мы определим **эрмитово сопряжённый вектор**

$$\overline{v}^\star=(v^i)^\star\ {}_i\overline{e}$$

**Определение 7.1.4.** Билинейное отображение

$$g:\overline{V}\times\overline{V}\rightarrow D$$

называется **эрмитовым скалярным произведением** в $D$-векторном пространстве $\overline{V}$, если ${}_{ii}g$ - эрмитово скалярное произведение в теле $D$ и ${}_{ij}g=0$ для $i\neq j$. □

## 7.2. Базис в евклидовом пространстве

Когда мы определяем метрику в центро-аффинном пространстве, мы получаем новую геометрию потому, что мы можем измерять расстояние и длину вектора. Если метрика положительно определена, мы будем называть пространство евклидовым $\mathcal{E}_n$, в противном случае мы будем называть пространство псевдоевклидовым $\mathcal{E}_{nm}$.

Преобразования, которые сохраняют длину, образуют группу Ли $SO(n)$ для евклидова пространства и группу Ли $SO(n, m)$ для псевдоэвклидова пространства, где $n$ и $m$ числа положительных и отрицательных слагаемых в метрике.

**Определение 7.2.1. Ортонормальный базис** $\overline{\overline{e}} = (\overline{e}_i)$ - это множество линейно независимых векторов $\overline{e}_i$ таких, что длина каждого вектора равна 1 и различные векторы ортогональны. □

**Определение 7.2.2. Многообразие базисов** $\mathcal{B}(\mathcal{E}_n)$ **евклидова пространства** - это множество ортонормальных базисов этого пространства. □

Мы будем называть активное преобразование **движением**. Мы будем называть пассивное преобразование **квазидвижением**.

Глава 8

# Математический анализ в $D$-аффинном пространстве

## 8.1. Криволинейные координаты в $D$-аффинном пространстве

Пусть $\overrightarrow{A}$ - $D$-аффинное пространство над непрерывным тело $D$.[8.1]

Пусть $({}_1\overline{e}, ..., {}_n\overline{e}, O)$ - $D_*{}^*$-базис аффинного пространства. Для каждого $i$, $i = 1$, ..., $i = n$, определим отображение[8.2]

(8.1.1) $$\overline{x}_i : D \to \overset{\circ}{A} \quad v^i \to O + v^i\ {}_i\overline{e}$$

Для произвольной точки $M(v^i)$ мы можем записать вектор $\overrightarrow{OM}$ в виде линейной комбинации

$$\overrightarrow{OM} = \overline{x}_i(v^i)$$

Мы можем рассматривать множество координат аффинного пространства $\overset{\circ}{A}$ как гомеоморфизм

$$f : \overset{\circ}{A} \to D^n$$

Если гомеоморфизм

$$f : \overset{\circ}{A} \to D^n$$

является произвольным отображением открытого подмножества $D$-аффинного пространства в открытое подмножество множества $D^n$, то мы будем говорить, что образ точки $A$ при отображении $f$ является **криволинейными координатами точки** $A$, а отображение $f$ будем называть системой координат. Если даны две системы координат, $f_1$ и $f_2$, то отображение $g = f_1 f_2^{-1}$ является гомеоморфизмом

$$g : D^n \to D^n$$
$$x'^i = x'^i(x^1, ..., x^n)$$

и называется преобразованием координат. Преобразование координат, порождённое аффинным преобразованием, можно рассматривать как частный случай.

**Замечание 8.1.1.** Рассмотрим поверхность $S$ в аффинном пространстве, определённую уравнением

$$f : D^n \to D$$
$$f(x^1, ..., x^n) = 0$$

[8.1] В этом разделе я изучаю криволинейные координаты аналогично тому, как это сделано в [3], глава V.

[8.2] В равенстве (8.1.1) нет суммы по индексу $i$.

Пусть $M(x^i)$ - произвольная точка. Малое изменение координат точки

$$x'^i = x^i + \Delta x^i$$

приводит к малому изменению функции

$$\Delta f = \partial f(\overline{x})(\Delta\overline{x})$$

В частности, если бесконечно близкие точки $M(x^i)$ и $L(x^i + \Delta x^i)$ лежат на одной поверхности, то вектор приращения $\overrightarrow{ML} = \Delta\overline{x}$ удовлетворяет дифференциальному уравнению

$$\partial f(\overline{x})(\Delta\overline{x}) = 0$$

и называется вектором, касательным к поверхности $S$.

Так как производная Гато отображения $f$ в точке $X$ является линейным отображением, то эта производная является 1-$D$-формой. Рассмотрим множество 1-$D$-форм

$$dx^i = \frac{\partial x^i}{\partial\overline{x}}$$

Очевидно, что

$$dx^i(\overline{a}) = \frac{\partial x^i}{\partial\overline{x}}(\overline{a}) = (x^i + a^i) - x^i = a^i$$

Производную Гато произвольной функции можно представить в виде

$$\frac{\partial f(\overline{x})}{\partial\overline{x}} = \frac{\partial f(\overline{x})}{\partial x^i}\left(\frac{\partial x^i}{\partial\overline{x}}\right) = \frac{\partial f(\overline{x})}{\partial x^i}\left(dx^i\right)$$

$$\frac{\partial f(\overline{x})}{\partial\overline{x}}(\overline{a}) = \frac{\partial f(\overline{x})}{\partial x^i}\left(\frac{\partial x^i}{\partial\overline{x}}(\overline{a})\right) = \frac{\partial f(\overline{x})}{\partial x^i}\left(dx^i(\overline{a})\right) = \frac{\partial f(\overline{x})}{\partial x^i}(a^i) \tag{8.1.2}$$

Следовательно, производная Гато функции $f$ является линейной комбинацией 1-$D$-форм $dx^i$. Пусть

$$g : D^n \rightarrow D^n$$
$$x'^i = x'^i(x^1, ..., x^n)$$

преобразование координат аффинного пространства. Из правила дифференцирования сложной функции (теорема [9]-6.2.12) следует

$$dx'^i = \frac{\partial x'^i}{\partial\overline{x}} = \frac{\partial x'^i}{\partial x^j}\left(\frac{\partial x^j}{\partial\overline{x}}\right) = \frac{\partial x'^i}{\partial x^j}\left(dx^j\right)$$

Координатная линия - это кривая, вдоль которой изменяется только одна координата $x^i$. Мы будем рассматривать координатную кривую как кривую, имеющую касательный вектор

$$\overline{x}_i(\overline{x})(v) = \frac{\partial\overline{x}}{\partial x^i}(v) \tag{8.1.3}$$

Пусть точка $O$ имеет координаты $\overline{x}$. Множество векторов $\overline{x}_1$, ..., $\overline{x}_n$ является базисом аффинного пространства

$$(\overline{x}_1(\overline{x}), ..., \overline{x}_n(\overline{x}), \overline{x})$$

$$\overline{v} = \overline{x}_i(\overline{x})(v^i) \tag{8.1.4}$$

Пусть

$$g : D^n \rightarrow D^n$$
$$x'^i = x'^i(x^1, ..., x^n)$$

преобразование координат аффинного пространства. Из правила дифференцирования сложной функции (теорема [9]-6.2.12) следует

$$\overline{x}'_i(\overline{x}')(v'^i) = \frac{\partial \overline{x}}{\partial x'^i}(v'^i) = \frac{\partial \overline{x}}{\partial x^j}\left(\frac{\partial x^j}{\partial x'^i}(v'^i)\right) = \overline{x}_j(\overline{x})\left(\frac{\partial x^j}{\partial x'^i}(v'^i)\right) \tag{8.1.5}$$

Из сравнения равенств (8.1.4), (8.1.5) следует закон преобразования координат вектора

$$v^j = \frac{\partial x^j}{\partial x'^i}(v'^i) \tag{8.1.6}$$

**Пример 8.1.2.** Рассмотрим линейное преобразование координат

$$\begin{aligned} x'^1 &= ax^1b + ax^2c \\ x'^2 &= x^1 + x^2 \end{aligned} \tag{8.1.7}$$

Тогда

$$x^1 = x'^2 - x^2 \tag{8.1.8}$$

$$x'^1 = a(x'^2 - x^2)b + ax^2c = ax'^2b + ax^2(c-b)$$

$$ax^2(c-b) = x'^1 - ax'^2b \tag{8.1.9}$$

Из равенств (8.1.8), (8.1.9) следует обратное преобразование

$$\begin{aligned} x^2 &= a^{-1}x'^1(c-b)^{-1} - x'^2b(c-b)^{-1} \\ x^1 &= x'^2 - a^{-1}x'^1(c-b)^{-1} + x'^2b(c-b)^{-1} \\ &= -a^{-1}x'^1(c-b)^{-1} + x'^2(1 + b(c-b)^{-1}) \end{aligned} \tag{8.1.10}$$

Из равенств (8.1.7) следует закон преобразования 1-$D$-форм

$$dx'^1 = adx^1b + adx^2c$$

$$dx'^2 = dx^1 + dx^2$$

Из равенств (8.1.6), (8.1.10) следует закон преобразования координат вектора

$$v^1 = -a^{-1}v'^1(c-b)^{-1} + v'^2(1 + b(c-b)^{-1})$$

$$v^2 = a^{-1}v'^1(c-b)^{-1} - v'^2b(c-b)^{-1}$$

□

## 8.2. Параллельный перенос

Пусть в точке $M_0$ вектор $\overline{v}$ имеет координаты $v_0^i$. Пусть вектор $\overrightarrow{M_0M_1}$ является вектором бесконечно малого смещения. Пусть $N_0 = M_0 + \overline{v}$, $N_1 = M_1 + \overline{v}$. Из аксиомы 3 определения 6.1.1 координаты векторов $\overrightarrow{M_0N_0}$ и $\overrightarrow{M_1N_1}$ равны. Как это построение будет выглядеть в криволинейных координатах?

Мы полагаем функции $v^i(\overline{x})$ непрерывно дифференцируемыми в смысле Гато. Следовательно, векторы локального репера $\overline{x}_i$ и координаты вектора $v^i$ являются непрерывно дифференцируемыми в смысле Гато функциями. В точке $M(\overline{x})$ справедливо равенство

$$\overline{v} = \overline{x}_i(\overline{x})(v^i(\overline{x})) \tag{8.2.1}$$

Так как $\overline{v} = \mathrm{const}$, то, дифференцируя равенство (8.2.1) по $\overline{x}$, получим

$$0 = \frac{\partial \overline{x}_i(\overline{x})}{\partial \overline{x}}(v^i)(\overline{a}) + \overline{x}_i(\overline{x})\left(\frac{\partial v^i}{\partial \overline{x}}(\overline{a})\right) \tag{8.2.2}$$

Из равенств (8.1.2), (8.2.2) следует

$$0 = \frac{\partial \overline{x}_i(\overline{x})}{\partial x^j}(v^i)(a^j) + \overline{x}_i(\overline{x})\left(\frac{\partial v^i}{\partial x^j}(a^j)\right) \tag{8.2.3}$$

Выражение $\dfrac{\partial \overline{x}_i(\overline{x})}{\partial x^j}(v^i)(a^j)$ является вектором аффинного пространства и в локальном базисе $\overline{x}_i(\overline{x})$ это выражение имеет разложение

$$\frac{\partial \overline{x}_i(\overline{x})}{\partial x^j}(v^i)(a^j) = \overline{x}_k(\overline{x})(\Gamma^k_{ji}(v^i)(a^j)) \tag{8.2.4}$$

где $\Gamma^k_{ji}$ - билинейные отображения, называемые **коэффициентами связности в $D$-аффинном пространстве**. Из равенств (8.2.4), (8.2.3) следует

$$0 = \overline{x}_k(\overline{x})(\Gamma^k_{ji}(v^i)(a^j)) + \overline{x}_i(\overline{x})\left(\frac{\partial v^i}{\partial x^j}(a^j)\right) \tag{8.2.5}$$

Так как $\overline{x}_k(\overline{x})$ - локальный базис, то из равенства (8.2.5) следует

$$0 = \overline{x}_k(\overline{x})\left(\Gamma^k_{ji}(v^i)(a^j) + \frac{\partial v^k}{\partial x^j}(a^j)\right)$$

$$0 = \Gamma^k_{ji}(v^i)(a^j) + \frac{\partial v^k}{\partial x^j}(a^j) \tag{8.2.6}$$

Так как $\overline{a}$ - произвольный вектор, то равенство (8.2.6) можно записать в виде

$$0 = \Gamma^k_{ji}(v^i) + \frac{\partial v^k}{\partial x^j} \tag{8.2.7}$$

**Теорема 8.2.1.** *В аффинном пространстве*

$$\Gamma^k_{ji}(v^i)(a^j) = \Gamma^k_{ij}(v^i)(a^j) \tag{8.2.8}$$

Доказательство. Согласно равенству (8.1.3)

$$\frac{\partial \overline{x}_i(\overline{x})}{\partial x^j}(v^i)(a^j) = \frac{\partial^2 \overline{x}}{\partial x^j \partial x^i}(v^i; a^j) \tag{8.2.9}$$

Согласно теореме [9]-9.1.6

$$\frac{\partial^2 \overline{x}}{\partial x^j \partial x^i}(v^i; a^j) = \frac{\partial^2 \overline{x}}{\partial x^i \partial x^j}(v^i; a^j) \tag{8.2.10}$$

Равенство (8.2.8) следует из равенств (8.2.4), (8.2.10). □

Пусть

$$g : D^n \to D^n$$
$$x'^i = x'^i(x^1, ..., x^n)$$

преобразование координат аффинного пространства. Продифференцируем равенство (8.1.5) по $x'^j$

$$\begin{aligned}
&\frac{\partial\overline{x}'_i(\overline{x}')}{\partial x'^j}(v'^i)(a'^j)\\
&=\frac{\overline{x}_k(\overline{x})\left(\frac{\partial x^k}{\partial x'^i}(v'^i)\right)}{\partial x^m}\left(\frac{\partial x^m}{\partial x'^j}(a'^j)\right)+\overline{x}_k(\overline{x})\left(\frac{\partial^2 x^k}{\partial x'^j\partial x'^i}(v'^i;a'^j)\right)
\end{aligned}\tag{8.2.11}$$

Из равенств (8.2.4), (8.2.11) следует

$$\begin{aligned}
&\overline{x}'_p(\overline{x}')(\Gamma'^p_{ji}(v'^i)(a'^j))\\
&=\overline{x}_r(\overline{x})\left(\Gamma^r_{mk}\left(\frac{\partial x^k}{\partial x'^i}(v'^i)\right)\left(\frac{\partial x^m}{\partial x'^j}(a'^j)\right)\right)+\overline{x}_r(\overline{x})\left(\frac{\partial^2 x^r}{\partial x'^j\partial x'^i}(v'^i;a'^j)\right)\\
&=\overline{x}_r(\overline{x})\left(\Gamma^r_{mk}\left(\frac{\partial x^k}{\partial x'^i}(v'^i)\right)\left(\frac{\partial x^m}{\partial x'^j}(a'^j)\right)+\frac{\partial^2 x^r}{\partial x'^j\partial x'^i}(v'^i;a'^j)\right)
\end{aligned}\tag{8.2.12}$$

Из равенств (8.1.5), (8.2.12) следует

$$\begin{aligned}
&\overline{x}'_p(\overline{x}')(\Gamma'^p_{ji}(v'^i)(a'^j))\\
&=\overline{x}'_p(\overline{x})\left(\frac{\partial x'^p}{\partial x^r}\left(\Gamma^r_{mk}\left(\frac{\partial x^k}{\partial x'^i}(v'^i)\right)\left(\frac{\partial x^m}{\partial x'^j}(a'^j)\right)+\frac{\partial^2 x^r}{\partial x'^j\partial x'^i}(v'^i;a'^j)\right)\right)
\end{aligned}\tag{8.2.13}$$

Из равенства (8.2.13) следует

$$\begin{aligned}
&\Gamma'^p_{ji}(v'^i)(a'^j)\\
&=\frac{\partial x'^p}{\partial x^r}\left(\Gamma^r_{mk}\left(\frac{\partial x^k}{\partial x'^i}(v'^i)\right)\left(\frac{\partial x^m}{\partial x'^j}(a'^j)\right)+\frac{\partial^2 x^r}{\partial x'^j\partial x'^i}(v'^i;a'^j)\right)
\end{aligned}\tag{8.2.14}$$

Глава 9

# Многообразие $D$-аффинной связности

Основная задача этой главы - дать общее представление о тех задачах, которые необходимо будет решать в процессе изучения дифференциальной геометрии.

## 9.1. Многообразие $D$-аффинной связности

**Определение 9.1.1.** Расслоенное центральное $D$-аффинное пространство называется **многообразием $D$-аффинной связности**. □

Мы полагаем, что размерность базы равна размерности слоя. Это позволяет отождествить слой с касательным пространством к базе. Рассмотрим область $U$ многообразия, в которой определён гомеоморфизм

$$\overline{x}:D^n\to U$$

Так же как и в случае аффинного пространства мы в каждом слое определим локальный базис

$$\overline{x}_k:\frac{\partial\overline{x}}{\partial x^k}$$

Мы полагаем, что отображение между слоями является морфизмом центрального $D$-аффинного пространства

$$d\overline{v}=\Gamma(\overline{v})(d\overline{x})$$

где билинейное отображение $\Gamma$ называется **$D$-аффинной связностью**. В локальном базисе связность можно представить как множество билинейных функций $\Gamma^k_{ij}$

$$\Gamma(\overline{v})(\overline{a})=\overline{x}_k(\overline{x})(\Gamma^k_{ji}(v^i)(a^j))$$

которые называются **коэффициентами $D$-аффинной связности**. При переходе от системы координат $\overline{x}$ к системе координат $\overline{x}'$ коэффициенты $D$-аффинной связности преобразуются согласно правилу

$$\begin{aligned}&\Gamma'^p_{ji}(v'^i)(a'^j)\\=&\frac{\partial x'^p}{\partial x^r}\left(\Gamma^r_{mk}\left(\frac{\partial x^k}{\partial x'^i}(v'^i)\right)\left(\frac{\partial x^m}{\partial x'^j}(a'^j)\right)+\frac{\partial^2x^r}{\partial x'^j\partial x'^i}(v'^i;a'^j)\right)\end{aligned}$$

Ковариантная производная векторного поля $\overline{v}$ имеет вид

$$D(\overline{v})(d\overline{x})=\partial(\overline{v})(d\overline{x})-\Gamma(\overline{v})(d\overline{x})$$

Векторное поле $\overline{v}$ параллельно переносится в данном направлении, если

$$\partial(\overline{v})(d\overline{x})=\Gamma(\overline{v})(d\overline{x})$$

Если $\overline{x}=\overline{x}(t)$ - отображение поля действительных чисел в рассматриваемое пространство, то это отображение определяет геодезическую, если касательный вектор параллельно переносится вдоль геодезической

$$\partial(\partial\overline{x}(d\overline{x}))(d\overline{x})=\Gamma(d\overline{x})(d\overline{x}) \tag{9.1.1}$$

Уравнение (9.1.1) эквивалентно уравнению

$$\partial^2\overline{x}(d\overline{x};d\overline{x})=\Gamma(d\overline{x})(d\overline{x}) \tag{9.1.2}$$

Глава 10

# Список литературы

[1] Серж Ленг, Алгебра, М. Мир, 1968

[2] S. Burris, H.P. Sankappanavar, A Course in Universal Algebra, Springer-Verlag (March, 1982), eprint http://www.math.uwaterloo.ca/ snburris/htdocs/ualg.html (The Millennium Edition)

[3] П. К. Рашевский, Риманова геометрия и тензорный анализ, М., Наука, 1967

[4] А. Г. Курош, Курс высшей алгебры, М., Наука, 1968

[5] А. Г. Курош, Общая алгебра, (лекции 1969 - 70 учебного года), М., МГУ, 1970

[6] Александр Клейн, Лекции по линейной алгебре над телом, eprint arXiv:math.GM/0701238 (2010)

[7] Александр Клейн, Расслоенное соответствие, eprint arXiv:0707.2246 (2007)

[8] Александр Клейн, Морфизм $T\star$-представлений, eprint arXiv:0803.2620 (2008)

[9] Александр Клейн, Введение в математический анализ над телом, eprint arXiv:0812.4763 (2010)

[10] Aleks Kleyn, Representation Theory: Representation of Universal Algebra, Lambert Academic Publishing, 2011

[11] John C. Baez, An Introduction to n-Categories, eprint arXiv:q-alg/9705009 (1997)

[12] П. Кон, Универсальная алгебра, М., Мир, 1968

[13] Paul Bamberg, Shlomo Sternberg, A course in mathematics for students of physics, Cambridge University Press, 1991

[14] Alain Connes, Noncommutative Geometry, Academic Press, 1994

Глава 11

# Предметный указатель

Глава 12

# Специальные символы и обозначения